\documentclass[10pt, reqno]{amsart}

\usepackage{preamble}

\title[]{Prismatic Steenrod operations and \\ arithmetic duality on Brauer groups}

\author{Shachar Carmeli and Tony Feng}

\begin{document}  

\begin{abstract}
We construct and analyze the ``syntomic Steenrod algebra'', which acts on the mod $p$ syntomic cohomology (also known as \'etale-motivic cohomology) of algebraic varieties in characteristic $p$. We then apply the resulting theory to resolve the last open cases of a 1966 Conjecture of Tate, concerning the existence of a symplectic form on the Brauer groups of smooth proper surfaces over finite fields. More generally, we exhibit symplectic structure on the higher Brauer groups of even dimensional varieties over finite fields.

Although the applications are classical, our methods rely on recent advances in perfectoid geometry and prismatic cohomology, which we employ to define a theory of ``spectral syntomic cohomology'' with coefficients in motivic spectra. We then organize the resulting cohomology theories into a category of ``spectral prismatic $F$-gauges'', generalizing the prismatic $F$-gauges of Drinfeld and Bhatt--Lurie, for which we establish a ``spectral Serre duality'' extending classical coherent duality. These abstract constructions are leveraged to explicitly compute the syntomic Steenrod operations. 
\end{abstract}

\maketitle

\tableofcontents

\section{Introduction}

\subsection{Classical motivations}\label{ssec: intro classical}At the 1962 ICM, Tate announced his famous \emph{arithmetic duality theorems} \cite{Tate62}, which exhibited parallels between: 
\begin{itemize}
\item the \'etale cohomology of global fields and the singular cohomology of 3-manifolds, and 
\item the \'etale cohomology of local fields and the singular cohomology of 2-manifolds.
\end{itemize}
Using these results, Tate constructed a skew-symmetric perfect pairing on (the non-divisible quotient of) the \emph{Tate--Shafarevich group} of a principally polarized abelian variety over a global field, and proved in \cite[Theorem 3.3]{Tate62} that it is alternating under a technical hypothesis (which is always satisfied for elliptic curves). Therefore, under the hypothesis, on this non-divisible quotient group there is a symplectic form, now known as the \emph{Cassels--Tate pairing}, which forces the order of the group to be a perfect square\footnote{We are implicitly invoking the fact that a finite abelian group with a symplectic form must have order a perfect square.} -- a numerology which had been empirically discovered by Selmer and proved by Cassels for elliptic curves. In fact, Venkatesh has suggested that Tate's discovery of arithmetic duality was motivated by a desire to explain this very numerology via a symplectic duality, which he realized should come from an arithmetic analogue of Poincar\'e duality wherein local and global fields played the role of manifolds.

A few years later, Tate's 1966 Bourbaki seminar \cite{Tate66} generalized the Birch and Swinnerton-Dyer Conjecture to abelian varieties and formulated its geometric analogue -- now called the \emph{Artin--Tate Conjecture} -- for a smooth, proper, geometrically connected surface $X$ over a finite field $k$ of characteristic $p$. Moreover, Tate conjectured that there should be an analogous symplectic form on the \emph{Brauer group} $\Br(X) := \rH^2_{\et}(X; \G_m)$, a torsion abelian group, conjecturally finite, which is closely connected to the Tate--Shafarevich group.\footnote{For example, work of Milne \cite{Milne75, Milne86} and Kato--Trihan \cite{KT} shows that the finiteness of $\Br(X)$, for all $X$, is equivalent to finiteness of the Tate--Shafarevich group plus the Birch and Swinnerton-Dyer Conjecture for all Jacobian varieties over function fields.} This conjecture has been studied in many works over the decades since, notably \cite{Manin67, Manin86, Urabe96, LLR05, Feng20}, and in this paper we will settle its remaining open cases. 

\subsubsection{The order of the Brauer group} Let $X$ be a smooth, proper, geometrically connected surface over a finite field $k$ of characteristic $p$. According to the Artin--Tate Conjecture, $\# \Br(X)$ is the main invariant featuring in a special value formula for the leading order term of the zeta function of $X$ at the center of its functional equation. Motivated by parallel properties of Tate--Shafarevich groups, Tate made the following Conjecture about $\# \Br(X)$. 

\begin{conj}[Tate, \cite{Tate66}]\label{conj: size} The order of $\Br(X)$ is a (finite) perfect square. 
\end{conj}

Although the finiteness of $\Br(X)$ is unknown, we know unconditionally that its \emph{non-divisible quotient} $\Br(X)_{\nd}$ -- the quotient of $\Br(X)$ by its subgroup of divisible elements -- is finite, and is equal to $\Br(X)$ if the latter is finite. Therefore, an unconditional version of Conjecture \ref{conj: size} is:

\begin{conj}[Tate, \cite{Tate66}]\label{conj: order 2} The order of $\Br(X)_{\nd}$ is a perfect square. 
\end{conj}

Conjecture \ref{conj: order 2} has been the subject of considerable attention. Purported ``counterexamples'' were found by Manin, but later debunked by Urabe, who then went on to prove \Cref{conj: order 2} for $p \neq 2$ in \cite{Urabe96}. For all $p$, it was proved by Liu--Lorenzini--Raynaud \cite{LLR05} that \emph{if} $\Br(X)$ is finite, then its order is a perfect square. (The proof of \cite{LLR05} proceeds by comparing $\# \Br(X)$ to the order of a certain Tate--Shafarevich group, whose finiteness is equivalent to the Birch and Swinnerton-Dyer Conjecture in that instance.) On the one hand, this gives great confidence that the conjectures are true; but on the other hand, the finiteness assumption is very strong: it is equivalent to the BSD Conjecture for Jacobians over function fields. The \emph{unconditional} statement that $\#\Br(X)_{\nd}$ is a perfect square has remained open in characteristic $p=2$ until now, when we will finally prove: 

\begin{thm}\label{thm: perfect square}
Let $X$ be a smooth, proper, geometrically connected surface over a finite field of characteristic $2$. Then the order of $\Br(X)_{\nd}$ is a perfect square. 
\end{thm}

\subsubsection{Tate's Symplecticity Conjecture}
As mentioned, Conjecture \ref{conj: size} was motivated by an analogy to the Birch and Swinnerton-Dyer Conjecture. Under this analogy, the Brauer group $\Br(X)$ corresponds to the Tate--Shafarevich group $\Sha$ of a Jacobian variety. Motivated by his proof that the Cassels--Tate pairing on (the non-divisible quotient of) such a $\Sha$ is symplectic under the technical hypothesis in \cite[Theorem 3.3]{Tate62}\footnote{This technical hypothesis was mistakenly dropped in the discussion of \cite{Tate66}. This led to a long misconception that was finally resolved in work of Poonen--Stoll \cite{PS99}.}, Tate conjectured that there should be a natural symplectic form on $\Br(X)_{\nd}$ as well. 

Indeed, for every prime $\ell \neq p$, M. Artin and Tate 
defined a non-degenerate skew-symmetric pairing on $
\Br(X)_{\nd}[\ell^{\infty}]$, which was named the 
\emph{Artin--Tate pairing} in \cite{Feng20}.  The restriction 
$\ell \neq p$ was due to the lack of a suitable $p$-adic 
cohomology theory at that time, and motivated the search for 
a $p$-adic cohomology theory to play an analogous role to $
\ell$-adic cohomology, but at the defining characteristic. 
This problem was solved by Milne's discovery of the 
\emph{logarithmic de Rham--Witt cohomology}, which is now 
also known under the synonymous names of \emph{$p$-adic 
\'etale-motivic cohomology}, and \emph{syntomic cohomology}. 
Imitating Artin--Tate's construction, Milne
constructed in \cite{Milne75} a non-degenerate skew-symmetric pairing on $
\Br(X)_{\nd}[p^{\infty}]$ when $p \neq 2$. The restriction $p 
\neq 2$ came from deficiencies in $p$-adic cohomology theory 
for $p=2$ at the time, and was removed in \cite{Milne86} 
using Illusie's development of the de Rham--Witt complex 
\cite{Ill79}. The upshot is that we now have a non-degenerate 
skew-symmetric form on all of $\Br(X)_{\nd}$ in all 
characteristics, which we call the \emph{Milne--Artin--Tate 
pairing}.

We remind the reader that a pairing $\langle -, - \rangle$ is said to be
\[
\text{\emph{skew-symmetric} if } \langle u, v \rangle = -  \langle v,u \rangle \quad \text{ for all } u,v 
\] 
\[
\text{and \emph{alternating} if } \langle u, u \rangle = 0 \quad \text{ for all } u.
\]
Alternating implies skew-symmetric, but the converse can fail if the group has non-trivial 2-torsion. Finally, a pairing is \emph{symplectic} if it is both alternating and non-degenerate. 

The skew-symmetry of the Artin--Tate pairing and Milne--Artin--Tate pairing can be proved in a couple of lines, and together with non-degeneracy implies that the $\ell$-primary part $\Br(X)_{\nd}[\ell^{\infty}]$ has order a perfect square for $\ell \neq 2$. However, it is not enough to deduce that $\#\Br(X)_{\nd}[2^\infty]$ is a perfect square. We refer to the assertion that the Milne--Artin--Tate pairing on $\Br(X)_{\nd}$ is actually symplectic as Tate's \emph{Symplecticity Conjecture}.\footnote{Milne's extension of the Artin--Tate pairing did not exist yet at the time of \cite{Tate66}, but we interpret this extended version of the question as fulfilling Tate's conjecture in \cite{Tate66}.}

\begin{conj}[Tate's Symplecticity Conjecture, \cite{Tate66}]\label{conj: symplecticity} The (Milne-)Artin--Tate pairing on $\Br(X)_{\nd}$ is symplectic. 
\end{conj}

We emphasize again that Conjecture \ref{conj: symplecticity} implies Conjecture \ref{conj: order 2}. Conjecture \ref{conj: symplecticity} has a long and tortuous history, which is described in the introduction of \cite{Feng20}. Since the pairing is skew-symmetric, the content of the Conjecture is concentrated at the $2$-primary part of $\Br(X)_{\nd}$. When $p \neq 2$, it was finally resolved in the second author's thesis \cite{Feng20}. The methods of \cite{Feng20} were restricted to $p \neq 2$ for fundamental reasons, since they use aspects of $\ell$-adic cohomology theory for $\ell=2$, which are only applicable when $\ell \neq p$. The second author has been trying since 2018 to solve the $p=2$ counterpart. One of our main results finally settles this problem: 

\begin{thm}\label{thm: MAT alternating}
Let $X$ be a smooth, proper, geometrically connected surface over a finite field of any characteristic $p$, including $p=2$. Then the Milne--Artin--Tate pairing on $\Br(X)_{\nd}$ is symplectic. 
\end{thm}

\begin{remark}[Higher dimensional generalizations] We construct a generalization of the Milne--Artin--Tate pairing for any smooth, proper, geometrically connected variety of even dimension over $\F_p$, and prove that it is symplectic -- see \S \ref{ssec: higher brauer}.   
\end{remark}

\begin{remark}
Conjecture \ref{conj: symplecticity} was motivated by the symplecticity of the Cassels--Tate pairing, but ironically the Cassels--Tate pairing turned out \emph{not} to be symplectic in the generality originally envisioned by Tate in \cite{Tate66}, as was discovered by Poonen--Stoll \cite{PS99}. The Milne--Artin--Tate pairing is also analogous to a topological duality, the \emph{linking form on a 5-manifold}, which is always skew--symmetric but also turns out not to be alternating in general. Thus, Theorem \ref{thm: MAT alternating} affirms a rather special feature of Brauer groups, not witnessed in other closely analogous mathematical settings. 
\end{remark}

The proof of the $p \neq 2$ case of Conjecture \ref{conj: symplecticity} in \cite{Feng20} exploited the \emph{\'etale Steenrod operations}, which are subtle symmetries of mod $2$ \'etale cohomology. As has been mentioned, a historical difficulty in constructing the pairing at $p=2$ was the lack of an appropriate cohomology theory. The cohomology theory now exists and is called (among other names) ``syntomic cohomology'', but its \emph{symmetries} are not satisfactorily understood. In order to prove Theorem \ref{thm: MAT alternating}, we therefore need to build a theory of \emph{syntomic Steenrod operations} acting on syntomic cohomology, which we will describe next. We emphasize that we develop the theory of syntomic Steenrod operations for all primes $p$, even though the application to Theorem \ref{thm: MAT alternating} concentrates on the case $p=2$. That is, \emph{the bulk of the paper is devoted to developing a general theory, which is not specific to characteristic $2$.}

\subsection{The syntomic Steenrod algebra}

Geisser--Levine \cite{GL00} proved that for smooth varieties over finite fields, syntomic cohomology coincides with $p$-adic \'etale-motivic cohomology. This contextualizes the issue of constructing syntomic Steenrod operations in terms of a classic problem raised by Voevodsky in his manuscript \cite{Voe02} on ``Open problems in motivic homotopy theory'': \emph{developing a theory of motivic Steenrod operations at the defining characteristic.} 

We recall some historical context for this problem. Away from defining characteristic, the mod $p$ motivic Steenrod algebra for varieties in characteristic $0$ was studied by Voevodsky \cite{Voe03, Voe10}, and for varieties in positive characteristic $\ell \neq p$ by Hoyois--Kelly--\O stv\ae r \cite{HKO}. Thanks to their work, the mod $p$ motivic Steenrod algebra is now well-understood away from characteristic $p$, but in characteristic $p$ it is still highly mysterious. Voevodsky conjectured a description of it in \cite{Voe02}, which implies that it should have a Serre--Cartan basis of power operations over the motivic cohomology of the base field. This conjecture remains wide open, but partial evidence was given by Frankland--Spitzweck \cite{FS18} who showed that Voevodsky's ``expected answer'' for the \emph{dual} motivic Steenrod algebra does at least appear as a (module-theoretic) summand of the true answer. This allows one to define power operations on motivic cohomology, as was suggested already in \cite{FS18} and carried out by Primozic in \cite{Pri20}, but from this definition it seems infeasible to control the \emph{properties} of these operations, such as the following basic information (necessary for our calculations):
\begin{itemize}
\item A formula for the product of power operations (which should be given by motivic \emph{Adem relations}), 
\item A formula for the coproduct of power operations (which should be given by a motivic \emph{Cartan formula}).
\end{itemize}
In the special case of mod $p$ Chow groups, these formulas were proved by Primozic \cite{Pri20}, but this case is not enough for us, nor does the argument generalize. In this approach to motivic Steenrod operations, one can control the product and coproduct only up to the ``error'' summand intervening between the motivic Steenrod algebra and the submodule generated by power operations, which should vanish according to Voevodsky's conjectures (but again, this is wide open). In the special case of the $\rH^{2i,i}$ line considered in \cite{Pri20}, the relevant ``error'' terms vanish by general vanishing properties of motivic cohomology. 

For the aforementioned reasons, as well as subtler ones that will be mentioned later, we introduce a new approach to mod $p$ (\'etale-)motivic Steenrod operations in characteristic $p$. Very recently, Annala--Elmanto \cite{AE} independently gave a new construction of motivic Steenrod operations in defining characteristic, for which they can prove motivic Adem relations and a motivic Cartan formula. Their approach is somewhat similar to ours, as will be discussed further in \Cref{rem:AE}, but our applications require still more information about the operations, which is only yielded by our method.

\subsubsection{Structure of the syntomic Steenrod algebra} We write $\Hsyn(-)$ for syntomic cohomology (cf. \S \ref{ssec: syntomic cohomology}), and 
\[
\Hsyn^{i,j}(X) := \Hsyn^i(X; \F_p(j)) \quad \text{and} \quad \Hsyn^{*,*}(X) := \bigoplus_{i,j \in \Z} \Hsyn^{i,j}(X).
\]
The index $i$ is the \emph{degree} and the index $j$ is the \emph{weight}. The cup product makes $\Hsyn^{*,*}(X)$ into a graded $\F_p$-algebra. 

 We construct and analyze the
\emph{syntomic Steenrod algebra} $
\Asyn^{*,*}$, which acts by natural transformations on the syntomic 
cohomology of algebraic varieties in 
characteristic $p$. Here is a summary of what we will prove about the syntomic Steenrod algebra.

\begin{thm}\label{thm: intro Steenrod}
Let $p$ be any prime and let $k$ be a field of characteristic $p$. Then there is a cocommutative Hopf algebra $\Asyn^{*,*}$ over $\Hsynpt$, which we call the syntomic Steenrod algebra, equipped with an algebra homomorphism to the ring of natural endomorphisms of syntomic cohomology $\Hsyn^{*,*}(-)$ viewed as a functor from varieties over $k$ to $\F_p$-vector spaces. Moreover, $\Asyn^{*,*}$ has the following properties. 
\begin{itemize}
\item $\Asyn^{*,*}$ is generated as an $\Hsyn^{*,*}(\Spec k)$-algebra by power operations $\Ps^i, \beta \Ps^i$ for $i \geq 0$, which change degree and weight as indicated: 
\begin{align}\label{eq: intro motivic weight change}
\Ps^i  &\co \Hsyn^{a,b}(-) \rightarrow \Hsyn^{a+2i(p-1), b+i(p-1)}(-) \\
\beta\Ps^i  &\co \Hsyn^{a,b}(-) \rightarrow \Hsyn^{a+2i(p-1)+1, b+i(p-1)}(-) 
\end{align}
\item A basis of $\Asyn^{*,*}$ over $\Hsyn^{*,*}(\Spec k)$ is given by 
\[
\Ps^\alpha :=  \beta^{\epsilon_r}  \Ps^{i_r} \ldots \beta^{\epsilon_1} \Ps^{i_1} \beta^{\epsilon_0} 
\]
as $\alpha$ ranges over elements of the set 
 \[
 \sI := \{(r, \epsilon_r, i_r, \ldots, \epsilon_1, i_1, \epsilon_0) \mid r \geq 0, i_j >0, \epsilon_j \in \{0,1\} , i_{j+1} \geq p i_j+ \epsilon_j\}.
 \]
\item The product on $\Asyn^{*,*}$ is given by explicit Adem relations (\S \ref{sssec: Adem}). 

\item The coproduct on $\Asyn^{*,*}$ is given by an explicit Cartan formula (\S \ref{sssec: Cartan formula}). 
\end{itemize}	
\end{thm}

Theorem \ref{thm: intro Steenrod} is proved in \S \ref{sec: syntomic Steenrod operations} for $k = \F_p$. It then follows for any extension $k/\F_p$ by tensoring over $\Hsyn^{*,*}(\Spec \F_p)$ with $\Hsyn^{*,*}(\Spec k)$. 

By the aforementioned work of Geisser--Levine, Theorem \ref{thm: intro Steenrod} can be summarized as the construction of an \'etale-motivic Steenrod algebra acting on mod $p$ \'etale-motivic cohomology, which has exactly the structure predicted by Voevodsky's conjectures. Note however that we do \emph{not} claim that $\Asyn^{*,*}$ is the full ring of stable cohomology operations on syntomic cohomology, so we are not proving (the \'etale localization of) Voevodsky's conjectures. Rather, we are showing that there is a Hopf sub-algebra of the true (\'etale-)motivic Steenrod algebra that behaves ``correctly'' in all aspects.

\emph{We also emphasize that the results discussed here apply to all $p$}, even though only the case $p=2$ is invoked for the applications to Brauer groups in \S \ref{ssec: intro classical}.

\begin{remark}\label{rem:AE}
As we were completing an initial draft of this paper, Annala--Elmanto communicated to us their independent work \cite{AE}, which provides another construction of motivic Steenrod operations in defining characteristic. They are able to establish the motivic Cartan formula and motivic Adem relations in general, improving upon \cite{Pri20}. Moreover, \'etale sheafifying their work gives another proof of Theorem \ref{thm: intro Steenrod}, and we expect that their operations recover ours in this way. 

The strategy of \cite{AE} follows a similar initial path to ours, both specializing from characteristic zero via infinitely ramified mixed-characteristic rings, though the implementation is different (even after \'etale sheafification); according to our understanding, the approaches had a common origin in ideas of Lurie (to be sketched below), who credits them to Akhil Mathew, who in turn credits inspiration to Niziol. 

However, for our applications we would not be able to get away with using the operations of \cite{AE} as a black box. For example, we also need the key compatibility statements in Theorem \ref{thm: intro compare operations} and Theorem \ref{thm: intro equivariance} below, which are bound up with our approach, and responsible for the bulk of this paper. 
\end{remark}

\subsubsection{A hint of the construction}
The construction of the syntomic Steenrod algebra $\Asyn^{*,*}$ documented here was explained to us by Jacob Lurie. It borrows elements of his vision for ``prismatic stable homotopy theory'', an extension of prismatic cohomology to extraordinary coefficients. Indeed, the classical Steenrod algebra can be thought of as the (derived) endomorphism algebra of the Eilenberg--MacLane spectrum $\F_p$ over the sphere spectrum $\sph$, and we will ultimately realize the syntomic Steenrod algebras as endomorphisms of syntomic cohomology over a certain ``syntomic sphere spectrum''.

We set up some language needed to articulate our strategy more precisely. For a scheme $S$, let $\SH_S$ be Morel--Voevodsky's \emph{$p$-complete motivic stable homotopy category} of $\A^1$-invariant cohomology theories over $S$. Thus $\SH_S$ contains an object $\MHFp$ representing mod $p$ motivic cohomology. Annala--Hoyois--Iwasa \cite{AHI1} constructed the category of $p$-complete \emph{motivic spectra} $\MS_S$, an enlargement of $\SH_S$ which includes non-$\A^1$-invariant cohomology theories. In particular, mod $p$ syntomic cohomology (which is not $\A^1$-invariant) promotes to an object $\MSFp \in \MS_S$. The objects of $\MS_S$ give rise to Nisnevich sheaves of ($p$-complete) \emph{spectra} on smooth schemes over $S$. If $S = \Spec R$, we will also denote $\MS_R:= \MS_S$ and $\SH_R := \SH_S$. 

The first instinct is to simply define the syntomic Steenrod algebra to be the algebra of (derived) endomorphisms of $\MSFp$ over the symmetric monoidal unit of $\MS_k$. This indeed gives the universal algebra of operations, but we would not be able to control the structure of this algebra. To define $\Asyn^{*,*}$ in a way that gives us a handle on the desired properties, we crucially use \emph{perfectoid geometry} (which bridges characteristic $0$ and characteristic $p$) to make a construction that can be ``controlled'' in terms of characteristic 0 objects which are already understood. 

We will summarize the construction below, but in a slightly oversimplified way that elides some technical issues. We start by considering the integral perfectoid ring $\Zpcyc = \Z_p[\mu_{p^\infty}]^{\wedge}_p$ obtained by adjoining all $p$-power roots of unity to $\Z_p$ and then $p$-adically completing. Its generic fiber will be denoted $\Qpcyc$ and its special fiber will be denoted $k$. Following a construction explained to us by Jacob Lurie (and credited by him to Akhil Mathew), we define a ``perfectoid nearby cycles'' functor $\psi \co \SH_K \rightarrow \MS_k$. The \emph{key calculation} is that $\psi$ carries the motivic cohomology spectrum $\MHFp \in \SH_K$ to the syntomic cohomology spectrum $\MSFp \in \MS_k$. This fact ultimately allows us to transmute information about the motivic Steenrod algebra in characteristic $0$, which was explicated by Voevodsky in \cite{Voe03, Voe10}, into information about the syntomic Steenrod algebra in characteristic $p$. Indeed, the motivic Steenrod algebra over $K$ can be defined as 
\[
\cA_{\mot}^{*,*} := \Ext^{*,*}_{\Sph}(\MHFp, \MHFp)
\]
where $\Sph \in \SH_K$ is the symmetric monoidal unit, which is the \emph{$p$-complete motivic sphere spectrum}. Taking inspiration from the key calculation mentioned above, we may regard $\psi(\Sph)$ as ``the syntomic sphere spectrum'' over $k$, and then try to define our syntomic Steenrod algebra as
\[
\Ext^{*,*}_{\psi(\Sph)}(\psi(\MHFp), \psi(\MHFp)).
\]
Our actual definition of $\Asyn^{*,*}$ is a technical variation on this idea.

\subsection{Application to arithmetic duality}
We will describe how the syntomic Steenrod algebra is applied to prove Theorem \ref{thm: MAT alternating}. 

\subsubsection{$\EE_\infty$ Steenrod operations} In fact, there is \emph{another} flavor of Steenrod operations acting on $\Hsyn^{*,*}(-)$, thanks to the realization of $\Hsyn^{*,*}(-)$ as the cohomology ring of a cochain complex $\RGamma_{\syn}(-)$ which has a natural $\EE_\infty$-structure. This leads to an action of $\EE_\infty$ Steenrod operations $\Pe^i$ and $\beta\Pe^i$ on $\Hsyn^{*,*}(-)$.

On motivic cohomology, there would also be the action of two types of Steenrod operations: the motivic Steenrod operations, and the $\EE_\infty$ Steenrod operations. Their interaction is invisible in topological cohomology and \'etale $\ell$-adic cohomology, where they essentially coincide in either of those settings. In other words, the two types of operations are collapsed onto each other under the realization from motivic cohomology to Betti or $\ell$-adic cohomology. However, they are distinct on syntomic cohomology, as can already be seen from the fact that they affect the degree and weights differently from the syntomic operations in \eqref{eq: intro motivic weight change},
\begin{align}\label{eq: intro E weight change}
\Pe^i  &\co \Hsyn^{a,b}(-) \rightarrow \Hsyn^{a+2i(p-1), pb}(-) \\
\beta\Pe^i  &\co \Hsyn^{a,b}(-) \rightarrow \Hsyn^{a+2i(p-1)+1, pb}(-) .
\end{align}
In our story, the interaction of these two flavors (motivic and $\EE_\infty$) of Steenrod operations plays a key role. Each individual flavor is by itself insufficient for the desired applications, but when combined they exactly supplement each other's deficiencies. Concretely, we can access the Milne--Artin--Tate pairing in terms of $\EE_\infty$-Steenrod operations, but then we cannot compute these operations. On the other hand, we can compute some syntomic Steenrod operations in terms of characteristic classes\footnote{Actually, this computation itself also uses certain properties of syntomic Steenrod operations that we prove using the eventual relation to $\EE_\infty$ Steenrod operations.}, since they are of motivic nature, but this is not useful a priori for understanding the pairing. We therefore also need a comparison theorem to combine the two types of operations. We will proceed to describe the situation, and its resolution, more precisely.

\subsubsection{Connection to the Milne--Artin--Tate pairing} The connection between $\EE_\infty$ Steenrod operations and the Milne--Artin--Tate pairing $\langle -, - \rangle_{\MAT}$ comes from a formula
\begin{equation}\label{eq: intro Bockstein formula}
\langle u, u \rangle_{\MAT} = \int_X \Pe^1 (\beta u)  \text{ for all } u \in \Br(X)[2]
\end{equation}
if $X$ is a smooth, proper, geometrically connected surface over a finite field $k$ of characteristic $p=2$. (There is a generalization of this formula to varieties of higher dimension, in Theorem \ref{thm: MAT form}.) We want to prove that the Milne--Artin--Tate pairing is alternating, so we want to show that $
\langle u, u \rangle_{\MAT} = 0$; the formula \eqref{eq: intro Bockstein formula} allows us to translate this into a problem of calculating the effect of certain $\EE_\infty$ Steenrod operations. There is a generalization of \eqref{eq: intro Bockstein formula} to $u \in \Br(X)[2^n]$, and also to elements of the ``higher Brauer groups'' of higher dimensional varieties, in Theorem \ref{thm: MAT form}, which we need but do not describe here.  

For $u \in \Br(X)[2]$ where $p \neq 2$, an analogous result was established in \cite{Feng20}, and was a crucial part of the strategy used there to show that the Artin--Tate pairing is alternating. The proof in \emph{loc. cit.} reduced to topological statements via \'etale homotopy theory, hence does not generalize immediately to the present situation. More seriously, however, \emph{the rest of the strategy in \cite{Feng20} definitively fails} when $p=2$. 

\subsubsection{Calculation of Steenrod operations} At this point, it may be helpful for the reader to refer to the Introduction of \cite{Feng20} for a summary of the proof. In brief, it draws inspiration from classical topology of manifolds to calculate the relevant Steenrod operations in \eqref{eq: intro Bockstein formula} in terms of characteristic classes, using an arithmetic analogue of Wu's theorem relating Steenrod operations and Stiefel--Whitney classes. 

But in the context of syntomic cohomology, simple weight considerations (using for example that Chern classes must live in the motivic line $\rH^{2i,i}$) reveal that the analogous formulas are only plausible for the syntomic Steenrod operations, rather than the $\EE_\infty$ operations.\footnote{This issue does not arise in the $p \neq 2$ situation considered in \cite{Feng20} because the \'etale sheaf $\mu_2$ is isomorphic to the \'etale sheaf $\Z/2$, so the weight index is negligible; this makes it plausible that the $\EE_\infty$ and \'etale motivic Steenrod operations essentially coincide in that case, which turns out to be true.} So at this point, the flavor of Steenrod operations that we can calculate is not the relevant one for the Milne--Artin--Tate pairing, presenting a major gap for our strategy. What saves us is a comparison theorem mediating between \emph{certain} Steenrod operations of different flavors, evaluated on \emph{certain} syntomic cohomology groups, and we discuss this next.

\subsubsection{The comparison theorem}
As has been mentioned, the syntomic and $\EE_\infty$ operations cannot agree in general, since they have different effects on weights: compare the codomains of \eqref{eq: intro motivic weight change} and \eqref{eq: intro E weight change}. Note, however, that the codomains \emph{sometimes} agree. For example, a crucial instance for Theorem \ref{thm: MAT alternating} is the case $p=2$, $a=3$, and $b=1$, where $\Pe^1$ and $\Ps^1$ coincidentally both take the form $\Hsyn^{3,1} \rightarrow \Hsyn^{5,2}$. We prove the following comparison theorem asserting that the two flavors of operations agree \emph{whenever} they have the same domain and codomain. 

\begin{thm}\label{thm: intro compare operations}
If $b = i$, so that the two maps 
\[
\Pe^i \co \Hsyn^{a,b} \rightarrow \Hsyn^{a+2i(p-1), pb} \quad \text{and} \quad 
\Ps^i \co \Hsyn^{a,b} \rightarrow \Hsyn^{a+2i(p-1), b + i(p-1)}
\]
have the same source and target, then they agree. 
\end{thm}

In fact, we prove a more complete statement, \Cref{thm:comparison}, which determines the relationship between $\Pe^i$ and $\Ps^i$ in all cases. To be clear, the main content comes from an analogous result of Bachmann--Hopkins \cite{BH25} for motivic cohomology, in characteristic 0, which we bootstrap to characteristic $p$ using our perfectoid nearby cycles functor. 

Theorem \ref{thm: intro compare operations} allows us to bridge the gap between the two halves of the strategy adapted from \cite{Feng20}, the first about calculating the Milne--Artin--Tate pairing in terms of $\EE_\infty$ Steenrod operations, and the second about calculating motivic Steenrod operations in terms of characteristic classes. It turns out, however, that the second half is much more difficult in our present setting, and we shall discuss this next.

\subsection{Spectral prismatization}\label{intro: spectral prismatization} It turns out that our strategy requires a certain subtle compatibility of the syntomic Steenrod operations with Poincar\'e duality (for the proof of the \emph{Arithmetic Wu formula}, \Cref{thm: syntomic Wu theorem}). The statement is elementary, so we give it below. Perhaps more interestingly, the proof leads us to develop an apparatus which is likely of deeper importance: a generalization of \emph{prismatization} in the sense of Drinfeld and Bhatt--Lurie, for the syntomic sphere spectrum, and an attendant ``spectral Serre duality''. In particular, we construct an approximation to Lurie's conjectural \emph{prismatic stable homotopy category} over $k$, which is an extension of the category of prismatic $F$-gauges over $k$ in the spirit of stable homotopy theory. To be clear, the idea for how to do this was explained to us by Lurie. We emphasize that we develop this story \emph{for general $p$}, even though only the case $p=2$ is invoked for the results on Brauer groups in \S \ref{ssec: intro classical}.

\subsubsection{Arithmetic duality and Steenrod equivariance} Let us formulate the key compatibility that we require. The Poincar\'e duality theorems for syntomic cohomology over finite fields were established by Milne in \cite{Milne76, Milne86}. Indeed, suppose that $X$ is a smooth, proper, geometrically connected variety of dimension $d$ over a characteristic $p$ finite field $k$. Then there is a trace map 
\[
\int_X \co \Hsyn^{2d+1,d}(X) \xrightarrow{\sim} \F_p
\]
and Milne proved that the pairing 
\[
\Hsyn^{a,b}(X) \times \Hsyn^{2d+1-a,d-b}(X) \rightarrow \Hsyn^{2d+1,d}(X) \xrightarrow{\int_X} \F_p
\]
is perfect for every $a,b \in \Z$. Dualizing the cup product
\[
\Hsyn^{*,*}(X) \otimes_{\F_p} \Hsyn^{*,*}(X) \rightarrow \Hsyn^{*,*}(X \times_{k} X)
\]
and then applying Poincar\'{e} duality, we obtain a commutative diagram 
\begin{equation}\label{diag: intro varphi}
\begin{tikzcd}
\Hsyn^{*,*}(X)^\vee \otimes_{\F_p} \Hsyn^{*,*}(X)^\vee \ar[d, "\wr"]  & \ar[l] \Hsyn^{*,*}(X \times_{k} X)^\vee \ar[d, "\wr"]  \\
\Hsyn^{*,*}(X) \otimes_{\F_p} \Hsyn^{*,*}(X) & \ar[l, "\varphi_*"] \Hsyn^{*,*}(X \times_{k} X)
\end{tikzcd}
\end{equation}
Note that the horizontal map $\varphi_*$ in the bottom row increases cohomological degree by $+1$. Both its source and target have natural actions of the syntomic power operations, the source because it is the cohomology of a variety, and the target by the coproduct on the syntomic Steenrod algebra (see \S \ref{ssec: coproduct}), which is given by the Cartan formula. Then our entire strategy hinges upon the following result.  

\begin{thm}\label{thm: intro equivariance} The map $\varphi_*$ is equivariant for the action of $\Asyn^{*,*}$.
\end{thm}

The main difficulty in the proof of Theorem \ref{thm: intro equivariance} comes from the \emph{arithmetic} nature of the duality on $\Hsyn^{*,*}(X)$, which combines the arithmetic duality on the ground field with the geometric duality on $X$. (See the discussion at the beginning of Part \ref{part: characteristic classes} for more about this.)

\begin{remark}
We will give a toy metaphor for Theorem \ref{thm: intro equivariance}. Heuristically, $\varphi_*$ behaves like a pushforward map\footnote{also sometimes called ``Umkehr map'' or ``wrong-way map''} on cohomology, hence the notation (even though there is no actual map $\varphi$ which induces it). Indeed, $\Spec k$ is topologically analogous to $S^1$, so that $X$ is topologically analogous to a manifold $M$ fibered over $S^1$. Then $X \times_{k} X$ is analogous to $M \times_{S^1} M$ and $\Hsyn^{*,*}(X) \otimes_{\F_p} \Hsyn^{*,*}(X)$ is analogous to the cohomology of $M \times M$. In these terms, $\varphi_*$ would be analogous to the pushforward map on cohomology associated to $M \times_{S^1} M \rightarrow M \times M$. Typically, pushforward maps are \emph{not} equivariant with respect to Steenrod operations; the failure of equivariance should be measured by the Stiefel--Whitney classes of the relative normal bundle. But in this particular situation, because $S^1$ has trivial tangent bundle, one would expect equivariance. Of course, this discussion is all heuristic: in our actual situation, there is not even a geometric object whose cohomology realizes $\Hsyn^{*,*}(X) \otimes_{\F_p} \Hsyn^{*,*}(X)$. 
\end{remark}

\subsubsection{Spectral prismatic $F$-gauges} The authors were stymied by this point for a long time, before stumbling upon a lifeline in the recent works of Drinfeld \cite{Drin24} and Bhatt--Lurie \cite{Bha22}. These works lift syntomic cohomology to \emph{prismatic $F$-gauges} (in the terminology of \cite{Bha22}), which are quasicoherent sheaves on certain stacks called \emph{prismatizations}. In particular, for each smooth and proper $f \co X \rightarrow \Spec k$, there is a perfect complex $\cH^X/p$ of quasicoherent sheaves on a stack $\FSyn{\F_p}$, such that $\RGamma(\FSyn{\F_p}; \cH^X/p)$ canonically identifies with $\RGamma_{\syn}(X; \F_p)$. The category of mod $p$ prismatic $F$-gauges over $k$ is\footnote{In the body of the paper, we take a very different definition as our starting point, and show that it is equivalent to this one.} 
\[
 \FGauge{\FF_p} := \QCoh(\FSyn{\F_p}).
 \]

Our proof of Theorem \ref{thm: intro equivariance} exploits prismatization in an essential way. While the full proof is too complicated to explain here (the entirety of \Cref{part: spectral prismatization} is devoted to it), we can hint at the ideas. 

Firstly, we may ``prismatize'' the problem by prismatizing the Steenrod algebra itself. This first involves constructing (following ideas of Lurie) a category of \emph{spectral prismatic $F$-gauges} $\pFGauge{\mathbb{S}}$ that extends $\FGauge{\FF_p}$ analogously to how the usual stable homotopy category extends the derived category of abelian groups. This category $\pFGauge{\mathbb{S}}$ is closely related to the ``prismatic stable homotopy category'' (over $k$)  envisioned by Lurie; the superscript ``pre'' reflects that it is only an approximation to the latter, which however is an equivalence on all the objects that we consider in this paper (see Remark \ref{rem:true-fgauge} for an explanation of the precise meaning of this statement). 


There is an adjunction 
\[
\iota^*\colon \pFGauge{\mbb{S}} \adj \FGauge{\FF_p} : \iota_*
\]
which should morally be thought of as coming from an embedding $\iota$ of $\FSyn{\F_p}$ into a ``spectral prismatization stack'' $\FSyn{\sph}$. Let $\one$ be the unit of $\FGauge{\FF_p}$, which in concrete terms corresponds to the structure sheaf $\cO_{\FSyn{\F_p}}$. We then define the internal Hom algebra in $\pFGauge{\mathbb{S}}$, 
\[
\sAsyn := \cRHom_{\pFGauge{\mathbb{S}}}(\iota_* \one , \iota_* \one ),
\]
as the ``prismatization of the syntomic Steenrod algebra''; it recovers $\Asyn^{*,*}$ by taking global sections in $\pFGauge{\mathbb{S}}$. By construction, all objects of $\FGauge{\FF_p}$ which come via $\iota^*$ from $\FSyn{\sph}$ are equipped with a tautological action of $\sAsyn$. 

Now let us explain how this helps for Theorem \ref{thm: intro equivariance}. It eventually allows us to ``localize'' the problem onto $\FSyn{\F_p}$. After such localization, there is a new perspective on Milne's Poincar\'e duality results due to Bhatt--Lurie \cite[\S 4]{Bha22}, which dissects arithmetic Poincar\'e duality for syntomic cohomology into two more elemental phenomena: 
\begin{enumerate}
\item ``geometric Poincar\'e duality'' (proved by Longke Tang \cite{Tang22}) for the association $X \mapsto \sfrac{\cH^X}{p} \in \FGauge{\FF_p}$, and 
\item Serre duality for $\FGauge{\FF_p}$.  
\end{enumerate}
Ultimately, this enables us to distill a generalization of Theorem \ref{thm: intro equivariance} which is local on the stack $\FSyn{\F_p}$. It says that a certain map of quasicoherent sheaves $\varphi^{\prism}$, which may be viewed as the ``prismatization of $\varphi_*$'', is compatible with the $\sAsyn$-action; this recovers Theorem \ref{thm: intro equivariance} after applying global sections. In turn, this compatibility then has a conceptual explanation: \emph{if the map $\varphi^{\prism}$ can be lifted to the spectral prismatization $\pFGauge{\mbb{S}}$, then it would be compatible with the $\sAsyn$-action for formal reasons}. We will discuss this lifting next.

\subsubsection{Spectral Serre duality} Serre duality is the key input for defining  the map $\varphi^{\prism}$; correspondingly, the key input for lifting $\varphi^{\prism}$ to $\FSyn{\sph}$ is to lift Serre duality to $\FSyn{\sph}$. For this, the key is to find a good candidate for the dualizing sheaf, and here we take our cue from the classical theory of \emph{Brown--Comenetz duality} \cite{BC76}, which is a generalization of Pontrjagin duality to spectra. 

To be precise, let $\Sp$ be the category of $p$-complete spectra and $\bI \in \Sp$ be the $p$-completion of the Brown--Comenetz spectrum. We ``pull back'' (in a suitable sense) $\bI$ from $\Sp$ to define a dualizing sheaf in $\pFGauge{\mbb{S}}$, which we show fits into a good theory of ``spectral Serre duality'' on $\pFGauge{\mbb{S}}$, compatible with the classical theory of coherent duality in $\FGauge{\FF_p}$. We moreover show that this spectral Serre duality plays well with the prismatized Steenrod algebra $\sAsyn$, which supplies the key ingredient to the proof of Theorem \ref{thm: intro equivariance}.

\begin{remark}
The lengths just described may seem fairly outrageous for proving such elementary statements as Conjecture \ref{conj: order 2} and Conjecture \ref{conj: symplecticity}. While the authors sympathize with this sentiment, they can assure the reader that they have tried many less arduous routes over the past decade, and have resorted to this one by necessity (and indeed, desperation). As noted already, these problems have resisted solution for almost sixty years. 
\end{remark}

\subsection{Outline of the paper} 
\begin{figure}[!h]
    \centering
    \includegraphics[scale=.55]{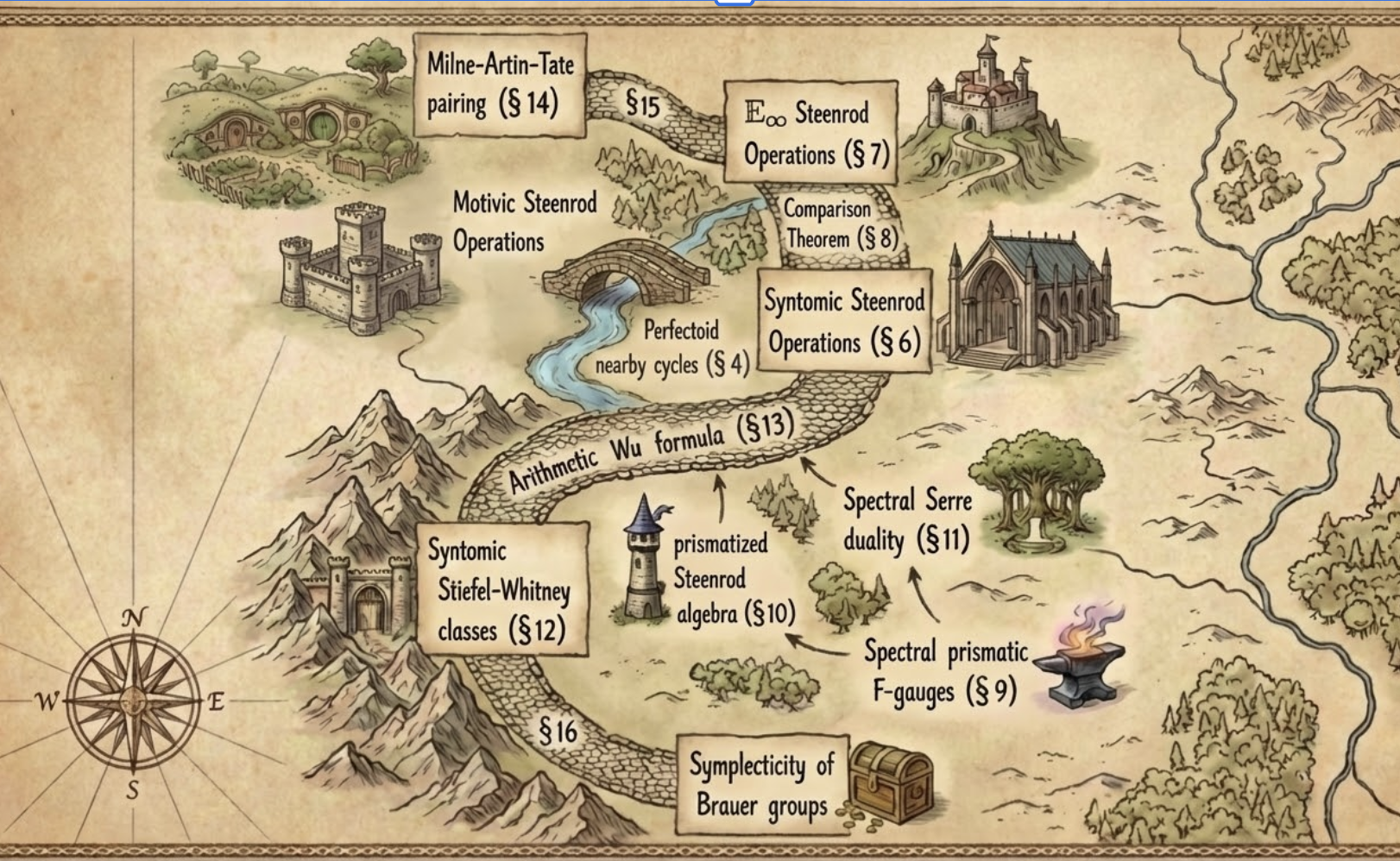}
    \caption{This roadmap depicts the long, winding journey to the symplectic arithmetic duality on higher Brauer groups.}\label{fig:roadmap}
\end{figure}

This paper is divided into five Parts. The contents of each Part are summarized in more detail where it appears. Figure \ref{fig:roadmap} depicts the logical flow to the classical applications to Brauer groups. 

\Cref{part: spectral syntomic cohomology} is concerned with syntomic cohomology. It constructs what we call \emph{spectral syntomic cohomology}, which is short for ``syntomic cohomology with coefficients in motivic spectra''. Furthermore, it defines certain module categories of ``coefficients'' for these spectral cohomology theories, which we call \emph{syntomic spectra}. It also introduces the perfectoid nearby cycles functor. 

\Cref{part: syntomic Steenrod algebra} is concerned with Steenrod operations. It constructs the syntomic Steenrod algebra and establishes its formal properties (e.g., \Cref{thm: intro Steenrod}). It also defines the $\EE_\infty$ Steenrod operations and establishes the comparison with syntomic Steenrod operations (e.g., \Cref{thm: intro compare operations}).

\Cref{part: spectral prismatization} is about prismatization of spectral syntomic cohomology. It defines the category $\pFGauge{\mbb{S}}$, studies a prismatized version of the syntomic Steenrod algebra, and establishes spectral Serre duality and, finally, Theorem \ref{thm: intro equivariance}. 

\Cref{part: characteristic classes} develops a theory of certain characteristic classes in syntomic cohomology, namely Stiefel--Whitney classes and Wu classes, and establishes their relation to the theory of Chern classes from \cite{BL22a}. It also proves the \emph{arithmetic Wu formula} relating syntomic Stiefel--Whitney classes and syntomic Wu classes; it is for this proof that we need \Cref{thm: intro equivariance} and the entire theory of spectral prismatization. This yoga of characteristic classes is used to eventually compute syntomic Steenrod operations. 

Finally, \Cref{part: Brauer groups} contains the applications to the main theorems on (higher) Brauer groups, including  \Cref{thm: perfect square}, \Cref{thm: MAT alternating}, and their higher dimensional generalizations.

\subsection{Acknowledgments} This paper owes a special thanks to Jacob Lurie, who generously taught us his ideas about spectral prismatization and the construction of syntomic Steenrod operations. In addition to guiding us with his high-level vision for prismatic stable homotopy theory, Lurie also supplied some technical arguments used in this paper, yet declined to be named a coauthor. 

We also thank Tom Bachmann, for patiently explaining the results of \cite{BH25}, and then writing it up for our use; and Matthew Morrow, for explaining several of the arguments used in \S \ref{sec: perfectoid nearby cycles}. 

We are grateful to Ben Antieau, Aravind Asok, Bhargav Bhatt, Marc Hoyois, Ryomei Iwasa, Akhil Mathew, Matthew Morrow, Leor Neuhauser, Maxime Ramzi, Longke Tang, John Voight, and Zhouli Xu for helpful correspondence and comments on an earlier draft. 

T.F.~was supported by the NSF (grants DMS-2302520 and DMS-2441922), the Simons Foundation, and the Alfred P. Sloan Foundation. 

S.C. was partially supported by the Danish National Research Foundation through the Copenhagen Centre for
Geometry and Topology (DNRF151), by a research grant from the Center
for New Scientists at the Weizmann Institute of Science, and would like to thank the Azrieli
Foundation for their support through an Early Career Faculty Fellowship.

\subsection{Notation}\label{ssec: notation}

We fix a prime number $p$.  

\subsubsection{Conventions on spectra}
We denote by $\Sp$ the $\infty$-category of \textbf{$p$-complete} spectra. We denote by $\sph$ the unit of $\Sp$, which goes by the name of \emph{$p$-complete sphere spectrum}. This boldface notation is to be contrasted with the blackboard-bold notation used for \emph{motivic spectra}: see Remark \ref{rem: motivic spectra notation}.

\subsubsection{Conventions on coefficients} 
In this paper, we will consider two p-adic ``directions'': one is the \emph{base}, which might (for example) be $\Q_p^{\cyc}, \Z_p^{\cyc}, \F_p$, and the other is the \emph{coefficients}, which might (for example) be $\F_p, \Z_p, \sph$, etc. 

We will study cohomology theories defined on schemes over the base, which are modules over the coefficients. We allow the coefficients to be ring spectra, while our ``test'' schemes are classical.

To make these directions more visually distinct, we write $\sO$ for a $p$-adic ring of integers, $K:= \sO[1/p]$ for its fraction field, and $k$ for its residue field, \emph{when they are being thought of as the ``base'' direction}. Starting with Part \ref{part: syntomic Steenrod algebra}, we will specifically take $\sO := \Z_p^{\cyc}$ (the $p$-adic completion of $\Z_p[\mu_{p^\infty}]$), so $K = \Q_p^{\cyc}$, and $k = \F_p$.  

\subsubsection{$\infty$-categories} We use the formalism of $\infty$-categories following \cite{HTT,HA}. A phrase that does not appear in \emph{loc. cit.}, but which we invoke frequently, is ``presentably symmetric monoidal $\infty$-category''. By this we mean a symmetric monoidal $\infty$-category $\cC$ which is presentable and whose tensor bifunctor $\otimes \co \cC \times \cC \rightarrow \cC$ commutes with colimits in each variable (separately).

\subsubsection{Sheaf and presheaf categories}\label{sssec:not-sheaf-categories}
 For a (small) $\infty$-category $\cC$ and a presentable $\infty$-category $\cD$, we denote by $\PShv(\cC; \cD)$ the presentable $\infty$-category of functors $\cC^{\op} \rightarrow \cD$. We think of this as the category of ``presheaves on $\cC$ valued in $\cD$''. If $\cC$ is equipped with a Grothendieck topology $\tau$, we denote 
\[
\PShv_{\tau}(\cC; \cD)  = \Shv(\cC_{\tau}; \cD)
\]
for the presentable $\infty$-category of sheaves on the site $\cC_\tau$. 

For a stack $S$, we denote by $\cD(S)$ the $\infty$-category of quasicoherent sheaves on $S$, and $\Perf(S) \subset \cD(S)$ the full subcategory of perfect complexes. For a commutative ring $R$, we abbreviate $\cD(R) := \cD(\Spec R)$.

\part{Spectral syntomic cohomology}\label{part: spectral syntomic cohomology}

This Part is concerned with the construction and properties of $p$-adic cohomology theories which include and generalize \emph{syntomic cohomology}. We first review what we mean by (classical) syntomic cohomology in \S\ref{sec: syntomic cohomology}, and establish some technical properties for later use. 

One of the goals of this Part is to expand the notion of syntomic cohomology to more general coefficients, parallel to how the theory of \emph{spectra} expands singular cohomology. We therefore refer to this expansion as \emph{spectral syntomic cohomology}, in order to evoke the image of ``syntomic cohomology with coefficients in spectra''. 

With this in mind, we review in \S\ref{sec: motivic spectra} the theory of \emph{motivic spectra}, which will be used as an ambient category to house our cohomology theories and constructions. 

The heart of this Part is \S\ref{sec: perfectoid nearby cycles}, which defines and analyzes a certain ``perfectoid nearby cycles'' functor $\psi$ from the motivic stable homotopy category over characteristic zero perfectoid fields $K$ (such as $\Q_p^{\cyc}$) to motivic spectra over their residue fields ($k  = \F_p$ in the example). The main theorem of \S\ref{sec: perfectoid nearby cycles} is that $\psi$ carries motivic cohomology over $K$ to syntomic cohomology over $k$. This calculation has to do with the idiosyncratic behavior of algebraic K-theory and motivic cohomology of smooth algebras over perfectoid valuation rings.\footnote{For example, an inspiration (which we do not logically use) is the main theorem of \cite{AMM22}, which establishes that the algebraic $K$-theory of a smooth algebra over a perfectoid valuation ring \emph{coincides} with the algebraic $K$-theory of its generic fiber.}

For a more general motivic spectrum $E \in \MS_{K}$, this justifies considering $\psi(E) \in \MS_{k}$ as the associated syntomic cohomology theory. Finally, in \S \ref{sec: syntomic spectra} we consider the case where $E = \Sph$ is the \emph{motivic sphere spectrum}, and define an associated module category, which we call the category of ``syntomic spectra''.

\section{Syntomic cohomology}\label{sec: syntomic cohomology}
In this section, we recall \emph{syntomic cohomology} and related notions, and prove some vanishing properties that will be needed later. 

Historically, the concept of syntomic cohomology has been studied in varying levels of generality by various different approaches. In this paper, we follow the recent perspective of Bhatt--Lurie \cite{BL22a}, building on work of Bhatt--Morrow--Scholze \cite{BMS2}. The advantage of this approach over earlier ones (e.g., \cite{FM87, Kat87, Tsu99, CN17}) is that it works for all weights, and over highly ramified (perfectoid) rings.

\subsection{Syntomic cohomology of schemes}\label{ssec:syntomic-cohomology}
Bhatt--Morrow--Scholze \cite{BMS2} defined the syntomic cohomology of $p$-adic \emph{formal} schemes. For each $n \in \Z$, and a formal scheme $\cX$ over $\Spf \Z_p$, they define complexes $\Z_p(n)^{\syn}_{\cX}$ on the quasisyntomic site of $\cX$, whose cohomology is the \emph{syntomic cohomology} $\RGamma_{\syn}(\cX; \Z_p(n))$.	

Bhatt--Lurie \cite[\S 8]{BL22a} defined a ``decompleted'' version of the theory, extending syntomic cohomology to schemes over $\Z_p$. Let us describe it for an affine scheme $X = \Spec R$ over $\Z_p$, the general case being deduced from this one by Zariski descent. Let $\wh{R}$ be the $p$-adic completion (in the derived sense) of $R$. Bhatt--Lurie constructed in \cite[\S 8.3]{BL22a} an \emph{\'etale comparison map}
\begin{equation}\label{eq:etale-comparison}
\gamma_{\et}\BK{n} \co  \RGamma_{\syn}(\Spf \wh{R}; \Z_p(n)) \rightarrow \RGamma_{\et}(\Spec \wh{R}[1/p]; \Z_p(n))
\end{equation}
and defined in \cite[Construction 8.4.1]{BL22a} the ``decompleted'' syntomic cohomology $\RGamma_{\syn}(\Spec R; \Z_p(n))$ as the homotopy fiber product in the derived category $\cD(\Z_p)$ of $p$-complete abelian groups,
\begin{equation}\label{eq: decompleted syntomic cohomology}
\begin{tikzcd}
\RGamma_{\syn}(\Spec R; \Z_p(n)) \ar[r] \ar[d] & \RGamma_{\et}(\Spec R[1/p]; \Z_p(n)) \ar[d] \\
\RGamma_{\syn}(\Spf \wh{R}; \Z_p(n)) \ar[r, "\gamma_{\et}\BK{n}"] & \RGamma_{\et}(\Spec \wh{R}[1/p]; \Z_p(n))
\end{tikzcd}
\end{equation}

The resulting cohomology theory $R \mapsto \RGamma_{\syn}(\Spec R; \Z_p(n))$ satisfies \'etale descent, and we denote by $\Z_p(n)_X^{\syn}$ the corresponding \'etale complex on $X = \Spec R$, so that 
\[
\RGamma_{\syn}(X; \Z_p(n)) \cong \RGamma_{\et}(X; \Z_p(n)_X^{\syn}).
\]

Now suppose $R$ is an algebra over $\Z_p^{\cyc}$. We record a consequence of this definition for later technical use. As in \cite[\S 8.5]{BL22a}, we let $\varepsilon := (1, \zeta_p, \zeta_{p^2}, \ldots)$ be a choice of compatible system of primitive $p$-power roots of unity in $\Zpcyc$. We may view $\varepsilon \in \Hsyn^0(\Spec \Zpcyc; \Z_p(1))$. Multiplication by $\varepsilon$ then induces a map 
\[
\RGamma_{\syn}(\Spec R; \Z_p(n)) \xrightarrow{\varepsilon} \RGamma_{\syn}(\Spec R; \Z_p(n+1)). 
\]

\begin{lemma}\label{lem: epsilon comparison}
Suppose $R$ is a commutative algebra over $\Z_p^{\cyc}$. Then the commutative diagram
\begin{equation}\label{eq: epsilon comparison}
\begin{tikzcd}
\RGamma_{\syn}(\Spec R; \Z_p(n)) \ar[d]  \ar[r, "\varepsilon"] &  \RGamma_{\syn}(\Spec R; \Z_p(n+1)) \ar[d] \\
\RGamma_{\syn}(\Spf \wh{R}; \Z_p(n))  \ar[r, "\varepsilon"]  & \RGamma_{\syn}(\Spf \wh{R}; \Z_p(n+1)) 
\end{tikzcd}
\end{equation}
is derived Cartesian. 
\end{lemma}

\begin{proof}
By the Cartesian nature of \eqref{eq: decompleted syntomic cohomology}, it suffices to prove that the analogous commutative square 
\[
\begin{tikzcd}
\RGamma_{\et}(\Spec R[1/p]; \Z_p(n)) \ar[d]  \ar[r, "\varepsilon"] &  \RGamma_{\et}(\Spec R[1/p]; \Z_p(n+1)) \ar[d] \\
\RGamma_{\et}(\Spec \wh{R}[1/p]; \Z_p(n))  \ar[r, "\varepsilon"]  & \RGamma_{\et}(\Spec \wh{R}[1/p]; \Z_p(n+1)) 
\end{tikzcd}
\] 
is Cartesian. Indeed, since $\varepsilon$ induces an isomorphism $\Z_p(n) \cong \Z_p(n+1)$ over $\Z_p^{\cyc}[1/p] = \Q_p^{\cyc}$, multiplication by $\varepsilon$ is in fact an isomorphism in both rows. 
\end{proof}

\subsection{Notation for syntomic cohomology}\label{ssec: syntomic cohomology}

Let $X$ be a scheme over $\Z_p$. For $n \geq 1$ and $b \in \Z$, we write $\Z/p^n(b)^{\syn}_X$ for the cone of multiplication by $p^n$ on $\Z_p(b)^{\syn}_X$. When $n=1$, we also write $\F_p(b)^{\syn}_X := \Z/p(b)^{\syn}_X$. We abbreviate the cohomology of these sheaves by 
\[
\Hsyn^{a,b}(X; \Z/p^n) := 	\Hsyn^a(X; \Z/p^n(b)) := \Het^a(X; \Z/p^n(b)^{\syn}_X), \text{ etc.}
\]
We use the notation $\Z_p(b)^{\syn}_X = \Z_p^{\syn}(b)_X$ interchangeably.

For integers $a,b \in \Z$, we write 
\[
\Hsyn^{a,b}(X) := \Hsyn^a(X; \F_p(b)).
\]
We will use the abbreviation
\[
\Hsyn^{*,*}(X) := \bigoplus_{a,b \in \Z} \Hsyn^a(X; \F_p(b))
\]
and define $\Hsyn^{a,*}(X)$ and $\Hsyn^{*,b}(X)$ similarly. In all the ``nice'' situations that we consider, these are supported on only finitely many $a$ and $b$; for example, Proposition \ref{prop: vanishing range syntomic cohomology} provides an estimate for the support of the cohomological degrees when $X$ is a smooth qcqs scheme over $\Z_p^{\cyc}$. 

There is a map
\[
\F_p(b)_X^{\syn} \otimes  \F_p(b')_X^{\syn} \rightarrow \F_p(b+b')_X^{\syn},
\]
which equips $\Hsyn^{*,*}(X)$ with a ring structure (and similarly with coefficients in $(\Z_p)^{\syn}_X$ or $(\Z/p^n)^{\syn}_X$ instead).

\subsection{Vanishing estimates for syntomic cohomology of formal schemes}

We will establish some vanishing estimates on syntomic cohomology, which are needed later for technical purposes. Below, for a (formal) scheme $X$ over $\Z_p$ we make use of the \emph{absolute prismatic cohomology} $\RGamma_{\prism}(X)$ defined in \cite{BL22a}, as well as the \emph{relative prismatic cohomology} $\RGamma_{\prism}(X/A)$ of Bhatt--Scholze \cite{BS22} when $X$ is defined over $A/I$ for a prism $(A,I)$.

\begin{lemma}\label{lem: prismatic vanishing}
Let $X := \Spec R$ be a smooth affine scheme of relative dimension $d$ over a perfectoid ring $\sO$, $\wh{R}$ be the $p$-adic completion of $R$, and $\wh{X} := \Spf \wh{R}$. Then $\rH^i_{\prism}(\wh{X})$
vanishes for $i \notin [0,d]$. 
\end{lemma}

\begin{proof}
Since $\sO$ is perfectoid, it corresponds to a perfect prism $(A,I)$; in particular, $\sO = A/I$. Hence the absolute prismatic cohomology $\RGamma_{\prism}(\wh{X})$ is canonically identified by \cite[Proposition 4.4.12]{BL22a} with the (relative) prismatic cohomology $\RGamma_{\prism}(\wh{X}/A)$ in the sense of Bhatt--Scholze \cite{BS22}. The prismatic cohomology $\RGamma_{\prism}(\wh{X}/A)$ is equipped with a complete filtration whose associated graded pieces are Breuil-Kisin twists of $\RGamma_{\ol{\prism}}(\wh{X}/A)$. Hence it suffices to show that $\rH^i_{\ol{\prism}}(\wh{X}/A) = 0$ for $i \notin [0,d]$. 

To this end, the \emph{Hodge-Tate comparison} \cite[Theorem 4.11]{BS22} identifies 
\[
\rH^i_{\ol{\prism}}(\wh{X}/A) \cong \Omega^i_{\wh{X}/\sO} \BK{-i}.
\]
Since $X$ is smooth of relative dimension $d$ over $\sO$, $\Omega^i_{\wh{X}/\sO} = 0$ for $i \notin [0,d]$, as desired.
\end{proof}

\begin{prop}\label{prop: vanishing range formal syntomic cohomology}
Let $\wh X$ be a smooth affine formal scheme of relative dimension $d$ over a perfectoid ring $\sO$. Then for all $n \in \Z$, $\Hsyn^i(\wh{X}; \Z_p(n))$ vanishes for $i \notin [0,d+1]$.
\end{prop}

\begin{proof}
By definition \cite[\S 7.4]{BMS2} (and in the notation of \cite[\S 1.6]{BL22a}), the syntomic cohomology of $\wh X$ is defined as the (derived) fiber
\[
\RGamma_{\syn}(\wh X; \Z_p(n)) := \Fib \Big(\Fil^n_{\cN} \RGamma_{\wh \prism}(\wh X) \BK{n} \xrightarrow{\varphi\BK{n} -1} \RGamma_{\wh \prism}(\wh X) \BK{n} \Big)
\]
where $\varphi$ is the Frobenius, $\RGamma_{\wh \prism}$ is the Nygaard-completed prismatic cohomology \cite[\S 1.5]{BL22a}, and $\Fil^n_{\cN} \RGamma_{\wh \prism}$ is the $n$th filtrant of the \emph{Nygaard filtration}. We note that we obtain an equivalent definition using prismatic cohomology without Nygaard completion (which is the starting point of \cite{BL22a})
\[
\RGamma_{\syn}(\wh X; \Z_p(n)) \cong \Fib \Big(\Fil^n_{\cN} \RGamma_{\prism}(\wh X) \BK{n} \xrightarrow{\varphi\BK{n} -1} \RGamma_{\prism}(\wh X) \BK{n} \Big).
\]
Arguing as in the proof of Lemma \ref{lem: prismatic vanishing}: the associated graded of the Nygaard filtration is identified by \cite[Theorem 1.16(3)]{BS22} up to Frobenius twist as truncations of Breuil--Kisin twists of $\RGamma_{\ol{\prism}}(\wh{X}/A)$, which are concentrated in degrees $[0,d]$. Since $\RGamma_{\wh \prism}(\wh X) $ is tautologically complete with respect to the Nygaard filtration, we deduce that 
$\Fil^n_{\cN} \RGamma_{\wh \prism}(\wh X) $ is also concentrated in degrees $[0,d]$. Hence $\RGamma_{\syn}(\wh X; \Z_p(n))$, being the derived fiber of a map between complexes concentrated in degrees $[0,d]$, is itself concentrated in degrees $[0,d+1]$. 
\end{proof}

\begin{prop}\label{prop: vanishing range syntomic cohomology}
Let $X := \Spec R$ be a smooth affine scheme of relative dimension $d$ over $ \Z_p^{\cyc}$. Then for all $n \in \Z$, $\Hsyn^i(X; \Z_p(n))$ vanishes for $i \notin [0,d+2]$.\footnote{We suspect that this bound can be improved to $[0,d+1]$.}
\end{prop}

\begin{proof}
Let $\varepsilon := (1, \zeta_p, \zeta_{p^2}, \ldots)$ be a choice of compatible system of primitive $p$-power roots of unity in $\Z_p^{\cyc}$. We may view $\varepsilon \in \Hsyn^0(\Spec  \Z_p^{\cyc}; \Z_p(1))$. By derived Nakayama, it suffices to show the vanishing modulo $p$, i.e., $\Hsyn^i(X; \F_p(n))=0$ for $i \notin [0, d+2]$. By \Cref{prop: vanishing range formal syntomic cohomology}, the bottom row of \eqref{eq: epsilon comparison} vanishes on cohomology in degree $i > d+1$. Then the Cartesian nature of \eqref{eq: epsilon comparison} implies that the top row is an isomorphism on cohomology in degrees $i >d+2$. Therefore, for $i>d+2$ we have 
\[
\Hsyn^{i,*}(\Spec R) \cong \colim \left(
\Hsyn^{i,*}(\Spec R) \xrightarrow{\varepsilon} 
\Hsyn^{i,*+1}(\Spec R) \xrightarrow{\varepsilon}  \ldots \right) \cong \Hsyn^{i,*}(\Spec R)[\varepsilon^{-1}].
\]
By \cite[Theorem 8.5.1]{BL22a}, the \'etale comparison map
\[
\Hsyn^{i,*}(\Spec R)[\varepsilon^{-1}] \rightarrow \Het^{i,*}(\Spec R[1/p]) := \bigoplus_{n \in \Z} \Het^i(\Spec R[1/p]; \F_p(n))
\]
is an isomorphism. It therefore suffices to show that $\Het^i(\Spec R[1/p]; \F_p(n)) = 0$ for $i > d+2$. Since $ \Q_p^{\cyc}$ has $p$-cohomological dimension $1$, and $\Spec R[1/p]$ has dimension $d$, this follows from the Hochschild--Serre spectral sequence (where $\CC_p$ is the completed algebraic closure of $\Q_p$)
\[
\rH^i(\Q_p^{\cyc}; \Het^j(\Spec R[1/p] \otimes_{\Q_p^{\cyc}}  \CC_p; \F_p(n))) \implies \Het^{i+j}(\Spec R[1/p]; \F_p(n))
\]
and the vanishing of $\Het^j(\Spec R[1/p] \otimes_{\Q_p^{\cyc}} \CC_p ; \F_p(n))$ for $j > d$ (by the Artin Vanishing Theorem \cite[Corollaire XIV.3.2]{SGA4}). 
\end{proof}

\subsection{Syntomic cohomology over finite fields} A case of particular interest to us is when $X$ is a smooth variety over a finite field $k$. In this case, we recall some older and more concrete approaches to syntomic cohomology. 

\subsubsection{Logarithmic de Rham--Witt sheaves} For smooth $X$ over $k$, the cohomology groups $\Hsyn^{*,*}(X)$ were first constructed by Milne \cite{Milne75, Milne86}, from the perspective of what he called \emph{logarithmic de Rham--Witt cohomology}. For a smooth variety $X$, Milne defined sheaves $\nu_n(b)$ in terms of the de Rham--Witt complex, which are shown in \cite[Corollary 8.21 and Remark 8.22]{BMS2} to be isomorphic to the (derived) pushforward of $\Z/p^n(b)[b]$ from the quasisyntomic site to the \'etale site of $X$.

\begin{example}\label{ex: nu sheaves examples}
For each $n  \geq 1$: 
\begin{itemize}
\item $\nu_n(0)$ is isomorphic to the constant sheaf $\Z/p^n \Z$.
\item  $\nu_n(1) \cong \rR \lambda_* \mu_{p^n}[1]$ where $\lambda$ is the map from the fppf site to the \'etale site.\footnote{For $b>1$, we note that $\nu_n(b)$ is \emph{not} in general isomorphic to $\lambda_* \mu_{p^n}^{\otimes b}$.}
\end{itemize}
\end{example}

\subsubsection{Mod $p$ motivic cohomology}\label{sssec: motivic}

Let $X$ be a smooth scheme over a field. For $b \in \Z$, there is a motivic complex $\Z(b)^{\mot}_X$ on $X$, which is a Nisnevich sheaf on $X$. The \emph{motivic cohomology} of $X$ is
\[
\RGamma_{\mot}(X; \Z(b)) := \RGamma_{\Nis}(X; \Z(b)^{\mot}_X).
\]

We denote by $\Z(b)^{\et}_X$ the \'etale sheafification of $\Z(b)^{\mot}_X$. The \emph{\'{e}tale-motivic cohomology} of $X$ is $\RGamma_{\et}(X; \Z(b)^{\et}_X)$. If $\ell$ is a prime number which is invertible in the base field, then it follows from a result of Suslin--Voevodsky \cite[Proposition 6.7]{SV98} that the \'{e}tale motivic complex $\Z/\ell^n(b)^{\et}_X := \Z(b)^{\et}_X \otimes_{\Z} \Z/\ell^n$ identifies with the Tate twist $\Z/\ell^n(b) := \mu_{\ell^n}^{\otimes b}$. See \cite{Gei05} for a survey of motivic cohomology, which provides references for these facts. 

On the other hand, if $\ell = p$ is the characteristic of the base field, then Geisser--Levine proved \cite[Theorem 8.3]{GL00} that $\Z(b)^{\et}_X \otimes_{\Z} \Z/p^n$ is isomorphic to Milne's $\nu_n(b)[-b]$, which as discussed above is identified with syntomic cohomology.

\section{Motivic spectra}\label{sec: motivic spectra}

In this section, we review \emph{motivic spectra}, and establish some machinery for promoting cohomology theories of interest, such as motivic and syntomic cohomology, to motivic spectra. 

To put this in context, recall that for a scheme $S$, Morel--Voevodsky defined the \emph{motivic stable homotopy category} over $S$: a universal category for \emph{$\A^1$-invariant} cohomology theories on schemes over $S$. But many of the cohomology theories that we care about (such as prismatic and syntomic cohomology) are not $\A^1$-invariant. Annala--Hoyois--Iwasa \cite{AHI1} defined the category of \emph{motivic spectra over $S$}, an enlargement of the motivic stable homotopy category which encompasses such cohomology theories. This seems to provide the ``right'' ambient context in which to construct and manipulate our spectral syntomic cohomology theories.

\subsection{Symmetric spectra objects}\label{ssec:symmetric_spectra_objects}
Given a presentably symmetric monoidal $\infty$-category $\cC$ and an object $c\in \cC$, one can form the $\infty$-category $\cC[c^{-1}]$, obtained by universally inverting the object $c$. 
The goal of this subsection is to recall an explicit construction of $\cC[c^{-1}]$ using ``symmetric $c$-spectra objects'' and related constructions. 

\subsubsection{$\Sigma$-objects}Let $\Sigma$ be the free commutative monoid on one element in the $\infty$-category of spaces. Concretely, $\Sigma$ is equivalent to the category of finite sets where maps are isomorphisms, so 
\[
\Sigma \cong \coprod_n B\Sigma_n 
\]
where $\Sigma_n$ is the symmetric group on $n$ elements. 

For an $\infty$-category $\cC$ we let $\cC^{\Sigma} := \Fun(\Sigma, \cC)$ be the $\infty$-category of functors from $\Sigma$ to $\cC$; informally speaking, $\cC^\Sigma$ is the $\infty$-category of tuples $(Y_0,Y_1,\dots)$ where each $Y_n\in \cC$ is endowed with an action of $\Sigma_n$. 

A presentably symmetric monoidal structure on $\cC$ induces one on $\cC^{\Sigma}$, via \emph{Day convolution}: for $Y = (Y_n)_{n\in \N}$ and $Z = (Z_n)_{n\in \N}$ in $\cC^{\Sigma}$, their Day convolution is $Y \otimes Z$ defined by\footnote{In the formula below, $\bigoplus$ means the coproduct, so that it applies even if $\cC$ is not additive. However, we shall only consider it for additive categories.}  
\[
(Y\otimes Z)_m = \bigoplus_{a+b = m} \Ind_{\Sigma_a\times \Sigma_b} ^{\Sigma_{m}}(Y_a \otimes Z_b).
\]

\begin{example}\label{ex:free-symmetric-algebra}
In this situation, for each $X\in \cC$ there is a canonical commutative algebra object 
\[
\Sym(X) := (\one,X,X^{\otimes 2}, \dots)\in \calg(\cC^{\Sigma})
\]
where $\one$ is the unit of $\cC$. 
Here $\Sigma_n$ acts on $X^{\otimes n}$ by permuting the tensor factors. This $\Sym(X)$ can be characterized as the free $\EE_\infty$-algebra on the object $X\tw{1}:= (0,X,0,\dots) \in \cC^\Sigma$, where $0$ denotes the initial object of $\cC$. 
\end{example}

 The construction of Example \ref{ex:free-symmetric-algebra} upgrades to a functor 
\begin{equation}\label{eq:free-symmetric-algebra}
\Sym\colon \cC \to \calg(\cC^{\Sigma}).
\end{equation}

\subsubsection{Stabilization} We set up a general framework for ``stabilization'' with respect to tensoring with an object. 

\begin{defn}\label{defn:stabilization}
Let $\cC$ be a presentably symmetric monoidal $\infty$-category and let $c\in \cC$. Then we define 
\[
\Sp_c^\lax(\cC) := \Mod_{\Sym(c)}(\cC^\Sigma).
\]
An object of $\Sp_c^\lax(\cC)$ consists of, in particular, the data of a sequence $(X_0,X_1,\dots)$ of objects of $\cC$, and the module structure over $\Sym(c)$ gives maps $c\otimes X_i \to X_{i+1}$ with adjoints $\sigma_i \colon X_i \to \hom_{\cC}(c,X_{i+1})$. We let
\begin{equation}\label{eq:sp-lax-inclusion}
\Sp_c(\cC) \subseteq \Sp_c^\lax(\cC) 
\end{equation}
to be the full subcategory spanned by the objects $(X_0,X_1,\dots)$ for which $\sigma_i\colon X_i \to \hom_{\cC}(c,X_{i+1})$ are all isomorphisms. The inclusion \eqref{eq:sp-lax-inclusion} admits a left adjoint 
\begin{equation}\label{eq:tau-C}
\tau_{\cC} \co \Sp_c^\lax(\cC)  \rightarrow \Sp_c(\cC)
\end{equation}
realizing $\Sp_c(\cC)$ as a localization of $\Sp_c^\lax(\cC) $. 
\end{defn}
By \cite[Proposition 1.3.14]{AHI1}, there is a canonical colimit preserving symmetric monoidal functor $\Sigma^\infty_c \colon \cC \to \Sp_c(\cC)$ which identifies $\Sp_c(\cC)$ with $\cC[c^{-1}]$, i.e., it is obtained from $\cC$ by universally inverting the object $c$. 

\subsubsection{Functoriality}\label{sssec:Functoriality_Symmetric_Spectra} For a symmetric monoidal colimit preserving functor $f^*\colon \cC \to \cD$ with right adjoint $f_*$, we obtain an adjunction 
\[
(f^*)^\lax\colon \Sp_c^\lax(\cC) \adj \Sp_{f^*c}^\lax(\cD) : (f_*)^\lax
\]
which are given by applying $f^*$ and $f_*$ pointwise. 

These induce an adjunction 
\begin{equation}\label{eq:sp_c-adjunction}
f^* \co \Sp_c(\cC) \adj \Sp_{f^*c}(\cD) \co f_*
\end{equation}
in which the right adjoint is the restriction of $f_*^\lax$ and the left adjoint is the composition of $(f^*)^\lax$ and the localization functor $\tau_{\cD} \co \Sp_{f^*c}^\lax(\cD) \to \Sp_{f^*c}(\cD)$. 
We abusively denote these functors again by $f_*$ and $f^*$ respectively. The following examples illustrate their calculation in more concrete terms. 

\begin{example}\label{ex:pull_push_lax_explicit}
The functor $f_*^\lax$---hence also the functor $f_*$ from \eqref{eq:sp_c-adjunction}---sends a $c$-spectrum $ Y = (Y_0, Y_1, \ldots) \in \cD^\Sigma$ to
$f_*  Y = (f_*Y_0,f_*Y_1,\dots) \in \cC^\Sigma$, and the structure maps are given by the composition 
\begin{align*}
f_*Y_n \oto{f_*\sigma_n(Y)} & f_*\cHom_{\cD}(f^*c,Y_{n+1})  \cong \cHom_{\cC}(c,f_*Y_{n+1}).
\end{align*}
\end{example}

\begin{example}\label{ex:pullback-MS}
The functor $(f^*)^{\lax}$ has a similar description. It is given on symmetric $c$-spectra objects by 
\[
(f^*)^{\lax}(Y_0,Y_1,\dots) = (f^*Y_0,f^*Y_1,\dots). 
\]
The structure maps $\sigma_n((f^*)^{\lax}Y)$ are given by 
\begin{equation}\label{eq:pullback_lax_tras_maps}
f^*Y_n \oto{f^*\sigma_n(Y)} f^*\cHom_{\cC}(c,Y_{n+1})
 \too  \cHom_{\cD}(f^*c,f^*Y_{n+1}) 
\end{equation}
Since the second map in this composition is not an isomorphism in general, $(f^*)^{\lax}$ does not necessarily carry the full subcategory $\Sp_c(\cC) \subset \Sp_c^{\lax}(\cC)$ to $\Sp_{f^*c}(\cD) \subset \Sp^{\lax}_{f^*c}(\cD)$. In particular, $f^* \co \Sp_c(\cC)  \rightarrow \Sp_{f^*c}(\cD)$ is not simply given by restricting $f^*$. Instead, it is the composition
\begin{equation}\label{eq:pullback_MS}
f^*\colon \Sp_c(\cC) \hookrightarrow \Sp_c^{\lax}(\cC) \oto{(f^*)^{\lax}} \Sp_{f^*c}^{\lax}(\cD) \oto{\tau_\cD} \Sp_{f^*c}(\cD),
\end{equation}
where $\tau_{\cD}$ is the left adjoint to the fully faithful embedding $\Sp_{f^*c}(\cD)\subset \Sp_{f^*c}^{\lax}(\cD)$.
\end{example}

\subsection{Motivic spectra}\label{sssec:motivic-spectra}
Let $S$ be a scheme. Following work of Annala--Hoyois--Iwasa \cite{AHI1}, we will define a presentably symmetric monoidal stable $\infty$-category $\MS_S$ of \emph{$p$-complete motivic spectra over $S$}.

Recall that we are denoting by $\Sp$ the $
\infty$-category of \emph{$p$-complete} spectra. We denote by $\Sm_S$ the category of smooth schemes over $S$. We will now consider the constructions of \S\ref{ssec:symmetric_spectra_objects} for the $\infty$-category 
\[
\cC_S:= \PShv_{\Nis,\ebu}(\Sm_S;\Sp),
\]
of Nisnevich sheaves of spectra on $\Sm_S$ satisfying elementary blowup excision. We denote by $\hom_{\cC_S}(-,-)$ the internal Hom of $\cC_S$. 

We denote the Yoneda embedding of $\Sm_S$ into $\cC_S$ as  
\begin{equation}
\Sigma^\infty_+ \co \Sm_S \rightarrow \cC_S.
\end{equation}
It factors as a composition 
\[
 \Sm_S \rightarrow (\Sm_S)_\star \xrightarrow{\Sigma^\infty} \cC_S
\]
where $(\Sm_S)_\star$ is the category of $S$-pointed smooth schemes over $S$, and the first functor sends $X$ to $X_+ := X \sqcup S$ with the obvious section. We view $\bP^1_S \in (\Sm_S)_\star$ as pointed by the constant section at $\infty$, and we abbreviate its image under $\Sigma^\infty$ by $\rPP^1_S$, or if $S$ is clear from the context, by $\rPP^1$. 

\begin{defn}
As a special case of \Cref{defn:stabilization}, we define the $\infty$-categories 
\[
\MS_S^{\lax} := \Sp^\lax_{\rPP^1}(\cC_S) \quad \text{and} \quad \MS_S := \Sp_{\rPP^1}(\cC_S).
\] 
\end{defn}

To unpack the definition: an object of $\MS^\lax_S$ can be represented by a sequence $Y := (Y_0,Y_1,\dots) \in \cC_S^\Sigma$ equipped with a module structure over $\Sym(\rPP^1)$. This structure provides, in particular, canonical maps (in the category $\cC_S$) 
\begin{equation}
\sigma_n= \sigma_n(Y) \colon Y_{n} \to \cHom_{\cC_S}(\rPP^1,Y_{n+1}),
\end{equation}
 $\MS_S$ is defined as the full subcategory of $\MS_S^{\lax}$ spanned by the objects $ Y$ for which the morphisms $\sigma_n$ are all isomorphisms.
 
 We write 
 \begin{equation}\label{eq:tau-S} 
 \tau_S \co \MS_S^{\lax} \rightarrow \MS_S
 \end{equation}
 for the left adjoint \eqref{eq:tau-C} to the inclusion. This exhibits $\MS_S$ as a symmetric monoidal localization of $\MS^\lax_S$, so it inherits a symmetric monoidal structure.

\begin{defn}[Motives of schemes]\label{defn:motives-of-schemes} There is a tautological colimit-preserving symmetric monoidal functor $\PShv(\Sm_S;\Sp) \rightarrow \MS_S$. Abusing notation, we also let $\Sigma^\infty_+$ be the composition 
\[
\Sigma^\infty_+ \co \Sm_S \xrightarrow{\Sigma^\infty_+} \PShv(\Sm_S;\Sp) \to \MS_S
\]
relying on context to make the target category clear (from now on, it will almost always be $\MS_S$). For $X \in \Sm_S$, we refer to $\Sigma^\infty_+  X$ as the ``motive of $X$'' in $\MS_S$. 
\end{defn}

We have the usual Tate twisting functor in $\MS_S$, defined as follows. 

\begin{defn}[Tate twist] 
Let $Y\in \MS_S$. For any $n \in \Z$, we denote 
\[
Y(n) := Y \otimes (\rPP^1)^{\otimes n} [-2n]
\]
and refer to $Y(n)$ as the ($n$-th) \emph{Tate twist} of $Y$. Note that since (tensoring with) $\bP^1$ is formally inverted in $\MS_S$, this makes sense also for negative values of $n$. 
\end{defn}

\subsubsection{Remarks on $p$-completion}\label{sssec:p-completion}

In the literature, $\MS_S$ usually refers to a version without $p$-complete coefficients; in other words, repeating the definition using spectra instead of $p$-complete spectra. Since this paper is concerned with $p$-adic cohomology theories, it is convenient for us to use $p$-complete coefficients everywhere. There is a standard yoga for porting over results from the non-completed situation to the $p$-complete situation, which we will describe now. 

The \emph{Lurie tensor product} of presentable $\infty$-categories \cite[\S 4.8.1]{HA} gives a purely categorical definition of $p$-completing a category: namely, tensoring with the $\infty$ category of $p$-complete spectra $\Sp$ (over the usual $\infty$-category of spectra). In particular, our $p$-complete version of $\MS_S$ is obtained from the one of Annala--Hoyois--Iwasa via this operation. The analogous remarks apply when discussing the motivic stable homotopy category $\SH_S$ below. 

For a presentable stable $\infty$-category $\cC$, let us denote by $\cC^\wedge_p$ its $p$-completion, i.e., its tensor product with the category of $p$-complete spectra. This is equipped with a canonical $p$-completion functor $\cC \rightarrow \cC^\wedge_p$, which is left adjoint to a fully faithful embedding $\cC^\wedge_p \inj \cC$. This construction has the following properties that we shall repeatedly use when importing results regarding $\MS$ and other categories under consideration. Let $F\colon \cC \to \cD$ be a colimit-preserving functor between presentable stable $\infty$-categories. 
\begin{enumerate}
\item There is an associated colimit preserving functor $F^\wedge_p \colon \cC^\wedge_p \to \cD^\wedge_p$, given by the composition 
\begin{equation}\label{eq:p-completion-of-functor}
F^\wedge_p \colon  \cC^\wedge_p \subseteq \cC \oto{F} \cD \to \cD^\wedge_p, 
\end{equation}
where the last functor is the $p$-completion functor. 

\item The right adjoint $(F^\wedge_p)^R$ of $F^\wedge_p \colon \cC^\wedge_p \to \cD^\wedge_p$ is the restriction of the right adjoint $F^R$ to the $p$-complete objects. 

\item If $F^R$ preserves colimits, then so does $(F^\wedge_p)^R$. Furthermore, if $F$ admits a left adjoint then so does $F^\wedge_p$, its left adjoint can be described as the $p$-completion of the left adjoint $F^L$ in the sense of \eqref{eq:p-completion-of-functor}. 
\end{enumerate}

\subsubsection{Relation to the motivic stable homotopy category}\label{sssec:relation-to-SH} 
Following work of Morel--Voevodsky, we let $\SH_S$ be the \emph{$p$-complete stable homotopy category}. It is a presentably symmetric monoidal stable $\infty$-category, given by the formula
\[
\SH_S = \Big(\rL_{\A^1} \PShv_{\Nis}(\Sm_S;\Sp)\Big)[(\bP^1)^{-1}].
\]
Here $\PShv_{\Nis}(\Sm_S;\Sp)$ is the $\infty$-category of Nisnevich sheaves on $\Sm_S$ valued in $\Sp$, then $\rL_{\A^1}$ is its $\A^1$-localization, and finally $[(\bP^1)^{-1}]$ inverts tensoring with the Tate motive. We emphasize that we demand $p$-completeness in the definitions of $\SH_S$, in contrast to the usual conventions; results from the usual situation will be ported over to ours via  the $p$-completion yoga discussed in \S \ref{sssec:p-completion}.

There is a fully faithful embedding $\SH_S \inj \MS_S$, whose essential image is spanned by objects which are $\A^1$-invariant. This has both a left and a right adjoint, hence preserves all colimits and limits. The left adjoint is $\A^1$-localization, and the right adjoint comes from \cite[Proposition 6.1]{AHI2}, by the general procedure of \S \ref{sssec:p-completion}. The left adjoint is symmetric monoidal, which gives a lax symmetric monoidal structure to the inclusion $\SH_S \inj \MS_S$.

We denote by $\Sph_S$ (or just $\Sph$ when context makes $S$ clear) the unit of $\SH_S$; this goes by the name of ``$p$-complete ($\A^1$-invariant) motivic sphere spectrum''. 

\subsubsection{Comparison with sheaves of spectra}
By construction, there is a colimit-preserving symmetric monoidal functor $\PShv_{\Nis}(\Sm_S;\Sp) \to \MS_S$, with right adjoint 
\begin{equation}\label{eq:loop-space}
\coh{(-)}:= \Omega^\infty_{\bP^1}\colon \MS_S \to \PShv_{\Nis}(\Sm_S;\Sp). 
\end{equation}
In terms of symmetric $\bP^1$-spectra, this construction extracts the $0$-th object of the symmetric spectrum, justifying the notation.
We think of $\coh{Y}$ as the cohomology theory on smooth $S$-schemes associated with the motivic spectrum $Y$. 

We turn to analyze the functor $Y\mapsto \coh{Y}$. For this, we need the following feature of Nisnevich sheaves, which we will repeatedly use. 

\begin{lemma}\label{lem:Nis_Finitary}
Let $S$ be a qcqs scheme. Then the fully faithful embeddings (cf. \S\ref{sssec:not-sheaf-categories} for notation)
\[
\PShv_{\Nis,\ebu}(\Sm_S;\Sp) \subseteq \PShv_{\Nis}(\Sm_S;\Sp) \subseteq \PShv(\Sm_S;\Sp)
\]
are both colimit-preserving. 
\end{lemma} 

\begin{remark}
By definition of $\PShv(\Sm_S;\Sp)$, the collection of evaluation functors 
\[
\{\PShv(\Sm_S;\Sp) \xrightarrow{\sF \mapsto \sF(X) } \Sp\}_{X \in \Sm_S}
\]
is colimit preserving and jointly conservative. Hence Lemma \ref{lem:Nis_Finitary} is equivalent to the statement that for every $X\in \Sm_S$, the functor 
\[
\cC_S:=\PShv_{\Nis,\ebu}(\Sm_S;\Sp) \xrightarrow{\sF \mapsto \sF(X) } \Sp 
\]
is colimit-preserving, and similarly for $\PShv_{\Nis}(\Sm_S;\Sp)$. 
\end{remark}

\begin{proof}[Proof of \Cref{lem:Nis_Finitary}]
Since these are fully faithful embeddings, we need to check that each subcategory is closed under colimits in the next. For the first inclusion, this is because the elementary blowup excision condition involves only pullback squares, so it is clearly closed under colimits (as we work with stable $\infty$-categories). 
For the second inclusion, Nisnevich descent is equivalent to satisfying Nisnevich excision, which is again a condition involving only pullback squares (e.g., as in \cite[Theorem 2.9]{DAGXI}). So this condition is also closed under colimits among all presheaves valued in stable $\infty$-categories. 
\end{proof}

\begin{prop}\label{forg_MS_sheaves}
Let $S$ be a qcqs scheme. For any $n \in \Z$, the functor 
\[
Y\mapsto \coh{Y(n)} \colon \MS_S \to \PShv_{\Nis}(\Sm_S;\Sp)
\]
is colimit-preserving, and the collection of all such functors for all $n \in \Z$ is jointly conservative. 
\end{prop}

\begin{proof}
For the colimit-preservation: since both the source and target are stable $\infty$-categories, and the functor $Y\mapsto \coh{Y(n)}$ preserves limits as it is a right adjoint, it suffices to show that the functor in question preserves \emph{filtered} colimits. Since the target is $p$-complete, we can check this after reducing modulo $p$. 

We claim that $\rPP^1/p \in \cC_S$ -- the cofiber of multiplication by $p$ on $\Sigma^\infty \bP^1$ -- is a compact object. Given the claim, it follows from \cite[Lemma 1.5.2]{Kuniversal} that $Y\mapsto \coh{Y(n)}$ preserves filtered colimits. In order to prove the claim, we observe that $\rPP^1$ is compact in the (not $p$-completed) version of $\MS_S$ from \cite{AHI1} because $\Map(\rPP^1,-)$ is given by evaluation at $\bP^1$, which preserves filtered colimits by the argument of \Cref{lem:Nis_Finitary}. This then implies that $\rPP^1/p$ is compact in the $p$-complete version of $\MS_S$ that we are considering here. 

For the joint conservativity,
since the shifts by $-2n$ are invertible, it is equivalent to show that the collection of functors 
\[
Y\mapsto\coh{(Y\otimes (\rPP^1)^{\otimes n})}
\]
is jointly conservative. This follows directly from \cite[Lemma 1.6.2]{Kuniversal}.
\end{proof}

Another consequence is the compact generation of $\MS_S$. 

\begin{prop}\label{prop:MS_compact_generators}
Let $S$ be a qcqs scheme. The category $\MS_S$ is compactly generated by the collection of objects (cf. \Cref{defn:motives-of-schemes}) $\{\Sigma^\infty_+ X / p(n)\}$ for $X\in \Sm_S$ and $n \in \Z$.
\end{prop}

\begin{proof}
First we argue that such objects generate $\MS_S$. Since our categories are $p$-complete, these objects generate in particular the objects (where we have not reduced mod $p$) $\Sigma^\infty_+ X (n)$ for $X\in \Sm_S$ and $n \in \Z$. So it suffices to show that the collection of functors $\{\Hom_{\MS_S}(\Sigma^\infty_+ X(n),-)\}$ is conservative. For $\cF \in \MS_S$, we have
\[
\Hom_{\MS_S}(\Sigma^\infty_+ X(n), \cF) = \cF(-n)(X)
\]
factors as the composition  $\cF \mapsto \cF(-n)_0$ from \Cref{forg_MS_sheaves}, followed by evaluation at $X$. The collection of functors $\{\cF \mapsto \cF(-n)_0\}_{n \in \Z}$ was proved to be conservative in \Cref{forg_MS_sheaves}. It therefore suffices to see that collection of functors $\{\cF \mapsto \cF(X)\}_{X \in \Sm_S}$ on $\PShv_{\Nis}(\Sm_S;\Sp)$ is conservative, but this is tautological. 

Next we argue that $\Sigma^\infty_+X/p$ is compact for all $X \in \Sm_S$. Note that the canonical functor $\PShv(\Sm_S;\Sp) \to \MS_S$ has a colimit-preserving right adjoint by the combination of \Cref{forg_MS_sheaves} and \Cref{lem:Nis_Finitary}, hence sends compact objects to compact objects. Since the Yoneda embedding of $X$ in $\PShv(\Sm_S; \mathsf{Spectra})$ is compact (where $\mathsf{Spectra}$ stands for the category of not necessarily $p$-complete spectra), its reduction mod $p$ is compact in $\PShv(\Sm_S;\Sp)$, and then its image $\Sigma^\infty_+X/p$ in $\MS_S$ is compact. 
\end{proof}

\subsubsection{Functoriality in the base} 
Let $f\colon S\to T$ be a morphism of qcqs schemes. We investigate the induced functors on the categories of motivic spectra. 

We have a symmetric monoidal adjunction\footnote{By definition, this means that the left adjoint is symmetric monoidal, so the right adjoint is lax symmetric monoidal.}
\begin{equation}\label{eq:cs-adjunction}
f^*\colon  \PShv_{\Nis,\ebu}(\Sm_T;\Sp) \adj \PShv_{\Nis,\ebu}(\Sm_S;\Sp) : f_*. 
\end{equation}
Moreover, the functor $f^*$ carries $\rPP^1_T$ to $\rPP^1_S$. 
Hence,  by the constructions of \Cref{sssec:Functoriality_Symmetric_Spectra}, we obtain symmetric monoidal adjunction 
\begin{equation}\label{eq:ms-lax-adj}
(f^*)^{\lax} \colon \MS_T^\lax \adj \MS_S^{\lax} : f_*^{\lax}
\end{equation}
of lax symmetric $\bP^1$-spectra objects, and a corresponding symmetric monoidal adjunction 
\begin{equation}\label{eq: MS adjunction}
f^* \colon \MS_T \adj \MS_S : f_*
\end{equation}
between the corresponding full subcategories of symmetric $\bP^1$-spectra objects. The concrete descriptions of these functors are summarized in \Cref{ex:pull_push_lax_explicit}.


\subsubsection{Properties} We turn to some of the basic properties of these functors. 
First, we note that $f^* \co \MS_T \rightarrow \MS_S$ commutes with Tate twist:
\[
f^*(Y(n))\cong f^*(Y\otimes (\rPP^1_T)^{\otimes n})[-2n] \cong  (f^*Y) \otimes (\rPP^1_S)^{\otimes n}[-2n] \cong (f^*Y)(n).
\]
Since the Tate twist functor $Y\mapsto Y(n)$ is invertible with inverse $Y\mapsto Y(-n)$, passing to the right adjoints we deduce that $f_* \co \MS_S \rightarrow \MS_T$ also commutes with Tate twist: 
\[
f_*(Y(-n))\cong (f_*Y)(-n).
\]

Next, observe that by construction the functor $\coh{(-)}$ from \eqref{eq:loop-space} commutes with $f_*$ from \eqref{eq: MS adjunction}. 
From this, we deduce the following.
\begin{prop} \label{prop:push_MS_colim}
    Let $f\colon S\to T$ be a morphism of qcqs schemes. Then the functor $f_*\colon \MS_S \to \MS_{T}$ is colimit-preserving.  
\end{prop}

\begin{proof}
Since the functors $Y\mapsto \coh{Y(n)}$ are colimit-preserving and jointly conservative by \Cref{forg_MS_sheaves}, it suffices to show that the functors
\[
Y\mapsto \coh{(f_*Y(n))}
\]
are colimit-preserving for all $n \in \Z$. The Tate twist functor is an equivalence, so it suffices to consider the case $n=0$. Since $\coh{(f_*Y)} \cong f_*(\coh{Y})$, it suffices to show that 
\[
f_*\colon \PShv_{\Nis}(\Sm_S;\Sp) \to \PShv_{\Nis}(\Sm_{T};\Sp) 
\]
is colimit-preserving.
Consider the commutative square
\[
\xymatrix{
\PShv_{\Nis}(\Sm_S;\Sp)\ar^{f_*}[r] \ar@{^{(}->}[d] & \PShv_{\Nis}(\Sm_T;\Sp) \ar@{^{(}->}[d] \\ 
\PShv(\Sm_S;\Sp) \ar^{f_*}[r] & \PShv(\Sm_T;\Sp).  
}
\]
By \Cref{lem:Nis_Finitary}, the vertical fully faithful functors are colimit-preserving. Therefore, in order to show that the upper horizontal functor is colimit-preserving, it suffices to show that the lower horizontal one
\[
f_*\colon \PShv(\Sm_S;\Sp) \to \PShv(\Sm_{T};\Sp) 
\]
is colimit-preserving, which is clear.  
\end{proof}




Note that in general we have a map $f^*(\coh{Y})\to \coh{(f^*Y)}$, but this map may not be an isomorphism. We next record a simple criterion for when it is an isomorphism. 

\begin{cor}\label{cor:sheaf_MS_pullback}
Let $f\colon S\to T$ be a morphism of qcqs schemes, and let $Y= (Y_0,Y_1,\dots) \in \MS_T$. If for all $n \in \N$ the maps
\begin{equation}\label{eq:sheaf_MS_pullback-hypothesis}
f^*\cHom_{\cC_T}(\rPP^1_T,Y_n) \to \cHom_{\cC_S}(\rPP^1_S,f^*Y_n) 
\end{equation}
are isomorphisms, then the natural map $f^*(\coh{Y}) \rightarrow \coh{(f^*Y)}$ is an isomorphism.
\end{cor}

\begin{proof}
By \eqref{eq:pullback_MS}, it would suffice to show that $(f^*)^{\lax}(Y) \in \MS_S$. In view of the description of $\sigma_n((f^*)^{\lax}Y)$ from \eqref{eq:pullback_lax_tras_maps}, this follows immediately from the hypothesis that \eqref{eq:sheaf_MS_pullback-hypothesis} is an isomorphism for all $n \in \N$. 
\end{proof}

\subsection{\'Etale sheafification}\label{ssec: etale sheafification}
There is an \'etale version of motivic spectra, which is obtained by repeating the constructions of \S \ref{sssec:motivic-spectra} with the \'etale topology replacing the Nisnevich topology.

There is an obvious forgetful functor $\MS_S^{\et} \to \MS_S$, with left adjoint the ``\'etale sheafification'' functor $\MS_S \to \MS_S^{\et}$. We denote the composition $\MS_S \to \MS_S^{\et} \to \MS_S$ by $L_{\et}$. 

\begin{remark}
The functor $L_{\et}$ can also be described more or less explicitly. 
Namely, \'etale sheafification (followed by the forgetful functor) induces an analogous functor $L_{\et}\colon \PShv_{\Nis,\ebu}(\Sm_S;\Sp) \to \PShv_{\Nis,\ebu}(\Sm_S;\Sp)$, which in turn induces
\[
L_{\et}^{\lax} \colon \MS_S^{\lax} \to \MS_S^{\lax} 
\]
given by $L_{\et}^{\lax}(Y_0,Y_1,\dots) = (L_{\et} Y_0,L_{\et} Y_1,\dots)$.  
Then $L_{\et} \co \MS_S \rightarrow \MS_S$ is given by the composition 
\begin{equation}\label{eq:etale_sheafification_MS}
L_{\et}\colon \MS_S \hookrightarrow \MS_S^{\lax} \oto{L_{\et}^{\lax}} \MS_S^{\lax} \oto{\tau_S} \MS_S,
\end{equation}
where $\tau_S$ is as in \eqref{eq:tau-S}.
\end{remark}

\subsection{Oriented graded algebras}\label{subsubsec:Oriented_Graded_Algebras}

In \cite[\S 6]{AHI1} the authors develop a useful technique to produce commutative algebras in $\MS_S$. We recall a simplified version of their construction here. It will be convenient to develop the theory in full generality for symmetric monoidal $\infty$-categories. 

\begin{defn}
Let $\cC$ be a presentably symmetric monoidal $\infty$-category and let $c\in \cC$. A \emph{$c$-preoriented graded algebra} in $\cC$ is a graded commutative algebra $E_\bullet\in \calg(\cC^\N)$ together with a map $\omega \colon c\to E_1$. Equivalently, $\omega$ can be regarded as a map $c\langle1\rangle \to E$ in $\cC^\N$ of graded objects, where $c\langle 1\rangle := (0,c,0,\dots)$ is the object $c$ concentrated in degree $1$. The 
$c$-preoriented graded algebras in $\cC$ organize into a presentably symmetric monoidal $\infty$-category $\calg^\por_c(\cC^\N)$. 
\end{defn}

\begin{defn}
Let $\cC$ be a presentably symmetric monoidal $\infty$-category and let $c\in \cC$ and let $(E_\bullet,\omega)$ be a $c$-preoriented graded algebra in $\cC$. 
For each $i\in \N$, the map $\omega \colon c\to E_1$ gives a map $c\otimes E_i \to E_1\otimes E_i \to E_{i+1}$ which adjoints to a map $\sigma_i \colon E_i \to \hom(c,E_{i+1})$. 
We say that $E$ is \emph{$c$-oriented} if $\sigma_i$ is an isomorphism for all $i$. 
We let $\calg^\ori_c(\cC^\N)\subseteq \calg^\por_c(\cC^\N)$ be the full subcategory spanned by the oriented graded algebras. 
\end{defn}

\subsubsection{Turning oriented graded algebras into symmetric spectra objects}\label{sssec:from_oriented_to_spectra}
There is a truncation functor $\pi \colon \Sigma \to \N$ which induces a lax symmetric monoidal functor $\pi^*\colon \cC^\N \to \cC^\Sigma$. Recall from \Cref{ex:free-symmetric-algebra} that $\Sym(c) \cong \Free_{\EE_\infty}(\pi^*c\langle 1\rangle)$ is the free $\EE_\infty$-algebra on $\pi^*c\langle 1\rangle$, so the data of a map from $\pi^*c\langle 1\rangle$ into a commutative algebra in $\cC^\Sigma$ is the same as that of a $\Sym(c)$-algebra. We thus obtain a functor 
\[
\nu:= \pi^*\colon \calg^\por_c(\cC^\N) \to \calg(\Sp_c^\lax(\cC)).
\]
In fact, as explained in \cite[\S 6]{AHI2}, by ``doubling the grading'' this functor canonically promotes to a functor into the category of graded algebras $\calg(\Sp_c^\lax(\cC)^\N)$, and we denote this promotion abusively by the same notation.  

\begin{remark}\label{rem:por_or_nu_pullback}
Unwinding the definitions, we see that the map $\nu$ takes a graded algebra $(E_0,E_1,\dots)$ to the symmetric sequence $(E_0,E_1,\dots)$ in which $\Sigma_n$ acts trivially on $E_n$, and the maps $\sigma_i \colon E_i \to \hom(c,E_{i+1})$ induced from the oriented graded algebra and lax spectrum structures agree. It follows immediately that $\nu$ carries $\calg^{\ori}_c(\cC^\N)$ into $\calg(\Sp_c(\cC))$; in fact, we have a pullback square of $\infty$-categories 
\[
\xymatrix{
\calg^{\ori}_c(\cC^\N)\ar[r]\ar^\nu[d] & \calg^{\por}_c(\cC^\N)\ar^\nu[d] \\
\calg(\Sp_c(\cC))\ar[r] & \calg(\Sp_c^\lax(\cC))
}
\]
\end{remark}

\begin{example}\label{ex:oriented_graded_c_invertible}
Assume that $c$ is already invertible in $\cC$, so that $\cC \xrightarrow{\sim} \Sp_c(\cC)$. Then, an oriented graded algebra in $\cC$ is a graded algebra $(E_0,E_1,\dots)$ together with a map $c\to E_1$ such that $E_i \iso \hom(c,E_{i+1})$. Since $c$ is invertible, this is equivalent to the original map $c\otimes E_i \to E_{i+1}$ being an isomorphism, so that $E_\bullet$ assumes the form $(E_0, E_0 \otimes c,E_0 \otimes c^{\otimes 2},\dots)$. Hence, an object of $\calg_c^{\ori}(\cC^\N)$ is the data of 
\begin{itemize}
\item a commutative algebra $E_0 \in \calg(\cC)$, together with 
\item a structure of graded commutative algebra on $\bigoplus_{n \in \N} (E_0\otimes c^{\otimes n})$.
\end{itemize}
This latter datum is equivalent to what is called a \emph{strict structure} on $E_0\otimes c$ (see for example \cite{carmeli2023strict}): a factorization of the map 
\[
\sph \oto{c} \Pic(\cC) \oto{(-)\otimes E_0} \Pic(E_0)
\]
from the sphere spectrum $\sph$ through the truncation map of spectra $\sph \to \Z$.  
In this situation, the functor 
\[
\nu\colon \calg_c^{\ori}(\cC^\N) \to \calg(\Sp_c(\cC)^\N)\cong \calg(\cC^\N)
\]
(in which the second functor is induced by the equivalence $\Omega^\infty_c\colon \Sp_c(\cC) \to \cC$)
agrees with the forgetful functor $\calg_c^{\ori}(\cC^\N) \to \calg(\cC^\N)$ that forgets the orientation. 
\end{example}

\subsubsection{Naturality of $\nu$}\label{sssec:functoriality_naturality_nu}

Let $f^*\colon \cC \to \cD$ be a colimit preserving symmetric monoidal functor with right adjoint $f_*$. They induce an adjunction  
\[
(f^*)^\por\colon \calg^\por_c(\cC^\N) \adj \calg^\por_{f^*c}(\cD^\N) : (f_*)^\por
\]
in which the functor $(f^*)^\por$ carries $(E_\bullet, \omega \colon c\to E_1)$ to $(f^*E_\bullet,f^*\omega \colon f^*c \to f^*E_1)$, and the right adjoint $f_*^\por$ carries $(E_\bullet, \omega \colon f^*c\to E_1)$ to $(f_*E_\bullet,c \oto{\text{unit}} f_*f^*c \oto{f_* \omega} f_*E_1)$. 

Since the functor $\pi^*$ is given by pre-composition with a symmetric monoidal functor $\Sigma \to \N$, it is compatible with post-composition with arbitrary colimit preserving symmetric monoidal functors $f^*\colon \cC \to \cD$. We immediately deduce the following.

\begin{prop}\label{prop:pull_realization}
Let $f^*\colon \cC \to \cD$ be a colimit preserving symmetric monoidal functor between presentably symmetric monoidal $\infty$-categories. Then the square 
\[
\xymatrix{
\calg^{\por}_c(\cC^\N)\ar^{\nu}[d]\ar^{(f^*)^{\por}}[r] & \calg^{\por}_{f^*c}(\cD^\N)\ar^\nu[d]\\
\calg(\Sp_c^\lax(\cC))\ar^{(f^*)^\lax}[r] & \calg(\Sp_{f^*c}^\lax(\cD))\\
}
\]
canonically commutes.
\end{prop}

Since the horizontal functors in this square have right adjoints, we obtain a Beck--Chevalley transformation $\nu f_*^{\por} \to  f_*^{\lax}\nu$. 

\begin{prop}\label{prop:realization_push}
Let $f^*\colon \cC \to \cD$ be a colimit preserving symmetric monoidal functor between presentably symmetric monoidal $\infty$-categories. Then the Beck--Chevalley transformation $\nu f_*^{\por} \to  f_*^{\lax}\nu$ is an isomorphism. Hence, we obtain a canonical homotopy rendering the diagram 
\[
\xymatrix{
\calg^{\por}_c(\cC^\N)\ar^{\nu}[d] & \calg^{\por}_{f^*c}(\cD^\N)\ar^\nu[d] \ar^{(f_*)^{\por}}[l]\\
\calg(\Sp_c^\lax(\cC)) & \calg(\Sp_{f^*c}^\lax(\cD)) \ar^{f_*^\lax}[l]\\
}
\]
commutative. 
\end{prop}

\begin{remark}
The analogous compatibilities hold for the factorization of $\nu$ through $\N$-graded motivic spectra, which was mentioned in the first paragraph of \S\ref{sssec:from_oriented_to_spectra}, by essentially the same arguments. 
\end{remark}

\begin{proof}Both $\nu f_*^{\por} (E_\bullet,\omega)$ and $f_*^{\lax}\nu  (E_\bullet,\omega)$ have underlying symmetric sequence $(f_*E_0,f_*E_1,\dots)$ by \Cref{ex:pull_push_lax_explicit} and its analogue for pre-oriented graded algebras. Moreover, from the definition of $\nu$ we see that the Beck--Chevalley map restricts to the identity on the underlying symmetric sequences, hence is an isomorphism. 
\end{proof}

\begin{cor}
With the same settings as in \Cref{prop:realization_push}, the functor $f_*^{\por}$ restricts to a functor $f_*^{\ori}\colon \calg_d^{\ori}(\cD^\N) \to \calg_{f_* d}^{\ori}(\cC^\N)$. 
\end{cor}

\begin{proof}
We have to check that if $(E_\bullet,\omega)\in \calg_d^\ori(\cD^\N)$ then $f_*\omega$ is an orientation of $f_*E_\bullet$. This is equivalent to the claim that $\nu f_*^{\por}(E_\bullet,\omega) \in \Sp_{f_*d}(\cC)$, which follows from $\nu f_*^{\por}(E_\bullet,\omega)\cong f_*^\lax\nu(E_\bullet,\omega)$ since $f_*^\lax$ carries spectra objects to spectra objects (as discussed at the beginning of \Cref{sssec:Functoriality_Symmetric_Spectra}). 
\end{proof}
Note that
the functor $f_*^{\ori}$ has a left adjoint $(f^*)^{\ori}$ given by the composition 
\[
\calg_c^{\ori}(\cC^\N) \to \calg_c^{\por}(\cC^\N) \oto{(f^*)^\por}\calg_{f^*c}^{\por}(\cD^\N) \to \calg_{f^*c}^{\ori}(\cD^\N)
\]
in which the last functor is the left adjoint to the inclusion of the oriented algebras into the preoriented ones.

\begin{cor}\label{cor:realization_reflective}
Let $\cC$ be a presentably symmetric monoidal $\infty$-category and let $c\in \cC$. Let $\Sigma^\infty_c \colon \cC \to \Sp_c(\cC)$ be the $c$-stabilization functor, with right adjoint $\Omega^\infty_c$. Then, for $(E_\bullet,\omega)\in \calg^{\ori}_c(\Sp_c(\cC)^\N)$ there is a natural isomorphism
\[
\nu((\Omega^\infty_c)^{\por} (E_\bullet,\omega) ) \cong E_\bullet 
\]
in $\calg(\Sp_c(\cC)^\N)$.
\end{cor}

\begin{proof}
Let $\tilde{c} \in \Sp_c(\cC)$ be the image of $c$ under $\Sigma_c^\infty$. Applying \Cref{prop:realization_push} to $f^* = \Sigma_c^\infty$, and using the observation at the end of \Cref{ex:oriented_graded_c_invertible}, we obtain
\[
\nu((\Omega^\infty_c)^{\por} (E_\bullet,\omega)) \cong (\Omega^\infty_{\tilde{c}})^\lax \nu(E_\bullet,\omega) \cong E_\bullet,
\]
as desired. 
\end{proof}

\subsubsection{Oriented graded algebras in $\cC_S$}

Recall that we defined $\cC_S := \PShv_{\Nis,\ebu}(\Sm_S;\Sp)$.

\begin{defn} Specializing the construction of oriented graded algebras to the case $\cC = \cC_S$ and $c = \rPP^1_S$, we can form the $\infty$-category of \emph{$\bP^1$-pre-oriented graded algebras} in $\cC_S$, which we denote  
\[
\calg^{\por}_{\bP^1}(\cC_S^\N):= \calg^{\por}_{\rPP^1}(\cC_S^\N)
\]
We further have the full subcategory 
\[
\calg^{\ori}_{\bP^1}(\cC_S^\N):= \calg^{\ori}_{\rPP^1}(\cC_S^\N)
\]
spanned by the $\bP^1$-oriented graded algebras. 
Recall that these are the pre-oriented graded algebras for which the resulting maps
\begin{equation}\label{eq: pre-oriented graded transition}
\sigma_n \colon E_n \to \hom_{\cC_S}(\rPP^1,E_{n+1})
\end{equation}
are isomorphisms. 

We denote the $\infty$-categories of pre-oriented and oriented graded algebras in $\cC_S$ by 
\[
\calg^{\ori}_{\bP^1}(\cC_S^\N)\subseteq \calg^{\por}_{\bP^1}(\cC_S^\N).
\]
\end{defn}

There is also a variant of this definition which will be useful later. 

\begin{variant}\label{var:pic_or}
Let $\pic = \rB\G_m$ be the Nisnevich sheaf of (connective, 1-truncated) spectra on schemes, assigning to a scheme $X$ its Picard stack $\pic(X)$, so that $\pi_0\pic(X)$ is the Picard group of $X$ and $\pi_1\pic(X) = \sO(X)^\times$. We implicitly restrict $\pic$ to smooth $S$-schemes. 
A \emph{$\pic$-pre-oriented graded algebra} in $\cC_S$ is a graded algebra $E\in \calg(\cC_S^\N)$ together with a map $\Sigma^\infty\pic\langle{1}\rangle \to E$, or equivalently, a map $\Sigma^\infty \pic \to E_1$. We denote by $\calg^{\por}_{\pic}(\cC_S^\N)$ the $\infty$-category of $\pic$-pre-oriented graded algebras in $\cC_S$.    

The canonical (pointed) map $\bP^1 \to \pic$ induces a forgetful functor $\calg^{\por}_{\pic}(\cC_S^{\N}) \to \calg^{\por}_{\bP^1}(\cC_S^{\N})$. We define the $\infty$-category of \emph{$\pic$-oriented graded algebras} $\calg^{\ori}_{\pic}(\cC_S^{\N})$ to be the preimage of $\calg^{\ori}_{\bP^1}(\cC_S^{\N})$ along this forgetful functor. In other words, we say that a $\pic$-pre-oriented graded algebra is \emph{oriented} if the underlying $\bP^1$-pre-oriented algebra is oriented. 
\end{variant}

\begin{remark}
Roughly speaking, a $\pic$-orientation of $E\in \calg(\cC_S^\N)$ is responsible for the existence of Thom isomorphisms for $E$-valued cohomology. 
\end{remark}

\subsubsection{Constructing motivic spectra from oriented graded algebras} 

Applying the construction from \Cref{sssec:from_oriented_to_spectra} to the case where $\cC = \cC_S$ and $c = \rPP^1$, we obtain a functor
\begin{equation}\label{eq: por to mslax}
\nu \co \calg^\por_{\bP^1}(\cC_S^{\N}) \to \calg(\MS_S^\lax).
\end{equation}
By \Cref{rem:por_or_nu_pullback}, this functor restricts to a functor from $\bP^1$-oriented  algebras to $\calg(\MS_S)$, namely, we have a commutative diagram
\begin{equation}\label{eq: or to ms}
\xymatrix{
\calg^\ori_{\bP^1}(\cC_S^{\N})\ar@{..>}^{\nu}[r]\ar@{^{(}->}[d] & \calg(\MS_S)\ar@{^{(}->}[d] \\	
\calg^\por_{\bP^1}(\cC_S^{\N}) \ar_\nu[r] & \calg(\MS_S^\lax)
}
\end{equation}


\subsubsection{}\label{sssec:tau-ori} The fully faithful inclusion  $\calg^\ori_{\bP^1}(\cC_S^{\N}) \subset \calg^\por_{\bP^1}(\cC_S^{\N}) $ admits a left adjoint which we call 
\begin{equation}
\tau_S \co \calg^\por_{\bP^1}(\cC_S^{\N}) \rightarrow \calg^\ori_{\bP^1}(\cC_S^{\N}).
\end{equation}
It is compatible with the $\tau_S$ from \eqref{eq:tau-S} under the realization functor $\nu$ \eqref{eq: por to mslax}, in the sense that the following diagram commutes
\begin{equation}\label{eq:tau-por-ori}
\begin{tikzcd}
\calg^\por_{\bP^1}(\cC_S^{\N}) \ar[r, "\tau_S"] \ar[d, "\nu"] &  \calg^\ori_{\bP^1}(\cC_S^{\N})\ar[d, "\nu"]  \\
\calg(\MS_S^\lax) \ar[r, "\tau_S"] & \calg(\MS_S) 
\end{tikzcd}
\end{equation}
and analogous remarks apply to the $\pic$-oriented variant.

\subsubsection{Functoriality in the base}
Let $f\colon S \to T$ be a morphism of qcqs schemes. The functor $f^* \colon \cC_T \to \cC_S$ carries the object $\rPP^1_T$ to the object $\rPP^1_S$, hence gives a functor 
\begin{equation}\label{eq:f-pullback-por}
(f^*)^{\por}\colon \calg^{\por}_{\bP^1}(\cC_T^\N) \to \calg^{\por}_{\bP^1}(\cC_S^\N). 
\end{equation}

By \Cref{prop:pull_realization}
it fits into a commutative square
\begin{equation}\label{eq: por algebra pullback}
\xymatrix{
\calg^\por_{\bP^1}(\cC_T^{\N})\ar^{\nu}[d] \ar^{(f^*)^{\por}}[r] & \calg^\por_{\bP^1}(\cC_S^{\N}) \ar^{\nu}[d] \\
\calg(\MS_T^\lax)  \ar^{(f^*)^\lax}[r]      & \calg(\MS_S^\lax)
}
\end{equation} 

As in \Cref{sssec:functoriality_naturality_nu} the functor $(f^*)^\por$ admits a right adjoint 
\[
(f_*)^\por \colon \calg^\por_{\bP^1}(\cC_S^{\N}) \to  \calg^\por_{\bP^1}(\cC_T^{\N}),
\]
which by \Cref{prop:realization_push} fits into a commutative diagram 

\[
\begin{tikzcd}
\calg^\por_{\bP^1}(\cC_S^{\N})\ar[d, "\nu"] \ar[r, "f_*^\por"] & \calg^\por_{\bP^1}(\cC_T^{\N}) \ar[d, "\nu"] \\
\calg(\MS_S^\lax)  \ar[r, "f_*^{\lax}"]      & \calg(\MS_T^\lax)
\end{tikzcd}
\]
and restricts to a functor between the full subcategories of $\bP^1$-oriented algebras, as expressed in the following Corollary.

\begin{cor} \label{cor:pushforward_oriented}
Let $f\colon S\to T$ be a morphism of qcqs schemes.
The functor 
\[
f_*^\por \colon  \calg^\por_{\bP^1}(\cC_S^{\N}) \to \calg^\por_{\bP^1}(\cC_T^{\N})
\]
carries $\bP^1$-oriented algebras to $\bP^1$-oriented algebras.
\end{cor}


\Cref{cor:pushforward_oriented} says in other words that $(f_*)^\por $ restricts to a functor 
\begin{equation}
f_* \co \calg^\ori_{\bP^1}(\cC_S^{\N}) \rightarrow \calg^\ori_{\bP^1}(\cC_T^{\N}).
\end{equation}
This $f_*$ fits into a commutative diagram 
\begin{equation}\label{eq:pushforward_pre-oriented_beck_chevalley}
\begin{tikzcd}
\calg^\ori_{\bP^1}(\cC_S^{\N})\ar[d, "\nu"] \ar[r, "f_*"] & \calg^\ori_{\bP^1}(\cC_T^{\N}) \ar[d, "\nu"] \\
\calg(\MS_S)  \ar[r, "f_*"]      & \calg(\MS_T)
\end{tikzcd}
\end{equation}

Finally, $f_*$ admits a left adjoint $f^* \co \calg^\ori_{\bP^1}(\cC_T^{\N}) \rightarrow \calg^\ori_{\bP^1}(\cC_S^{\N})$, defined as the composition
\begin{equation}\label{eq:f-pullback-ori}
f^* \co \calg^\ori_{\bP^1}(\cC_T^{\N}) \inj \calg^\por_{\bP^1}(\cC_T^{\N}) \oto{(f^*)^{\por}} \calg^\por_{\bP^1}(\cC_S^{\N}) \xrightarrow{\tau_S} \calg^\ori_{\bP^1}(\cC_S^{\N}).
\end{equation}

\begin{remark}\label{rem:etale_sheafification_por_compatibility}
Since the \'{e}tale sheafification functor is a composition of a symmetric monoidal left adjoint and a right adjoint, similar considerations apply to it. Namely, for every qcqs scheme $S$, we have a functor 
\[
L_{\et}^\por \colon \calg^\por_{\bP^1}(\cC_S^{\N})\to \calg^\por_{\bP^1}(\cC_S^{\N})
\]
fitting into a
commutative square
\begin{equation}\label{eq:etale_sheafification_por_compatibility}
\xymatrix{
\calg^\por_{\bP^1}(\cC_S^{\N})\ar^{\nu}[d] \ar^{L_{\et}^\por}[r] & \calg^\por_{\bP^1}(\cC_S^{\N}) \ar^{\nu}[d] \\
\calg(\MS_S^\lax)  \ar^{L_{\et}^\lax}[r]      & \calg(\MS_S^\lax)
}
\end{equation} 
Composing $L_{\et}^\por$ with the functor $\tau_S$ from \eqref{eq:tau-por-ori}, we obtain the oriented version of \'etale sheafification: 
\begin{equation}
L_{\et} \colon \calg^\ori_{\bP^1}(\cC_S^{\N}) \oto{L_{\et}^\por} \calg^\por_{\bP^1}(\cC_S^{\N}) \oto{\tau_S} \calg^\ori_{\bP^1}(\cC_S^{\N}).
\end{equation}
From the definition of \'etale sheafification and \eqref{eq:tau-por-ori}, it is clear that these functors are compatible with \'etale sheafification of (lax) motivic spectra, in the sense of the commutative diagrams
\[
\begin{tikzcd}
\calg^\ori_{\bP^1}(\cC_S^{\N}) \ar[r,"L_{\et}"] \ar[d, "\nu"] &  \calg^\ori_{\bP^1}(\cC_S^{\N})  \ar[d, "\nu"]  \\
\calg(\MS_S) \ar[r, "L_{\et}"] & \calg(\MS_S)
\end{tikzcd} \hspace{1cm} \begin{tikzcd}
\calg^\por_{\bP^1}(\cC_S^{\N}) \ar[r,"L_{\et}^\por"] \ar[d, "\nu"] &  \calg^\por_{\bP^1}(\cC_S^{\N})  \ar[d, "\nu"]  \\
\calg(\MS_S^{\lax}) \ar[r, "L_{\et}^{\lax}"] & \calg(\MS_S^{\lax})
\end{tikzcd}
\]
\end{remark}

\subsection{Promoting motivic and syntomic cohomology to motivic spectra}\label{ssec: promoting}
The language of oriented graded algebras allows us to lift syntomic cohomology into a motivic spectrum. Apart from that, it will be a convenient language to compare syntomic cohomology with motivic and \'{e}tale cohomology, so we shall now explain how all three theories organize into $\bP^1$-oriented graded algebras, resulting in motivic spectra.

\subsubsection{Syntomic cohomology} For $n \in \Z$, let $\SZp(n)$ denote the absolute syntomic cohomology functor of \cite[\S 8.4]{BL22a}. Taking the direct sum over $n \in \N$, these assemble into a functor $\SZp(\bu) := \bigoplus_{n \in \N} \SZp(n)$ from schemes to graded commutative algebras in $\cD(\Z_p)$, the \emph{$p$-completed derived category} of abelian groups. Since it satisfies Nisnevich descent and elementary blowup excision, we may regard 
\[
\SZp(\bu) = (\SZp(n))_{n \in \N} \in \calg(\cC_S^\N).
\]
Then we can \emph{shear} it to get a new graded commutative algebra\footnote{Note that the commutative algebra structure on the shearing uses the fact that $\SZp(\bu)$ lands in $\Z$-algebras. In general, such a shearing of a graded commutative algebra is only an $\EE_2$-algebra.} 
\[
\SZp(\bu)[2\bu]  = (\SZp(n)[2n])_{n \in \N} \in \calg(\cC_S^\N).
\]

In \cite[Proposition 7.5.2]{BL22a}, Bhatt--Lurie construct the \emph{syntomic first Chern class} 
\[
c_1^{\syn} \colon \pic \to \SZp(1)[2].
\]
Forgetting the spectrum structure, we may regard $\pic$ as sheaf of pointed spaces together with a map of sheaves of spectra $\Sigma^\infty \pic \to \SZp(1)[2]$. Composing with it, we obtain a $\pic$-pre-orientation
\[
\xi^{\syn} \colon \Sigma^\infty \pic \to \pic \oto{c_1^{\syn}} \SZp(1)[2]. 
\]
\begin{prop} \label{prop:pic_orientation_syntomic}
The $\pic$-pre-orientation $\xi^{\syn}$ is an orientation, so that 
\[
(\SZp(\bu)[2\bu]_S,\xi^{\syn}) \in \calg^{\ori}_{\pic}(\cC_S^\N)
\]
for every qcqs scheme $S$. 
 \end{prop}

\begin{proof}
By definition of orientation, we have to show that the map 
\[
\SZp(n) \to \cHom_{\cC_S}(\Sigma^\infty\bP^1,\SZp(n+1)[2])
\]
induced from $\xi^{\syn}$
is an isomorphism. This can be checked after taking sections over an arbitrary smooth $S$-scheme $X$, so it suffices to show that $\xi^{\syn}$ and evaluation at the basepoint of $\bP^1$ together induce an isomorphism
\[
\Hsyn^{*+2,*+1}(X\times \bP^1;\Z_p) \iso  \Hsyn^{*+2,*+1}(X;\Z_p) \oplus \Hsyn^{*,*}(X;\Z_p).
\]
This is a special case of the projective bundle formula of \cite[Theorem 9.1.1]{BL22a}.
\end{proof}
This leads to the following construction.
 
\begin{defn}[Syntomic cohomology as a motivic spectrum]\label{def:mot_syn_coh}
Let $S$ be a qcqs scheme. We define the commutative algebra 
\[
(\MSZp)_S = \MSZp(0)_S \in \calg(\MS_S)
\]
to be the commutative algebra in motivic spectra corresponding to $(\SZp(\bu)[2\bu]_S,\xi^{\syn})$, or in other words
\[
(\MSZp)_S := \nu(\SZp(\bu)[2\bu]_S,\xi^{\syn}). 
\]
We further define $\MSZp(n)_S := (\MSZp)_S(n)$. When $S$ is clear from the context, we shall omit it from the notation.
We also denote $\MSFp:= \MSZp \otimes_{\Z_p} \F_p$ (the tensor product being derived, of course). 
\end{defn}

\begin{remark}\label{rem: motivic spectra notation}
One can view ``syntomic cohomology'' as a cohomology theory or as a motivic spectrum. These are morally similar but the latter is a more refined piece of structure, and we distinguish them by the font, with boldface indicating syntomic sheaves and blackboard-bold indicating the corresponding motivic spectra.  
\end{remark}

The most important property of $\MSZp(n)_S$ from our perspective is that it comes equipped with a natural isomorphism
\[
\coh{(\MSZp(n)_S)} \cong \SZp(n)_S;
\]
in other words, the cohomology theory represented by $\MSZp(n)_S$ is the restriction of absolute syntomic cohomology with $n^{\mrm{th}}$ Tate twist to smooth $S$-schemes.

Syntomic cohomology is an ``absolute'' cohomology theory for schemes, not restricted to any specific base. This assertion has refinement to motivic spectra:
$\MSZp(n)_{S}$ is the $S$-component of an ``absolute motivic spectrum''. 
We shall only consider the incarnation of this fact for individual morphisms.

\begin{prop}\label{prop:syntomic_absolute}
Let $f\colon S \to T$ be a morphism of qcqs schemes. Then there is a canonical isomorphism
\[
f^*(\MSZp(n)_T) \cong \MSZp(n)_S.
\]
\end{prop}

\begin{proof}
Since $f^*$ intertwines the Tate twists, it suffices to treat the case $n=0$. 
Let $(\SZp(\bu)[2\bu]_T,\xi^{\syn}_T)$ be the $\pic$-oriented graded algebra from \Cref{def:mot_syn_coh} and similarly for $S$.
First, we shall construct a canonical isomorphism
\begin{equation}\label{eq: syntomic pullback spectrum} 
(f^*)^{\por}(\SZp(\bu)[2\bu]_T,\xi^{\syn}_T) \cong (\SZp(\bu)[2\bu]_S,\xi^{\syn}_S).
\end{equation}
By definition, $f^*\SZp(n)_T$ is constructed by applying the following sequence of operations:
\begin{enumerate}
\item Left Kan extension from smooth $T$-schemes to all (qcqs) $T$-schemes. Let us denote the inclusion of smooth $T$-schemes into all qcqs $T$-schemes by $\mu\colon \Sm_T \to \Sch_T$, and the relevant left Kan extension of presheaves by $\mu_\sharp$, with right adjoint $\mu^*$.  
\item Restriction to smooth $S$-schemes; note that every smooth $S$-scheme acquires the structure of a qcqs $T$-scheme via $f$. 
\item Nisnevich sheafification.
\end{enumerate}
Let $\SZp(\bu)[2\bu]$ be the absolute syntomic cohomology complex, regarded as a sheaf of commutative graded algebras on qcqs schemes. Then $\SZp(\bu)_S[2\bu]$ is obtained from $\SZp(\bu)[2\bu]$ by restriction to smooth $S$-schemes, and similarl for $T$. The counit of the symmetric monoidal adjunction $\mu_\sharp \dashv \mu^*$ is a  natural transformation
\begin{equation}\label{eq:syntomic-kan-extension}
\mu_\sharp \SZp(\bu)_T[2\bu]\cong \mu_\sharp \mu^* (\SZp(\bu)[2\bu]|_{\Sch_T}) \to \SZp(\bu)[2\bu]|_{\Sch_T}, 
\end{equation}
which restricts to a map 
\begin{equation}\label{eq:syntomic-kan-pullback}
(\mu_\sharp \SZp(\bu)_T[2\bu])|_{\Sm_S} \to \SZp(\bu)[2\bu]|_{\Sm_S} = \SZp(\bu)_S[2\bu]. 
\end{equation}
We want to show that \eqref{eq: syntomic pullback spectrum} is an isomorphism; since both sides are Nisnevich sheaves, it suffices to assume that $T,S$ are affine and that \eqref{eq:syntomic-kan-pullback} is an isomorphism. As that arises as restriction from \eqref{eq:syntomic-kan-extension}, it suffices to show that the latter map is an isomorphism, which means in other words that syntomic cohomology for qcqs $T$-schemes is the left Kan extension of its restriction to smooth $T$-schemes. For $T = \Spec \Z$, this is \cite[Proposition 8.4.10]{BL22a}, but the proof applies to an arbitrary affine base. 

To construct \eqref{eq: syntomic pullback spectrum}, it remains to verify that $f^*$ carries the $\pic$-pre-orientation $\xi^{\syn}_T$ to $\xi^{\syn}_S$; this would follow from the functoriality of the first Chern class $c_1^{\syn}\colon \pic \to \SZp(1)[2]$ in arbitrary maps of qcqs schemes. In turn, this functoriality follows by Zariski descent from the affine case, where it is part of the definition of $c_1^{\syn}$ in \cite[\S 7]{BL22a}. This completes the construction of \eqref{eq: syntomic pullback spectrum}. 

Now, recalling that $\tau_S \colon\MS_S^\lax \to \MS_S$ is the left adjoint to the inclusion, we have
\begin{align*}
f^*(\MSZp)_{T} \stackrel{(1)}\cong & \tau_S (f^*)^{\lax} \nu (\SZp(\bu)[2\bu]_T,\xi^{\syn}_T) \\ 
\stackrel{(2)}\cong & \tau_S \nu (f^*)^{\por}  (\SZp(\bu)[2\bu]_T,\xi^{\syn}_T)\\ \stackrel{(3)}\cong & \tau_S \nu (\SZp(\bu)[2\bu]_S,\xi^{\syn}_S) \\
\stackrel{(4)}\cong & \tau_S (\MSZp)_S \stackrel{(5)}\cong (\MSZp)_S,   
\end{align*}
where the isomorphisms are explained as follows:
\begin{enumerate}
\item is direct from the definitions of $(\MSZp)_T$ and $f^*$ (cf. \eqref{eq:pullback_MS}),
\item is the commutativity \eqref{eq: por algebra pullback}, 
\item is \eqref{eq: syntomic pullback spectrum}, 
\item is by definition of $(\MSZp)_S$, and
\item follows from the fact that $(\MSZp)_S$ already belongs to $\MS_S\subseteq \MS_S^\lax$.  
\end{enumerate}
This completes the proof. 
\end{proof}

\subsubsection{Motivic cohomology}\label{ssec:motivic-cohomology}

Our next goal is to explain how the motivic cohomology spectrum arises from an oriented graded algebra. This will be done by running the machinery of \Cref{sssec:from_oriented_to_spectra} in reverse. The task is easier over a field of characteristic zero, and since this is the only case we need, we only consider this generality. 

Let $K$ be a field of characteristic $0$. Voevodsky defined the Eilenberg--MacLane spectra in \cite{Voe10}, thanks to which we have the $p$-adic motivic cohomology spectrum $\MHZp \in \SH_K$. We want to upgrade it to a $\pic$-oriented graded commutative algebra $(\MHZp(\bu)[2\bu], \xi^{\mot}) \in \calg^{\ori}_{\pic}(\cC_K^{\N})$. To achieve this, we will promote $\MHZp$ into a $\pic$-oriented graded algebra in $\SH_K\subseteq \MS_K$ in the following sense (compare \Cref{var:pic_or}). 

\begin{defn}
A \emph{$\pic$-oriented graded algebra in $\MS_K$} is a graded algebra $E_\bullet \in \calg(\MS_K^\N)$ together with a map $\Sigma^\infty \pic \to E_1$ such that its restriction along $\Sigma^\infty \bP^1 \to \Sigma^\infty \pic$ is a $\bP^1$-orientation of $E_\bullet$ in $\MS_K$.  
\end{defn}
The forgetful functor $(-)_0\colon \MS_K \to \cC_K$ induces a functor 
\[
(-)_0^{\ori}\colon \calg^{\ori}_{\pic}(\MS_K^\N) \to \calg^{\ori}_{\pic}(\cC_K^\N). 
\]
Our strategy is to construct the $\pic$-oriented graded algebra $\HZp(\bullet)[2\bullet]$ as the image under $(-)_0^{\ori}$ of a $\pic$-oriented graded algebra $\MHZp(\bullet)[2\bullet]$ in $\SH_K \subseteq \MS_K$. 
 
 We will realize $\MHZp(\bullet)[2\bullet]$ as the associated graded for Voevodsky's \emph{slice filtration} \cite[\S 2]{Voe02} on the motivic $K$-theory spectrum. 
First, we need to recall some aspects of the formalism of the slice filtration. Define $\fslice^{\geq 0} \SH_K \subset \SH_K$ to be the full subcategory generated under colimits and shifts from objects of the form $\Sigma^\infty_+ X$ for $X\in \Sm_\Qpcyc$. Then for all $n \in \Z$,  
\[
\fslice^{\geq n} \SH_K = (\rPP^1)^{\otimes n} \otimes \fslice^{\geq 0} \SH_K.
\]
This defines a $\Z$-indexed filtration by full stable subcategories closed under colimits and compatible with the symmetric monoidal structure, i.e., 
\[
\fslice^{\geq i} \SH_K \otimes \fslice^{\geq j} \SH_K \rightarrow \fslice^{\geq i+j} \SH_K.
\]
It follows formally that there is a symmetric monoidal functor $f^{\ge \bullet}\colon \SH_K \to \Fil(\SH_K)$ (the target being the $\infty$-category of $\Z$-filtered objects in $\SH_K$) such that $f^{\ge i}$ is the coreflection of $\SH_K$ onto $\fslice^{\geq i} \SH_K$, i.e., $f^{\ge i}$ is right adjoint to the embedding $\fslice^{\geq i} \SH_K \subseteq \SH_K$. For $X\in \SH_K$, we further denote by $f^{\le i} X \in \SH_K$ the cofiber of the counit map $f^{\ge i+1} X \to X$, and by $f^i X \in \SH_K$ the cofiber of the map $f^{\ge i+1} X \to f^{\ge i} X$. Hence for the associated graded functor $\gr \co \Fil(\SH_K) \rightarrow (\SH_K)^{\Z}$, we have 
\[
\gr(f^{\ge \bullet} X) \cong \bigoplus_{n \in \Z} f^n X 
\in \SH_K^{\Z}. 
\]

Let $\KGL \in \calg(\SH_K)$ be the $p$-completion of the $\A^1$-invariant algebraic $K$-theory spectrum \cite{MotSpec}. Recall that $\Sph$ denotes the motivic $p$-adic sphere spectrum. Voevodsky proved in  \cite{voevodsky2003zeroslicespherespectrum} (in the setting that $K$ has characteristic zero) that 
\begin{equation}\label{eq:zero-slice}
f^0\Sph \iso f^0\KGL \cong  \MHZp
\end{equation}
as commutative algebras in $\SH_K$. Using the Bott periodicity isomorphism $\KGL(1)[2]\cong \KGL$, we obtain that the associated graded algebra of $f^{\ge \bullet}\KGL$ is given by   
\begin{equation}
\gr(f^{\ge \bullet} \KGL) \cong \MHZp(\bullet)[2\bullet] \in (\SH_K)^{\Z}.
\end{equation}
Note that this identification holds for all $\bullet\in \Z$. Restricting to the non-negative part of the grading, we obtain an identification 
\begin{equation} \label{eq:motivic_gr_kgl}
\gr(f^{\ge \bullet} \KGL) \cong \MHZp(\bullet)[2\bullet] \in (\SH_K)^{\N}
\end{equation}
of $\N$-graded motivic spectra. The LHS has a tautological commutative algebra structure induced by that on $\KGL$, which restricts to the isomorphism \eqref{eq:zero-slice} of commutative algebras in degree $0$. Henceforth, by $\gr(f^{\ge \bullet} \KGL)$ we mean this $\N$-graded commutative algebra. We may therefore transfer the commutative algebra structure to the RHS, to view $\MHZp(\bullet)[2\bullet] \in \calg((\SH_K)^{\N})$. 

As explained in \cite[\S 6]{AHI2}, there is a canonical map $\beta \colon \Sigma^\infty \pic \to f^{\ge 1}\KGL$, which turns $f^{\ge \bullet} \KGL$ into a $\pic$-oriented \emph{filtered} algebra in $\SH_K$ in the sense that after restriction along the map $\bP^1 \to \pic$, the resulting maps 
\[
f^{\ge i}\KGL \to \hom_{\SH_K}(\Sigma^\infty \bP^1,f^{\ge i+1}\KGL) \cong f^{\ge i+1}\KGL(-1)[-2]
\]
(where $\Sigma^\infty \bP^1$ refers to the image in $\SH_\Qpcyc$)
are all isomorphisms. Passing to the associated graded algebra of the filtration, we thus obtain a $\pic$-orientation $\xi \colon \Sigma^\infty \pic \to \MHZp(1)[2]$ of $\MHZp(\bullet)[2\bullet]$. 

\begin{defn}[Motivic cohomology as an Oriented Algebra]\label{defn: motivic cohomology}
Let $\Qpcyc$ be a field of characteristic zero. We define the oriented graded algebra  $(\HZp(\bu)[2\bu],\xi^\mot)\in \calg_{\pic}^{\ori}(\cC_\Qpcyc^\N)$ to be the image of $(\MHZp(\bu)[2\bu],\xi)$ under the functor 
$(-)_0^{\ori}\colon \calg_{ \pic}^{\ori}(\SH_{\Qpcyc}^\N) \to \calg_{\pic}^{\ori}(\cC_{\Qpcyc}^{\N})$. 
\end{defn}

By design, this oriented graded algebra models motivic cohomology. 

\begin{prop}
We have $\nu(\HZp(\bu)[2\bu],\xi^\mot) \cong \MHZp$ in $\calg(\SH_\Qpcyc)$.
\end{prop}

\begin{proof}
Since $(\HZp(\bu)[2\bu],\xi^\mot) = (\MHZp(\bu)[2\bu],\xi)_0^{\ori}$, the claim follows from \Cref{cor:realization_reflective} (where $(-)_0$ is denoted $\Omega^\infty_c$). 
\end{proof}
\begin{remark}
The map $\xi^\mot \colon \Sigma^\infty \pic \to \HZp(1)[2]$ decomposes as  $\Sigma^\infty \pic \to \pic \to \HZp(1)[2]$ where the first one is the canonical map coming from the fact that $\pic$ is a sheaf of spectra and the second map is the $p$-completion and shift of the classical isomorphism $\Z^{\mot}(1)[1]\cong \G_m$.  
\end{remark}

\section{Perfectoid nearby cycles}\label{sec: perfectoid nearby cycles}
Let $\sO$ be a rank-one $p$-adic perfectoid valuation ring, which we moreover assume\footnote{This assumption is likely superfluous, and should be removable using forthcoming work of Bouis--Kundu \cite{BK25} and Bachmann--Elmanto--Morrow \cite{BEM}} can be written as a $p$-completed filtered colimit $\sO = (\colim_n \sO_n)^\wedge_p$ where $\sO_n$ is a local finite $\Z_p$-algebra. Let $K := \sO[1/p]$ be the fraction field of $\sO$ and $k$ be the residue field of $\sO$. 

\begin{example}
The main example of interest to us is $\sO = \Z_p^{\cyc}$, in which case we may for example take $\sO_n = \Z_p[\mu_{p^n}]$. 
\end{example}

In this section we define a functor $\psi \co \SH_{K} \rightarrow \MS_k$, whose form is reminiscent of the nearby cycles functor. (This $\psi$ in fact arises as the restriction of a functor $\MS_K \rightarrow \MS_k$ to $\SH_K \subset \MS_K$.) We will see that, thanks to certain miraculous properties of perfectoid rings, this functor is well-behaved, and carries motivic cohomology (as a motivic spectrum over $K$) to syntomic cohomology (as a motivic spectrum over $k$). One significance of this fact is that \emph{it will eventually allow us to transport knowledge about the motivic Steenrod algebra from characteristic $0$ to knowledge about the syntomic Steenrod algebra in characteristic $p$}, the latter of which is a priori mysterious. In characteristic $0$, the motivic Steenrod algebra was calculated by Voevodsky \cite{Voe03, Voe10}, and through $\psi$ this will give our initial traction on the syntomic Steenrod algebra.

\begin{remark}
Some of our results in this section are similar to some results discovered independently by Bouis--Kundu in the recent paper \cite{BK25}. The technical arguments appear to share some similarities, which we have not attempted to analyze in detail. 

 By combining \cite{BK25} with forthcoming results of Bachmann--Elmanto--Morrow on motivic cohomology in mixed characteristic, it should be possible to ``descend'' our results to the motivic level. A more precise statement is made in \Cref{rem:motivic-descent} below. 
\end{remark}


\subsection{The functor}
We define the (lax symmetric monoidal) functor 
\begin{equation}\label{eq: Psi}
\Psi \co \SH_{\Qpcyc} \rightarrow \MS_{\Zpcyc}
\end{equation}
as the composition of the (lax symmetric monoidal) functors
\[
\Psi\colon \SH_K \inj \MS_K \oto{j_*} \MS_\sO \oto{L_{\et}} \MS_\sO
\]
from \S \ref{sssec:relation-to-SH}, \eqref{eq: MS adjunction}, and \S \ref{ssec: etale sheafification} respectively. 

\begin{remark}\label{rem:Let}
We emphasize that although $L_{\et}$ factors through $\MS_{\sO}^{\et}$, we are forgetting back down to $\MS_{\sO}$ in the definition of $\Psi$. (Later, we use the notation $\cL_{\et}$ for \'etale sheafification without forgetting back down to the Nisnevich topology.)
\end{remark}

\begin{remark}
The definition of $\Psi$ could have been made starting over a general extension $K/\Q_p$, but it would not behave well in general. The highly ramified nature of the perfectoid field $K$ is needed for $\Psi$ to have good properties, a phenomenon that can be traced back to observations of Niziol in \cite{Niz98}. We do expect the analogous definition to behave well whenever $\sO$ is a perfectoid valuation ring over $\Z_p$. 
\end{remark}

Consider the closed embedding  $i \co \Spec k \inj \Spec \Zpcyc$. We define the (lax symmetric monoidal) functor 
\begin{equation}\label{eq: psi}
\psi \co \SH_{\Qpcyc} \rightarrow \MS_{k}
\end{equation}
as the composition of $\Psi$ with $i^* \co  \MS_{\Zpcyc} \rightarrow \MS_{k}$. 

\begin{remark}
The definition of $\psi$ resembles that of the $p$-adic nearby cycles, except we have not base changed to an algebraic closure of $\Qpcyc$ as we would do for the formation of nearby cycles. That this construction behaves well relies crucially on the perfectoid nature of $\Qpcyc$ (and that we are taking $p$-adic coefficients). A related phenomenon appears in \cite[Theorem 10.1]{BMS2}\footnote{Although the result there is stated with $C$ algebraically closed, it is not necessary for the proof.}.
\end{remark}

Over the rest of the section, we will investigate the properties of this functor when evaluated on motivic cohomology (viewed as a motivic spectrum via \Cref{defn: motivic cohomology}). In particular, the goal of this section is to construct an isomorphism 
\[
\psi((\MHZp)_\Qpcyc) \cong (\MSZp)_k \in \MS_k.
\]

\begin{remark}\label{rem: tate twists compatible} The functors $\Psi$ and $\psi$ are each compositions of the functors considered in \S \ref{sec: motivic spectra}, all of which are compatible with the Tate twist in the obvious way. Hence we have (compatible) isomorphisms of functors 
\[
\psi( - (n)) \cong \psi(-)(n) \co \SH_K \rightarrow \MS_k
\]
for all $n \in \Z$, and similarly for $\Psi$. 
\end{remark}

\subsection{``Local constancy'' of syntomic cohomology}

Our first goal is to show that the functor $\Psi$ carries the motivic cohomology spectrum over $\Qpcyc$ to the syntomic cohomology spectrum over $\Zpcyc$, or more precisely to construct an isomorphism 
\begin{equation}\label{eq: Psi motivic}
\Psi((\MHZp)_\Qpcyc) \cong (\MSZp)_{\Zpcyc} \in \MS_{\Zpcyc}.
\end{equation}
We will do this by comparing the oriented graded algebras $(\SZp(\bu)[2\bu]_\Zpcyc,\xi^{\syn})$ and $(\HZp(\bu)[2\bu]_\Qpcyc,\xi^{\mot})$. In fact, rather than comparing them directly, we shall compare each with the oriented graded algebra computing $p$-adic \'{e}tale cohomology. 

\subsubsection{\'Etale comparison for syntomic cohomology}
Let $\EZp(n)_K = \limit_N \mu_{p^N}^{\otimes n}$, regarded as an \'etale sheaf on smooth schemes over $\Spec \Qpcyc$. There is an obvious commutative algebra structure on $\bigoplus_{n \in \N} \EZp(n)_K$, which we shear to obtain a commutative algebra 
\[
\EZp(\bu)[2\bu]_K = (\EZp(n)[2n]_K)_{n \in \N}  \in \calg(\PShv_{\et}(\Sm_{\Qpcyc}; \Sp)^{\N}).
\]
The first Chern class for \'etale cohomology of $\rB\G_m$ equips $\EZp(\bu)[2\bu]_K $ with an orientation, promoting $\EZp(\bu)[2\bu]_K$ to an object of $\calg^{\ori}_{\pic}(\PShv_{\et}(\Sm_{\Qpcyc}; \Sp)^{\N})$. 

\begin{prop}[Bhatt--Lurie]
There exists a canonical morphism of \'etale $\pic$-oriented graded algebras over $\Zpcyc$, 
\[
\gamma_{\et} \colon (\SZp(\bu)[2\bu]_\Zpcyc, \xi^{\syn}) \to j_*(\EZp(\bu)[2\bu]_\Qpcyc, \xi^{\et})   \in \calg^{\ori}_{\pic}(\PShv_{\et}(\Sm_{\Zpcyc}; \Sp)^{\N}).
\]
\end{prop}

\begin{proof}This is the \'etale comparison morphism of Bhatt--Lurie \cite[Theorem 8.3.1]{BL22a}.
\end{proof}


\subsubsection{\'Etale sheafification of motivic cohomology}
Let us write 
\begin{equation}
\cL_{\et} \co \PShv_{\Nis}(\Sm_\Qpcyc;\Sp) \rightarrow \PShv_{\et}(\Sm_\Qpcyc;\Sp)
\end{equation} for the \'etale sheafification functor. Note the distinction from $L_{\et}$, which is the composition
\[
\begin{tikzcd}
\PShv_{\Nis}(\Sm_\Qpcyc;\Sp) \ar[r, "\cL_{\et} "] \ar[dr, "L_{\et}"'] &  \PShv_{\et}(\Sm_\Qpcyc;\Sp) \ar[d, "\text{forget}"'] \\
 &\PShv_{\Nis}(\Sm_\Qpcyc;\Sp)
\end{tikzcd}
\]

Next we use that the \'etale sheafification of $p$-adic motivic cohomology is $p$-adic \'etale cohomology. Such statements were initially proved by Suslin--Voevodsky \cite{SV98}, and reproved by Geisser--Levine in \cite{GL01}, but we want to invoke a more structured version involving the commutative algebra structure on $\bigoplus_{n \in \Z} \HZp(n)[2n]_K$ constructed in \Cref{ssec:motivic-cohomology}, that does not seem to have been established in the literature until \cite{BEM}. More precisely, \cite[Proposition 6.5]{BEM} implies that
\begin{equation}\label{eq: gabber-suslin-structured}
\cL_{\et} (\HZp(\bu)[2\bu]_K,\xi^{\mot}) \cong (\EZp(\bu)[2\bu]_K, \xi^{\et}) \qin \calg^{\ori}_{\pic}(\PShv_{\et}(\Sm_\Qpcyc;\Sp)^{\N}).
\end{equation}

 
Now, we have a Beck--Chevalley comparison map 
\[
\delta \colon \cL_{\et} j_*(\HZp(\bu)[2\bu]_{\Qpcyc},\xi^{\mot}) \to  j_* \cL_{\et}(\HZp(\bu)[2\bu]_{\Qpcyc},\xi^{\mot}) \stackrel{\eqref{eq: gabber-suslin-structured}}\cong j_*(\EZp(\bu)[2\bu]_{\Qpcyc},\xi^{\et}).
\]. 

\subsubsection{Upshot} In summary, we have  a diagram of solid arrows in $\calg^{\ori}_{\pic}(\PShv_{\et}(\Sm_{\Zpcyc}; \Sp))$, 
\begin{equation}\label{eq:diagram_mot_et_syn}
\xymatrix{
 & \cL_{\et} j_*(\HZp(\bu)[2\bu]_{\Qpcyc},\xi^{\mot})\ar^{\delta}[d] \\ 
(\SZp(\bu)[2\bu]_{\sO},\xi^{\syn}) \ar^{\gamma_{\et}}[r]\ar@{..>}^{\sim}[ru] & j_*(\EZp(\bu)[2\bu]_\Qpcyc,\xi^{\et}) 
}
\end{equation}
Our next goal is to construct the dashed isomorphism in \eqref{eq:diagram_mot_et_syn}, making the triangle commute. This will require some preparations.

\subsubsection{Beilinson's t-structure on graded objects} Recall that if $\cC$ is a stably symmetric monoidal $\infty$-category with a compatible\footnote{meaning that the tensor product of connective objects is connective} t-structure $\tau^{\le *}$,  
then the $\infty$-category $\cC^\N$ (of $\N$-graded objects in $\cC$) inherits a symmetric monoidal t-structure $\tau^{\le *}_\gr$ with 
\[
\tau^{\le n}_\gr(X_0,X_1,\dots) = (\tau^{\le n} X_0, \tau^{\le n-1} X_1,\dots, \tau^{\le n-m}X_m,\dots).
\]
This is a graded version of Beilinson's $t$-structure \cite{Bei87}, and we shall therefore refer to it as the \emph{Beilinson $t$-structure on $\cC^\N$}.

\begin{example}\label{ex:beilinson-on-etale-sheaves}
For a qcqs scheme $S$, the $\infty$-category $\PShv_{\et}(\Sm_S;\Sp)$ of \'{e}tale sheaves of $p$-complete spectra admits a canonical $t$-structure for which $\sF$ is connective (i.e., concentrated in non-positive degrees) if and only if $\sF/p$ is connective as an ordinary sheaf of spectra. This induces a Beilinson t-structure on $\PShv_{\et}(\Sm_S;\Sp)^{\N}$. 
\end{example}

We will start by constructing the dashed arrow in \eqref{eq:diagram_mot_et_syn} at the level of commutative graded algebras, ignoring the pre-orientations. At this level, we may regard the objects as lying in $\PShv_{\et}(\Sm_S;\Sp)^{\N}$ with the t-structure of \Cref{ex:beilinson-on-etale-sheaves}, and we will show that both maps $\gamma_{\et}$ and $\delta$ from \eqref{eq:diagram_mot_et_syn} identify their respective sources as the \emph{connective cover} (i.e. the truncation $\tau^{\le 0}$) of $\EZp(\bu)[2\bu]$; this will show in particular that they are canonically isomorphic to each other. 

\subsubsection{The map $\gamma_{\et}$ is a connective cover}
We begin by focusing on the map $\gamma_{\et}$ from \eqref{eq:diagram_mot_et_syn}. 

\begin{prop}\label{prop:comparison_syntomic_etale}
The map 
\[\gamma_{\et}\colon \SZp(\bu)[2\bu]_{\Zpcyc} \to j_* \EZp(\bu)[2\bu]_{\Qpcyc}
\]
exhibits $\SZp(\bu)[2\bu]_{\Zpcyc}$ as the connective cover of $j_* \EZp(\bu)[2\bu]_{\Qpcyc}$ in $\PShv_{\et}(\Sm_{\Zpcyc};\Sp)^{\N}$, with respect to the t-structure of \Cref{ex:beilinson-on-etale-sheaves}. In other words, for every $n \in \N$ the \'etale comparison map induces an isomorphism
\begin{equation}\label{eq:comparison_syntomic_etale}
\gamma_{\et}\BK{n} \colon \SZp(n)_{\Zpcyc} \iso \tau^{\le n} j_* \EZp(n)_{\Qpcyc} \in \PShv_{\et}(\Sm_{\Zpcyc};\Sp).
\end{equation}

\end{prop}

\begin{proof}
Using the definition of the Beilinson $t$-structure and the universal coefficients cofiber sequence, the assertion can be checked after reduction modulo $p$. Hence, it suffices to prove the analogous question with $\F_p$ coefficients instead: the \'etale comparison map  $
\gamma_{\et}\BK{n} \colon \SFp(n)_\Zpcyc \to j_* (\EFp(n)_{\Qpcyc})$ induces an isomorphism $\SFp(n)_\Zpcyc \iso \tau^{\le n} j_* (\EFp(n)_{\Qpcyc})$. 

Since all the sheaves in question are cohomologically bounded below, hence hypercomplete, we can check both properties on the stalks. Let $X\in \Sm_{\Zpcyc}$, and let $x\in X$ with strictly Henselian local ring $R := \mathcal{O}_{X,x}^{\rm{sh}}$.  
Note that by \eqref{eq: gabber-suslin-structured}, we have 
\[
j_* (\EFp(n)_K)(\Spec R)\cong \RGamma_{\et}(\Spec R[1/p];\mu_p^{\otimes n}).
\]
Then it would suffice to show that the \'etale comparison map 
\begin{equation}\label{eq: gamma_n}
\gamma_{\et}\BK{n} \colon \RGamma_{\syn}(\Spec R;\F_p(n))\to \RGamma_{\et}(\Spec R[1/p];\mu_p^{\otimes n})
\end{equation}
induces an isomorphism
\[
\RGamma_{\syn}(\Spec R;\F_p(n))\iso \tau^{\le n}\RGamma_{\et}(\Spec R[1/p];\mu_p^{\otimes n}).
\]

If $x$ lies over the generic point $\Spec \Qpcyc$ of $\Spec \Zpcyc$, then the result follows from the observations that 
\begin{itemize}
\item $\gamma_{\et}\BK{n} $ is (tautologically) an isomorphism for schemes over $\Q_p$, and 
\item $\RGamma_{\et}(\Spec R[1/p];\mu_p^{\otimes n})$ is concentrated in degree $0$ (as $R[1/p] = R$ is strictly Henselian).
\end{itemize}

Assume next that $x$ lies over the special point $\Spec k$. By \cite[Theorem G]{AMMN22}, the cohomology groups of  $\RGamma_{\syn}(\Spec R; \F_p(n))$
are supported in cohomological degrees $\leq n$. By \cite[Theorem 1.8]{BM23}, the map \eqref{eq: gamma_n} induces an isomorphism on cohomology groups in degrees $\leq n-1$, and an injection in degree $n$. It therefore suffices to show that it is surjective in degree $n$. By the Bloch--Kato type isomorphism of \cite[Theorem 3.5]{LM23}, the target is generated by products of cohomology classes in degree $1$, so we are reduced to the case $n=1$, where we want to show that the map
\begin{equation}\label{eq: symbol degree 1}
\gamma_{\et}\BK{1} \co \Hsyn^{1,1}(\Spec R;\F_p) \rightarrow \rH^1_{\et}(\Spec R[1/p]; \mu_p)
\end{equation}
is an isomorphism. By the Kummer sequence and the definition of $\gamma_{\et}\BK{1}$, \eqref{eq: symbol degree 1} identifies with the map
\begin{equation}
R^{\times} \otimes_{\Z} \Z/p \rightarrow R[1/p]^{\times} \otimes_{\Z} \Z/p.
\end{equation}
induced by $R \subset R[1/p]$. By construction, $R$ is a filtered colimit of smooth $\Zpcyc$-algebras along \'etale maps, or in other words $R$ is \emph{essentially smooth} over $\Zpcyc$, so we reduce to the following Lemma (taking $\varpi=p$). 
\end{proof}

\begin{lemma}
Let $R$ be an essentially smooth local algebra over a perfectoid valuation ring $\Zpcyc$, with pseudo-uniformizer $\varpi$. Then the map 
\[
R^\times \otimes_{\Z} \Z/p \rightarrow R[1/\varpi]^\times  \otimes_{\Z} \Z/p
\]
is surjective. 
\end{lemma}

\begin{proof}
Let $\mf p$ be the maximal ideal of $\Zpcyc$; we shall use the same notation for its extension to a prime ideal of $R$.

\emph{Claim:} the localization $R_{\mf p}$ is a valuation ring, with the same value group as $\Zpcyc$. Granting the claim for now, let us see how to finish the proof. Since $R$ is essentially smooth, it is normal, so we have $R = R_{\mf p} \cap R[1/\varpi]$, hence
\[
R^\times   = R_{\mf p}^\times \cap R[1/\varpi]^\times.
\]
This implies that the obvious map 
\begin{equation}\label{eq: value group}
\frac{R[1/\varpi]^\times}{R^\times} \rightarrow \frac{R_{\mf p}[1/\varpi]^\times}{R_{\mf p}^\times }
\end{equation}
is injective, and it is surjective by the claim, so it is an isomorphism. In other words, the cokernel of the map
$R^\times \rightarrow R[1/\varpi]^\times$ is isomorphic to the value group of $R_{\mf p}$. But the claim asserts that this value group is isomorphic to the value group of $\Zpcyc$, and the latter is $p$-divisible since $\Zpcyc$ is perfectoid \cite[Lemma 3.2]{Sch12}, so it vanishes modulo $p$. Thanks to \eqref{eq: value group} being an isomorphism, this completes the proof up to establishing the claim, which we do next.

For this, we present $R_{\mf p} = \colim_i R_i$ as a filtered colimit of localizations of smooth $\Zpcyc$-algebras along \'etale transition maps (with each $R_i$ being local). By the structure theorem for smooth maps, $R_i$ is a localization of an \'etale algebra over a polynomial ring $\Zpcyc[T_1, \ldots, T_n]_{\mf p}$. If the claim is granted for $\Zpcyc[T_1, \ldots, T_n]_{\mf p}$, then by \cite[Tag 0ASF, specifically Tag 0ASJ]{stacks-project}, each $R_i$ is a valuation ring and the transition maps induce isomorphisms on value groups. Hence we have reduced to the case $R_{\mf p} = \Zpcyc[T_1, \ldots, T_n]_{\mf p}$. In that case, from direct inspection of the definitions we see that $R_{\mf p}$ is a valuation ring for the valuation induced by the ``Gauss norm'' on $\sO[T_1, \ldots, T_n]$, 
\[
||\sum_I a_I T^I||  = \max_I |a_I|,
\]
which clearly has the same value group as $\Zpcyc$. This establishes the claim, which completes the proof. 
\end{proof}

\subsubsection{The map $\delta$ is a connective cover} 
Next, we consider the comparison map $\delta$ from \eqref{eq:diagram_mot_et_syn}. We shall similarly show that it is a connective cover. 

\begin{prop}\label{prop:comparison_motivic_etale}
The map 
\[
\delta \colon \cL_{\et} j_* (\HZp(\bu)[2\bu]_{\Qpcyc}) \to j_*(\EZp(\bu)[2\bu]_{\Qpcyc})
\]
exhibits $\cL_{\et}j_*(\HZp(\bu)[2\bu]_{\Qpcyc})$ as the connective cover of $j_*(\EZp(\bu)[2\bu]_{\Qpcyc})$ in $\PShv_{\et}(\Sm_{\sO};\Sp)^{\N}$, with respect to the t-structure of \Cref{ex:beilinson-on-etale-sheaves}. In other words, for every $n \in \N$ it induces an isomorphism 
\begin{equation}\label{eq:comparison_motivic_etale}
\cL_{\et} j_* (\HZp(n)_{\Qpcyc}) \xrightarrow{\sim} \tau^{\le n}j_*(\EZp(n)_{\Qpcyc} ) \in \PShv_{\et}(\Sm_S;\Sp).
\end{equation}
\end{prop}

\begin{proof}
According to the Beilinson--Lichtenbaum Conjecture proved by Voevodsky, for smooth $Y/\Qpcyc$ we have a natural isomorphism (cf. \cite[\S 1.4]{HW19})
\begin{equation}\label{eq: bloch-kato}
\HFp(n)|_{Y_{\Zar}} \cong \tau^{\leq n} (\rR \nu_* \mu_p^{\otimes n})
\end{equation}
where $\nu \co Y_{\et} \rightarrow Y_{\Zar}$ is the change-of-topology map between the small \'{e}tale and small Zariski sites of $Y$. 

As in the proof of \Cref{prop:comparison_syntomic_etale},
in order to prove that \eqref{eq:comparison_motivic_etale} is an isomorphism we may reduce to $\F_p$-coefficients and check the claim at stalks along $R:=\cO_{X,x}^{\mrm{sh}}$ for all $X \in \Sm_{\Zpcyc}$ and $x \in X$.

By \eqref{eq: bloch-kato} $j_*(\HFp(n)_{X})$ is the sheaf given by
\begin{equation}\label{eq: bloch-kato-U}
(j_*\HFp(n)_K)(U) \cong \RGamma_{\Zar}(U[1/p]; \tau^{\leq n} (\rR\nu_* \mu_p^{\otimes n})) \quad \text{ for all } U\in X_{\Zar}.
\end{equation}
Taking the colimit of the isomorphisms \eqref{eq: bloch-kato-U} over \'{e}tale neighborhoods $U$ of $x\in X$, we obtain a commutative triangle 
\begin{equation}\label{eq:bloch-kato-triangle}
\xymatrix{
(j_*\HFp(n)_{\Qpcyc})_{(X,x)} \ar^{\delta_x}[rd] \ar^{\sim}[rr] & &
 \RGamma_{\Zar}(\Spec R[1/p];\tau^{\le n}\rR\nu_*\mu_p^{\otimes n}) \ar[ld]
 \\ 
 & \RGamma_{\et}(\Spec R[1/p];\mu_p^{\otimes n}).
}
\end{equation}
The map from right-to-bottom in \eqref{eq:bloch-kato-triangle} tautologically becomes an isomorphism after applying the truncation functor $\tau^{\le n}$ to both sides, and we want to show that it induces an isomorphism 
\begin{equation}\label{eq: delta truncation}
\RGamma_{\Zar}(\Spec R[1/p] ;\tau^{\le n}\rR\nu_*\mu_p^{\otimes n}) \xrightarrow{\sim} \tau^{\le n}\RGamma_{\et}(\Spec R[1/p] ;\mu_p^{\otimes n}).
\end{equation}
For this, it suffices to show that the source of the map is already cohomologically bounded by $n$. By the hypercohomology spectral sequence 
\[
\rH^i_{\Zar}(\Spec R[1/p] ; \rR^j\nu_*\mu_p^{\otimes n}) \Rightarrow \rH^{i+j}_{\Zar}(\Spec R[1/p] ; \rR\nu_*\mu_p^{\otimes n}) 
\]
this follows from the next Lemma. 
\end{proof}

\begin{lemma}\label{lem:zariski-coh-vanish}
Let $X$ be a smooth $\Zpcyc$-scheme and let $x\in X$ with strict Henselization $R := \cO_{X,x}^{\mrm{sh}}$. Then for all $j \leq n$ and all $i >0$ we have
 \begin{equation}\label{eq: coh vanishing 2}
 \rH^i_{\Zar} (\Spec R[1/p]; \rR^j \nu_* \mu_p^{\otimes n}) = 0. 
\end{equation}
\end{lemma}

\begin{proof}
The proof will be carried by reduction to the case of smooth schemes over a discrete valuation ring rather than over a perfectoid base. 
Recall that we assumed that $\Zpcyc$ could be presented as a filtered colimit $\Zpcyc = (\colim_n \Zpcyc_n)^\wedge_p$, where $\Zpcyc_n$ is a finite extension of $\Z_p$. Let $\Zpcyc':=\bigcup_n \Zpcyc_n$. By \cite[Corollary 1.4]{Tang24}, the map $\Zpcyc' \rightarrow \Zpcyc$ is ind-smooth. Hence we may write any smooth $R/\Zpcyc$ as a filtered colimit
\[
R = \colim_m R_m, \quad \text{where each $R_m$ is an essentially smooth local algebra over $\Zpcyc_m$.}
\]
This presents $\Spec R$ as a cofiltered inverse limit of qcqs schemes with affine transition maps, hence by \cite[Tag 03Q4]{stacks-project} the natural map 
\[
\colim_{m} \rH^i_{\Zar}(\Spec R_m[1/p];\rR^j\nu_*\mu_p^{\otimes n}) \rightarrow \rH^i_{\Zar}(\Spec R[1/p] ;\rR^j\nu_*\mu_p^{\otimes n})  
\]
is an isomorphism. Since $R_m$ is an essentially smooth local algebra over the discrete valuation ring $\Zpcyc_m$, \cite[Lemma 4.2(i)]{LM23} applies to say that for $i>0$ and $j\le n$ each of the terms in the filtered system on the LHS vanish. Hence their colimit vanishes, completing the proof. 
\end{proof}

\begin{cor}\label{cor:comparison_syntomic_motivic_without_or}
There is a unique\footnote{Up to contractible space of choices, as usual.} isomorphism $\SZp(\bu)[2\bu]_\Zpcyc \iso \cL_{\et} j_*(\HZp(\bu)[2\bu]_\Qpcyc)$ of commutative graded algebras in $\PShv_{\et}(\Sm_\Zpcyc;\Sp)^{\N}$, together with a commutative triangle 
\begin{equation}\label{eq:diagram_mot_et_syn_without_or}
\xymatrix{
 & \cL_{\et} j_*(\HZp(\bu)[2\bu]_{\Qpcyc})\ar^{\delta}[d] \\ 
\SZp(\bu)[2\bu]_{\Zpcyc} \ar^{\gamma_{\et}}[r]\ar@{..>}^{\sim}[ru] & j_* (\EZp(\bu)[2\bu]_\Qpcyc) 
}.
\end{equation}

\end{cor}

\begin{proof}
We have already constructed $\delta$ and $\gamma$ as maps of commutative graded algebras in $\calg(\PShv_{\et}(\Sm_\Zpcyc;\Sp)^{\N})$. By \Cref{prop:comparison_syntomic_etale} and \Cref{prop:comparison_motivic_etale}, we have identifications of both $\SZp(\bu)[2\bu]_\Zpcyc$ and $\cL_{\et} j_*(\HZp(\bu)[2\bu]_\Qpcyc)$ with the connective cover of $j_*(\EZp(\bu)[2\bu]_\Qpcyc)$ for the Beilinson t-structure of \Cref{ex:beilinson-on-etale-sheaves}, at the level of underlying $\N$-graded sheaves of spectra (i.e., forgetting the commutative algebra structure). The connective cover of a commutative algebra inherits a unique commutative algebra structure, making it universal among maps from connective commutative algebras. Thus, $\delta$ and $\gamma_{\et}$ automatically identify their respective sources with this connective cover, as objects of $\calg(\PShv_{\et}(\Sm_\Zpcyc;\Sp)^{\N})$. 
\end{proof}

\subsection{Matching the $\pic$-orientations}
By pushing the isomorphism of Corollary \ref{cor:comparison_syntomic_motivic_without_or}
down to the Nisnevich site, we obtain an isomorphism 
\begin{equation}\label{eq:comparison_syntomic_motivic_without_or}
\SZp(\bu)[2\bu]\iso L_{\et}j_* (\HZp(\bu)[2\bu]_{\Qpcyc}) \in \calg(\PShv_{\Nis}(\Sm_{\Zpcyc}; \Sp)^{\N}).
\end{equation}
In order to promote this isomorphism to the level of motivic spectra, it remains (by the construction of $\MSZp$ and $\MHZp$ as motivic spectra in \S \ref{ssec: promoting}) to show that this isomorphism is compatible with the respective $\bP^1$-orientations. In fact, we even have compatibility with the $\pic$-orientations. 

\begin{prop}\label{prop:comparison_etale_motivic_or_to_or}
The isomorphism \eqref{eq:comparison_syntomic_motivic_without_or} carries the $\pic$-orientation $\xi^{\syn}$ of $\SZp(\bu)[2\bu]_{\Zpcyc}$ to the $\pic$-orientation $\xi^{\mot}$ of $L_{\et}j_*(\HZp(\bu)[2\bu]_\Qpcyc)$; in other words, it gives the desired commutative triangle \eqref{eq:diagram_mot_et_syn}. 
\end{prop}

\begin{proof}
Let $\xi'$ be the image of $\xi^{\syn}$ under \eqref{eq:comparison_syntomic_motivic_without_or}, so that we wish to prove that $\xi' = \xi^{\mot}$. By the construction of $\xi^{\et}$ we have $\delta(\xi^{\mot}) = \xi^{\et}$ while by the defining property of the Bhatt-Lurie comparison map $\gamma_{\et}$ we have $\delta(\xi') = \gamma_{\et}(\xi^{\syn}) = \xi^{\et}$. Hence the images of the $\pic$-orientations $\xi'$ and $\xi^{\mot}$ under the comparison map $\delta \colon L_{\et}j_*(\HZp(\bu)[2\bu]_{\Qpcyc}) \to j_*(\EZp(\bu)[2\bu]_{\Qpcyc})$ 
agree. To conclude, it remains to show that the map $\delta$ is injective on homotopy classes of pre-orientations, i.e., we have an injection 
\[
 \Hom_{\cC_\Zpcyc}(\Sigma^\infty\pic,L_{\et}j_*(\HZp(1)[2]_{\Qpcyc})) \to \Hom_{\cC_\Qpcyc}(\Sigma^\infty\pic,\EZp(1)[2]_{\Qpcyc}). 
\]
Using \Cref{cor:comparison_syntomic_motivic_without_or}, we can identify this map with the map 
\[
\gamma_{\et}\BK{1} \colon  \rHsyn^{2,1}(\pic_{\Zpcyc};\Z_p) \to \rHet^{2,1}(\pic_{\Qpcyc};\Z_p)
\]
where the tilde indicates \emph{reduced} cohomology. Since the map $\gamma_{\et}$ is a map of $\pic$-oriented theories, we have a commutative square 
\[
\xymatrix{
\Hsyn^{0,0}(\Spec \Zpcyc ;\Z_p) \ar^{\gamma_{\et}\BK{0}}[d] \ar[r] & \rHsyn^{2,1}(\pic_{\Zpcyc};\Z_p) \ar^{\gamma_{\et}\BK{1}}[d]  \\ 
\Het^{0,0}(\Spec \Qpcyc ;\Z_p) \ar[r] & \rHet^{2,1}(\pic_\Qpcyc;\Z_p).  
}
\]
Since $\gamma_{\et}\BK{0}$ is clearly an isomorphism, it would suffice to show that the horizontal maps are isomorphisms.

To see these, note that by \cite[Theorem 9.3.1]{BL22a} (applied to $\mrm{BGL}_1 = \pic$) we have 
\[
\rHsyn^{2,1}(\pic_\Zpcyc;\Z_p) \cong \bigoplus_{n = 1}^\infty \Hsyn^{2-2n,1-n}(\Spec \Zpcyc ;\Z_p),
\]
and it is classical (and follows) that a similar result holds for \'{e}tale cohomology.  
All the terms with $n>1$ vanish for degree reasons, and the projection to the first summands are precisely the horizontal maps in the square above, proving the result. 
\end{proof}

Thanks to \Cref{prop:comparison_etale_motivic_or_to_or}, the isomorphism \eqref{eq:comparison_syntomic_motivic_without_or} promotes to an isomorphism of $\pic$-pre-oriented (hence a fortiori also of $\bP^1$-pre-oriented) commutative algebras, 
\begin{equation}\label{eq:comparison_syntomic_motivic_por}
(\SZp(\bu)[2\bu], \xi^{\syn}) \iso L_{\et}^{\por} j_* (\HZp(\bu)[2\bu]_{\Qpcyc}, \xi^{\mot})  \in \calg^{\por}_{\pic}(\cC_{\Zpcyc}^{\N}).
\end{equation}

\begin{cor}\label{cor:Psi_syntomic}
The functor $\Psi$ carries the motivic cohomology object $(\MHZp)_{\Qpcyc} \in \calg(\SH_{\Qpcyc})$ to the syntomic cohomology object $(\MSZp)_{\Zpcyc} \in \calg(\MS_{\Zpcyc})$: 
\[
\Psi((\MHZp)_{\Qpcyc}) \cong (\MSZp)_{\Zpcyc}  \in \calg(\MS_{\Zpcyc}).
\]
\end{cor}

\begin{remark}
Note that this implies that $\Psi(\MHZp(n)_\Qpcyc)\cong \MSZp(n)_\Zpcyc$ for all $n\in \Z$, by the compatibility of $\Psi$ with Tate twists (Remark \ref{rem: tate twists compatible}).
\end{remark}
 
\begin{proof} Straight from the definitions of $\Psi$ and of motivic cohomology as a motivic spectrum (cf. \Cref{defn: motivic cohomology}), we have 
\begin{equation}\label{eq:Psi-MHZp}
\Psi((\MHZp)_\Qpcyc)  :=  L_{\et}j_*\nu (\HZp(\bu)[2\bu]_K,\xi^{\mot})  :
= \tau_{\Zpcyc} L_{\et}^{\lax} j_* \nu (\HZp(\bu)[2\bu]_K,\xi^{\mot})
\end{equation}
where $\tau_{\sO} \colon \MS_{\sO}^\lax \to \MS_{\sO}$ is from \eqref{eq:tau-S} and $\nu \colon \calg^{\ori}_{\bP^1}(\cC_{S}^\N) \to \calg(\MS_S)$ is from \eqref{eq: or to ms}. Then we rewrite this as 
\begin{align*}
\Psi((\MHZp)_\Qpcyc) = \tau_{\Zpcyc} L_{\et}^{\lax} j_* \nu (\HZp(\bu)[2\bu]_K,\xi^{\mot}) & \stackrel{(1)}\cong \tau_{\Zpcyc} \nu L_{\et}^{\por} j_*(\HZp(\bu)[2\bu]_K,\xi^{\mot}) \\ 
&\stackrel{(2)}\cong \tau_{\Zpcyc} \nu  (\SZp(\bu)[2\bu]_{\Zpcyc},\xi^{\syn}) \\ &\stackrel{(3)}\cong \tau_{\Zpcyc} (\MSZp)_{\Zpcyc} \\
&\stackrel{(4)}= (\MSZp)_{\Zpcyc}.  
\end{align*}
where: 
\begin{enumerate}
\item is the compatibility of $\nu$ with pushforward (see \eqref{eq:pushforward_pre-oriented_beck_chevalley}) and \'etale sheafification (see \Cref{rem:etale_sheafification_por_compatibility} and \eqref{eq:etale_sheafification_por_compatibility}). 
\item is \eqref{eq:comparison_syntomic_motivic_por}. 
\item is the definition of $\MSZp$ as a motivic spectrum (cf. \Cref{def:mot_syn_coh}). 
\item is because $\MSZp$ is already in $\MS_{\Zpcyc} \subseteq \MS_{\Zpcyc}^{\lax}$. 
\end{enumerate}
\end{proof}

\begin{remark}\label{rem:motivic-descent}
We can contemplate the ``motivic descent'' of the functor $\Psi$, defined similarly but \emph{without} \'etale sheafification: 
\[
\Psi^{\mot}  \co \SH_K \inj \MS_K \xrightarrow{j_*} \MS_{\sO}.
\]
By combining \cite[Corollary D]{BK25} with \cite{Bou24} and forthcoming results of Bachmann--Elmanto--Morrow on motivic cohomology in mixed characteristic, it should be possible to prove that $\Psi^{\mot}  ((\MHZp)_K) \cong (\MHZp)_{\sO} \in \calg(\MS_{\sO})$, and that \'etale sheafifying this identity recovers \Cref{cor:Psi_syntomic}. At present, even the meaning of the object ``$(\MHZp)_{\sO}$'' is somewhat ambiguous, as there are multiple approaches to motivic cohomology which do not obviously agree in this (non-noetherian, mixed characteristic) setting; we understand that this is one of the issues which is addressed by \cite{Bou24} and \cite{BEM}. 
\end{remark}

\subsection{Restriction to the special fiber}
Consider the closed embedding  $i \co \Spec k \inj \Spec \Zpcyc$. Recall that we defined the (lax symmetric monoidal) functor 
\begin{equation}
\psi \co \SH_{\Qpcyc} \rightarrow \MS_{k}
\end{equation}
as the composition of $\Psi$ with $i^* \co  \MS_{\Zpcyc} \rightarrow \MS_{k}$. 

\begin{cor}\label{cor: psi Z_p}The functor $\psi$ carries the motivic cohomology object $(\MHZp)_{\Qpcyc} \in \calg(\SH_{\Qpcyc})$ to the syntomic cohomology object $(\MSZp)_{k} \in \calg(\MS_{k})$: 
\begin{equation}\label{eq: psi Z_p 1}
\psi((\MHZp)_{\Qpcyc}) \cong (\MSZp)_{k}  \in \calg(\MS_{k}).
\end{equation}
Furthermore, for all $n \in \Z$, we have 
\begin{equation}\label{eq: psi Z_p 2}
\psi((\MHZp(n))_\Qpcyc) \cong (\MSZp(n))_{k} \in \MS_k.
\end{equation}
\end{cor}

\begin{proof} \Cref{cor:Psi_syntomic} gives an isomorphism
\[
\psi((\MHZp)_\Qpcyc) \cong  i^*\Psi((\MHZp)_{\Qpcyc}) \cong i^*(\MSZp)_{\Zpcyc} \in \calg(\MS_{k})
\]   
and \Cref{prop:syntomic_absolute} gives an isomorphism 
\[
 i^*(\MSZp)_{\Zpcyc} \cong (\MSZp)_{k} \in \calg(\MS_{k}).
\]
Composing these gives the isomorphism \eqref{eq: psi Z_p 1}. Then \eqref{eq: psi Z_p 2} follows from compatibility with twisting, Remark \ref{rem: tate twists compatible}.
\end{proof}

\subsection{Technical properties}

We record the following technical result, which will be needed later.

\begin{lemma}\label{lem: psi direct sum}
Suppose $\{a_m\}$ and $\{b_m\}$ are two sequences of integers with $\lim_{m \rightarrow \infty} a_m = \infty$. Then, with $\MHFp \in \SH_{\Qpcyc}$ denoting the mod $p$ motivic cohomology spectrum over $\Qpcyc$, the natural assembly map
\begin{equation}\label{eq: Psi direct sum comparison}
 \bigoplus_m \Psi(\MHFp)[a_m](b_m) \rightarrow \Psi\left(\bigoplus_m \MHFp[a_m](b_m) \right)
\end{equation}
is an isomorphism in $\MS_{\Zpcyc}$, and the natural assembly map
\begin{equation}\label{eq: direct sum comparison}
\bigoplus_m \psi(\MHFp)[a_m](b_m) \rightarrow \psi\left(\bigoplus_m \MHFp[a_m](b_m) \right)
\end{equation}
is an isomorphism in $\MS_k$.
\end{lemma}

\begin{proof}
By definition, $\psi$ is the composition $\SH_{\Qpcyc} \xrightarrow{\Psi} \MS_{\Zpcyc} \xrightarrow{i^*} \MS_k$. Since $i^*$ is a left adjoint, it preserves all colimits. Therefore, \eqref{eq: direct sum comparison} being an isomorphism follows from \eqref{eq: Psi direct sum comparison} being an isomorphism.

In turn, $\Psi$ is defined as a composition $\SH_{\Qpcyc} \inj \MS_{\Qpcyc} \xrightarrow{j_*} \MS_{\Zpcyc} \xrightarrow{L_{\et}} \MS_{\Zpcyc}$. The functor $\SH_{\Qpcyc} \inj \MS_{\Qpcyc}$ preserves all colimits (\S \ref{sssec:relation-to-SH}), as does $j_*$ (\Cref{prop:push_MS_colim}), so the crux is to control $L_{\et}$.\footnote{Here it matters that we defined $L_{\et}$ as (1) \'etale sheafification followed by (2) forgetting back down to the Nisnevich topology, cf. Remark \ref{rem:Let}. Step (1) preserves all colimits, being a left adjoint, but step (2) does not.} In general, for a collection of motivic spectra $\{E_m \in \MS_{\Zpcyc}\}$, there is an assembly map 
\begin{equation}\label{eq: assembly}
\bigoplus_m L_{\et} (E_m) \rightarrow L_{\et} \left( \bigoplus_m E_m \right) \in \MS_{\Zpcyc},
\end{equation}
exhibiting the RHS as the \'etale sheafification of the LHS. In the case at hand, we have $E_m := j_* \MHFp[a_m](b_m)$, so the LHS of \eqref{eq: assembly} is identified by \Cref{cor:Psi_syntomic} with 
\begin{equation}\label{eq: direct sum of Psi}
\bigoplus_m L_{\et} (E_m) \cong \bigoplus_m \MSFp[a_m](b_m)_{\Zpcyc}
\end{equation}
To show that \eqref{eq: assembly} is an isomorphism, we will argue that \eqref{eq: direct sum of Psi} already satisfies \'etale descent. For this, it suffices to see that the natural map
\begin{equation}\label{eq: sum to product}
\bigoplus_m  \MSFp[a_m](b_m)_{\Zpcyc} \rightarrow \prod_m \MSFp[a_m](b_m)_{\Zpcyc}
\end{equation}
is an isomorphism, because the RHS satisfies \'etale descent (as limits preserve the sheaf property). By \Cref{forg_MS_sheaves}, it suffices to show that the map \eqref{eq: sum to product} is an isomorphism when evaluated on any smooth scheme $X$ over $\Zpcyc$ (as any additional Tate twist can be absorbed into the formulation by shifting the $b_m$'s). Since affine schemes form a basis for the Nisnevich topology, it suffices to consider the case where $X$ is affine. Then the assertion that it is an isomorphism can be checked in each fixed cohomological degree, where (using the fact that evaluation of Nisnevich sheaves on $X$ preserves infinite direct sums by \Cref{lem:Nis_Finitary}) it becomes the statement that
\[
\bigoplus_m \Hsyn^{i+a_m, b_m}(X) \rightarrow  \prod_m \Hsyn^{i+a_m, b_m}(X)
\]
is an isomorphism for every $i \in \Z$. But by  \Cref{prop: vanishing range syntomic cohomology} and the assumption that $\lim_{m \rightarrow \infty} a_m = \infty$, for each fixed $i$ all but finitely many of the factors $\Hsyn^{i+a_m, b_m}(X)$ vanish, so this is clear. 
\end{proof}

\section{Categories of syntomic spectra}\label{sec: syntomic spectra}

One of the motivations for Drinfeld's and Bhatt--Lurie's ``stacky'' approach to prismatic cohomology was to define appropriate categories of \emph{modules} for prismatic (resp. syntomic) cohomology. We are interested in the generalization of this problem for spectral syntomic cohomology. The approach of Drinfeld and Bhatt--Lurie is by first constructing stacks, and then taking their categories of quasicoherent sheaves. We will develop a different approach that accesses the desired categories without constructing stacks. In this section, which constitutes the first step of our approach, we will apply some general categorical constructions to $\MS_S$ in order to define certain module categories for spectral syntomic cohomology, or what we call \emph{syntomic spectra} in short. 

\subsection{Module categories for cosimplicial commutative algebras}
Let $\cC$ be a symmetric monoidal $\infty$-category. For $R\in \calg(\cC)$, we can form the $\infty$-category of modules $\Mod_R(\cC)$. 

\subsubsection{Cosimplicial modules} When $\cC$ is presentable\footnote{or more generally, when $\cC$ has geometric realizations of simplicial objects that distribute over its tensor product.}, there is a straightforward generalization of $\Mod_R(\cC)$ to diagrams of commutative algebras in $\cC$, and in particular for cosimplicial diagrams. 

\begin{defn}\label{def: simplicial module cat}
Let $\Delta$ be the simplex category, with objects $[n]:= \{0, 1, \ldots, n\}$ for $n \in \N$. Let $\cC$ be a presentably symmetric monoidal $\infty$-category and let $R_\bullet\in \calg(\cC)^{\Delta}$ be a cosimplicial commutative algebra in $\cC$. We denote 
\[
\Mod_{R_\bullet}(\cC):= \limit_{[n]\in \Delta} \Mod_{R_n}(\cC). 
\]

Concretely, an object of $\Mod_{R_\bullet}(\cC)$ consists of a system of modules $\{M_n \in \Mod_{R_n}(\cC)\}_{[n]\in \Delta}$ together with compatible isomorphisms $M_n\otimes_{R_n} R_{n+1} \cong M_{n+1}$ for each of the face maps $[n]\to [n+1]$. 
\end{defn}

\begin{example}
If $R_\bullet$ is the constant cosimplicial diagram on $R$, then we obtain a natural equivalence $\Mod_R(\cC) = \Mod_{R_\bullet}(\cC)$.
\end{example}

\subsubsection{Lax cosimplicial modules} We now define a lax variant of cosimplicial module categories. 

\begin{defn}\label{def:compare_module_lim_module_lax}
With $\Delta$, $\cC$, and $R_{\bu}$ as in \Cref{def: simplicial module cat}, let $\cC^\Delta:= \Fun(\Delta,\cC)$, endowed with the pointwise symmetric monoidal structure. Since $R_\bullet$ is a commutative algebra in $\cC^\Delta$, we can form the module category $\Mod_{R_\bullet}(\cC^{\Delta})$. In fact, this module category is the \emph{lax limit} of the diagram $[n]\mapsto \Mod_{R_n}(\cC)$.

Concretely, objects of $\Mod_{R_\bullet}(\cC^{\Delta})$ consist of systems of modules $\{M_n \in \Mod_{R_n}(\cC)\}_{[n]\in \Delta}$ together with compatible comparison maps $M_n\otimes_{R_n} R_{n+1} \to M_{n+1}$ which are \emph{not} required to be isomorphisms. 

\end{defn}

\begin{remark}\label{rem:compare_module_lim_module_lax}
Since the morphisms in the diagram $[n]\mapsto \Mod_{R_n}(\cC)$ are all colimit-preserving, the limit is a lax symmetric monoidal colocalization of the lax limit, hence we have a lax symmetric monoidal limit-preserving functor 
\begin{equation}\label{eq:U_R}
U_R\colon \Mod_{R_\bullet}(\cC^\Delta) \to \Mod_{R_\bullet}(\cC)  
\end{equation}
whose left adjoint is the embedding of the objects for which the comparison maps $M_n\otimes_{R_n} R_{n+1} \to M_{n+1}$ are all isomorphisms.  
\end{remark}

We are interested in the following type of examples.

\begin{example}\label{ex: cocech}
For $R\in \calg(\cC)$, we obtain a cosimplicial commutative algebra $\Am{R}\in \calg(\cC)^\Delta$ as the coČech nerve of the map $\one_\cC \to R$: 
\[
\Am{R}:=
\left(
\begin{tikzcd}
    R
    \ar[r,yshift=1.5pt,->]
    \ar[r,yshift=0pt,<-]
    \ar[r,yshift=-1.5pt,->]
    &
    R \otimes R
    \ar[r,yshift=3pt,->]
    \ar[r,yshift=1.5pt,<-]
    \ar[r,yshift=0pt,->]
    \ar[r,yshift=-1.5pt,<-]
    \ar[r,yshift=-3pt,->]
    &
   R \otimes R \otimes R
    \ar[r,yshift=4.5pt,->]
    \ar[r,yshift=3pt,<-]
    \ar[r,yshift=1.5pt,->]
    \ar[r,yshift=0pt,<-]
    \ar[r,yshift=-1.5pt,->]
    \ar[r,yshift=-3pt,<-]
    \ar[r,yshift=-4.5pt,->]
    &
    \cdots 
\end{tikzcd} 
\right)
\]
\end{example}

\begin{remark}
The category of cosimplicial modules $\Mod_{R^{\otimes \bullet}}(\cC^\Delta)$ is closely related to the construction of \emph{synthetic $R$-modules} (as developed in \cite{pstrkagowski2023synthetic}). Correspondingly, all our applications below could be seen as a ``synthetic'' approach to the syntomic Steenrod algebra, but we will make no direct reference to this point of view.
\end{remark}

\subsubsection{Functoriality}

For a symmetric monoidal colimit-preserving functor $\phi \colon \cC \to \cD$ and $R_\bullet\in \calg(\cC)^\Delta$, the functor $\phi$ induces a canonical symmetric monoidal colimit-preserving functor $\cC^\Delta \to \cD^\Delta$, and then a functor  
\[
\phi \colon  
\Mod_{R_\bullet}(\cC^\Delta) \to \Mod_{\phi(R_\bullet)}(\cD^\Delta),
\]
which restricts to a functor 
\begin{equation}\label{eq:sym-mon-wt-phi}
\wt{\phi} \colon  
\Mod_{R_\bullet}(\cC) \to \Mod_{\phi(R_\bullet)}(\cC).
\end{equation}
Indeed, since $\phi$ commutes with relative tensor products, the comparison maps for $\phi(M_\bullet)$ are identified with the image of the comparison maps for $M_\bullet$ under $\phi$, hence remain isomorphisms. 

Now suppose that $\phi$ is only lax symmetric monoidal (and not necessarily colimit-preserving). Since module categories are natural in lax symmetric monoidal functors, we still have an induced functor 
\[
\phi \co  \Mod_{R_\bullet}(\cC^\Delta) \to \Mod_{\phi(R_\bullet)}(\cD^\Delta) 
\]
between the lax limits, and we can now define:
\begin{equation}\label{eq:wt-phi}
\wt{\phi}\colon \Mod_{R_\bullet}(\cC) \into \Mod_{R_\bullet}(\cC^\Delta) \oto{\phi} \Mod_{\phi(R_\bullet)}(\cD^\Delta) \oto{U_{\phi(R_\bullet)}}\Mod_{\phi(R_\bullet)}(\cD),
\end{equation}
recalling the last functor from \eqref{eq:U_R}. This $\wt{\phi}$ agrees with \eqref{eq:sym-mon-wt-phi} in case $\phi$ is colimit preserving and symmetric monoidal.

\begin{remark}\label{rem:condition_for_psi_enh_be_pointwise}
Note, however, that if $\phi$ is not colimit-preserving and symmetric monoidal, then $\tilde{\phi}$ may not be computed levelwise; in other words, for $M_\bullet\in \Mod_{R_\bullet}(\cC)$ there is a natural map  $\phi(M_n)\to \tilde{\phi}(M)_n$ which is, in general, not an isomorphism. Unwinding the definition, we see that it is an isomorphism for all $n$ precisely when $[n]\mapsto \phi(M_n)$ satisfies the base change condition:
\begin{equation}\label{eq:condition_for_psi_enh_be_pointwise}
\phi(M_{n})\otimes_{\phi(R_n)} \phi(R_{n+1}) \iso \phi(M_{n}\otimes_{R_n} R_{n+1}) \quad \text{ for all $n$}.
\end{equation}
\end{remark}

\subsection{Motivic Adams approximation}

Recall that we denote by $\Sph \in \SH_{\Qpcyc} \subset \MS_{\Qpcyc}$ the unit of $\SH_{\Qpcyc}$. Let $\Sph \rightarrow \MHFp $ be the tautological map of ring spectra in $\SH_{\Qpcyc}$. Taking the coCech nerve (\Cref{ex: cocech}), we obtain the cosimplicial commutative algebra in $\calg(\SH_{\Qpcyc})$ 
\[
\SphFp := \left( \begin{tikzcd}
    \MHFp
    \ar[r,yshift=1.5pt,->]
    \ar[r,yshift=0pt,<-]
    \ar[r,yshift=-1.5pt,->]
    &
    \MHFp \otimes \MHFp
    \ar[r,yshift=3pt,->]
    \ar[r,yshift=1.5pt,<-]
    \ar[r,yshift=0pt,->]
    \ar[r,yshift=-1.5pt,<-]
    \ar[r,yshift=-3pt,->]
    &
   \MHFp  \otimes \MHFp  \otimes \MHFp 
    \ar[r,yshift=4.5pt,->]
    \ar[r,yshift=3pt,<-]
    \ar[r,yshift=1.5pt,->]
    \ar[r,yshift=0pt,<-]
    \ar[r,yshift=-1.5pt,->]
    \ar[r,yshift=-3pt,<-]
    \ar[r,yshift=-4.5pt,->]
    &
    \cdots
\end{tikzcd} \right)
\]
\begin{remark}
Informally speaking, we think of $\SphFp$ as the ``completion of $\Sph$ along $\Sph \rightarrow \MHFp $''. The necessity of working with $\SphFp $ instead of $\Sph$ comes from the well-known issue that the \emph{motivic Adams spectral sequence does not converge} (unlike its classical counterpart). Instead, we manually replace the motivic sphere spectrum by its Adams spectral sequence, in an appropriate sense.
\end{remark}
\subsubsection{} Applying $\Psi$, we obtain (since $\Psi$ is lax symmetric monoidal) a cosimplicial commutative algebra in $\MS_{\Zpcyc}$, 
\begin{equation}\label{eq: Psi cosimplicial diagram}
    \Psi(\SphFp)  = \left( \begin{tikzcd} \Psi(\MHFp )
    \ar[r,yshift=1.5pt,->]
    \ar[r,yshift=0pt,<-]
    \ar[r,yshift=-1.5pt,->]
    &
    \Psi(\MHFp  \otimes  \MHFp )
    \ar[r,yshift=3pt,->]
    \ar[r,yshift=1.5pt,<-]
    \ar[r,yshift=0pt,->]
    \ar[r,yshift=-1.5pt,<-]
    \ar[r,yshift=-3pt,->]
    &
   \Psi(\MHFp  \otimes  \MHFp  \otimes  \MHFp )
    \ar[r,yshift=4.5pt,->]
    \ar[r,yshift=3pt,<-]
    \ar[r,yshift=1.5pt,->]
    \ar[r,yshift=0pt,<-]
    \ar[r,yshift=-1.5pt,->]
    \ar[r,yshift=-3pt,<-]
    \ar[r,yshift=-4.5pt,->]
    &
    \cdots 
\end{tikzcd}\right)
\end{equation}
and then applying $i^*$ gives the cosimplicial commutative algebra in $\MS_k$, 
\begin{equation}\label{eq: psi cosimplicial diagram}
\psi(\SphFp) = \left(\begin{tikzcd}\psi(\MHFp )
    \ar[r,yshift=1.5pt,->]
    \ar[r,yshift=0pt,<-]
    \ar[r,yshift=-1.5pt,->]
    &
    \psi(\MHFp  \otimes  \MHFp )
    \ar[r,yshift=3pt,->]
    \ar[r,yshift=1.5pt,<-]
    \ar[r,yshift=0pt,->]
    \ar[r,yshift=-1.5pt,<-]
    \ar[r,yshift=-3pt,->]
    &
   \psi(\MHFp  \otimes  \MHFp  \otimes  \MHFp )
    \ar[r,yshift=4.5pt,->]
    \ar[r,yshift=3pt,<-]
    \ar[r,yshift=1.5pt,->]
    \ar[r,yshift=0pt,<-]
    \ar[r,yshift=-1.5pt,->]
    \ar[r,yshift=-3pt,<-]
    \ar[r,yshift=-4.5pt,->]
    &
    \cdots
\end{tikzcd}\right) 
\end{equation}

\subsubsection{} More generally, any $E \in \calg(\SH_{\Qpcyc})$ is an algebra over $\Sph$ in a canonical way, hence can be tensored with $\MHFp$ to produce a cosimplicial commutative algebra in $\SH_{\Qpcyc}$,
\begin{equation}\label{eq: Psi E cosimplicial diagram}
E_p^{\bu} := \left( \begin{tikzcd}
    E \otimes  \MHFp 
    \ar[r,yshift=1.5pt,->]
    \ar[r,yshift=0pt,<-]
    \ar[r,yshift=-1.5pt,->]
    &
    E \otimes  \MHFp  \otimes \MHFp 
    \ar[r,yshift=3pt,->]
    \ar[r,yshift=1.5pt,<-]
    \ar[r,yshift=0pt,->]
    \ar[r,yshift=-1.5pt,<-]
    \ar[r,yshift=-3pt,->]
    &
    E \otimes \MHFp  \otimes  \MHFp  \otimes \MHFp
    \ar[r,yshift=4.5pt,->]
    \ar[r,yshift=3pt,<-]
    \ar[r,yshift=1.5pt,->]
    \ar[r,yshift=0pt,<-]
    \ar[r,yshift=-1.5pt,->]
    \ar[r,yshift=-3pt,<-]
    \ar[r,yshift=-4.5pt,->]
    &
    \cdots
\end{tikzcd}\right)
\end{equation}
We then obtain cosimplicial commutative algebras $\Psi(E_p^{\bu}) \in \MS_{\Zpcyc}$ and $\psi(E_p^{\bu}) \in \MS_{k}$ as above.

\subsection{Spectral syntomic module categories}\label{ssec: spectral module categories}
Consider the module categories in the sense of \Cref{def: simplicial module cat},
\[
\Mod_{\Psi(\SphFp)}(\MS_{\Zpcyc}) \quad \text{and} \quad \Mod_{\psi(\SphFp)}(\MS_{k}).
\]

\subsubsection{Enhanced functors} 

We can now upgrade the functors $\Psi$ and $\psi$ with variants that track the motivic Adams approximations. First, recall that the recipe of \eqref{eq:wt-phi} produces $\wt \Psi \co \Mod_{\SphFp}(\SH_{\Qpcyc}) \rightarrow \Mod_{\Psi(\SphFp)}(\MS_{\Zpcyc})$.  

\begin{defn}
We let $\psienh$ and $\Psienh$ be the compositions
\[
\Psienh \colon \SH_{\Qpcyc} \oto{(-)\otimes \SphFp} \Mod_{\SphFp}(\SH_{\Qpcyc}) \oto{\tilde{\Psi}} \Mod_{\Psi(\SphFp)}(\MS_{\Zpcyc})  
\]
and 
\[
\psienh := i^*\circ \Psienh \colon \SH_K \to \Mod_{\psi(\SphFp)}(\MS_k).
\]
\end{defn}


\subsubsection{Evaluation on motivic cohomology}\label{sssec:dual-motivic-steenrod}

Our next goal is to analyze the object $\Psienh(\MHFp(n)_\Qpcyc)$ of $\Mod_{\Psi(\SphFp)}(\MS_\Zpcyc)$. This will rely on certain facts about the dual motivic Steenrod algebra, which will also be of crucial importance later. Define the set 
\begin{equation}\label{eq: sI}
\sI := \{(r, \epsilon_r, i_r, \ldots, \epsilon_1, i_1, \epsilon_0) \mid r \geq 0, i_j >0, \epsilon_j \in \{0,1\} , i_{j+1} \geq p i_j+ \epsilon_j\}.
\end{equation}
By \cite[Theorem 1.1]{HKO}, which in this characteristic zero situation goes back to Voevodsky \cite{Voe03, Voe10}, we have a decomposition in $\SH_K$, 
\begin{equation}\label{eq: dual steenrod}
\MHFp  \otimes_{\Sph} \MHFp \cong \bigoplus_{\alpha \in \sI}  \MHFp[p_\alpha](q_\alpha) 
\end{equation}
where
\begin{equation}\label{eq: bidegree}
p_\alpha = \sum_{j=0}^r \epsilon_j + \sum_{j=1}^r 2 i_j(p-1) \quad \text{ and } \quad 
q_\alpha = \sum_{j=1}^r i_j(p-1).
\end{equation}

Consequently, for any $n \geq 1$ we learn that $(\MHFp)^{\otimes n+1}$ is a direct sum of shifts of $(\MHFp)^{\otimes n}$ as a module over $(\MHFp)^{\otimes n}$ along any of the $n+1$ maps $(\MHFp)^{\otimes n} \rightarrow (\MHFp)^{\otimes n+1}$ induced by the simplex category. 

As explained in \Cref{rem:condition_for_psi_enh_be_pointwise}, for any $E\in \SH_\Qpcyc$ there is a natural transformation  
\begin{equation}\label{eq:psienh-comparison}
\Psi((-)\otimes  \SphFp) \to \Psienh(-) \co \SH_K \rightarrow  \Mod_{\Psi(\SphFp)}(\MS_\Zpcyc^\Delta),
\end{equation}
but it is not a natural isomorphism. 

\begin{prop} \label{prop:psienh-on-Fp} The natural transformation \eqref{eq:psienh-comparison} evaluates to an isomorphism on the object $\MHFp(n)_\Qpcyc$ for every $n \in \Z$. That is, we have a canonical isomorphism
\begin{equation}\label{eq:psienh-on-Fp}
\Psi(\MHFp(n) \otimes \SphFp)  \iso  \Psienh(\MHFp(n)) \in \Mod_{\Psi(\SphFp)}(\MS_{\Zpcyc}) .
\end{equation}

\begin{proof} Since all functors are compatible with Tate twisting, it suffices to treat the case $n=0$. 

By \Cref{rem:condition_for_psi_enh_be_pointwise}, it would suffice to show that for  $E= \MHFp$ and for each of the coface maps $(\MHFp)^{\otimes i} \to (\MHFp)^{\otimes {i+1}}$, the resulting comparison map 
\[
\Psi\Big(E \otimes (\MHFp)^{\otimes i}\Big) \otimes_{\Psi\big((\MHFp)^{\otimes i}\big)} \Psi\Big((\MHFp)^{\otimes i+1}\Big)  \to \Psi\Big(E\otimes (\MHFp)^{\otimes i+1}\Big) 
\]
is an isomorphism. Using the computation of $\MHFp\otimes_{\Sph} \MHFp$ in \eqref{eq: dual steenrod}, we can express this map in the form 
\[
\Psi\Big(\bigoplus_{\alpha \in \sI} (\MHFp)^{\otimes i}[p_\alpha](q_\alpha) \Big)\otimes_{\Psi\big((\MHFp)^{\otimes i}\big)} \Psi\Big((\MHFp)^{\otimes i+1}\Big) \to \Psi\Big(\bigoplus_{\alpha \in \sI} (\MHFp)^{\otimes i+1}[p_\alpha](q_\alpha)\Big), 
\]
and it fits into a commutative square 
\[
\xymatrix{
\bigoplus_\alpha\Psi\Big((\MHFp)^{\otimes i}\Big)\otimes_{\Psi\big((\MHFp)^{\otimes i}\big)} \Psi\Big((\MHFp)^{\otimes i+1}\Big) [p_\alpha] (q_\alpha) \ar[r]\ar[d]& \bigoplus_\alpha \Psi( (\MHFp)^{\otimes i+1})[p_\alpha](q_\alpha)\ar[d]  \\
\Psi\Big(\bigoplus_\alpha (\MHFp)^{\otimes i}[p_\alpha](q_\alpha) \Big)\otimes_{\Psi\big((\MHFp)^{\otimes i}\big)} \Psi\Big((\MHFp)^{\otimes i+1}\Big) \ar[r] & \Psi\Big(\bigoplus_\alpha (\MHFp)^{\otimes i+1}[p_\alpha](q_\alpha)\Big), 
}
\]
Since the upper horizontal map is clearly an isomorphism, it remains to show that the vertical maps are both isomorphisms. From \eqref{eq: dual steenrod} we see that for every $i \geq 1$, the arguments of $\Psi$ in the bottom row are sums of shifted and twisted copies of $\MHFp$ with shifts that grow to $\infty$. Hence \Cref{lem: psi direct sum} implies that $\Psi$ commutes with the direct sum decompositions in the diagram above, giving the result.

\end{proof}

\end{prop}

The main advantage of $\Psienh(\MHFp)$ over $\Psi(\MHFp) = \MSFp$ is that we can compute its version of the ``dual Steenrod algebra''. More precisely, while it is unclear how to compute $\MSFp \otimes \MSFp$ within $\MS_{\Zpcyc}$ (or even in $\MS_k$), for $\Psienh(\MHFp)$ we have the following Proposition. 

\begin{prop}\label{prop:Psi_enh_on_Fp^2}
In $\Mod_{\Psi(\SphFp)}(\MS_{\Zpcyc})$, we have 
\[
\Psienh(\MHFp) \otimes_{\Psi(\SphFp)} \Psienh(\MHFp) \iso \Psienh(\MHFp \otimes_{\Sph}  \MHFp) \cong \bigoplus_{\alpha \in \sI} \Psienh(\MHFp)[p_\alpha](q_\alpha). 
\]
\end{prop}

\begin{proof}
By \Cref{prop:psienh-on-Fp} we have $\Psienh(\MHFp)\cong \Psi(\MHFp\otimes \SphFp)$. The tensor product in a limit of $\infty$-categories is computed coordinate-wise, so it suffices to show that for every $i \geq 1$, the canonical map 
\[
\Psi\Big(\MHFp \otimes  (\MHFp)^{\otimes i} \Big) \otimes_{\Psi\big((\MHFp)^{\otimes i}\big)} \Psi\Big(\MHFp \otimes  (\MHFp)^{\otimes i} \Big) \to \Psi\Big(\MHFp\otimes  (\MHFp)^{\otimes i} \otimes  \MHFp\Big) 
\]
is an isomorphism. 
Using the computation of $\MHFp \otimes \MHFp$ in \eqref{eq: dual steenrod}, we can rewrite this map as 
\begin{align*}
& \Psi\Big(\bigoplus_{\alpha\in \sI} (\MHFp)^{\otimes i}[p_\alpha](q_\alpha)\Big) \otimes_{\Psi\big((\MHFp)^{\otimes i}\big)} \Psi\Big(\bigoplus_{\alpha \in \sI}(\MHFp)^{\otimes i}[p_\alpha](q_\alpha)\Big) \\
&\hspace{1in} \to \Psi\Big(\bigoplus_{(\alpha,\alpha')\in \sI \times \sI} (\MHFp)^{\otimes i}[p_\alpha + p_{\alpha'}](q_\alpha + q_{\alpha'})\Big).
\end{align*}
From \eqref{eq: dual steenrod} we see that for every $i \geq 1$, the arguments of $\Psi$ involving infinite direct sums can be written as infinite direct sums of $\MHFp$ with shifts that tend to $\infty$. Hence $\Psi$ commutes with the direct sum decompositions by \Cref{lem: psi direct sum}. Then the result follows from the distributivity of the tensor product over direct sums. 
\end{proof}

By applying the colimit-preserving, symmetric monoidal functor $i^*$, and using similar considerations, we obtain similar results for $\psi$ and $\psienh$.

\begin{prop}\label{prop:psi_enh_on_Fp^2} We have
\[
\psi((\MHFp)_\Qpcyc \otimes \SphFp) \cong  \psienh((\MHFp)_\Qpcyc) \in \Mod_{\psi(\SphFp)}(\MS_{k})
\]
and
\[
 \psienh((\MHFp)_\Qpcyc) \otimes_{\psi(\SphFp)} \psienh((\MHFp)_\Qpcyc) \cong \psienh((\MHFp\otimes_{\Sph} \MHFp)_\Qpcyc)  \in \Mod_{\psi(\SphFp)}(\MS_{k}).
\]
\end{prop}

\subsubsection{Module categories over motivic ring spectra} Let $E\in \calg(\MS_\Qpcyc)$. We can then form the cosimplicial commutative algebra $E_p^\bu := E\otimes \SphFp$, and then the module categories $\Mod_{\Psi(E_p^\bu)}(\MS_\Zpcyc)$ and $\Mod_{\psi(E_p^\bu)}(\MS_k)$.
There is a canonical map $E \rightarrow E_p^{\bu}$, so by functoriality $\Psi$ induces a functor 
\begin{equation}\label{eq: Psi E comparison}
\Mod_{\Psi(E)}(\MS_{\Zpcyc})\rightarrow \Mod_{\Psi(E_p^{\bu})}(\MS_{\Zpcyc})
\end{equation}
and similarly $\psi$ induces a functor 
\begin{equation}\label{eq: psi E comparison}
\Mod_{\psi(E)}(\MS_{k}) \rightarrow \Mod_{\psi(E_p^{\bu})}(\MS_{k}).
\end{equation}

\begin{lemma}\label{lem: enhanced Psi category}
If $E \in \SH_{\Qpcyc}$ is an $\MHFp $-algebra, then both functors \eqref{eq: Psi E comparison} and \eqref{eq: psi E comparison} are equivalences. 
\end{lemma}

\begin{proof}
We will prove that \eqref{eq: Psi E comparison} is an equivalence, the case of \eqref{eq: psi E comparison} being completely analogous. 
Since $E$ is an $\MHFp$-algebra, we may tautologically present $E \cong E \otimes_{\MHFp} \MHFp$. This induces an isomorphism
\[
E \otimes_{\Sph} (\MHFp)^{\otimes n} \cong  E \otimes_{\MHFp} (\MHFp)^{\otimes n+1}.
\]
These isomorphisms fit together into an isomorphism of cosimplicial ring spectra
\begin{equation}\label{eq: E extra codegeneracy}
\begin{tikzcd}
    E \otimes  \MHFp
    \ar[r,yshift=1.5pt,->]
    \ar[r,yshift=0pt,<-]
    \ar[r,yshift=-1.5pt,->] \ar[d, "\wr"] 
    &
    E \otimes  (\MHFp)^{\otimes 2}
    \ar[r,yshift=3pt,->]
    \ar[r,yshift=1.5pt,<-]
    \ar[r,yshift=0pt,->]
    \ar[r,yshift=-1.5pt,<-]
    \ar[r,yshift=-3pt,->] \ar[d, "\wr"] 
    &
   E \otimes (\MHFp)^{\otimes 3}
    \ar[r,yshift=4.5pt,->]
    \ar[r,yshift=3pt,<-]
    \ar[r,yshift=1.5pt,->]
    \ar[r,yshift=0pt,<-]
    \ar[r,yshift=-1.5pt,->]
    \ar[r,yshift=-3pt,<-]
    \ar[r,yshift=-4.5pt,->] \ar[d, "\wr"] 
    &
    \cdots \ar[d, "\wr"]  \\
    E \otimes_{\MHFp} (\MHFp)^{\otimes 2}
    \ar[r,yshift=1.5pt,->]
    \ar[r,yshift=0pt,<-]
    \ar[r,yshift=-1.5pt,->] \ar[u]
    &
    E \otimes_{\MHFp} (\MHFp)^{\otimes 3}
    \ar[r,yshift=3pt,->]
    \ar[r,yshift=1.5pt,<-]
    \ar[r,yshift=0pt,->]
    \ar[r,yshift=-1.5pt,<-]
    \ar[r,yshift=-3pt,->] \ar[u]
    &
   E \otimes_{\MHFp} (\MHFp)^{\otimes 4}
    \ar[r,yshift=4.5pt,->]
    \ar[r,yshift=3pt,<-]
    \ar[r,yshift=1.5pt,->]
    \ar[r,yshift=0pt,<-]
    \ar[r,yshift=-1.5pt,->]
    \ar[r,yshift=-3pt,<-]
    \ar[r,yshift=-4.5pt,->] \ar[u]
    &
    \cdots \ar[u]
\end{tikzcd}
\end{equation}
But the bottom cosimplicial diagram has a contraction to $E$ by the extra codegeneracy argument, since it can be prolonged to a diagram 
\[
\begin{tikzcd}
E \ar[r] & 
    E \otimes_{\MHFp} (\MHFp)^{\otimes 2}
    \ar[r,yshift=1.5pt,->]
    \ar[r,yshift=0pt,<-]
    \ar[r,yshift=-1.5pt,->] 
    &
    E \otimes_{\MHFp} (\MHFp)^{\otimes 3}
    \ar[r,yshift=3pt,->]
    \ar[r,yshift=1.5pt,<-]
    \ar[r,yshift=0pt,->]
    \ar[r,yshift=-1.5pt,<-]
    \ar[r,yshift=-3pt,->] 
    &
   E \otimes_{\MHFp} (\MHFp)^{\otimes 4}
    \ar[r,yshift=4.5pt,->]
    \ar[r,yshift=3pt,<-]
    \ar[r,yshift=1.5pt,->]
    \ar[r,yshift=0pt,<-]
    \ar[r,yshift=-1.5pt,->]
    \ar[r,yshift=-3pt,<-]
    \ar[r,yshift=-4.5pt,->] 
    &
    \cdots 
    \end{tikzcd}
    \]
The image of the diagram $E \rightarrow \eqref{eq: E extra codegeneracy}$ under any functor will still have a contraction, so the natural map induces an equivalence
\[
\Mod_{\Psi(E)}(\MS_{\Zpcyc}) \xrightarrow{\sim} \limit_{\eqref{eq: E extra codegeneracy}} \Mod_{\Psi( E  \otimes_{\MHFp}(\MHFp)^{\otimes n})}(\MS_{\Zpcyc}) \cong  \Mod_{\Psi(E_p^{\bu})}(\MS_{\Zpcyc}),
\]
as desired. 
\end{proof}

\begin{cor}\label{cor:module-category}
Regarding $\Psienh(\MHFp)$ as a commutative algebra in $\Mod_{\Psi(\SphFp)}(\MS_\Zpcyc)$, there is a natural equivalence of categories 
\[
\Mod_{\Psienh(\MHFp)}\left(\Mod_{\Psi(\SphFp)}(\MS_\Zpcyc)\right) \cong \Mod_{\MSFp}(\MS_\Zpcyc),
\]
and similarly for $\psi$ and $\MS_k$. 
\end{cor}

\begin{proof}
In general, if $A$ is an $R$-algebra in a symmetric monoidal $\infty$-category $\cC$, then we have a natural equivalence $\Mod_A(\Mod_R(\cC)) \cong \Mod_{A}(\cC)$. Using the isomorphism  $\Psienh(\MHFp) \cong \Psi(\MHFp\otimes \SphFp)$ established in \Cref{prop:psienh-on-Fp}, we can apply this levelwise to $\cC := \MS_\Zpcyc$, $R = \Psi(\SphFp)$ and $A = \Psi(\MHFp \otimes \SphFp)$ to deduce that 
\[
\Mod_{\Psienh(\MHFp)}\left(\Mod_{\Psi(\SphFp)}(\MS_\Zpcyc)\right) \cong \Mod_{\Psi((\MHFp)_p^{\bu})}(\MS_\Zpcyc).
\]
Then \Cref{lem: enhanced Psi category} identifies the latter category with $\Mod_{\Psi(\MHFp)}(\MS_\Zpcyc)$, which is finally identified with $\Mod_{\MSFp}(\MS_\Zpcyc)$ using \Cref{cor:Psi_syntomic}.
\end{proof}

This has the following more concrete consequence.
Recall that objects of $\Mod_{\Psi(\SphFp)}(\MS_{\Zpcyc})$ are, in particular, cosimplicial diagrams in $\MS_{\Zpcyc}$, so we can take their limit to obtain an object of $\MS_{\Zpcyc}$. A similar discussion applies with $k$ in place of $\Zpcyc$, giving lax symmetric monoidal functors 
\begin{equation}\label{eq:simplicial-limit}
\lim_\Delta \colon \Mod_{\Psi(\SphFp)}(\MS_{\Zpcyc}) \to \MS_{\Zpcyc} \hspace{1cm} \text{and} \hspace{1cm} 
\lim_\Delta \colon \Mod_{\Psi(\SphFp)}(\MS_{k}) \to \MS_{k}.
\end{equation}

\begin{cor}\label{cor:psi_enh_on_Fp}
For all $n \in \Z$, there is a canonical isomorphism 
\[
\lim_\Delta \Psienh(\MHFp(n)_\Qpcyc) \cong \MSFp(n)_\Zpcyc,
\]
and similarly 
\[
\lim_\Delta \psienh(\MHFp(n)_\Qpcyc) \cong \MSFp(n)_k.
\]
\end{cor}

\begin{proof}
We prove the claim over $\Zpcyc$, the proof for $k$ being analogous. Using the compatibility of all the constructions involved with Tate twists, it is enough to prove the case $n=0$.

Unwinding the definitions, the functor $\lim_\Delta$ fits into a commutative diagram of lax symmetric monoidal functors 
\[
\xymatrix{
\Mod_{\Psienh(\MHFp)}(\Mod_{\Psi(\SphFp)}(\MS_{\Zpcyc})) \ar^-\sim[r]\ar^{\text{forget}}[d] & \Mod_{\MSFp}(\MS_{\Zpcyc}) \ar^{\text{forget}}[d] \\ 
\Mod_{\Psi(\SphFp)}(\MS_{\Zpcyc}) \ar^{\lim_\Delta}[r]  & \MS_{\Zpcyc} \\ 
}
\]
The upper horizontal functor is the symmetric monoidal equivalence of \Cref{cor:module-category}, hence carries the unit $\Psienh(\MHFp)$ of the upper left category to the unit $(\MSFp)_\Zpcyc$ of the upper right category. Comparing the values of the two paths in the diagram on $\Psienh(\MHFp)$, we obtain the desired identification. 
\end{proof}

\part{Steenrod operations}\label{part: syntomic Steenrod algebra} As mentioned in the Introduction, there are two different flavors of Steenrod operations acting on syntomic cohomology: the \emph{syntomic Steenrod operations}, and the \emph{$\EE_\infty$ Steenrod operations}. In this Part, we define and study these operations and their interaction.  

Henceforth, we make the specific choice $\sO := \Z_p^{\cyc}$, $K := \Q_p^{\cyc}$, and $k := \F_p$. We consider the ``perfectoid nearby cycles'' functor $\psi$ from \S\ref{sec: perfectoid nearby cycles}, using this choice. In \S \ref{sec: syntomic Steenrod operations} we construct the syntomic Steenrod algebra as the Ext algebra of (derived) endomorphisms of $\MSFp$ over the ``syntomic sphere spectrum'' $\psi(\SphFp)$. 

In \S \ref{sec: E_infty operations}, we define the $\EE_\infty$ operations. Their existence is well-known, but we will utilize aspects of their specific construction via the Tate Frobenius, so we take the opportunity to document the foundations. 

Then in \S \ref{sec: comparing operations}, we formulate and prove the Comparison Theorem which determines the precise relationship between the two types of operations. 

\section{Syntomic Steenrod operations}\label{sec: syntomic Steenrod operations}

In this section we define the \emph{syntomic Steenrod algebra $\Asyn^{*,*}$ over $k$} and over $\Zpcyc$. By analyzing the perfectoid nearby cycles functor $\psi$, we construct power operations in $\Asyn^{*,*}$ imported from Voevodsky's work on the motivic Steenrod algebra. We show that these power operations freely generate $\Asyn^{*,*}$ over $\Hsynpt$, and establish Adem relations in \S \ref{sssec: Adem} and a Cartan formula in \S \ref{sssec: Cartan formula}, which together describe the Hopf algebra structure of $\Asyn^{*,*}$ completely. 


\subsection{Syntomic Steenrod algebra}
Let $\MHFp \in \SH_K \subset \MS_K$ be the ($\A^1$-invariant) mod $p$ motivic cohomology spectrum over $K$ (\Cref{defn: motivic cohomology}). The \emph{motivic Steenrod algebra over $K$} is\footnote{Here, and in the rest of the section, the bigrading on the Ext groups is given by the cohomological degree and the weight.} 
\[
\Amot^{*,*} := \Ext_{\SH_K}^{*,*}(\MHFp, \MHFp).
\]
Its structure was studied by Voevodsky, and we will recall his results below.

According to \Cref{cor:psi_enh_on_Fp}, we have a canonical lift  of the syntomic cohomology motivic spectrum $(\MSFp)_{\sO} \in \MS_{\sO}$ to $\Psienh((\MHFp)_K) \in \Mod_{\Psi(\SphFp)}(\MS_{\sO})$ under the map $\Mod_{\Psi(\SphFp)}(\MS_{\sO}) \rightarrow \MS_{\sO}$ given by formation of limits:
\begin{equation}
\begin{tikzcd}
(\MHFp)_K \in \SH_K  \ar[r, "{\Psienh}"] \ar[dr, "\Psi"'] & \Mod_{\Psi(\SphFp)}(\MS_{\sO}) \ar[d, "\lim_\Delta"] \\
&  (\MSFp)_{\sO} \in  \MS_{\sO}
\end{tikzcd} \hspace{1cm}
\begin{tikzcd}
(\MHFp)_K \in \SH_K  \ar[r, "{\psienh}"] \ar[dr, "\psi"'] & \Mod_{\psi(\SphFp)}(\MS_{k}) \ar[d, "\lim_\Delta"] \\
&  (\MSFp)_{k} \in  \MS_{k}
\end{tikzcd}
\end{equation}
where $\lim_{\Delta}$ is as in \eqref{eq:simplicial-limit}. We will abuse notation and write $(\MSFp)_{\sO} \in \Mod_{\Psi(\SphFp)}(\MS_{\sO})$ to denote this lift, relying on context to disambiguate the ambient category. Similarly, we will denote the canonical lift of $(\MSFp)_k \in \MS_k$ to $\Mod_{\psi(\SphFp)}(\MS_{k})$ by $(\MSFp)_k$.

\begin{defn}[Syntomic Steenrod algebra]
We define the \emph{syntomic Steenrod algebra} over $k$ to be 
\[
 \Asyn^{*,*}  = \cA_{\syn, k}^{*,*} := \Ext^{*,*}_{\Mod_{\psi(\SphFp)}(\MS_{k})}((\MSFp)_k, (\MSFp)_k).
\]
Here the bigrading $*,*$ is by cohomological degree and weight, respectively. We are using the abuse of notation discussed just above by denoting $(\MSFp)_k$ for the object $\psienh((\MHFp)_\Qpcyc)$.

Similarly, we define the syntomic Steenrod algebra over $\Zpcyc$ to be 
\[
\cA_{\syn, \Zpcyc}^{*,*} := \Ext^{*,*}_{\Mod_{\Psi(\SphFp)}(\MS_{\Zpcyc})}((\MSFp)_{\Zpcyc}, (\MSFp)_{\Zpcyc}).
\]
\end{defn}

\begin{remark}[Cohomology operations] Since $\MSFp$ represents syntomic cohomology, there is a tautological homomorphism from $\Asyn^{p_\alpha,q_\alpha}$ to natural transformations $\Hsyn^{*,*}(-) \rightarrow \Hsyn^{*+p_\alpha, *+q_\alpha}(-)$ of syntomic cohomology on the category of schemes over $k$, and similarly for $\sO$. 
\end{remark}

By functoriality, $\psienh$ induces a homomorphism
\begin{equation}\label{eq: psi to Steenrod algebra}
\psi \co \Amot^{*,*} \rightarrow  \Asyn^{*,*} 
\end{equation} 
and similarly $\Psienh$ induces a homomorphism
\begin{equation}
\Psi \co \Amot^{*,*} \rightarrow  \cA_{\syn, \Zpcyc}^{*,*}.
\end{equation}

\subsubsection{Voevodsky's power operations} For each $i \in \N$, Voevodsky constructed in \cite{Voe03} a \emph{motivic power operation} in $\SH_{\Qpcyc}$ (and more generally in the stable homotopy category over any characteristic zero field)
\[
\Pm^i \in \Amot^{2i(p-1),i(p-1)}.
\]

We denote by $\beta$ the \emph{Bockstein operation} 
\[
\beta \in \Amot^{1,0}
\]
induced by the exact triangle of motivic cohomology complexes,
\[
\MHFp \rightarrow (\ZZ/p^2)^{\mot} \rightarrow \MHFp,
\]
and 
\[
\Bm^i  = \beta \Pm^i \in \Amot^{2i(p-1)+1,i(p-1)}.
\]

Recall the set \eqref{eq: sI} 
\[
\sI := \{(r, \epsilon_r, i_r, \ldots, \epsilon_1, i_1, \epsilon_0) \mid r \geq 0, i_j >0, \epsilon_j \in \{0,1\} , i_{j+1} \geq p i_j+ \epsilon_j\}.
\]
To each $\alpha \in \sI$, Voevodsky defined the motivic power operation 
\begin{equation}\label{eq: char 0 power operation}
\Pm^\alpha :=  \beta^{\epsilon_r} \Pm^{i_r} \ldots \beta^{\epsilon_1} \Pm^{i_1} \beta^{\epsilon_0}  \in \Amot^{p_\alpha, q_\alpha},
\end{equation}
where $p_\alpha$ and $q_\alpha$ are defined in \eqref{eq: bidegree}. Thus $\Pm^\alpha$ is a cohomology operation of bi-degree $(p_\alpha, q_\alpha)$. 

\subsubsection{Syntomic Steenrod operations} For each $\alpha \in \sI$ we denote by $\Ps^\alpha \in \Asyn^{*,*}$ the image of Voevodsky's operation $\Pm^\alpha$ under the map \eqref{eq: psi to Steenrod algebra}. 

\begin{remark}\label{rem: operations as psi}
Thanks to \Cref{cor:psi_enh_on_Fp}, we know that under the forgetful functor $\Mod_{\psi(\SphFp)}(\MS_{k}) \rightarrow \MS_k$ the operation $\Ps^\alpha$ is simply sent to $\psi(\Pm^\alpha) \in \Ext_{\MS_k}^{p_\alpha, q_\alpha}(\MSFp, \MSFp)$, and a similar remark applies in $\MS_{\Zpcyc}$. 
\end{remark}

\begin{notation}We write
\[
\Bs^i := \beta \circ \Ps^i \in \Asyn^{2i(p-1)+1, i(p-1)}.
\]
If $p=2$, then we denote 
\[
\Sqs^{2i} := \Ps^{i} \in \Asyn^{2i,i} \quad \text{and} \quad \Sqs^{2i+1} := \beta \circ \Ps^i \in \Asyn^{2i+1,i}.
\]
\end{notation}

\subsubsection{Adem relations}\label{sssec: Adem} We may now state the Adem relations for our syntomic Steenrod operations. For $p=2$, the integral variants of the Adem relations, which will not be used in this paper, involve the element
$\tau \in \Hsyn^{0,1}(S) \cong \mu_p$ corresponding to $-1 \in \G_m(S)$ for a scheme $S$ over $\sO$. The proofs also reference $\rho \in \Hsyn^{1,1}(\sO)$, the Kummer image of $-1$, which vanishes in our situation because $\sO$ contains $p$-power roots of $-1$ by assumption. 

\begin{prop}[Adem relations: odd $p$]\label{prop: adem p odd} 
Let $p$ be odd. For $0<a<pb$, we have 
\[
\Ps^a \Ps^b = \sum_{i=0}^{\lfloor a/p\rfloor} (-1)^{a+i} \binom{(p-1)(b-i)-1}{a-pi} \Ps^{a+b-i} \Ps^i \in \Asyn^{*,*},
\]
and for $0<a \leq pb$, we have 
\begin{align*}
\Ps^a \Bs^b & = \sum_{i=0}^{\lfloor a/p\rfloor} (-1)^{a+i} \binom{(p-1)(b-i)}{a-pi}\Bs^{a+b-i} \Ps^i \\
& \hspace{1cm} + \sum_{i=0}^{\lfloor (a-1)/p \rfloor} (-1)^{a+i-1} \binom{(p-1)(b-i)-1}{a-pi-1} \Ps^{a+b-i} \Bs^i \in \Asyn^{*,*}.
\end{align*}

The same relations hold in $\cA_{\syn, \sO}^{*,*}$. 
\end{prop}

\begin{proof}
Since \eqref{eq: psi to Steenrod algebra} is an algebra homomorphism, the relations in $\Asyn^{*,*}$ follow by applying $\psienh$ to the analogous statement in $\SH_K$ for $\Pm^a$ and $\Bm^b$, which is \cite[Theorem 10.3]{Voe03}. For $\cA_{\syn, \sO}^{*,*}$, the same argument applies using $\Psienh$ instead. 
\end{proof}

\begin{prop}[Adem relations: $p=2$]\label{prop: adem p=2}
Let $p=2$ and $0<a<2b$. 
\begin{enumerate}
\item[(i)] If $a \equiv b \equiv 0 \pmod{2}$, then we have 
\[
\Sqs^a \Sqs^b = \sum_{\substack{i=0 \\ \text{$i$ even}}}^{\lfloor a/2\rfloor} \binom{b-i-1}{a-2i} \Sqs^{a+b-i} \Sqs^i \in \Asyn^{*,*}.
\]

\item[(ii)] If $a \equiv 0 \pmod{2}$ and $b \equiv 1 \pmod{2}$, then we have 
\begin{align*}
\Sqs^a \Sqs^b & = \sum_{i=0}^{\lfloor a/2\rfloor}  \binom{b-i-1}{a-2i}\Sqs^{a+b-i} \Sqs^i  \in \Asyn^{*,*}.
\end{align*}

\item[(iii)] If $a \equiv 1 \pmod{2}$ and $b \equiv 0 \pmod{2}$, then we have 
\begin{align*}
\Sqs^a \Sqs^b & = \sum_{\substack{i=0 \\ i \text{ even}}}^{\lfloor a/2 \rfloor} \binom{b-i-1}{a-2i} \Sqs^{a+b-i} \Sqs^i \in \Asyn^{*,*}. 
\end{align*}

\item[(iv)] Finally, if $a \equiv b \equiv 1 \pmod{2}$, then we have 
\begin{align*}
\Sqs^a \Sqs^b & = \sum_{\substack{i=0 \\ i \text{ odd}}}^{\lfloor a/2 \rfloor} \binom{b-i-1}{a-2i} \Sqs^{a+b-i} \Sqs^i \in \Asyn^{*,*}.
\end{align*}
\end{enumerate}
Over $\sO$, (i) is replaced by
\[
\Sqs^a \Sqs^b = \sum_{i=0}^{\lfloor a/2\rfloor} \tau^{i \text{ mod } 2} \binom{b-i-1}{a-2i} \Sqs^{a+b-i} \Sqs^i \in \cA^{*,*}_{\syn, \sO}
\]
while relations (ii), (iii), and (iv) are the same as above. 
\end{prop}

\begin{proof}
The same argument as for Proposition \ref{prop: adem p odd} applies here, using instead \cite[Theorem 10.2]{Voe03} (with typos corrected as in \cite[Theorem 5.1]{HKO}) and noting that $\rho=0 \in \Hsyn^{1,1}(\Spec k)$. 
\end{proof}

\subsection{Freeness of the syntomic Steenrod algebra}\label{ssec:syntomic-Steenrod-basis}
We will show that $\Asyn^{*,*}$ is free over $\Hsynpt$, with basis given by the power operations. 

\subsubsection{Dual basis}\label{sssec:dual-syntomic-steenrod-basis} Recall from \Cref{prop:psi_enh_on_Fp^2} that the natural map
\begin{equation}\label{eq: dual Steenrod 1}
\psienh(\MHFp) \otimes_{\psienh(\Sph)} \psienh(\MHFp)  \rightarrow \psienh(\MHFp \otimes_{\Sph}  \MHFp)\in \Mod_{\psi(\SphFp)}(\MS_{k}),
\end{equation}
coming from the lax symmetric monoidality of $\psienh$, is an isomorphism. In \eqref{eq: dual steenrod}, there is a distinguished generator $\xi_\alpha$ of the $\alpha$ summand $\MHFp[p_\alpha](q_\alpha)$, pinned down by the property that it is dual to the power operation $\Pm^\alpha$ from \eqref{eq: char 0 power operation} in the following sense. Rewrite \eqref{eq: dual steenrod} as
\begin{equation}\label{eq: dual steenrod xi}
\MHFp  \otimes_{\Sph} \MHFp \cong \bigoplus_{\alpha \in \sI} \MHFp  \, \xi_\alpha.
\end{equation}
Working in the category $\SH_{\Qpcyc}$, we have 
\begin{equation}\label{eq: hom-tensor K}
\Ext^{a,b}_{\SH_{\Qpcyc}}(\MHFp, \MHFp) \cong \Ext^{a,b}_{\MHFp}(\MHFp \otimes_{\Sph}  \MHFp, \MHFp).
\end{equation}
From \eqref{eq: dual steenrod xi}, we have
\begin{equation}\label{eq: hom-tensor split}
\Ext^{a,b}_{\MHFp}(\MHFp \otimes_{\Sph} \MHFp, \MHFp) \cong  \prod_{\alpha \in \sI} \Ext^{a,b}_{\MHFp}(\MHFp  \, \xi_\alpha, \MHFp).
\end{equation}
Then the generator $\xi_\alpha$ for $\MHFp \xi_\alpha$ induces an isomorphism 
\[
\MHFp \ \xi_\alpha  \cong \MHFp[p_\alpha](q_\alpha)
\]
which is pinned down by the characterization that the corresponding map $\MHFp \rightarrow  \MHFp[p_\alpha](q_\alpha)$ in the left side of \eqref{eq: hom-tensor K} is the power operation $\Pm^\alpha$.

 Applying Proposition \ref{prop:Psi_enh_on_Fp^2} to \eqref{eq: dual steenrod} and using Lemma \ref{lem: psi direct sum}, we obtain a decomposition 
\begin{equation}\label{eq: tensor split}
\psienh(\MHFp \otimes_{\Sph} \MHFp) \cong \bigoplus_{\alpha \in \sI} \psienh(\MHFp)  \, \xi_\alpha  \cong \bigoplus_{\alpha \in \sI}  \psienh(\MHFp) [p_\alpha](q_\alpha)
\end{equation}
where $\xi_\alpha$ is dual to the syntomic Steenrod operation $\Ps^\alpha  \in \Asyn^{p_\alpha, q_\alpha}$ in the analogous sense. 

\subsubsection{Basis of the syntomic Steenrod algebra}\label{sssec:syntomic-Steenrod-basis} In particular, by \Cref{prop:Psi_enh_on_Fp^2} and the splitting \eqref{eq: tensor split}, we have that 
\begin{align*}
\Asyn^{*,*} &= \Ext^{*,*}_{\Mod_{\psi(\SphFp)}(\MS_{k})}(\psienh(\MHFp), \psienh(\MHFp)) \\
&\cong \Ext^{*,*}_{\psienh(\MHFp)}(\psienh(\MHFp) \otimes_{\psienh(\Sph)}  \psienh(\MHFp), \psienh(\MHFp)) \\ 
& \cong \prod_{\alpha \in \sI} \Ext^{*, *}_{\MSFp}(\MSFp[p_\alpha](q_\alpha), \MSFp)
\end{align*}
where in the last isomorphism we used Lemma \ref{lem: enhanced Psi category} to simplify the simplicial structure, and Corollary \ref{cor: psi Z_p}. Since for each fixed cohomological degree the produce is supported on only finitely many $\alpha$, this graded product agrees with the graded direct sum, hence is free as a graded module over $\Ext^{*,*}_{\MSFp}(\MSFp, \MSFp) =  \Hsyn^{*,*}(\Spec k)$.

\subsection{Dual syntomic Steenrod algebra}\label{ssec: syntomic dual Steenrod algebra} We will now investigate the structure of the ``dual syntomic Steenrod algebra'' over $k$. Everything we say below applies verbatim over $\sO$ instead of $k$, with the same arguments; we will omit the statements to keep the narrative from sprawling.

\begin{defn}\label{defn:dual-sAsyn}
Consider the object
\begin{equation}\label{eq:dual-sAsyn}
\dsAsyn   := \psienh(\MHFp) \otimes_{\psienh(\Sph)} \psienh(\MHFp) \in \Mod_{\psi(\SphFp)}(\MS_k).
\end{equation}

A priori, $\dsAsyn$ naturally has the structure of a commutative Hopf algebroid over $\psienh(\MHFp)$ in $\Mod_{\psi(\SphFp)}(\MS_{k})$ (this is a general pattern in algebra; for example, see \cite{SHA} which was a motivation for our construction here). 
\end{defn}

\begin{lemma}\label{lem: dual Hopf algebroid is algebra}
The Hopf algebroid structure on $\dsAsyn$ is actually a Hopf algebra structure, meaning that the two tautological maps in \[
\Hom_{\Mod_{\psi(\SphFp)}(\MS_{k})}(\psienh(\MHFp), \dsAsyn)
\]
coincide. 
\end{lemma}

\begin{proof}
The Hopf algebroid structure on $\MHFp \otimes_{\Sph} \MHFp \in \SH_{\Qpcyc}$ is described in \cite[\S 5.1, Theorem 5.6]{HKO}. From the formulas there we see that the left and right unit maps in 
\begin{equation}\label{eq: Qpcyc left/right units}
\Hom_{\SH_{\Qpcyc}}(\MHFp,  \MHFp \otimes_{\Sph} \MHFp) 
\end{equation} agree if $p$ is odd, and differ by a multiple of $\rho$ if $p=2$. Thanks to \Cref{prop:Psi_enh_on_Fp^2}, the left and right units maps in 
\begin{equation}\label{eq: k left/right units}
\Hom(\psienh(\MHFp), \psienh(\MHFp) \otimes_{\psienh(\Sph)} \psienh(\MHFp))
\end{equation}
are obtained by applying $\psienh$ to those in \eqref{eq: Qpcyc left/right units}. Hence their difference immediately vanishes if $p$ is odd, and if $p=2$ it also vanishes because $\rho = 0$; indeed, we have $\rho \in \Hsyn^{1,1}(k) = 0$.   
\end{proof}

This allows us to view $\dsAsyn$ unambiguously as an $\psienh(\MHFp)$-algebra. The image of $\dsAsyn$ under the functor $\lim_\Delta$ from \Cref{cor:psi_enh_on_Fp} becomes an $\MSFp$-algebra in $\Mod_{\MSFp}(\MS_k)$, which we also denote by $\dsAsyn$. The two versions of $\dsAsyn$ are identified under \Cref{cor:module-category}, so this abuse of notation should cause little confusion. 

\begin{defn}[Dual syntomic Steenrod algebra]
We define the \emph{dual syntomic Steenrod algebra} to be 
\[
\dAsyn_{*,*} := \Ext^{-*,-*}_{\Mod_{\MSFp}(\MS_k)}(\MSFp, \dsAsyn).
\]
\end{defn}

\begin{lemma}
As an $\Hsynpt$-module, $\dAsyn_{*,*}$ is free and reflexive. 
\end{lemma}

\begin{proof}
According to \eqref{eq: tensor split} and \Cref{prop:psi_enh_on_Fp^2}, we have an isomorphism 
\begin{equation}\label{eq:dsAsyn-splitting}
\dsAsyn  \cong \bigoplus_{\alpha \in \sI} \MSFp \xi_\alpha  \cong \bigoplus_{\alpha \in \sI} \MSFp [p_\alpha] (q_\alpha).
\end{equation}
This exhibits $\dsAsyn$ as a sum of shifts and twists of $\MSFp$, such that for any given $N$, there are only finitely many $\alpha$ such that $p_\alpha < N$. In particular, $\dAsyn_{*,*}$ is finite-dimensional in any degree. In fact, this shows that $\dsAsyn$ is even free and reflexive over $\MSFp$. The freeness and reflexivity of $\dAsyn_{*,*} $ over $\Hsynpt$ follows immediately. 

\end{proof}

By base change, we have (as already used in \S \ref{sssec:syntomic-Steenrod-basis}) 
\[
\hom_{\Mod_{\MSFp}(\MS_k)}(\dsAsyn, \MSFp) \cong \hom_{\Mod_{\psi(\SphFp)}(\MS_{k})}(\MSFp, \MSFp).
\]
Thanks to the freeness and reflexivity of $\dsAsyn$ over $\MSFp$ from \eqref{eq:dsAsyn-splitting}, we obtain an isomorphism 
\begin{equation}\label{eq:Asyn-dual-of-dAsyn}
\Hom_{\Hsynpt}(\dAsyn_{*,*}, \Hsynpt) \cong \Asyn^{-*,-*},
\end{equation}
and double-dualizing (and using reflexivity) identifies
\begin{equation}
\dAsyn_{*,*} \cong  \Hom_{\Hsynpt}(\Asyn^{-*,-*}, \Hsynpt).
\end{equation}
Furthermore, it is clear from tracking definitions that this is compatible with the respective Hopf algebra structures. In summary, we have established the following.

\begin{cor}\label{cor:dual-syntomic-steenrod}
The dual syntomic Steenrod algebra $\dAsyn_{*,*}$ is a commutative graded Hopf algebra over $\Hsynpt$, which is free and reflexive as a graded $\Hsynpt$-module. Its graded dual is identified as a graded Hopf algebra over $\Hsynpt$ with the syntomic Steenrod algebra $\Asyn^{*,*}$, which is also free and reflexive as a graded $\Hsynpt$-module. 
\end{cor}


\subsection{Cartan formula}\label{sssec: Cartan formula} We are now able to identify the coproduct on the syntomic Steenrod algebra explicitly.

\begin{prop}\label{prop: Cartan} Assume $p \neq 2$. Then the comultiplication $\Delta \co \Asyn^{*,*} \rightarrow \Asyn^{*,*} \otimes_{\Hsynpt} \Asyn^{*,*}$ satisfies 
\[
\Delta(\Ps^i) = \sum_{j=0}^i \Ps^j \otimes \Ps^{i-j}
\]
and 
\[
\Delta(\Bs^i) = \sum_{j=0}^i (\Bs^j \otimes \Ps^{i-j} +  \Ps^j \otimes \Bs^{i-j}).
\]
The same formulas hold for $\Delta \co \cA^{*,*}_{\syn, \sO} \rightarrow \cA^{*,*}_{\syn, \sO} \otimes_{\Hsyn^{*,*}(\sO)} \cA^{*,*}_{\syn, \sO}$.

Assume $p=2$. Then 
\[
\Delta(\Sqs^{2i}) = \sum_{j=0}^i \Sqs^{2j} \otimes \Sqs^{2i-2j}  \in \Asyn^{*,*} \otimes_{\Hsynpt}  \Asyn^{*,*},
\]
while over $\sO$ we have instead
\[
\Delta(\Sqs^{2i})  = \sum_{j=0}^i \Sqs^{2j} \otimes \Sqs^{2i-2j}  + \tau \sum_{j=0}^{i-1} (\Sqs^{2j+1}) \otimes (\Sqs^{2i-2j-1}) \in  \cA^{*,*}_{\syn, \sO} \otimes_{\Hsyn^{*,*}(\sO)} \cA^{*,*}_{\syn, \sO}. 
\]
The formula 
\[
\Delta (\Sqs^{2i+1}) = \sum_{j=0}^i \left( \Sqs^{2j+1} \otimes \Sqs^{2i-2j} + \Sqs^{2j} \otimes  \Sqs^{2i-2j+1}\right)
\]
holds for both $\Asyn^{*,*}$ and $\cA^{*,*}_{\syn, \sO}$.
\end{prop}

\begin{proof} 
The comultiplication on $\Asyn^{*,*}$ is induced by the multiplication on $\dsAsyn$. By Proposition \ref{prop:Psi_enh_on_Fp^2}, the vertical maps are isomorphisms in the natural diagram   
\[
\begin{tikzcd}\dsAsyn  \otimes_{\psienh(\MHFp)}  \dsAsyn \ar[r] \ar[d, "\sim"]  & \dsAsyn \ar[d, "\sim"]  \\
\psienh(\MHFp \otimes  \MHFp) \otimes_{\psienh(\MHFp)} \psienh(\MHFp \otimes  \MHFp)   \ar[r] &  \psienh(\MHFp \otimes  \MHFp)
\end{tikzcd}
\]
where the horizontal arrows are the multiplication maps. Invoking Proposition \ref{prop:Psi_enh_on_Fp^2} again,
we have an identification 
\begin{align*}
\psienh(\MHFp \otimes  \MHFp) \otimes_{\psienh(\MHFp)} \psienh(\MHFp \otimes  \MHFp)   & \cong  \psienh(\MHFp \otimes  \MHFp \otimes \MHFp) \\
& \cong \psienh((\MHFp \otimes \MHFp) \otimes_{\MHFp} (\MHFp \otimes \MHFp) ),
\end{align*}
Moreover, this fits into the natural commutative diagram
\begin{equation}\label{eq: dual steenrod multiplication 2}
\adjustbox{scale = 0.9, center}{\begin{tikzcd}
\psienh(\MHFp \otimes   \MHFp) \otimes_{\psienh(\MHFp)} \psienh(\MHFp \otimes  \MHFp)  \ar[r] \ar[d, "\wr"] &  \psienh(\MHFp \otimes  \MHFp) \ar[d, equals] \\
\psienh\left((\MHFp \otimes  \MHFp) \otimes_{\MHFp} (\MHFp \otimes \MHFp) \right) \ar[r] &  \psienh(\MHFp \otimes  \MHFp)  
\end{tikzcd}}
\end{equation}
where the top horizontal arrow is the multiplication map for the commutative algebra $\psienh(\MHFp \otimes  \MHFp)$, and the bottom horizontal arrow is $\psienh$ applied to the multiplication map for Voevodsky's dual motivic Steenrod algebra $\MHFp \otimes  \MHFp \in \calg(\SH_{\Qpcyc})$,
\begin{equation}\label{eq: VV dual product}
(\MHFp \otimes \MHFp) \otimes_{\MHFp} (\MHFp \otimes  \MHFp) \rightarrow (\MHFp \otimes  \MHFp)  \in \SH_{\Qpcyc}. 
\end{equation}
The latter map \eqref{eq: VV dual product} is identified explicitly by \cite[Proposition 9.7]{Voe03}, and is given by the asserted formula if $p$ is odd. If $p=2$, 
then \cite[Proposition 9.7]{Voe03} says that
\[
\Delta(\Sqm^{2i}) = \sum_{j=0}^i \Sqm^{2j} \otimes \Sqm^{2i-2j} + \tau \sum_{j=0}^{i-1} (\Sqm^{2j+1}) \otimes (\Sqm^{2i-2j-1})
\]
and 
\[
\Delta (\Sqm^{2i+1}) = \sum_{j=0}^i \left( \Sqm^{2j+1} \otimes \Sqm^{2i-2j} + \Sqm^{2j} \otimes  \Sqm^{2i-2j+1}\right) + \rho \sum_{j=0}^{i-1} \Sqm^{2j+1} \otimes \Sqm^{2i-2j-1}.
\]
Applying the functor $\Psienh$ to these relations, and noting that $\rho = 0$, we obtain the asserted formulas for $\cA^{*,*}_{\syn, \sO}$. Instead applying $\psienh$ to the above relations, and further noting that $\tau = 0$ in $\Hsyn^{*,*}(\Spec k)$, we obtain the asserted formulas for $\Asyn^{*,*}$. 
\end{proof}

\subsection{Structure of the syntomic Steenrod algebra}

By Lemma \ref{lem: dual Hopf algebroid is algebra}, $\Asyn^{*,*}$ is a Hopf algebra over $\Hsynpt$. We note that this latter ring is explicit and simple: $\Hsyn^*(\Spec k; \F_p(i)) = 0$ for $i \neq 0$, so we have 
\begin{equation}\label{eq: Hsynpt}
\Hsynpt \cong \Het^*(\Spec k; \F_p) \cong \F_p[\epsilon]/\epsilon^2
\end{equation}
where $\epsilon$ is a generator of $\Het^1(\Spec k; \F_p) \cong \F_p$. Note that $\Ps^i$ and $\Bs^i$ kill $\Hsynpt$ for degree and weight reasons whenever $i>0$, so $\Asyn^{*,*}$ acts trivially on $\Hsynpt$. 

We now summarize this section's results on the structure of the syntomic Steenrod algebra. 

\begin{thm}\label{thm: syntomic structure}
The syntomic Steenrod algebra $\Asyn^{*,*}$ is a free and reflexive cocommutative Hopf algebra over $\Hsynpt$, with a $\Hsynpt$-basis consisting of the $\Ps^\alpha$ for $\alpha \in \sI$. The algebra structure is then determined by the Adem relations (\S \ref{sssec: Adem}), and coproduct is determined by the Cartan formula (\S \ref{sssec: Cartan formula}). 
\end{thm}

Notice that the product and coproduct of power operations $\Ps^i$ and $\Bs^i \in \Asyn^{*,*}$ do not involve the element $\epsilon$ from \eqref{eq: Hsynpt}, and as remarked above act trivially on $\epsilon$. This allows us to ``descend'' $\Asyn^{*,*}$ from $\Hsynpt$ to $\F_p$. 

\begin{defn} We define the \emph{reduced syntomic Steenrod algebra} $\rAsyn^{*,*} \subset \Asyn^{*,*}$ to be the $\F_p$-subalgebra generated by all the power operations $\Ps^i, \Bs^i \in \Asyn^{*,*}$ for $i \in \Z_{\geq 0}$. We define the \emph{reduced dual syntomic Steenrod algebra} $\drAsyn_{*,*} \subset \dAsyn_{*,*}$ analogously.  
\end{defn}

Then $\rAsyn^{*,*}$ is a cocommutative Hopf algebra over $\F_p$, equipped with a natural isomorphism  
\[
\rAsyn^{*,*} \otimes_{\F_p} \Hsynpt \cong \Asyn^{*,*}
\]
of Hopf algebras over $\Hsynpt$. From Theorem \ref{thm: syntomic structure}, we see that
\begin{itemize}
\item An $\F_p$-basis of $\rAsyn^{*,*}$ is given by $\Ps^\alpha$ for $\alpha \in \sI$. 
\item The algebra structure is given by the Adem relations. 
\item The coproduct is given by the Cartan formula. 
\end{itemize} 
Naturally, we also have the dual statements for $\drAsyn_{*,*}$.

\section{$\EE_{\infty}$ Steenrod operations}\label{sec: E_infty operations}
The cohomology ring of any $\EE_\infty$-$\F_p$-algebra is equipped with power operations, an observation that goes back at least to May \cite{May70}. Although this is well-known, it is treated from different perspectives and in different languages in the literature, so we take this section to set up definitive foundations for use in this paper. Our presentation is guided by later considerations, e.g., we need to set up the $\EE_\infty$ power operations on syntomic cohomology as morphisms of \emph{spectra} (which recover the familiar operations upon taking cohomology). Our formulation is based on the Tate Frobenius; we learned this perspective on $\EE_\infty$ Steenrod operations from Lurie \cite[\S 2]{DAGXIII} and Nikolaus--Scholze \cite[\S IV.1]{NS18}. 


\subsection{The Tate-valued Frobenius} Let $C_p \cong \Z/p$ be the cyclic group of order $p$ together with a chosen generator. Recall that for a semiadditive $\infty$-category $\cC$ with all limits and colimits and $X\in \cC^{\rB C_p}$ (in other words, $X$ is an object of $\cC$ with an action of $C_p$), we have a canonical ``norm map'' 
\[
\Nm\colon X_{hC_p} \to X^{hC_p}
\]
from the homotopy orbits to the homotopy fixed points of the action. The cofiber of this map is, by definition, the \emph{Tate construction} 
\[
X^{tC_p} := \mathrm{Cofib}(X_{hC_p} \to X^{hC_p}). 
\]
In \cite[Definition IV.1.1]{NS18}, the authors construct a natural ``Frobenius'' map 
\begin{equation}\label{eq: spectral Tate Frob}
\Fr\colon R \to R^{tC_p}
\end{equation}
(where the Tate construction is with respect to the trivial action) for every $\EE_\infty$-algebra $R$ in $\Sp$, equipped with the trivial action of $C_p$. This construction can be easily generalized from spectra to sheaves of spectra, as follows. 

\begin{defn}[Tate-valued Frobenius map]
Let $\cC$ be a site. For a sheaf of commutative ring spectra $\cR\in \calg(\Shv(\cC;\Sp))$, equipped with the trivial action of $C_p$, we define the \emph{Tate-valued Frobenius map}
\[
\Fr\colon \cR\to \cR^{tC_p}
\]
to be the sheafification of the levelwise Frobenius map 
\[
U\in \cC \mapsto (\Fr\colon \cR(U) \to \cR(U)^{tC_p})
\]
from \eqref{eq: spectral Tate Frob}. 
\end{defn}

----------------------
\shachar{revised version here}

\subsection{$\EE_\infty$ power operations} \label{sec:E_infty_ops_def} 

We first spell out some canonical (co)homology classes of $\rB C_p$. 
There is a canonical generator $x \in \rH^1(\rB C_p;\F_p)  \cong \Hom(C_p, \F_p)$ corresponding to the fixed isomorphism $C_p \cong \F_p$. 
We can interpret $x\in \rH^1(\rB C_p;\F_p) \cong \Hom(\rB C_p, \F_p[1])$ as an $\F_p$-module map in spectra
\begin{equation}\label{eq:c_i}
    x \colon \F_p \otimes \rB C_p \to \F_p[1]. 
\end{equation}
If $p=2$, then $\rH^*(\rB C_p;\F_p)$ is a polynomial ring on $x$. If $p \neq 2$, then the Bockstein map $\beta$ identifies $\rH^1(\rB C_p; \F_p) \xrightarrow{\sim} \rH^2(\rB C_p, \F_p)$, and $\rH^*(\rB C_p;\F_p)$ is the graded symmetric algebra on $\rH^1(\rB C_p; \F_p) \oplus \rH^2(\rB C_p; \F_p)$. 
In both cases, let $t:= \beta(x)$ be the resulting generator of $\rH^2(BC_p,\F_p)$ (so that when $p=2$ we have $t = x^2$). 
Thus, the monomials $t^n$ and $xt^n$ for $n \ge 0$ form a graded basis of $\rH^*(BC_p;\F_p)$ and we let $\{c_j\}_{j\ge 0}$ be a dual basis of $\rH_*(BC_p;\F_p)$, with $c_j$ dual to the generator of degree $j$.   

By abuse of notation, we denote again by
$t$ the image of $t$ Tate cohomology $t\in \rH^2_{\mathrm{Tate}}(C_p;\F_p) := \rH^2(\F_p^{tC_p})$. 
Crucially, note that in the Tate cohomology the element $t$ becomes invertible.

Let $\cX$ be a sheaf of $\F_p$-modules. Then, since the Tate construction is lax symmetric monoidal \cite[Theorem I.3.1]{NS18}, the sheaf $\cX^{tC_p}$ promotes to a sheaf of modules over the ring spectrum $\F_p^{tC_p}$. Hence, the powers of $t$ define isomorphisms 
\[
t^i\colon \cX^{tC_p} \iso \cX^{tC_p}[2i]
\]

\begin{defn}
Let $\cC$ be a site and let $\cR\in \calg(\Shv(\cC;\F_p))$, equipped with the trivial $C_p$-action. For $i\in \Z$ we let $\Pe^i \colon \cR \to \cR[2i(p-1)]$ be the composition 
\[
\Pe^i \co \cR\oto{\Fr} \cR^{tC_p} \oto{t^{i(p-1)-1}} \cR^{tC_p}[2i(p-1)-2] \to \cR_{hC_p}[2i(p-1)-1] \oto{x} \cR[2i(p-1)] 
\]
where the middle map is the boundary map of the cofiber sequence 
\[
\cR_{hC_p} \oto{\Nm} \cR^{hC_p} \to \cR^{tC_p}.
\]
\end{defn}

\begin{remark}
In fact, using similar considerations, one can show that the norm map $\cR_{hC_p}\to \cR^{hC_p}$ is the zero map, so that we have a direct sum decomposition $\cR^{tC_p}\cong \cR^{hC_p} \oplus \cR_{hC_p}[1]$. Via this decomposition, we can re-write the appropriate shift of the map $\Pe^i$ as 
\[
\cR[-2i(p-1)] \oto{\Frob} \cR^{tC_p}[-2i(p-1)] \oto{t^{i(p-1)}} \cR^{tC_p} \to \cR^{hC_p} \to \cR
\]
where the last map is the natural inclusion of the fixed points and the one before is the projection of the Tate construction to the fixed points from above. While this description will not be used, it is perhaps more familiar. 
\end{remark}

Thus on cohomology groups, $\Pe^i$ induces a map 
\[
\Pe^i \co \rH^{a}(U;\cR) \to \rH^{a}(U;\cR[2i(p-1)]) \cong \rH^{a+2i(p-1)}(U;\cR).
\]

\begin{remark} \label{rem:power_op_no_sheafification}
    Note that the construction of $\Pe^i$ does not involve sheafification. In other words, the map $\Pe^i$ is the restriction of the corresponding map on presheaves (i.e., for $\cC$ with the indiscrete Grothendieck topology). For presheaves, the construction is ``pointwise'' in the sense that it is given by applying the ordinary power operation $\Pe^i$ for $\EE_\infty$-$\F_p$-algebra in spectra at each object:
    \[
    \Pe^i(U) \co 
    \cR(U) \oto{\Fr} \cR(U)^{tC_p} \oto{t^{i(p-1)-1}} \cR(U)^{tC_p}[2i(p-1)-2] \to \cR(U)_{hC_p}[2i(p-1)-1] \oto{x} \cR(U)[2i(p-1)]. 
    \]
\end{remark}

\begin{remark}
Assume that $\cR\cong \bigoplus_{n\in \Z} \cR_n$ is a $\Z$-graded algebra (with trivial action of $C_p$). Then the Frobenius map becomes a graded map after rescaling the grading of $\cR$ by $p$. Indeed, we may write $\Fr$ as a composition of the Tate diagonal
\[
\cR\to (\cR^{\otimes p})^{tC_p}
\]
followed by the multiplication map 
\[
(\cR^{\otimes p})^{tC_p} \to \cR^{tC_p}.
\]
The second map clearly respects the grading, while the first map has a graded refinement provided by \cite[Example A.10 and preceding discussion]{AMMN22} if the grading on the source is scaled by a factor of $p$. 

In particular, the map $\Pe^i$ carries the graded piece $\cR_n$ to the graded piece $\cR_{pn}$, hence decomposes into homogeneous pieces of the form 
\[
\Pe^i\colon \cR_n \to \cR_{pn}[2i(p-1)].
\]
In particular, using the notation $\rH^{a,b}(U; \cR):= \rH^a(U; \cR_b)$, we obtain on cohomology a map of signature 
\[
\Pe^i\colon \rH^{a,b}(U;\cR) \to \rH^{a + 2i(p-1),pb}(U;\cR) \quad \text{for all }  U\in \cC.
\]

\end{remark}

\begin{defn} Let $S$ be a qcqs scheme. Applying this construction to the graded algebra 
\[
\SFp(\bu)_S = \bigoplus_{n \in \Z} \SFp(n)_S \quad \text{in} \quad \PShv_{\Nis}(\Sm_S;\Sp),
\]
we obtain the \emph{$\EE_\infty$-power operations}
\begin{equation}\label{eq:E_infty_operations_syntomic}
\Pe^i\colon \Hsyn^{a,b}(X) \to \Hsyn^{a+2i(p-1),pb}(X)
\end{equation}
for any scheme $X/S$. 
\end{defn}

\subsection{Functoriality along geometric  morphisms}

Let $\phi\colon \cC \to \cD$ be a morphism of sites, so that we obtain a symmetric monoidal adjunction between their categories of $\Sp$-valued sheaves,
\[
\phi^*\co \Shv(\cD;\Sp) \adj \Shv(\cC;\Sp) : \phi_*. 
\] 
We will show that, in this situation, the constructions $\Pe^i$ are compatible with the functors $\phi^*$ and $\phi_*$. To avoid notational ambiguity, in this subsection we shall temporarily denote the map $\Pe^i\colon \cR\to \cR[2i(p-1)]$ by $\Pe^i(\cR)$. (Our results will imply that we never need to worry about this again, as all power operations of this shape are compatible in the obvious way.)

\begin{prop} \label{prop:power_op_push}
Let $\phi \colon \cC \to \cD$ be a morphism of sites, and let $\cR\in \calg(\Shv(\cC;\F_p))$. Then the two maps 
\[
\Pe^i(\phi_*\cR)\colon \phi_*\cR \to \phi_*\cR[2i(p-1)] 
\]
and 
\[
\phi_*\Pe^i(\cR)\colon \phi_*\cR \to \phi_*\cR[2i(p-1)] 
\]
are naturally homotopic. 

Similarly, if $\cR' \in \calg(\Shv(\cD;\F_p))$ then the two maps 
\[
\phi^*\Pe^i(\cR')\colon \phi^*\cR' \to \phi^*\cR'[2i(p-1)] 
\]
and 
\[
\Pe^i(\phi^*\cR') \colon \phi^*\cR' \to \phi^*\cR'[2i(p-1)]  
\]
are naturally homotopic.
\end{prop}

\begin{proof}
We start by showing the claim for $\phi_*$. Using \Cref{rem:power_op_no_sheafification}, we reduce to the case of presheaf categories, where $\phi_*$ is given by pre-composition with $\phi$. Then the claim follows immediately from the fact that $\Pe^i$ is computed pointwise, as discussed in \Cref{rem:power_op_no_sheafification}. 

We turn to the claim regarding $\phi^*$. For this, abbreviate $m := 2i(p-1)$ and consider the following diagram
\[
\begin{tikzcd}[column sep = huge]
\phi^*\mathcal{R}' \arrow[r, "\phi^* \mathrm{coev}"] 
  \arrow[d, "\phi^*\Pe^i(\mathcal{R}')"'] 
  \arrow[rrr, "\mathrm{Id}", bend left=15] 
& \phi^*\phi_*\phi^*\mathcal{R}' 
  \arrow[d, "\phi^*\Pe^i(\phi_* \phi^* \mathcal{R}')" ] 
  \arrow[r, equals] 
& \phi^*\phi_*\phi^*\mathcal{R}' 
  \arrow[d, "\phi^*\phi_* \Pe^i(\phi^* \mathcal{R}')"] 
  \arrow[r, "\mathrm{ev} \phi^*"] 
& \phi^*\mathcal{R}' 
  \arrow[d, "\Pe^i(\phi^* \mathcal{R}')"] \\
\phi^*\mathcal{R}'[m] 
  \arrow[r, "\phi^* \mathrm{coev}"] 
  \arrow[rrr, "\mathrm{Id}", bend right=15] 
& \phi^*\phi_*\phi^*\mathcal{R}'[m] 
  \arrow[r, equals] 
& \phi^*\phi_*\phi^*\mathcal{R}'[m] 
  \arrow[r, "\mathrm{ev} \phi^*"] 
& \phi^*\mathcal{R}'[m]
\end{tikzcd}
\]
where $\coev$ and $\ev$ are the unit and counit of the adjunction $(\phi^*, \phi_*)$. We claim that it commutes:
\begin{itemize}
\item The left square commutes by the naturality of $\Pe^i$. 
\item The right square commutes by the naturality of the counit. 
\item The upper and lower regions commute by the triangle identities.
\item The middle square commutes by the claim for $\phi_*$.
\end{itemize}
The claim for $\phi^*$ now follows by comparing the two outer circuits from $\phi^*\cR'$ to $\phi^*\cR'[m]$ along the boundary of the diagram. 
\end{proof}

\subsection{$\EE_\infty$-operations for syntomic cohomology}

Let $S$ be a qcqs scheme. The formalism of power operations of sheaves of $\F_p$-algebras applies in particular to the syntomic cohomology $\SFp(\bu)_S$, giving maps
\[
\Pe^i\colon \SFp(n)_S \to \SFp(pn)_S[2i(p-1)].  
\]

Now specialize this discussion to the setup of \S \ref{sec: perfectoid nearby cycles}, with a perfectoid valuation ring $\Zpcyc$ with generic fiber $\Qpcyc$ and special fiber $k$. From Corollary \ref{cor:Psi_syntomic} and Corollary \ref{cor: psi Z_p}, we have identifications
\begin{equation}\label{eq: Psi psi Z_p}
\SFp(\bu)_\Zpcyc \cong L_{\et}j_*\HFp(\bu)_{\Qpcyc} \text{ and } \SFp(\bu)_{k} \cong i^*\SFp(\bu)_\Zpcyc.
\end{equation}

\begin{cor}\label{cor:E_infty_ops_psi}
Via the identifications in \eqref{eq: Psi psi Z_p}, the maps 
\[
\Pe^i \colon \SFp(n)_\Zpcyc \to \SFp(pn )_\Zpcyc[2i(p-1)] 
\]
are the images of the corresponding maps for $\HFp(n)_\Qpcyc$ under the functor $L_{\et} j_*$. Similarly, the maps 
\[
\Pe^i\colon \SFp(n)_k \to \SFp(pn)_k[2i(p-1)] 
\]
are the images of the corresponding maps for $\SFp(n)_\Zpcyc$ under the functor $i^*$.
\end{cor}

\begin{proof}
Both the functors $L_{\et} j_*$ and $i^*$ are compositions of push-forward and pullback maps along morphisms of sites. Hence, the results follow from \Cref{prop:power_op_push}. 
\end{proof}

As another application, we can now compute the power operations on the commutative $\F_p$-algebra  $\RGamma_{\syn}(X;\F_p(*))$ in terms of the ``local operations'' on the level of sheaves. 

\begin{prop}\label{prop:E_infty_op_global_sections}
Let $S$ be a qcqs scheme, so that we have a graded $\EE_\infty$-$\F_p$-algebra $\SFp(\bu)(S)$. The resulting maps of $\F_p$-spectra
\[
\Pe^i\colon \SFp(n)(S) \to \SFp(pn)(S)[2i(p-1)] 
\]
agree with the maps induced from the morphism of sheaves of spectra  
\[
\Pe^i \colon \SFp(n)_S \to \SFp(pn)_S[2i(p-1)] 
\]
by taking global sections over $S$. 
\end{prop}

\begin{proof}
This follows immediately from \Cref{prop:power_op_push} using the morphism of sites opposite to the inclusion 
$\{S\}\into \Sm_S$. 
\end{proof}

\subsection{Comparison with unstable construction}
The main advantage of the Tate-cohomology perspective on power operations is that it is done entirely within the framework of sheaves of spectra and ``linear'' maps between them. However, more classical constructions (e.g., \cite{May70}) involve non-linear maps, and linearity is proven a posteriori. For later computational purposes, it will be useful for us to reformulate the operations $\Pe^i$ also in this language.

\begin{defn}
Let $\cD$ be a presentably symmetric monoidal $\infty$-category. For $M\in \cD$, we denote
\[
D^p(M):= (M^{\otimes p})_{hC_p},
\] 
where $C_p$ acts on $M^{\otimes p}$ by cyclically permuting the tensor factors. 
\end{defn}

Recall that part of the data of an $\EE_\infty$-algebra $A$ is a $p$-fold multiplication map $D^p(A) \oto{\mrm{mult}} A$.  

\begin{defn}\label{def:total_power}
Let $\cC$ be a presentably symmetric monoidal $\infty$-category.
Given $M\in \cC$, $A\in \calg(\cC)$, and a map $\alpha \colon M\to A$ in $\cC$, we denote by $\PP^p(\alpha)\colon D^p(M)\to A$ the composition 
\[
\PP^p(\alpha)\colon D^p(M)\to D^p(A) \oto{\mrm{mult}} A. 
\]
We refer to $\PP^p(\alpha)$ as the \emph{total $p$-th power} of $\alpha$. 
\end{defn}

When $\cC$ as in \Cref{def:total_power} is linear over $\F_p$, we can now give a definition of power operations. We will rely on the theory of Picard spectra as explained, for example, in \cite[\S 3]{carmeli2023strict}. In particular, to $\cC$ we can associate a spectrum called its \emph{Picard spectrum} $\pic(\cC)$. Recall that an object $L\in \pic(\cC)$ is called \emph{strict} if it extends to a map of spectra $\Z \to \pic(\cC)$, or equivalently, to a symmetric monoidal exact functor $\Perf^{\gr}(\mathbb{S}) \to \cC$ from the category of graded finite spectra. In this case, and using the $\F_p$-linearity of $\cC$, we claim that the action of $C_p$ on $L[a]^{\otimes p}$ is trivial, so that $D^p(L[a])\cong L^{\otimes p} \otimes \rB C_p[pa]$. For the case $a=0$ this follows (for example) from \cite[Proposition 3.15]{carmeli2023strict}; for general $a$ it additionally uses the fact that, as a complex with action of the symmetric group $\Sigma_p$ we have $\F_p[a]^{\otimes p} \cong \F_p[pa]\otimes \rho^{\otimes a}$ where $\rho$ is the sign representation of $\Sigma_p$.

\begin{defn}\label{defn:norm-E-infty-ops}
Let $\cC$ be an $\F_p$-linear presentably symmetric monoidal $\infty$-category, and let $L\in \pic(\cC)$ be strict. Let $A\in \calg(\cC)$. We define the map  
\begin{equation}\label{eq:norm-E-infty-ops}
\Pe^i\colon \Map(L^{\otimes b} [-a],A) \to \Map(L^{\otimes {pb}} [-a-2i(p-1)],A)  
\end{equation}
to be the composite 
\begin{align*}
\Map(L^{\otimes {b}} [-a],A) \oto{D^p} 
&\Map(D^p(L^{\otimes b} [-a]),D^p(A)) \oto{\mathrm{mult} \circ (-)} \Map(D^p(L^{\otimes {b}} [-a]),A) \\
\cong & \Map(L^{\otimes pb} [-pa] \otimes \rB C_p,A) \oto{(-) \circ c_j} \Map(L^{\otimes pb} [-a - 2i(p-1)],A)
\end{align*}
where
\begin{equation}\label{eq:Steenrod-j}
j = pa - (a + 2i(p-1)) = (p-1)(a-2i)
\end{equation}
and $c_j$ is the generator of $\rH_j(\rB C_p;\F_p)$ 
pinned down at the beginning of \S\ref{sec:E_infty_ops_def}, which gives a map $\one_\cC[j] \to \one_\cC \otimes_{\F_p} \rB C_p$ by the $\F_p$-linearity of $\cC$. (By convention, if $j<0$ then we set $c_j= 0$.)
\end{defn}

See \cite{May70} for properties of the $\Pe^i$ thus defined, including Cartan relation and Adem relations. In this paper, we will only need the simple properties in the following Example. 

\begin{example}\label{ex:Pe-edge-cases} When $a=2i$, we see from \eqref{eq:Steenrod-j} and the definition that the induced map on cohomology 
\[
\Pe^i \co \pi_0  \Map(L^{\otimes b} [-a],A) \to \pi_0 \Map(L^{\otimes {pb}} [-pa],A)  
\]
is the $p$th power operation. Furthermore, if $2i>a$ so that \eqref{eq:Steenrod-j} is negative, then $\Pe^i$ vanishes.
\end{example}

Currently, we have a collision of notation between $\Pe^i$ from \Cref{defn:norm-E-infty-ops} and the operations defined in \S \ref{sec:E_infty_ops_def}. However, we will see that they agree when both are defined. To relate the two constructions, let $\cC$ be a site and let $\cD = \Shv(\cC;\Mod_{\F_p}^\gr(\Sp))$ be the $\infty$-category of graded sheaves of $\F_p$-module spectra over $\cC$. Then a graded commutative $\F_p$-algebra $A \in \calg^{\gr}(\Shv(\cC;\Sp))$ can be regarded as a commutative algebra $A \in \calg(\cD)$. Moreover, there is a strict Picard object $L\in \cD$ such that $\Map(L^{\otimes b},A)$ is the $b$-th graded piece of the sections of $A$. Applying the construction $\Pe^i$ above to $\pi_0\Map_U(L^{\otimes b}[-a],A):= \rH^{a,b}(U;A)$ for $U\in \cC$, we obtain maps 
\begin{equation}\label{eq:norm-E-infty-power-ops}
\Pe^i \colon \rH^{a,b}(U;A) \to \rH^{a+2i(p-1),pb}(U;A).
\end{equation}
\begin{prop}\label{prop:compare_Pe_defs}
In the situation above, the maps \eqref{eq:norm-E-infty-power-ops} above agree with the power operations defined in \S\ref{sec:E_infty_ops_def}.
\end{prop}

\begin{proof}
The proof when $A = \F_p$ and $\cC$ is the site of sheaves over a point, with trivial grading, is given in \cite[ Proposition IV.1.16]{NS18}. The same argument applies in the more general case we are considering. 
\end{proof}

\section{Comparing syntomic and $\EE_{\infty}$ Steenrod operations}\label{sec: comparing operations}
Over the past two sections, we have defined two different ``flavors'' of Steenrod operations on $\Hsyn^{*,*}(X)$ for a scheme $X/\Zpcyc$: the $\EE_\infty$-operations, and the syntomic Steenrod operations. These do \emph{not} agree (in general), as they have a different effect on the weights. Our goal in this section is to study their precise relationship. In particular, we will see that they \emph{do} agree in cases where their weights coincide.

\subsection{Formulation of the comparison}
Recall that we have fixed $\Zpcyc = \Z_p^{\cyc}$, with fraction field $\Qpcyc = \Q_p^{\cyc}$ and residue field $k$. Let $X$ be a scheme over $\Zpcyc$. Then we have defined two types of Steenrod operations on $\Hsyn^{*,*}(X)$:
\begin{itemize}
\item The $\EE_{\infty}$ Steenrod operations 
\[
\Pe^i \co \Hsyn^{a ,b}(X) \rightarrow \Hsyn^{a+2i(p-1),pb}(X).
\]
\item The syntomic Steenrod operations
\[
\Ps^i \co \Hsyn^{a,b}(X) \rightarrow \Hsyn^{a+2i(p-1),b+i(p-1)}(X).
\]
\end{itemize}

Choose a primitive $p$th root of unity $\zeta \in \mu_p(\Zpcyc)$. Then we obtain a corresponding element $\tau \in \Hsyn^{0,1}(\Zpcyc)$. Although $\tau$ depends on the choice of $\zeta$, the element $\tau^{p-1} \in \Hsyn^{0,p-1}(\Zpcyc)$ (often called $v_1$) is independent of the choice. All the formulas below refer to $\tau$ only through its $(p-1)^{\mrm{st}}$ power, hence are independent of the choice of $\zeta$. 

\begin{thm}\label{thm:comparison} Let $X$ be any scheme over $\Zpcyc$. Let $i\in \N$ and $a,b\in \Z$. 
\begin{itemize}
\item For $i>b$, we have  
\[
\Ps^i = \tau^{(p-1)(i-b)}\Pe^i \co \Hsyn^{a,b}(X) \rightarrow \Hsyn^{a+2i(p-1),b+i(p-1)}(X).
\]
\item For $b\ge i$ we have 
\[
\Pe^i = \tau^{(p-1)(b-i)}\Ps^i \co \Hsyn^{a ,b}(X) \rightarrow \Hsyn^{a+2i(p-1),bp}(X).
\]
\end{itemize}
\end{thm}

The main content behind \Cref{thm:comparison}  is work of Bachmann--Hopkins \cite{BH25}, which proves the analogous result for motivic cohomology of characteristic zero schemes. Our argument merely bootstraps their result using the perfectoid nearby cycles functor. 

For future calculations, the following special case of \Cref{thm:comparison} will be crucial. 

\begin{cor}\label{cor: operations agree}
For $i = b$ and any scheme $X/k$, the maps 
\begin{equation}\label{eq: compare Pe}
\Pe^i  \co \Hsyn^{a,b}(X)  \rightarrow \Hsyn^{a+2i(p-1), pb}(X)
\end{equation}
and 
\begin{equation}\label{eq: compare Pm}
\Ps^i \co \Hsyn^{a,b}(X) \rightarrow \Hsyn^{a+2i(p-1), b+(p-1)i}(X)
\end{equation}
agree. If $i<b$, then \eqref{eq: compare Pe} vanishes. If $i>b$, then \eqref{eq: compare Pm} vanishes. 
\end{cor}

\begin{proof}
Since $X$ is over $k$, multiplication by $\tau^{p-1} \in \Hsyn^{0,p-1}(\sO)$ factors over multiplication by its image in $\Hsyn^{0,p-1}(k) = 0$. 
\end{proof}

\begin{remark}
Analogous questions for mod $\ell$ motivic cohomology when $\ell \neq p$, and $k$ contains a primitive $\ell$th root of unity, were considered by Brosnan-Joshua in \cite{BJ15} and partially answered there. Their ``simplicial operations'' are what we call the $\EE_{\infty}$ Steenrod operations. Their methods are limited to $\ell \neq p$, but in any case the analogues of their results would not be enough for our purposes. Indeed, an important example for us is that the operations
\[
\Pe^1  \co \Hsyn^{3,1}(X) \rightarrow \Hsyn^{5,2}(X)
\]
and 
\[
\Ps^1 \co \Hsyn^{3,1}(X) \rightarrow \Hsyn^{5,2}(X)
\]
agree for $p=2$. If $X$ is over $\F_p$ with $p \neq 2$, the result of \cite{BJ15} for coefficients over $\Z/2$ says that these agree \emph{after further composing} with the map $\Het^{5}(X; \F_2(2)) \rightarrow \Het^{5}(X; \F_2(3))$ given by multiplication by the Bott element. In our situation of interest, where $X$ is a smooth proper surface over $k$, we actually have that $\Hsyn^{5,3}(X) =0$, so the analogous statement would be vacuous. 
\end{remark}

A priori, the operations $\Ps^i$ are hard to compute, even in the range of degrees where the classical Steenrod operations are trivial (e.g., on $\Hsyn^j$ where $j<i$). A consequence of \Cref{thm:comparison} is the computation of $\Ps^i$ in several important cases. 

\begin{cor}\label{cor: Sqm special cases}Let $X$ be a scheme over $\Zpcyc$. 

(1) The operation 
\begin{equation}\label{eq: Sqm 2i,i}
\Ps^{i} \co \Hsyn^{2i,i}(X) \rightarrow \Hsyn^{2pi,pi}(X)
\end{equation}
is given by raising to the $p$-th power. 

(2) Suppose $X$ is over $k$. Then the operation $\Ps^i$ vanishes on $\Hsyn^{a,b}(X)$ if $i > b$ for all $a$, or if $i =b$ and $2i \geq a$.
\end{cor}

\begin{proof}
(1) It was pointed out in \Cref{ex:Pe-edge-cases} that  
\[
\Pe^i \co \Hsyn^{2i,i}(X) \rightarrow \Hsyn^{2pi,pi}(X)
\]
is given by raising to the $p$th power. That map agrees with \eqref{eq: Sqm 2i,i} by \Cref{cor: operations agree}. 

(2) If $i > b$, then the vanishing follows from \Cref{cor: operations agree} and the observation that the restriction of $\tau^{p-1}$ to $\Hsyn^{0,p-1}(k)$ vanishes. If $i=b $, then \Cref{cor: operations agree} says that the operation agrees with 
\[
\Pe^i \co \Hsyn^{a,b}(X) \rightarrow \Hsyn^{a+2i(p-1),pb}(X),
\]
which then vanishes since $2i>a$ (as pointed out in \Cref{ex:Pe-edge-cases}).  

\end{proof}

\begin{remark}
For $X/k$, Annala--Elmanto have an alternate approach to Corollary \ref{cor: Sqm special cases} in \cite[Theorem 3.7(3) and Corollary 3.8]{AE}, which does not require the comparison to $\EE_\infty$ operations. However, \Cref{thm:comparison} is still needed crucially for the application to Brauer groups. 	

\end{remark}

\begin{question} It would be interesting to further investigate the computability of  the syntomic Steenrod operations. For example, our results do not (immediately) answer the question: Does $\Ps^i$ vanish on $\Hsyn^{a,b}(X)$ if $2i>a$ and $i < b$? 
\end{question}

\subsection{Proof of \Cref{thm:comparison}}

We will prove \Cref{thm:comparison} essentially by reduction to the case of a characteristic zero base, in which the analogous result for operations in motivic cohomology is proved in \cite{BH25}. To perform such a reduction, it is necessary to upgrade from statements about maps defined on cohomology (as in \Cref{thm:comparison}) to highly structured statements regarding maps of sheaves of spectra, so that we can apply the natural functoriality of such sheaf categories associated with $k,\Zpcyc,$ and $\Qpcyc$. 
Fortunately, all the maps under consideration have been constructed as maps of sheaves: 
\begin{itemize}
\item The map $\Pe^i\colon \Hsyn^{a,b}(X) \to \Hsyn^{a + 2(p-1)i,pb}(X)$ is obtained from the map of sheaves of spectra 
\[
\Pe^i\colon \SFp(b) \to \SFp(pb)[2(p-1)i]
\] 
by evaluation at $X$ and then taking $a^{\mrm{th}}$ cohomology groups. 

\item The map $\Ps^i \colon \Hsyn^{a,b}(X) \to \Hsyn^{a + 2(p-1)i,b + (p-1)i}(X)$ is obtained from the map of motivic spectra 
\[
\Ps^i \colon \MSFp \to \MSFp((p-1)i)[2(p-1)i]
\]
by twisting by $b$, evaluating at $X$ and then taking $a^{\mrm{th}}$ cohomology groups. As an intermediate step, the map $\Ps^i$ of motivic spectra restricts to maps of Nisnevich sheaves of spectra that we abusively denote again by 
\[
\Ps^i \colon \SFp(b) \to \SFp(b + i(p-1))[2i(p-1)].
\]
\end{itemize}

We will compare $\Pe^i$ and $\Ps^i$ in this incarnation: as maps of Nisnevich sheaves valued in $\Sp$. 

\begin{prop}\label{prop:comp_E_infty_motivic_O}
Let $i\in \N$ and $b\in \Z$. Then we have the following relations of morphisms in $\PShv_{\Nis}(\Sm_{\Zpcyc}; \Sp)$: 
\begin{itemize}
\item For $i>b$, we have  
\[
\Ps^i \cong \tau^{(p-1)(i-b)}\Pe^i \colon \SFp(b) \to \SFp(b + i(p-1))[2i(p-1)].
\]
\item For $b\ge i$, we have 
\[
\Pe^i \cong \tau^{(p-1)(b-i)} \Ps^i \co \SFp(b) \to \SFp(pb)[2(p-1)i].
\]
\end{itemize}
In particular, if $i=b$ then $\Pe^i \cong \Ps^i$. 
\end{prop}

\begin{proof}
By definition (cf. Remark \ref{rem: operations as psi}), the operation 
\[
\Ps^i \colon (\MSFp)_\Zpcyc \to \MSFp((p-1)i)_\Zpcyc[2(p-1)i]
\]
is obtained from Voevodsky's power operation $\Pm^i \colon (\MHFp)_{\Qpcyc} \to \MHFp((p-1)i)_\Qpcyc[2(p-1)i]$ by applying the functor $\Psienh$. Restricting to sheaves of spectra, we deduce that 
\[
\Ps^i \colon \SFp(b)_\Zpcyc \to \SFp(b+(p-1)i)_\Zpcyc[2(p-1)i]
\]
is obtained from the morphism 
\[
\Pm^i \colon \HFp(b)_\Qpcyc \to \HFp(b+(p-1)i)_\Qpcyc[2(p-1)i]
\]
by applying the functor 
$L_{\et} j_* \colon \PShv_{\Nis}(\Sm_\Qpcyc;\Sp) \to \PShv_{\Nis}(\Sm_\Zpcyc;\Sp)$. 
Similarly, the operation $\Pe^i$ for the syntomic complexes over $\Zpcyc$ is obtained from the operation $\Pe^i$ for the motivic cohomology complexes over $\Qpcyc$ by applying the functor $L_{\et} j_*$, thanks to \Cref{cor:E_infty_ops_psi}. It remains to show that the same-notated relations as in the theorem hold between the $\EE_\infty$-operations of the motivic complexes over $\Qpcyc$ and  Voevodsky's motivic power operations, which is the result \cite[Corollary 1.10]{BH25} of Bachmann--Hopkins. 
\end{proof}

\begin{proof}[Proof of \Cref{thm:comparison}] In view of \Cref{prop:comp_E_infty_motivic_O}, it remains to see that the operations $\Ps^i$ and $\Pe^i$ on syntomic cohomology groups are obtained from the operations on the underlying $\Sp$-valued Nisnevich sheaves from  \Cref{prop:comp_E_infty_motivic_O}, via evaluation at $X$ and passage to cohomology groups. For $\Ps^i$ this is tautological, and for $\Pe^i$ it is \Cref{prop:E_infty_op_global_sections}. 
\end{proof}

\part{Spectral prismatization}\label{part: spectral prismatization}

Prismatization over $k$, developed by Drinfeld \cite{Drin24} and Bhatt--Lurie \cite{BL22a}, lifts syntomic cohomology to quasicoherent sheaf theory on a stack $\FSyn{}$. In this Part, we carry out a generalization of this procedure for spectral syntomic cohomology, which we correspondingly call \emph{spectral prismatization}. In particular, in \S \ref{sec: spectral prismatization} we define a category $\pFGauge{\mbb{S}}$ of ``(pre)spectral prismatic $F$-gauges'' for $\Spec k$, which is a major step towards constructing Lurie's envisioned \emph{prismatic stable homotopy category} over $k$.

Then in \S \ref{sec: prismatization of Steenrod}, we prismatize the syntomic Steenrod algebra and its dual, lifting them to objects of $\pFGauge{\mbb{S}}$. This gets used to prove the compatibility statement between Poincar\'e duality and syntomic Steenrod operations, formulated as \Cref{thm: intro equivariance} in the Introduction, which is necessary for the eventual application to Brauer groups.

 In \cite[\S 4]{Bha22}, Bhatt--Lurie reinterpreted Poincar\'e duality for syntomic cohomology in terms of Serre duality on the stack $\FSyn{}$. The first step for proving \Cref{thm: intro equivariance} is to correspondingly lift the compatibility statement to the stack $\FSyn{}$, where it becomes an assertion about the compatibility of the prismatized syntomic Steenrod algebra and Serre duality. This compatibility is then explained by a theory of ``spectral Serre duality'' for $\pFGauge{\mbb{S}}$, which we develop in \S \ref{sec: Steenrod equivariance}.

\section{Spectral prismatic $F$-gauges}\label{sec: spectral prismatization}

In \S \ref{sec: syntomic spectra}, we defined the category of ``syntomic spectra'' $\Mod_{\psi(\SphFp)}(\MS_{k})$. In this section, we will define certain subcategories 
\[
\FGauge{\ZZ_p} \subset \Mod_{\MSZp}(\MS_{k}) 
\quad \text{and} \quad  \pFGauge{\mathbb{S}} \subset \Mod_{\psi(\SphFp)}(\MS_{k}),
\]
consisting of ``geometric objects''. We will eventually see that $\FGauge{\ZZ_p}$ can be identified with the category of ``prismatic $F$-gauges'' in the sense of Bhatt--Lurie, which justifies calling $\pFGauge{\mathbb{S}}$ the category of ``spectral prismatic $F$-gauges''. 

The philosophical significance of $\pFGauge{\mathbb{S}}$ is that it approximately matches Lurie's envisioned ``prismatic stable homotopy category'' over $k$. The practical significance, for our later applications, is that the categories $\Mod_{\MSFp}(\MS_k)$ and $\Mod_{\psi(\SphFp)}(\MS_{k})$ are ``too big'' to support a reasonable version of Serre duality; cutting down to the subcategories of (spectral) prismatic $F$-gauges will remedy this issue.

\subsection{Generation of classical prismatic $F$-gauges}\label{sssec:generation-prismatic-fgauge} Let $k$ be a finite field of characteristic $p$. Drinfeld and Bhatt--Lurie have defined a formal stack $\FSyn{} := k^{\Syn}$ over $\Spf \Z_p$ (following the notation of \cite[Chapter 4]{Bha22}). The derived category of quasicoherent sheaves $\cD(\FSyn{})$ is called the ``(derived) category of prismatic $F$-gauges over $k$''.  

The key preparation for defining the desired spectral enhancement $\pFGauge{\mathbb{S}}$ is to develop an alternative characterization of $\cD(\FSyn{})$. This subsection proves a conjecture of Bhatt \cite[Remark 4.4.6]{Bha22} about generation of $\cD(\FSyn{})$ by certain nice objects, which is the main content behind this alternative characterization.

\subsubsection{Geometry of $\FSyn{}$}\label{sssec:FSyn-geometry} We will need to invoke some explicit aspects of the construction in \cite[Definition 4.1.1]{Bha22}, which we review in terms that apply for any perfect field $k$ of characteristic $p$. 
\begin{enumerate}
\item The \emph{Nygaard filtered prismatization} of $\Spec k$ is the formal stack 
\[
(\Spec k)^{\cN}  = \Big[\sfrac{\Spf W(k)[u,t]/(tu-p)}{ \G_m} \Big].
\]
Here we need to specify the $\G_m$-action on $\Spf W(k)[u,t]/(tu-p)$, which amounts to specifying a grading on $W(k)[u,t]/(tu-p)$. Our convention is that $t$ is the \emph{Rees parameter}, so it has degree $-1$ action, and $u$ has degree $+1$. 

\item There are two open embeddings $j_{\mrm{HT}}, j_{\dR}\co  (\Spec k)^{\prism} \cong \Spf W(k) \rightrightarrows (\Spec k)^{\cN}$ of the \emph{prismatization} $(\Spec k)^{\prism}$ into $(\Spec k)^{\cN}$. The \emph{Hodge-Tate embedding} $j_{\mrm{HT}}$ is given by the locus $u \neq 0$, and the \emph{de Rham embedding} $j_{\dR}$ is given by Frobenius onto the locus $t \neq 0$. 
\end{enumerate}
Then $\FSyn{}$ is obtained by gluing $j_{\mrm{HT}}$ and $j_{\dR}$ along the obvious isomorphism. See Figure \ref{fig:ksyn} for a visualization.
\begin{figure}
\begin{tikzpicture}
    \draw[thick, blue, <->] (0,-3) -- (0,3) node[above] {Hodge-Tate};
    \draw[thick, red, <->] (-4,0) -- (4,0) node[right] {de Rham};
    
    \filldraw[purple] (0,0) circle (3pt) node[below right] {Hodge point};
    
    \foreach \r in {1, 1.5, 2, 2.5} {
        \draw[gray, thin] (\r, 0) arc (0:90:\r);
        \draw[gray, thin] (-\r, 0) arc (180:270:\r);
    }
\end{tikzpicture}
\caption{A visualization of $(\Spec k)^{\cN}$. There are two (open) embeddings of $(\Spec k)^{\prism}$, one corresponding to Hodge--Tate cohomology and the other corresponding to de Rham cohomology; their intersection corresponds to Hodge cohomology. The two axes are glued via Frobenius to form $\FSyn{}$.}\label{fig:ksyn}
\end{figure}

Let $\FSyn{\F_p}$ be the base change of $\FSyn{}$ to $\F_p$. While $\FSyn{}$ is a formal stack over $\Spf \Z_p$, the description above makes clear that $\FSyn{\F_p}$ is an algebraic stack over $\F_p$.

\subsubsection{Syntomification of schemes}\label{sssec:syntomification-of-schemes}
More generally, if $X$ is a scheme (or stack) over $k$, we can form its syntomification $X^{\Syn}$ as in \cite[\S 4.1]{Bha22}. In fact, the construction $X\mapsto X^{\Syn}$ is obtained via the ``transmutation'' procedure as in \cite{Bha22}, using the $k$-algebra stack $\G_{a,k}^{\Syn}$, so that 
\[
X^{\Syn}(R) = X(\G_a^{\Syn}(R)).
\]
In particular, this shows that the construction $X\mapsto X^{\Syn}$ is limit-preserving. 

\subsubsection{Prismatic $F$-gauges of schemes}

For each smooth $k$-scheme $f \co X \rightarrow \Spec k$, we have an object
\[
\cH^X := Rf^{\Syn}_* (\cO_{X^{\Syn}}) \in \cD(\FSyn{})
\]
discussed in \cite[\S 4.2]{Bha22}. When $f$ is smooth and proper, this is a perfect complex. We abbreviate 
\[
\ol \cH^X := \cH^X \otimes \F_p \in \cD(\FSyn{\F_p})
\]
which is a perfect complex (equivalently, a dualizable object) if $f$ is smooth and proper. These constructions organize into functors 
\begin{equation}\label{eq:classical-cohomological-motive}
\cH^{(-)} \co \Sm_k^{\op} \rightarrow \cD(\FSyn{}) \quad \text{and} \quad \ol \cH^{(-)} \co \Sm_k^{\op} \rightarrow \cD(\FSyn{\F_p}).
\end{equation}

\begin{thm}\label{thm: smooth proper motives generate}
The collection $\{\cH^X\}$, as $X$ ranges over smooth projective varieties over $\F_p$, compactly generates $\cD(\FSyn{})$. 
\end{thm}

\begin{remark}
The proof of \Cref{thm: smooth proper motives generate} will also show more generally that for any perfect field $k/\F_p$, the category $\cD(\FSyn{})$ is compactly generated by $\{\cH^X\}$ for smooth projective $X/k$, \emph{affirmatively} answering the question raised in \cite[Remark 4.4.6]{Bha22}. 
\end{remark}

\subsubsection{Preliminary observations}\label{sssec:preliminary-observations} Now we begin some technical preparations for the proof of \Cref{thm: smooth proper motives generate}. Let $\iota \co [\Spec k/\G_m] \inj \FSyn{} $ be the closed embedding of the special fiber of the ``Hodge point''. Let $\delta := \iota_* \cO_{[\Spec k/\G_m]}$ be the corresponding skyscraper sheaf. 
\begin{lemma}
Let $\cO$ be the structure sheaf of $\FSyn{}$ and $\delta$ as above. Then the category $\cD(\FSyn{})$ is generated under colimits by the objects $\{\cO[m]\BK{n}, \delta[m]\BK{n}\}_{m,n \in \Z}$. Here $\BK{n}$ denotes the Breuil-Kisin twist by $n$.
\end{lemma}

\begin{proof}

We note that since $\FSyn{}$ is a $p$-adic formal stack, $\cD(\FSyn{})$ is generated by $\cD(\FSyn{\F_p})$. Therefore, it suffices to show that if $\cF \in \cD(\FSyn{\F_p})$ is right-orthogonal to all shifts and Breuil-Kisin twists of $\cO$ and $\delta$, then $\cF = 0$.

Suppose that
\[
\Hom_{\FSyn{}}(\delta[m]\BK{n}, \cF) = 0 
\]
for every $m,n \in \Z$. By adjunction, this implies that $\iota^! \cF = 0$. Consider the commutative diagram 
\begin{equation}
\begin{tikzcd}
\Spec k \ar[d, "\pr'", twoheadrightarrow]\ar[r, "\iota'"] & \Spf W(k)[u,t]/(tu-p) \ar[d, "\pr", twoheadrightarrow] \\
(\Spec k)/\G_m \ar[r, "\iota"] & \FSyn{}
\end{tikzcd}
\end{equation}
The commutativity of the diagram implies that 
\begin{equation}\label{eq:shrieks-vanish}
\iota'^!  \pr^! (\cF) \cong \pr'^! \iota^! (\cF) = 0.
\end{equation}
Since $\iota'$ is a regular embedding cut out by $u=0$ and $t=0$, the functor $\iota'^!$ identifies up to shift with the derived mod $(u,t)$ reduction. More precisely, we have a natural isomorphism $\iota'^! \cong \iota'^*[-2]$ (all sheaf operations are derived). Then from \eqref{eq:shrieks-vanish} we see that $\iota'^* \pr^! (\cF) = 0$. Since $\pr$ is smooth, $\pr^!$ differs from $\pr^*$ by a shift and invertible twist, so we deduce that 
\[
0 = \iota'^* \pr^* (\cF) = \pr'^* \iota^* (\cF).
\]
Since $\pr'$ is faithfully flat, we conclude that $\iota^* \cF = 0$. 
 
 As we assumed that $\cF$ is supported on the special fiber $\FSyn{\F_p}$, this shows that $\cF = j_* \cF_0$ is pushed forward from the open complement of $\iota$ in $\FSyn{\F_p}$. But this open complement is isomorphic to $\Spec k$, so if $\cF_0$ is right-orthogonal to all shifts of the pullback of $\cO$, then it must be zero. 
\end{proof}

Let $\tw{\cH^X}_X \subset \cD(\FSyn{})$ be the full subcategory generated under colimits by objects $\cH^X$ for smooth projective $X/k$. Clearly we have $\{\cO[m]\BK{n}\}_{m, n\in \Z} \subset \tw{\cH^X}_X$. To prove Theorem \ref{thm: smooth proper motives generate}, it therefore suffices to argue that $\delta \in \tw{\cH^X}_X $. Let $k \rightarrow k'$ be a finite extension. In turn, it suffices to build $\delta'$, the pushforward of the skyscraper sheaf along $\iota' \co [\Spec k'/\G_m] \inj (\Spec k')^{\Syn} \rightarrow \FSyn{}$, since $\delta$ is a summand of the pushforward of $\delta'$. We will do this with a very explicit construction built out of supersingular elliptic curves.

\subsubsection{Supersingular elliptic curves}
Let $E$ be a supersingular elliptic curve over $\F_q$, for $q=p^2$, such that $\Frob_q$ acts as multiplication by $p$. This implies that all endomorphisms of $E$ are defined over $\F_q$, so that $\End_{\F_q}(E)$ is a maximal order in a division algebra ramified exactly at $p$ and $\infty$. Let $D := \End(E) \otimes_{\Z} \Z_p$ and $\varpi \in D$ be a uniformizer.

For the rest of the section, we take $k := \F_q$. The claim in question for any other $k'$ easily reduces to this case, by the observations in \S \ref{sssec:preliminary-observations}. 

For the perfect complex $\cH^E \in \cD(\FSyn{})$, let $H := \cH^1(\cH^E)$. We claim that $H$ is a vector bundle on $\FSyn{}$. As explained in the proof of \cite[Theorem 3.5.1]{Bha22}, this follows from the fact that 
\begin{itemize}
\item $E$ has $p$-torsionfree crystalline cohomology, and 
\item the Hodge--de Rham spectral sequence degenerates for $E$. 
\end{itemize}
Note that $\cH^0(\cH^E) \cong \cO$ and $\cH^2(\cH^E) \cong \cO\BK{-1}$, so that $H = \cH^1(\cH^E)$ indeed lies in $\tw{\cH^X}_X$.

\begin{lemma}\label{lem: Nygaard filtration on H}
The Nygaard filtration of $H$ is given by $\Fil^0_{\cN}(H) = H$ and $\Fil^n_{\cN}(H) = \varpi p^{n-1} H$ for $n \geq 1$. 
\end{lemma}

\begin{proof}
The statement $\Fil^0_{\cN}(H) = H$ is tautological from the construction of the Nygaard filtration. Since $E$ has dimension $1$, it is also immediate that $\Fil^{n+1}_{\cN}(H) = p \Fil^n_{\cN}(H)$ for $n \geq 1$. 

It only remains to see that $\Fil^1_{\cN}(H) = \varpi H$. We always have $H \supset \Fil^1_{\cN}(H) \supset p H$, and both inclusions are strict because $\Fil^1_{\cN}(H)/ p \Fil^0_{\cN}H$ is the Hodge filtration on the de Rham cohomology of $E$, which is non-trivial. Since $H^1(E)$ is free of rank 1 over $D$, the only $D$-module lying strictly between $H$ and $pH$ is $\varpi H$, so this must be identified with $\Fil^1_{\cN}(H)$. 
\end{proof}

\subsubsection{The prismatic $F$-gauge $M$}\label{sssec: M} Recall the Nygaard filtered prismatization of $k$, 
\[
(\Spec k)^{\cN}  = \Big[\sfrac{\Spf W(k)[u,t]/(tu-p)}{ \G_m}\Big].
\]
Here our convention is that $t$ is the \emph{Rees parameter}, so that it has weight $-1$ for the $\G_m$-action, and $u$ has weight $+1$. There is an \'etale covering $(\Spec k)^{\cN} \rightarrow \FSyn{}$ obtained by gluing the two open sections $j_{\mathrm{HT}}, j_{\dR} \co  \Spf W(k)  \cong (\Spec k)^{\prism} \rightrightarrows (\Spec k)^{\cN}$ along the Frobenius of $W(k)$.

From the definition of $(\Spec k)^{\cN}$, we see that quasicoherent sheaves on $(\Spec k)^{\cN}$ are identified with graded modules over the ring $W(k)[u,t] / (tu-p)$ which are $p$-adically complete in a suitable sense; see \cite[\S 3.3]{Bha22} for the precise formulation. 

\begin{notation}\label{not: gauge diagram}
We will depict a graded module $M_{\bu} \cong \bigoplus_n M_n$ over $W(k)[u,t]/(tu-p)$ as a diagram 
\[
\adjustbox{scale = 0.9, center}{\begin{tikzcd}
\deg & \ldots & -1 & 0 & 1  & \ldots \\
M_{\bu} & \cdots \ar[r, bend left, "u"]  & M_{-1} \ar[l, bend left, "t"] \ar[r, bend left, "u"] & M_0 \ar[l, bend left, "t"]  \ar[r, bend left, "u"] & M_1 \ar[r, bend left, "u"]  \ar[l, bend left, "t"]  & \cdots   \ar[l, bend left, "t"]
\end{tikzcd}}
\]

\end{notation}

\begin{example}\label{ex:nygaard-H}
The $W(k)[u,t]/(tu-p)$-module $H_{\bu}$ associated to $H|_{(\Spec k)^{\cN}}$ has $H_i$ being the $i$th Nygaard filtrant, with $t$ being the inclusion and $u$ being multiplication by $p$. Hence in terms of \Cref{not: gauge diagram}, $H_{\bu}$ has the form 
\begin{equation}\label{diag: H gauge}
\adjustbox{scale = 0.9, center}{\begin{tikzcd}
\ldots & -1 & 0 & 1 & 2 & \ldots \\
\cdots \ar[r, bend left, "u = p"] & H \ar[l, bend left, "t", "\sim"']  \ar[r, bend left, "u=p"] & H \ar[l, bend left, "t", "\sim"']  \ar[r, bend left, "u=p"]  & \varpi H  \ar[l, bend left, "t"]   \ar[r, bend left, "u=p", "\sim"']   & p \varpi H  \ar[r, bend left, "u=p", "\sim"'] \ar[l, bend left, "t"]  & \cdots  \ar[l, bend left, "t"]  
\end{tikzcd}}
\end{equation}

\end{example}

Let $M := H/ \varpi H$, considered as a coherent sheaf on $\FSyn{}$. By Lemma \ref{lem: Nygaard filtration on H}, the Nygaard filtration on $\varpi H$ is 
\[
\Fil^n_{\cN}(\varpi H) = \begin{cases} \varpi H & n = 0, \\  \varpi^2 p^{n-1} H & n \geq 1. \end{cases}
\]
Hence the underlying graded module of $ M|_{(\Spec k)^{\cN}}$ is $\bigoplus_{n \in \Z} (H/\varpi H) \cong \bigoplus_{n \in \Z} \F_q$. We see from \eqref{diag: H gauge} that as a graded $W(k)[u,t]/(tu-p)$-module, $M|_{(\Spec k)^{\cN}}$ is the direct sum of a class $v_0$ in degree $0$ which is annihilated by $u$ and a class $w_1$ in degree $1$ which is annihilated by $t$. 
In other words, in terms of Notation \ref{not: gauge diagram}, $M|_{(\Spec k)^{\cN}}$ looks like 
\begin{equation}\label{diag: M gauge}
\adjustbox{scale = 0.9, center}{\begin{tikzcd}
\ldots & -1 & 0 & 1 & 2 & \ldots \\
\cdots & H/\varpi \ar[l, bend left, "t", "\sim"'] & H/\varpi \ar[l, bend left, "t", "\sim"']   & \varpi H / \varpi^2 H \ar[r, bend left, "u", "\sim"']   & p \varpi H / p \varpi^2 H \ar[r, bend left, "u", "\sim"'] & \cdots  \\
\ldots & t v_0 \ar[l] & \ar[l] v_0 & w_1 \ar[r] &  u w_1 \ar[r] & \ldots 
\end{tikzcd}}
\end{equation}
where the maps $u$ and $t$ vanish when not depicted in the diagram. The following Lemma articulates the uniqueness of a prismatic $F$-gauge with this form. 

\begin{lemma}\label{lem: M unique}There is a unique isomorphism class of prismatic $F$-gauges whose pullback to $(\Spec k)^{\cN}$ is isomorphic to $M|_{(\Spec k)^{\cN}}$.
\end{lemma}

\begin{proof}Examining the explicit construction of $\FSyn{}$, we see that the additional datum required to descend $\cM|_{(\Spec k)^{\cN}}$ from $(\Spec k)^{\cN}$ to $\FSyn{}$ is that of a Frobenius-semilinear isomorphism between the localizations with respect to $u$ and $t$, which intertwines the $u$-action with the $t$-action. By the form of \eqref{diag: M gauge}, such an isomorphism must restrict to an isomorphism between $M_0$ and $M_1$, which are both 1-dimensional $\F_q$-vector spaces, and be determined by this restriction. By rescaling one of the generators if necessary, we can assume that this isomorphism carries the image of $v_0$ (in the localization with respect to $t$) to the image of $w_1$ (in the localization with respect to $u$), which now uniquely specifies the isomorphism class of the descent. 
\end{proof}

\subsubsection{Square of supersingular elliptic curve} 
Note that we have
\[
\cH^{E \times_k E} \cong \cH^E \otimes_{\FSyn{} } \cH^E \in \Perf(\FSyn{}).
\]
Therefore, $H \otimes_{\FSyn{} } H =: H^{\otimes 2} \in \Perf(\FSyn{})$ also lies in $\tw{\cH^X}_X$.

We analyze some related prismatic $F$-gauges. Note that by geometric Poincar\'e duality \cite{Tang22}, we have an isomorphism
\[
\End_{\FSyn{}} (H) \cong H \otimes H \BK{1}.
\]
We have a natural map $D \otimes_{\Z_p} \Z_q \rightarrow \End_{\FSyn{}}(H)$. This induces a map of rank 4 vector bundles over $\FSyn{\F_p}$,
\[
D \otimes_{\Z_p} \cO\BK{-1} \rightarrow 
\End_{\FSyn{}} (H) \BK{-1} \cong H \otimes_{\FSyn{} } H =: H^{\otimes 2}.
\]

By Lemma \ref{lem: Nygaard filtration on H}, as a graded module over $W(k)[u,t]/(tu-p)$, $H^{\otimes 2}|_{(\Spec k)^{\cN}}$ looks like 
\[
\adjustbox{scale = 0.9, center}{\begin{tikzcd}[row sep = tiny, column sep = tiny]
\ldots & -1 & 0 & 1  & 2  & \ldots  \\
\ldots & H^{\otimes 2} &  H^{\otimes 2} & (\varpi H) \otimes H &  p (\varpi H) \otimes H & \ldots \\
& & &  +  H \otimes (\varpi H) & +  p (H \otimes \varpi H) \\
& & & & + (\varpi H \otimes \varpi H)  
\end{tikzcd}}
\]
with the action of $t$ being the obvious inclusions, and the action of $u$ being multiplication by $p$. In particular, we have a nesting of coherent sheaves on $\FSyn{}$, 
\begin{equation}\label{eq: nesting}
p (H \otimes_{\FSyn{} } H)  \inj D \otimes_{\Z_p} \cO\BK{-1} \inj (H \otimes_{\FSyn{} } H).
\end{equation} 
By the standard structure theory of endomorphisms of supersingular elliptic curves \cite[Theorem 42.1.9, \S 23.1.3]{Voi21}, we can choose a $W(k)$-module basis of $H|_{(\Spec k)^{\cN}}$ such that 
\begin{itemize}
\item $\varpi$ acts as $\begin{pmatrix} 0 & 1 \\ p & 0 \end{pmatrix}$, and 
\item $D \otimes_{\Z_p} \Z_q \cong \Gamma_0(\Z_q)$ is the Eichler order of matrices with coefficients in $\Z_q$ which are upper-triangular modulo $p$. 
\end{itemize}
With this choice of basis, the maps in \eqref{eq: nesting} look as follows when pulled back to $(\Spec k)^{\cN}$, in terms of Notation \ref{not: gauge diagram}.
\begin{equation}\label{eq: M' gauge}
\adjustbox{scale =0.8, center}{\begin{tikzcd}[row sep = huge] 
\deg & \ldots & 0 & 1 & 2  &  \ldots  \\
\End(H)\BK{-1} &  \cdots \ar[r, "u", bend left]  & \begin{pmatrix} * & * \\ * & * \end{pmatrix} \ar[l, "t", bend left] \ar[r, "u", bend left] &  \begin{pmatrix} * & * \\ p & * \end{pmatrix}  \ar[l, "t", bend left] \ar[r, "u", bend left] &  \begin{pmatrix} p & * \\ p^2 & p \end{pmatrix}  \ar[l, "t", bend left] \ar[r, "u", bend left]  & \cdots \ar[l, "t", bend left] \\
D \otimes_{\Z_p} \cO\BK{-1} \ar[u, hook] & \cdots \ar[r, "u", bend left]  & \begin{pmatrix} * & * \\ p & * \end{pmatrix} \ar[u, hook]  \ar[l, "t", bend left] \ar[r, "u", bend left] & \begin{pmatrix} * & * \\ p & * \end{pmatrix} \ar[l, "t", bend left] \ar[r, "u", bend left] \ar[u, hook] &  \begin{pmatrix} p & p \\ p^2 & p \end{pmatrix} \ar[u, hook] \ar[l, "t", bend left] \ar[r, "u", bend left] & \cdots \ar[l, "t", bend left] \\
p \End(H)\BK{-1} \ar[u, hook]  & \cdots  \ar[r, "u", bend left]  & \begin{pmatrix} p & p \\ p & p \end{pmatrix} \ar[u, hook]  \ar[l, "t", bend left] \ar[r, "u", bend left] &  \begin{pmatrix} p & p \\ p^2 & p \end{pmatrix}  \ar[u, hook] \ar[l, "t", bend left] \ar[r, "u", bend left] &  \begin{pmatrix} p^2 & p \\ p^3 & p^2 \end{pmatrix} \ar[u, hook] \ar[l, "t", bend left] \ar[r, "u", bend left] & \cdots \ar[l, "t", bend left] \\
\end{tikzcd}}
\end{equation}

\subsubsection{The prismatic $F$-gauge $M'$}
There is a trace map 
\begin{equation}\label{eq: trace}
H^{\otimes 2} \rightarrow \cO\BK{-1}
\end{equation}
corresponding under geometric Poincar\'e duality to the trace map $\End_{\FSyn{}}(H) \rightarrow \cO_{\FSyn{}}$. Let us write $(H^{\otimes 2})_0$ for the kernel of this map, which is a rank 3 vector bundle on $\FSyn{}$.

The restriction of \eqref{eq: trace} to $D \otimes_{\Z_p} \cO\BK{-1}$ is the trace map $D \otimes_{\Z_p} \Z_q \rightarrow \Z_q$ tensored (over $\Z_q$) with $\cO\BK{-1}$. Hence the kernel of this restriction is $D_0 \otimes_{\Z_q} \cO\BK{-1}$ where $D_0  \subset D \otimes_{\Z_p} \Z_q$ is the subspace of trace 0 matrices. 

Then consider the prismatic $F$-gauge 
\[
\wt M := \coker\Big(p (H^{\otimes 2})_0 \rightarrow D_0 \otimes_{\Z_q} \cO\BK{-1}\Big) \in \Coh(\FSyn{}).
\]
From examining \eqref{eq: M' gauge}, we see that the restriction $\wt M|_{(\Spec k)^{\cN}}$ has the following description as a graded module over $\Z_q[u,t]/(tu-p)$: 
\[
\begin{tikzcd}[row sep = tiny]
\ldots & 0 & 1 & 2  &  \ldots  \\
 \F_q^2 & \F_q^2  &  \F_q^3  &  \F_q^2 & \F_q^2 
\end{tikzcd}
\]
More precisely, $\wt M|_{(\Spec k)^{\cN}}$ has 3 generators over $\Z_q[u,t]/(tu-p)$, in degree 1: 
\begin{enumerate}
\item one annihilated by $u$ (corresponding to the lower-left entry in \eqref{eq: M' gauge}) which we call $v_1'$, 
\item one annihilated by $t$ (corresponding to the upper-right entry in \eqref{eq: M' gauge}) which we call $w_1'$, and 
\item one not killed by any power of $u$ or $t$ (corresponding to anti-diagonal entry in \eqref{eq: M' gauge}) which we call $y_1'$. This generator spans the line bundle $\cO\BK{-1}$ on the special fiber $\wt M|_{(\Spec k)^{\cN}_{\F_p}}$. 
\end{enumerate}

\[
\adjustbox{scale=0.92,center}{
\begin{tikzcd}[row sep=1.4em, column sep=2.8em]
\deg & \cdots & -1 & 0 & 1 & 2 & 3 & \cdots \\
& \cdots
& t^2v_1'
& tv_1' \arrow[l, "t"'] 
& v_1' \arrow[l, "t"']
& 0
& 0
& \cdots \\
& \cdots
& t^2y_1'
& ty_1' \arrow[l, "t"'] 
& y_1' \arrow[l, "t"'] \arrow[r, "u"]
& uy_1' \arrow[r, "u"]
& u^2y_1'
& \cdots \\
& \cdots
& 0
& 0
& w_1' \arrow[r, "u"]
& uw_1' \arrow[r, "u"]
& u^2w_1'
& \cdots
\end{tikzcd}}
\]

Then the $(t,u)$-power torsion subsheaf of $\wt M$ admits the following description: it is the unique coherent sheaf on $\FSyn{\F_p}$ (up to isomorphism) which when pulled back to $(\Spec k)^{\cN}_{\F_p}$, is isomorphic to the graded $\F_q[u,t]/(tu-p)$-module with two generators $v_1'$ and $w_1'$ in degree 1, one annihilated by $u$ and the other annihilated by $t$. The uniqueness aspect here has the same meaning (and proof) as in Lemma \ref{lem: M unique}.

Let us denote by $M' \subset \wt M$ this $(t,u)$-power torsion subsheaf. Then the quotient $\wt M/M'$ is some line bundle $\cL$ on the special fiber $\FSyn{\F_p} \inj \FSyn{}$, since it pulls back to a line bundle on the \'etale cover $(\Spec k)^{\cN}_{\F_p}$. Thus we have an extension
\begin{equation}\label{eq: M' filtration}
0 \rightarrow M' \rightarrow  \wt M \rightarrow \cL \rightarrow 0.
\end{equation}
Moreover, the line bundle $\cL$ becomes isomorphic to $\cO\BK{-1}$ after restriction to $(\Spec k)^{\cN}_{\F_p}$. Then $\cL$ is determined by the gluing isomorphism, which in terms of the given trivialization is multiplication by an element $\lambda \in k^\times$. Since multiplying the trivialization by $\alpha$ multiplies the gluing isomorphism by $\phi(\alpha) \alpha^{-1}$, multiplying $\lambda$ by $\phi(\alpha) \alpha^{-1}$ has no effect on the isomorphism class of $\cL$. Hence, after making a finite base change $k \rightarrow k'$ if necessary so that $\lambda$ can be written in the form $\phi(\alpha) \alpha^{-1}$, we can trivialize $\cL$. Therefore, after making such an extension, we find that $\cL$, and then by \eqref{eq: M' filtration} also $M'$, lies in $\tw{\cH^X}_X$. 

\subsubsection{Finding $\delta$}\label{sssec: finding delta}
Recall the prismatic $F$-gauge $M$ from \S\ref{sssec: M}. Observe that there is a monomorphism of prismatic $F$-gauges $M \rightarrow M'$, which in the explicit presentations above is given by $v_0 \mapsto  t v_1'$ and $w_1 \mapsto w_1'$. (This description of the map on $(\Spec k)^{\cN}$ descends to $\FSyn{}$ since it can be arranged to be compatible with the gluing.) We claim that the cokernel $M'/M$ is isomorphic to $\delta'\BK{-1}$. This can be checked after pulling back to $(\Spec k)^{\cN}$. On this pullback, we see explicitly that $(M'/M)|_{(\Spec k)^\cN}$ is the graded $\Z_q[u,t]/(tu-p)$-module which is a 1-dimensional $\F_q$-vector space in degree 1 (generated by the image of $v_1'$), annihilated by both $t$ and $u$. This verifies the claim, and completes the proof of Theorem \ref{thm: smooth proper motives generate}. \qed 

\subsubsection{Shortcut for $p=2$}
If $p=2$, then we remark that there is a more direct way to build $\delta$ from \eqref{eq: M' filtration}, without having to pass to a finite extension. 

We equip $\End(H)$ with the following filtration. Observe that we have the tautological map $\cO_{\FSyn{}} \rightarrow \End(H)$ given by the identity element, and the trace map $\End(H) \rightarrow \cO_{\FSyn{}}$. These maps compose to multiplication by $\rank H = 2$, which vanishes on the special fiber. Thus we obtain a filtration on $\End(H)/p$, which translates to a filtration $F_{\bu} (H^{\otimes 2}/p)$, with 
\[
\begin{tikzcd}[row sep = tiny]
F_0  \ar[d, equals] & F_1  \ar[d, equals]  & F_2 \ar[d, equals]  & F_3 \ar[d, equals]\\
0 \ar[r, hook]  & \cO\BK{-1}/p \ar[r, hook] & \ker(H^{\otimes 2}/p \rightarrow \cO\BK{-1}/p) \ar[r, hook]  & (H^{\otimes 2}/p)
\end{tikzcd}
\]

The pullback of this filtration to $D \otimes_{\Z_p} \cO\BK{-1}/p$ has $F_1(D \otimes_{\Z_p} \cO\BK{-1}/p) \cong \cO\BK{-1}/p$ and $F_2(D \otimes_{\Z_p} \cO\BK{-1}/p) = D_0 \otimes \cO\BK{-1}/p$. Composing its inclusion into $D_0 \otimes_{\Z_q} \cO\BK{-1}/p$ with the projection modulo the image of $p( H^{\otimes 2})_0$, we obtain a map
\[
\cO \BK{-1}/p \rightarrow 
\coker\left( p (H^{\otimes 2})_0 \rightarrow D_0 \otimes_{\Z_q} \cO\BK{-1} \right)  = \wt{M}
\]
which splits the filtration \eqref{eq: M' filtration}. Hence it shows that $\cL = \cO\BK{-1}/p$ in the notation there, and we may proceed as in  \S \ref{sssec: finding delta} to conclude that $\delta \in \tw{\cH^X}_X$.

\subsection{The category of spectral prismatic $F$-gauges}\label{ssec:spectra-prismatic-f-gauges}

Recall that if $X \in \Sm_k$, then there is an associated object $\Sigma_+^\infty X \in \MS_k$ which we think of as the ``motive of $X$'' (cf. \Cref{defn:motives-of-schemes}).

\begin{defn}\label{defn:mod-p-FGauge}We define the $\infty$-category 
\[
\FGauge{\ZZ_p}\subseteq \Mod_{\psi(\MHZp)}(\MS_k) \cong \Mod_{\MSZp}(\MS_k)
\]
to be the full subcategory generated under colimits, twists, and shifts from objects of the form $\Sigma^\infty_+X\otimes\MSZp$ where $X$ is \emph{smooth and projective} over $k$. 

We define the $\infty$-category of \emph{mod $p$ prismatic $F$-gauges over $k$} 
\[
\FGauge{\FF_p}\subseteq \Mod_{\psi(\MHFp)}(\MS_k) \cong \Mod_{\MSFp}(\MS_k)
\]
similarly, with $\MSFp$ in place of $\MSZp$. 
\end{defn}

\begin{remark}We will see later (in \Cref{prop: classical QCoh}) that there are natural equivalences $\FGauge{\ZZ_p} \cong \cD(\FSyn{})$ and $\FGauge{\F_p} \cong \cD(\FSyn{\FF_p})$, justifying the names of our categories. 
\end{remark}

In particular, $\FGauge{\FF_p}$ is a compactly generated (compactness of $\Sigma^\infty_+ X \otimes \MSFp$ holding by \Cref{prop:MS_compact_generators}), symmetric monoidal presentable subcategory of $\Mod_{\MSFp}(\MS_k)$.

\subsubsection{The spectral enhancement} We will now define an enlargement of $\FGauge{\FF_p}$, whose relationship to $\FGauge{\FF_p}$ is analogous to the relationship between the $p$-complete stable homotopy category $\Sp$ and the derived category of $\F_p$-modules $\cD(\F_p)$. 

\begin{defn}
We define $\pFGauge{\mbb{S}}$ to be the fibered product of $\infty$-categories 
\[
\pFGauge{\mbb{S}} := \Mod_{\psi(\SphFp)}(\MS_k)\times_{\Mod_{\MSFp}(\MS_k)} \FGauge{\FF_p}.  
\]
In other words, $\pFGauge{\mbb{S}}$ is the full subcategory of $\Mod_{\psi(\SphFp)}(\MS_k)$ spanned by all objects whose image under the (colimit-preserving, symmetric monoidal) functor
\[
\Mod_{\psi(\SphFp)}(\MS_k) \xrightarrow{(-) \otimes_{\psi(\SphFp)} \psienh(\MHFp) }  \Mod_{\psienh(\MHFp)}(\Mod_{\psi(\SphFp)}(\MS_k) ) \cong \Mod_{\MSFp}(\MS_k),
\]
the last equivalence being \Cref{cor:module-category}, belongs to $\FGauge{\FF_p}$.
\end{defn}

\begin{remark}[Relation to prismatic stable homotopy theory]\label{rem:true-fgauge}
Lurie has conjectured the existence of a ``prismatic stable homotopy category'' over a base scheme, whose relationship to prismatic $F$-gauges should be analogous to the relationship between the stable homotopy category and the derived category of $\Z$-modules. One motivation for this hypothetical category is to capture algebraic invariants studied in \cite{BMS2} such as $\mathrm{THH}$, $\mathrm{TC}$, etc., which are not $\Z$-linear. 

For the base scheme $\Spec k$, this category can be constructed unconditionally (although it is not yet in the literature); let us denote it $\FGauge{\mbb{S}}$. It is closely related to $\pFGauge{\mbb{S}}$ but not quite the same: we expect that $\FGauge{\mbb{S}}$ is the right-completion of $\pFGauge{\mbb{S}}$ with respect to a certain ``prismatic t-structure''. In particular, there is a functor $\pFGauge{\mbb{S}} \rightarrow \FGauge{\mbb{S}}$ which induces an equivalence on the subcategories of eventually connective objects for the prismatic t-structure. It turns out that in this paper, we will only deal with eventually connective objects, so it suffices for our purposes to work in $\pFGauge{\mbb{S}}$. Furthermore, the definition of the prismatic t-structure, and the proof of the desired properties of $\FGauge{\mbb{S}}$, involve considerable technicalities; since this paper already feels technical enough, we defer the study of $\FGauge{\mbb{S}}$ to a future work. 
\end{remark}

\subsubsection{Modules over syntomic cohomology} Next we will show that the syntomic cohomology object $\MSFp$ admits a natural structure of a commutative algebra in $\pFGauge{\mbb{S}}$, whose module category recovers the usual mod $p$ prismatic $F$-gauges: 
\[
\Mod_{\MSFp}(\pFGauge{\mbb{S}})\cong \FGauge{\FF_p}.
\]
In geometric terms, this implies that the canonical forgetful functor $\FGauge{\FF_p} \to \pFGauge{\mbb{S}}$ behaves like a pushforward along an affine morphism of stacks. 

\begin{prop}\label{prop:iota_affine}
The symmetric monoidal adjunction 
\[
(-)\otimes_{\psi(\SphFp)} \MSFp  \colon \Mod_{\psi(\SphFp)}(\MS_k) \adj \Mod_{\MSFp}(\MS_k) : \textrm{forget}
\]
restricts to an adjunction
\[
\iota^*\colon \pFGauge{\mbb{S}} \adj \FGauge{\FF_p} : \iota_*.
\]
The functor $\iota_*$
is colimit-preserving, conservative, and linear over $\pFGauge{\mbb{S}}$; in particular, we have a canonical projection formula 
\[
\iota_*(M \otimes \iota^* N)\cong \iota_*(M)\otimes N \quad \text{ for all } M\in \FGauge{\FF_p}, N\in \pFGauge{\mbb{S}}. 
\]
\end{prop}

\begin{proof}
The functor $(-) \otimes_{\psi(\SphFp)} \MSFp$ carries $\pFGauge{\mbb{S}}$ into $\FGauge{\FF_p}$ by definition. To show that its right adjoint carries $\FGauge{\FF_p}$ to $\pFGauge{\mbb{S}}$, by the definition of $\pFGauge{\mbb{S}}$ it would suffice to show that the composition 
\[
\Mod_{\MSFp}(\MS_k) \xrightarrow{\text{forget}} \Mod_{\psi(\SphFp)}(\MS_k) \xrightarrow{(-) \otimes_{\psi(\SphFp)} \MSFp } \Mod_{\MSFp}(\MS_k)  
\]
carries any $M \in \FGauge{\FF_p}$ into $\FGauge{\FF_p}$. \Cref{prop:psi_enh_on_Fp^2} implies (together with \Cref{cor:psi_enh_on_Fp}) that 
\[
M\otimes_{\psi(\SphFp)} \MSFp \cong M\otimes_{\MSFp} \left(\MSFp \otimes_{\psi(\SphFp)} \MSFp \right) \cong \bigoplus_{\alpha \in \sI} M[p_\alpha] (q_\alpha)
\] 
(where $\sI, p_\alpha, q_\alpha$ are as defined in \S \ref{sssec:dual-motivic-steenrod}) which is a direct sum of shifts and twists of $M$. Since $\FGauge{\FF_p}$ is closed under shifts, twists, and infinite direct sums in $\Mod_{\MSFp}(\MS_k)$, we see that $M\otimes_{\psi(\SphFp)} \MSFp \in \FGauge{\FF_p}$, as desired. 

The assertions that $\iota_*$ is conservative, colimit-preserving, and linear over $\pFGauge{\mbb{S}}$ now follow from the similar properties of the forgetful functor $\Mod_{\MSFp}(\MS_k)\to \Mod_{\psi(\SphFp)}(\MS_k)$. These properties hold for any forgetful functor from modules over an algebra in an arbitrary presentably symmetric monoidal category, which applies here because of \Cref{cor:module-category}.
\end{proof}

Since $\iota_*$ is lax symmetric monoidal, it automatically factors as 
\begin{equation}\label{eq:iota-enhance}
\begin{tikzcd}
\FGauge{\FF_p} \ar[r, "\iota_*^{\enh}"] \ar[dr, "\iota_*"']  &  \Mod_{\iota_*\MSFp}(\pFGauge{\mbb{S}}) \ar[d, "\text{forget}"] \\
& \pFGauge{\mbb{S}}
\end{tikzcd}
\end{equation}

\begin{cor}\label{cor:fgauge-module-category}
The functor $\iota_*^{\enh}$ in \eqref{eq:iota-enhance} is an equivalence. 
\end{cor}
	
\begin{proof}
This follows from the colimit-preservation, conservativity, and linearity of $\iota_*$ established in \Cref{prop:iota_affine}; see for example \cite[Proposition 2.5]{ChromFourier}.
\end{proof}

\subsubsection{Compatibility with the syntomic Steenrod algebra}
Since the fully faithful embedding
\begin{equation}\label{eq:spectral-F-gauge-embedding}
\pFGauge{\mbb{S}} \into \Mod_{\psi(\SphFp)}(\MS_k)
\end{equation}
carries $\iota_*\MSFp$ to $\psienh(\MHFp)$, the endomorphism algebra $\Ext^{*,*}_{\pFGauge{\mbb{S}}}(\MSFp)$ identifies with the syntomic Steenrod algebra introduced previously in \S \ref{sec: syntomic Steenrod operations}.

\begin{cor}\label{cor:Ext_F_p_in_F_Gauge}
The fully faithful embedding \eqref{eq:spectral-F-gauge-embedding} induces an isomorphism
\[
\Ext^{*,*}_{\pFGauge{\mbb{S}}}(\iota_*\MSFp,\iota_*\MSFp) \xrightarrow{\sim} \Asyn^{*,*}.
\]
\end{cor}


\subsection{Comparing $\FGauge{\FF_p}$ with prismatic $F$-gauges}\label{ssec:comparison-F-gauges}

We will now compare the categories $\FGauge{\ZZ_p}$ and $\FGauge{\FF_p}$ with the derived categories of prismatic $F$-gauges arising in the work of Bhatt--Lurie \cite{Bha22}.

\subsubsection{The classical syntomification}\label{sssec:classical-syntomification}
Recall the syntomification $\FSyn{}$ from \S \ref{sssec:generation-prismatic-fgauge}, whose base change to $\F_p$ is $\FSyn{\F_p}$.

\begin{defn}To keep notation from becoming too unwieldy, we shall henceforth abbreviate
\[
\cS := \FSyn{\F_p} \quad \text{and} \quad \cD(\cS) := \cD(\FSyn{\F_p}).
\]
\end{defn}

\subsubsection{Motivic spectra valued in prismatic $F$-gauges}

The first step in relating $\FGauge{\FF_p}\subseteq \Mod_{\MSFp}(\MS_k)$ and $\cD(\cS)$ is to construct a comparison functor between them. The idea is that the functor should be an enhanced version of $X\mapsto \RGamma_{\syn}(X;\F_p)$, valued in the category $\cD(\cS)$ rather than $\cD(\F_p)$.  

Recall that we have constructed the commutative algebra $\MSFp \in \calg(\MS_k)$ using the formalism of \emph{oriented graded algebra} from \S\ref{subsubsec:Oriented_Graded_Algebras}. We need a generalization of this formalism with coefficients in a general target category. 

\begin{variant}
Let $\cD$ be a $p$-complete presentably symmetric monoidal stable $\infty$-category, let $S$ be a qcqs scheme, and let 
\[
\cC_S(\cD):= \PShv_{\Nis,\ebu}(\Sm_S;\cD) \cong \cC_S \otimes \cD
\]
be the category of $\cD$-valued Nisnevich sheaves on $\Sm_S$ satisfying elementary blowup excision. 

Furthermore, let $\MS_S(\cD):= \MS_S \otimes \cD$ be the category of $\cD$-valued motivic spectra over $S$. 

We can define the category $\calg^{\por}(\cC_S(\cD))$ as the category of $\N$-graded commutative algebras $E$ in $\cC_S(\cD)$ together with a map $\rPP^1 \otimes \one_\cD \to E_1$. As in the case where $\cD= \Sp$, such data determines maps $E_i \to \hom(\rPP^1 \otimes \one_\cD,E_{i+1})$ and we let $\calg^{\ori}(\cC_S(\cD)) \subset \calg^{\por}(\cC_S(\cD))$ be the full subcategory spanned by the $E$ for which these maps are isomorphisms. 
Then, using similar arguments to the ones in \S\ref{subsubsec:Oriented_Graded_Algebras}, there is a canonical functor 
\begin{equation}\label{eq:D-valued-nu}
\nu_*\colon \calg^{\ori}(\cC_S(\cD)) \to \calg(\MS_S(\cD))
\end{equation}
turning a $\cD$-valued oriented graded algebra into a commutative algebra in $\MS_S(\cD)$. 
\end{variant}

\subsubsection{Lifting syntomic cohomology to an oriented graded algebra}
We will lift syntomic cohomology to an object of $\calg^{\ori}(\cC_k(\cD(\cS)))$.

Recall that for each smooth proper $f \co X \rightarrow \Spec k$, there is a perfect complex
\[
\cH^X := Rf^{\Syn}_* (\cO_{X^{\Syn}}) \in \Perf(\FSyn{})
\]
and its mod $p$ reduction 
\[
\ol \cH^X := \cH^X \otimes_{\Z_p} \F_p \in \Perf(\cS).
\]

The assignment $X\mapsto X^{\Syn}$ gives a functor 
\[
(-)^{\Syn}\colon \Sm_k^{\op} \to \mathsf{Stacks}_{\FSyn{}}^{\op},
\]
which is symmetric monoidal for the op-Cartesian structure on source and target (this amounts to the fact that the construction $X \mapsto X^{\syn}$ preserves limits, as discussed in \S \ref{sssec:syntomification-of-schemes}). The assignment $X \mapsto \ol{\cH}(X)$ therefore induces a lax symmetric monoidal functor
\[
\ol \cH \colon \Sm_k^\op \to \cD(\cS) 
\]
which then promotes uniquely to a functor $\Sm_k^\op \to \calg(\cD(\cS))$, which we also denote $\ol{\cH}$. This functor satisfies Nisnevich descent and elementary blowup excision, hence gives an object of $\calg(\cC_k(\cD(\cS)))$. 

Let $\cO_{\cS}\BK{1}$ be the \emph{Breuil-Kisin line bundle} on $\cS$ (cf. \cite{Bha22}). As in \S \ref{ssec: promoting}, it induces a graded commutative algebra 
\[
\cO_{\cS}\BK{\bu}[2 \bu] := (\cO_{\cS}\BK{n}[2n])_{n \in \N} \in \calg(\cC_k(\cD(\cS))^\N),
\]
and by tensoring this with $\ol{\cH}$ we obtain the graded algebra 
\[
\ol \cH\BK{\bu}[2\bu]:= \ol\cH \otimes \cO_{\cS}\BK{\bu}[2 \bu] \in \calg(\cC_k(\cD(\cS))^\N). 
\]

Finally, in order to obtain a $\bP^1$-pre-orientation of $\ol \cH\BK{\bu}[2\bullet]$, note that the first Chern class $c_1^{\syn}$ provides a map 
\[
\ol \cH^{\bP^1_k}\to \ol \cH \oplus \ol \cH\BK{-1}[-2] 
\]
which is an isomorphism \cite[Remark 4.3.6]{Bha22}. Hence, we obtain a canonical lift of the pre-orientation of syntomic cohomology, yielding a promotion of $\ol \cH\BK{\bu}[2\bu]$ to an object of $\calg^{\ori}(\cC_k(\cD(\cS))^\N)$.

\subsubsection{Lifting syntomic cohomology to a $\cD(\cS)$-valued motivic spectrum} We can now apply the $\cD(\cS)$-valued version \eqref{eq:D-valued-nu} of the functor $\nu_*$  to $\ol \cH\BK{\bu}[2\bu] \in \calg^{\ori}(\cC_k(\cD(\cS))^\N)$, to obtain a commutative algebra 
\[
\ol\cmH_k := \nu_* \ol \cH\BK{\bu}[2\bu] \in \calg(\MS_k(\cD(\cS))),
\]
which lifts $(\MSFp)_k$ to an object of $\calg(\MS_k(\cD(\cS)))$. It represents the functor $\ol \cH$, in the sense that the mapping spectrum $\hom(-, \ol\cmH_k)$  naturally promotes to a $\cD(\cS)$-valued functor satisfying 
\[
\hom(\Sigma^\infty_+ X, \ol\cmH_k) \cong \ol \cH^X \in \cD(\cS) \quad \text{ for all } X\in \Sm_k.
\]

\subsubsection{The comparison functor} As a general feature of $p$-complete symmetric monoidal stable $\infty$-categories, we have a symmetric monoidal adjunction $\Sp \adj \MS_k$, where the left adjoint is the ``constant spectrum''. The right adjoint is colimit-preserving by \Cref{forg_MS_sheaves} and \Cref{prop:MS_compact_generators}. Hence we may view this adjunction as an adjunction in the 2-category of presentable $\infty$-categories. Tensoring with $\cD(\cS)$, we deduce that there is a symmetric monoidal adjunction 
\begin{equation}\label{eq:constant-gamma-adjunction}
\mathsf{constant} \co \cD(\cS) \adj \MS_k(\cD(\cS)) \co \Gamma
\end{equation}
in which the right adjoint $\Gamma \co \MS_k(\cD(\cS)) \to \cD(\cS)$ is colimit-preserving and lax symmetric monoidal.

The canonical functor $\Sp \rightarrow \cD(\cS)$ induces a functor 
\begin{equation}\label{eq:MS-to-MS(D(S))}
\MS_k = \MS_k(\Sp) \rightarrow \MS_k(\cD(\cS)).
\end{equation}
The image of $\MSFp$ under this functor admits a canonical commutative algebra map to $\ol \cmH_k$, adjoint to the tautological isomorphism to $\MSFp$ from the image of $\ol \cmH_k$ under the right adjoint of \eqref{eq:MS-to-MS(D(S))}. This induces a functor  
\begin{equation}\label{eq:tensor-to-H}
(-)\otimes_{\MSFp} \ol \cmH_k \co \Mod_{\MSFp}(\MS_k) \rightarrow \Mod_{\ol \cmH_k}(\MS_k(\cD(\cS))).
\end{equation}

\begin{defn}\label{defn:Upsilon-Fp}
We define the functor 
\[
\Upsilon\colon \Mod_{\MSFp}(\MS_k) \to \cD(\cS) 
\]
to be the composition 
\[
\Mod_{\MSFp}(\MS_k) \oto{(-)\otimes_{\MSFp} \ol \cmH_k} \Mod_{\ol \cmH_k}(\MS_k(\cD(\cS))) \oto{\Gamma} \cD(\cS),
\]
where the second functor is the restriction to $\ol \cmH_k$-modules of the right adjoint $\Gamma$ from \eqref{eq:constant-gamma-adjunction}. 
\end{defn}

\begin{remark}\label{rem:Upsilon-Zp}
A similar discussion with $\Z_p$ in place of $\F_p$ leads to a functor 
\begin{equation}
\Upsilon \co \Mod_{\MSZp}(\MS_k) \rightarrow \cD(\FSyn{}).
\end{equation}
\end{remark}

By construction, $\Upsilon$ is a composition of colimit-preserving lax symmetric monoidal functors, hence is itself colimit-preserving and lax symmetric monoidal. 
We will show that $\Upsilon$ restricts to a symmetric monoidal equivalence 
$\FGauge{\FF_p} \cong \cD(\cS)$.

\subsubsection{Consequences of Atiyah duality}

Below, for a dualizable object $c$ in a symmetric monoidal category $\cC$ we denote by $c^\vee$ its dual. 

\begin{thm}[Atiyah duality for motivic spectra, \cite{AHI2}]\label{thm:atiyah-duality}
Let $S$ be a qcqs scheme, and let $X$ be a smooth and projective $S$-scheme of relative dimension $n$ over $S$. Then $\Sigma^\infty_+ X$ is dualizable in $\MS_S$, with dual $\Sigma^\infty_+X\langle -T_X\rangle$, the motivic Thom spectrum of the negative tangent bundle on $X$. Furthermore, we have 
\[
\hom(\Sigma_+^\infty X, (\MSFp)_S) \cong (\Sigma_+^\infty X  \otimes (\MSFp)_S)^\vee \cong \Sigma_+^\infty X \otimes (\MSFp)_S[-2n](-n) .
\]
\end{thm}

\begin{proof}
The first part is the Atiyah duality isomorphism of \cite[Theorem 1.1]{AHI2}. The last assertion follows from the Thom isomorphism associated with the $\pic$-orientation of $\SFp(\bu)[2\bu]$, cf. \cite[\S 6]{AHI1}. 
\end{proof}

\begin{remark}\label{rem:upsilon-on-dualizable}
For $M$ dualizable in $\Mod_{\MSFp}(\MS_k)$, the value of $\Upsilon$ on the dual $M^\vee$ is given by
\[
\Upsilon(M^\vee) = \Gamma(M^\vee\otimes_{\MSFp} \ol \cmH_k) \cong \Gamma(\hom_{\MSFp}(M,\ol\cmH_k))  \in \cD(\cS). 
\]
In particular, if $X \in \Sm_k$ with $\Sigma^\infty_+ X$ dualizable in $\MS_k$ (which is satisfied when $X$ is smooth and projective, by \Cref{thm:atiyah-duality}) then  $(\Sigma^\infty_+ X) \otimes \MSFp$ is dualizable in $\Mod_{\MSFp}(\MS_k)$, and we have
\begin{equation}\label{eq:upsilon-on-dualizable}
\Upsilon((\Sigma^\infty_+ X \otimes \MSFp)^\vee ) \cong \ol \cH^X \in \cD(\cS).   
\end{equation}
\end{remark}

We now deduce that the restriction of $\Upsilon$ to $\FGauge{\FF_p}$ is fully faithful.

\begin{prop}\label{prop: classical QCoh}
The functor $\Upsilon$ (from \Cref{defn:Upsilon-Fp}) induces a symmetric monoidal equivalence 
\[
\FGauge{\FF_p} \cong \cD(\cS) = \QCoh(\FSyn{\F_p})
\]
and its $\Z_p$-variant (from \Cref{rem:Upsilon-Zp}) induces a symmetric monoidal equivalence 
\[
\FGauge{\ZZ_p} \cong \QCoh(\FSyn{})
\]
\end{prop}

\begin{proof}
The arguments for the two cases are completely analogous, so we will just write the first one. 

The category $\FGauge{\FF_p}$ is, by definition, generated under colimits by Tate twists of objects of the form $\Sigma^\infty_+X \otimes \MSFp$ for $X$ smooth and projective. Then by \Cref{thm: smooth proper motives generate} the image of $\FGauge{\FF_p}$ under $\Upsilon$ generates $\cD(\cS)$ under colimits. It remains to establish that the (colimit preserving) restriction of the functor $\Upsilon\colon \Mod_{\MSFp}(\MS_k) \to \cD(\cS)$ to the full subcategory $\FGauge{\FF_p}$ is symmetric monoidal and fully faithful.

First we show that $\Upsilon$ is symmetric monoidal. Since $\Upsilon$ preserves colimits and sends Tate twists to Breuil-Kisin twists, in order to check that its lax symmetric monoidal structure is actually symmetric monoidal, it suffices to check that the relevant map is an isomorphism on generators. Noting that by \Cref{thm:atiyah-duality}, $\FGauge{\FF_p}$ is generated by Tate twists of objects of the form $(\Sigma^\infty_+X \otimes \MSFp)^\vee$, what we want is to check that the map
\begin{align}\label{eq:K\"unneth-1} \Upsilon((\Sigma^\infty_+X \otimes \MSFp)^\vee )\otimes \Upsilon((\Sigma^\infty_+Y \otimes \MSFp)^\vee )  \rightarrow  \Upsilon((\Sigma^\infty_+X \otimes \MSFp)^\vee\otimes (\Sigma^\infty_+Y \otimes \MSFp)^\vee )  
\end{align}
is an isomorphism for all smooth projective schemes $X,Y$ over $k$. The source identifies with $\ol\cH^X\otimes \ol\cH^Y$ by  \eqref{eq:upsilon-on-dualizable}. Using  \eqref{eq:upsilon-on-dualizable} again -- in addition to the symmetric monoidality of each operation $\Sigma^\infty_+$, tensoring with $\MSFp$, and formation of the dual (of a dualizable object) -- the target identifies with $\ol\cH^{X\times Y}$. Under these identifications, \eqref{eq:K\"unneth-1} identifies with the K\"unneth isomorphism\footnote{which follows from the fact that $X \mapsto X^{\syn}$ preserves limits, so that 
\[
(X \times_k Y)^{\syn} \cong X^{\syn} \times_{(\Spec k)^{\syn}} Y^{\syn}.
\]} for $\cD(\cS)$-valued syntomic cohomology
\[
\ol\cH^{X\times Y} \cong \ol\cH^X\otimes \ol\cH^{Y}
\]
for smooth projective $X,Y$. 

Next we show that $\Upsilon$ is fully faithful. Since $\Upsilon$ is colimit-preserving and symmetric monoidal, and since our generators for $\FGauge{\FF_p}$ are compact, dualizable, and closed under the symmetric monoidal operation, it suffices to check that $\Upsilon$ induces an isomorphism on the Hom spaces from a generator to the unit. By the further compatibility of $\Upsilon$ with Tate twists, we reduce to checking that the map 
\begin{equation}\label{eq:upsilon-2}
\RHom((\Sigma^\infty_+X\otimes \MSFp),\MSFp) \to \RHom(\Upsilon(\Sigma^\infty_+X\otimes \MSFp),\Upsilon(\MSFp)) 
\end{equation}
is an isomorphism for all smooth projective $X/k$. The source identifies with the syntomic cohomology $\RGamma_{\syn}(X;\F_p)$ since $\MSFp$ is a motivic spectrum representing syntomic cohomology, and the target identifies with the global sections $\RGamma(\cS;\ol\cH^X)$, which again identifies with $\RGamma_{\syn}(X;\F_p)$ since the functor $\ol\cH$ was constructed to lift of mod $p$ syntomic cohomology to $\cD(\cS)$. Moreover, these identifications carry \eqref{eq:upsilon-2} to the identity map, by construction. 
\end{proof}

\subsection{Spectral prismatic $F$-gauges of schemes}\label{ssec:spectral-cohomological-motives}
An important desideratum for the category $\pFGauge{\mbb{S}}$ is the existence of a functor $X \mapsto \sH^X \co \Sm_k \rightarrow \pFGauge{\mbb{S}}$ that ``lifts'' the functor $X \mapsto \cH^X$ from \S \ref{sssec:generation-prismatic-fgauge}, in the sense of a commutative diagram 
\begin{equation}\label{eq:motive-lift}
\begin{tikzcd}
\Sm_k^{\op} \ar[r, "\sH^{(-)}"] \ar[dr, "\cH^{(-)}"']  &\pFGauge{\mbb{S}} \ar[d, "\wt \iota^*"] \\
& \FGauge{\ZZ_p}
\end{tikzcd}
\end{equation}
We will construct this commutative diagram (and in fact, a more general one encompassing qcqs schemes over $k$). 

\subsubsection{Construction of the functor $\sH^{(-)}$} Recall that $\pFGauge{\mbb{S}}$ was defined as a full subcategory of $\Mod_{\psi(\SphFp)}(\MS_k)$. Let 
\begin{equation}\label{eq:prismatic-colocalization-sph}
\cR_{\mbb{S}}\colon \Mod_{\psi(\SphFp)}(\MS_k) \to \pFGauge{\mbb{S}}
\end{equation}
be the right adjoint to the tautological inclusion. Similarly, let 
\begin{equation}\label{eq:prismatic-colocalization-Zp}
\cR_{\FF_p}\colon \Mod_{\MSFp}(\MS_k) \to \FGauge{\FF_p} \quad \text{and} \quad \cR_{\ZZ_p}\colon \Mod_{\MSZp}(\MS_k) \to \FGauge{\ZZ_p}
\end{equation}
be the right adjoints to the tautological inclusions.

Recall that the construction of the equivalence $\FGauge{\ZZ_p} \cong \QCoh(\FSyn{})$ in \S \ref{ssec:comparison-F-gauges} was arranged to intertwine the functor $X \mapsto \cH^X$ with the composite functor
\begin{equation}\label{eq:cH^X-formula}
\Sm_k^\op \oto{\hom_{\MS_k}(\Sigma^\infty_+(-),\MSZp)} \Mod_{\MSZp}(\MS_k) \oto{\cR_{\ZZ_p}} \FGauge{\ZZ_p} . 
\end{equation}
This motivates the following definition.

\begin{defn}
We define the functor $\sH^{(-)}\colon \Sm_k^{\op} \to \pFGauge{\mbb{S}}$ to be the composite functor 
\begin{equation}\label{eq:spherical-Fgauge}
\sH^{(-)}\colon \Sm_k^{\op} \oto{\hom_{\MS_k}(\Sigma^\infty_+(-),\psienh(\Sph))} \Mod_{\psi(\SphFp)}(\MS_k) \oto{\cR_{\mbb{S}}} \pFGauge{\mbb{S}}.
\end{equation}
\end{defn}

Note that $\sH^{(-)}$ is lax symmetric monoidal, being a composition of two lax symmetric monoidal functors. Later in Corollary \ref{cor:sH-symmetric-monoidal}, we will see that it is even symmetric monoidal. 

\begin{example}
If $X$ is smooth and projective over $k$, then $\Sigma^\infty_+ X$ is dualizable by Theorem \ref{thm:atiyah-duality}, and both $\Sigma^\infty_+ X \otimes \psienh(\Sph)$ and its dual in $\Mod_{\psi(\SphFp)}(\MS_k)$ belong to $\pFGauge{\mbb{S}}$. Therefore, in this case we have 
\[
\sH^X \cong (\Sigma^\infty_+ X)^\vee \otimes \psienh(\Sph). 
\]
\end{example}

\subsubsection{} In order to show that the functor $\sH^{(-)}$ fits into a commutative diagram \eqref{eq:motive-lift}, we will show that each of the constituents in the definition of $\sH^{(-)}$ is compatible with base change to $\MSZp$, in a suitable sense. 
The first part is handled by the following Proposition. 

\begin{prop}\label{prop:inner-hom-base-change}
Let $X$ be a smooth scheme over $k$. Then the canonical map 
\[
\hom_{\MS_k}(\Sigma^\infty_+X,\psienh(\Sph))\otimes_{\psienh(\Sph)}\MSZp \to \hom_{\MS_k}(\Sigma^\infty_+X,\MSZp) \in \Mod_{\MSZp}(\MS_k)
\]
is an isomorphism. 
\end{prop}

\begin{proof}Since our categories are $p$-complete, we can check that the map is an isomorphism after reducing mod $p$, where it becomes the natural map 
\begin{equation}\label{prop:inner-hom-base-change-mod-p}
\hom_{\MS_k}(\Sigma^\infty_+X,\psienh(\Sph))\otimes_{\psienh(\Sph)}\MSFp \to \hom_{\MS_k}(\Sigma^\infty_+X,\MSFp) \in \Mod_{\MSFp}(\MS_k).
\end{equation}
For $[n] \in \Delta$, abbreviate $A_n := \psi((\MHFp)^{\otimes (n+1)})$, and recall that by Definition (cf. \Cref{def: simplicial module cat})
\begin{equation}\label{eq:simplicial-limit-A}
\Mod_{\psi(\SphFp)}(\MS_k) := \lim_{[n] \in \Delta} \Mod_{A_n}(\MS_k).
\end{equation}
Accordingly, the functor $\Sigma^\infty_+ X\otimes(-)\colon \Mod_{\psi(\SphFp)}(\MS_k) \to \Mod_{\psi(\SphFp)}(\MS_k)$ can be presented as the limit of the vertical functors in the following diagram:
\begin{equation}
\begin{tikzcd}
    \Mod_{A_0}(\MS_k)
    \ar[r,yshift=1.5pt,->]
    \ar[r,yshift=0pt,<-]
    \ar[r,yshift=-1.5pt,->] \ar[d,"{\Sigma^\infty_+ X \otimes (-)}"] 
    &
    \Mod_{A_1}(\MS_k)
    \ar[r,yshift=3pt,->]
    \ar[r,yshift=1.5pt,<-]
    \ar[r,yshift=0pt,->]
    \ar[r,yshift=-1.5pt,<-]
    \ar[r,yshift=-3pt,->] \ar[d,"{\Sigma^\infty_+ X \otimes (-)}"] 
    &
   \Mod_{A_2}(\MS_k)
    \ar[r,yshift=4.5pt,->]
    \ar[r,yshift=3pt,<-]
    \ar[r,yshift=1.5pt,->]
    \ar[r,yshift=0pt,<-]
    \ar[r,yshift=-1.5pt,->]
    \ar[r,yshift=-3pt,<-]
    \ar[r,yshift=-4.5pt,->] \ar[d,"{\Sigma^\infty_+ X \otimes (-)}"] 
    &
    \cdots   \\
    \Mod_{A_0}(\MS_k)
    \ar[r,yshift=1.5pt,->]
    \ar[r,yshift=0pt,<-]
    \ar[r,yshift=-1.5pt,->] 
    &
    \Mod_{A_1}(\MS_k)
    \ar[r,yshift=3pt,->]
    \ar[r,yshift=1.5pt,<-]
    \ar[r,yshift=0pt,->]
    \ar[r,yshift=-1.5pt,<-]
    \ar[r,yshift=-3pt,->] 
    &
   \Mod_{A_2}(\MS_k)
    \ar[r,yshift=4.5pt,->]
    \ar[r,yshift=3pt,<-]
    \ar[r,yshift=1.5pt,->]
    \ar[r,yshift=0pt,<-]
    \ar[r,yshift=-1.5pt,->]
    \ar[r,yshift=-3pt,<-]
    \ar[r,yshift=-4.5pt,->] 
    &
    \cdots 
\end{tikzcd}
\end{equation}
In these terms, our goal is to compare the right adjoint of the limit of the vertical functors and the right adjoint of the leftmost vertical functor. More generally, we will show that the right adjoint of the limit is obtained from the right adjoints of the vertical functors in the diagram by passage to the limit (and, in particular, that they assemble to a natural transformation of diagrams). \tony{I think there is a disconnect between the previous two sentences and the next one. The previous sentences discuss some general situation while the next one jumps into some Beck--Chevalley condition} \shachar{The reference to the Beck-Chevalley condition is the content of the citation. I made its usage a bit more explicit (but you still have to go to the ref to see why the BC condition is what you'r after). I also fixed a mistake in the argument (reducing to the free rank 1 module implicitly) so this might have been another source of confusion...} To show this, by \cite[Proposition 2.1.7]{ACS19}, it would suffice to show that for each $n$ and each of the structure maps $A_{n} \to A_{n+1}$, the square 
\[
\xymatrix{
\Mod_{A_n}(\MS_k)\ar^{(-)\otimes_{A_n} A_{n+1}}[rr]\ar^{\Sigma^\infty_+ X\otimes (-)}[d] & & \Mod_{A_{n+1}}(\MS_k) \ar^{\Sigma^\infty_+ X\otimes (-)}[d] \\
\Mod_{A_n}(\MS_k)\ar^{(-)\otimes_{A_n}A_{n+1}}[rr] & & \Mod_{A_{n+1}}(\MS_k)
}
\]
is (vertically) right adjointable. In other words, we have to show that the induced 
Beck--Chevalley maps 
\begin{equation}\label{eq:hom-base-change-level-n}
\hom_{\MS_k}(\Sigma^\infty_+X,M)\otimes_{A_n} A_{n+1} \to \hom_{\MS_k}(\Sigma^\infty_+X,M\otimes_{A_n} A_{n+1}).
\end{equation}
obtained by passing to the right adjoints of the vertical functors in these squares  
are isomorphisms for all $M\in \Mod_{A_n}(\MS_k)$.
Tensoring the isomorphism of \Cref{prop:Psi_enh_on_Fp^2} with $\hom_{\MS_k}(\Sigma^\infty_+X,M)$ and $M$ respectively, we can identify \eqref{eq:hom-base-change-level-n} with the assembly map
\[
\bigoplus_{\alpha \in \sI} \hom_{\MS_k}(\Sigma^\infty_+X,M)[p_\alpha](q_\alpha) \to \hom_{\MS_k}\Big(\Sigma^\infty_+X,\bigoplus_{\alpha \in \sI} M[p_\alpha](q_\alpha)\Big).
\]
Hence, it would suffice to show that the functor $\hom_{\MS_k}(\Sigma^\infty_+X,-)\colon \MS_k \to \MS_k$ preserves colimits. Indeed, note that this functor factors as 
\[
\MS_k \oto{f^*}\MS_X \oto{f_*} \MS_k
\]
for $f\colon X\to \Spec k$ the structure morphism. Now, $f^*$ preserves colimits because it is a left adjoint, and $f_*$ preserves colimits by \Cref{prop:push_MS_colim}, so we are done. 
\end{proof}

\subsubsection{} Next we will establish compatibility of the right adjoints \eqref{eq:prismatic-colocalization-sph} and \eqref{eq:prismatic-colocalization-Zp} with base change. We have a commutative diagram 
\[
\begin{tikzcd}
\pFGauge{\mbb{S}} \ar[r, hook]\ar[d, "(-)\otimes_{\psienh(\Sph)} \MSZp"']  & \Mod_{\psi(\SphFp)}(\MS_k) \ar[d, "(-)\otimes_{\psienh(\Sph)} \MSZp"] \\
\FGauge{\ZZ_p}\ar[r, hook] & \Mod_{\MSZp}(\MS_k)
\end{tikzcd}
\]
which induces a Beck--Chevalley natural transformation 
\begin{equation}\label{eq:Beck--Chevalley-tensor}
\cR_{\mbb{S}}(M)\otimes_{\psienh(\Sph)}\MSZp \to \cR_{\ZZ_p}(M\otimes_{\psienh(\Sph)} \MSZp )
\end{equation}
of functors $\Mod_{\psi(\SphFp)}(\MS_k) \rightarrow  \FGauge{\ZZ_p}$. 

\begin{prop}\label{prop:right-adjoint-base-change}
The natural transformation \eqref{eq:Beck--Chevalley-tensor} is an isomorphism. 
\end{prop}
\begin{proof}Since our categories are $p$-complete, we can check that the map is an isomorphism after reducing mod $p$, where it becomes the analogous Beck--Chevalley map
\begin{equation}\label{eq:Beck--Chevalley-tensor-mod-p}
\cR_{\mbb{S}}(M)\otimes_{\psienh(\Sph)}\MSFp \to \cR_{\FF_p}(M\otimes_{\psienh(\Sph)} \MSFp ).
\end{equation}
As in the proof of \Cref{prop:inner-hom-base-change}, abbreviate $A_n := \psi((\MHFp)^{\otimes (n+1)})$ for each $[n] \in \Delta$. Define
\[
\FGauge{A_n} := \Mod_{A_n}(\pFGauge{\mbb{S}}).
\]
Then it follows immediately from the definition of $\pFGauge{\mbb{S}}$ that 
\[
\pFGauge{\mbb{S}} \cong \lim_{[n]\in \Delta} \FGauge{A_n}
\] 
and the functor $\pFGauge{\mbb{S}} \to \Mod_{\psi(\SphFp)}(\MS_k)$ is the limit of the map of simplicial diagrams (where the vertical functors are colimit-preserving)
\begin{equation}
\begin{tikzcd}
    \FGauge{\FF_p}
    \ar[r,yshift=1.5pt,->]
    \ar[r,yshift=0pt,<-]
    \ar[r,yshift=-1.5pt,->] \ar[d, hook] 
    &
    \FGauge{A_1}
    \ar[r,yshift=3pt,->]
    \ar[r,yshift=1.5pt,<-]
    \ar[r,yshift=0pt,->]
    \ar[r,yshift=-1.5pt,<-]
    \ar[r,yshift=-3pt,->] \ar[d, hook] 
    &
   \FGauge{A_2}
    \ar[r,yshift=4.5pt,->]
    \ar[r,yshift=3pt,<-]
    \ar[r,yshift=1.5pt,->]
    \ar[r,yshift=0pt,<-]
    \ar[r,yshift=-1.5pt,->]
    \ar[r,yshift=-3pt,<-]
    \ar[r,yshift=-4.5pt,->] \ar[d, hook] 
    &
    \cdots \ar[d, hook]  \\
    \Mod_{\MSFp}(\MS_k)
    \ar[r,yshift=1.5pt,->]
    \ar[r,yshift=0pt,<-]
    \ar[r,yshift=-1.5pt,->] 
    &
    \Mod_{A_1}(\MS_k)
    \ar[r,yshift=3pt,->]
    \ar[r,yshift=1.5pt,<-]
    \ar[r,yshift=0pt,->]
    \ar[r,yshift=-1.5pt,<-]
    \ar[r,yshift=-3pt,->] 
    &
   \Mod_{A_2}(\MS_k)
    \ar[r,yshift=4.5pt,->]
    \ar[r,yshift=3pt,<-]
    \ar[r,yshift=1.5pt,->]
    \ar[r,yshift=0pt,<-]
    \ar[r,yshift=-1.5pt,->]
    \ar[r,yshift=-3pt,<-]
    \ar[r,yshift=-4.5pt,->] 
    &
    \cdots 
\end{tikzcd}
\end{equation}

 In order to show that \eqref{eq:Beck--Chevalley-tensor-mod-p} is an isomorphism, it suffices by \cite[Proposition 2.1.7]{ACS19} to show that the levelwise Beck--Chevalley transformation is an isomorphism for each $[n] \in \Delta$. Using that $(-)\otimes_{\psienh(\Sph)} \MSFp \cong (-)\otimes_{A_n} A_{n+1}$ on $A_n$-modules (as in the proof of \Cref{prop:inner-hom-base-change}), we can write the map in question as 
\begin{equation}\label{eq:R-Beck--Chevalley-level-n}
\cR_{A_n}(M_n) \otimes_{A_n} A_{n+1} \to \cR_{A_{n+1}}(M_n\otimes_{A_n} A_{n+1}),
\end{equation}
where $M_n:= M\otimes A_n \in \Mod_{A_n}(\MS_k)$.  Note that since $A_{n}$ is an $\MSFp$-algebra, we also have the presentation 
\[
\FGauge{A_n} \cong \Mod_{A_n}(\FGauge{\FF_p}).
\]
By definition of $\FGauge{\FF_p}$, this implies that $\FGauge{A_{n+1}}$ is generated under colimits and shifts by the objects of the form $\Sigma^\infty_+ X \otimes A_{n+1}$ for projective $X\in \Sm_k$, so it suffices to show that \eqref{eq:R-Beck--Chevalley-level-n} is carried to an isomorphism by $\Hom_{A_{n+1}}(\Sigma^\infty_+X\otimes A_{n+1},-)$ for all such $X$. Here we abbreviate
\[	
\Hom_{A_{n+1}}(M,N) := \Hom_{\Mod_{A_{n+1}}(\MS_k)}(M,N).
\]
(If $M,N$ lie in $\FGauge{A_{n+1}}$ then we use the same notation $\Hom_{A_{n+1}}(M,N)$; since the embedding $\FGauge{A_{n+1}} \inj \Mod_{A_{n+1}}(\MS_k)$ is fully faithful, there is never risk of confusion.)

 In other words, we want to check that for all projective $X \in \Sm_k$, the map
\begin{equation}\label{eq:Beck--Chevalley-levelwise}
\Hom_{A_{n+1}}(\Sigma^\infty_+X \otimes A_{n+1},\cR_{A_n}(M_n) \otimes_{A_n} A_{n+1}) \to \Hom_{A_{n+1}}(\Sigma^\infty_+ X\otimes A_{n+1},\cR_{A_{n+1}}(M_n\otimes_{A_n} A_{n+1}))
\end{equation}
is an isomorphism. Using \Cref{prop:Psi_enh_on_Fp^2} to identify $(-)\otimes_{A_n} A_{n+1} \cong \bigoplus_{\alpha \in \sI} (-)[p_\alpha](q_\alpha)$, we may rewrite the LHS as 
\begin{align*}
\Hom_{A_{n+1}}(\Sigma^\infty_+X \otimes A_{n+1},\cR_{A_n}(M_n) \otimes_{A_n} A_{n+1})&\cong  \Hom_{A_{n}}(\Sigma^\infty_+X \otimes A_{n},\cR_{A_n}(M_n) \otimes_{A_n} A_{n+1}) \\ 
& \cong \Hom_{A_{n}}\Big(\Sigma^\infty_+X \otimes A_{n},\bigoplus_{\alpha \in \sI} \cR_{A_n}(M_n)[p_\alpha](q_\alpha) \Big) \\
&\stackrel{(\star)}\cong \bigoplus_{\alpha \in \sI} \Hom_{A_n}(\Sigma^\infty_+X\otimes A_n,\cR_{A_n}(M_n)[p_\alpha](q_\alpha)) \\
&\cong \bigoplus_{\alpha \in \sI} \Hom_{A_n}(\Sigma^\infty_+X\otimes A_n,M_n[p_\alpha](q_\alpha)) \\ &\stackrel{(\star \star)}\cong 
\Hom_{A_{n}}\Big(\Sigma^\infty_+X \otimes A_{n},\bigoplus_{\alpha \in \sI} M_n[p_\alpha](q_\alpha)\Big).
\end{align*}
The only non-formal steps are $(\star)$ and $(\star \star)$, where we commuted the formation of $\Hom$ with an infinite direct sum. The explanation is the same for both, so we will focus on ($\star$). To justify this, note that there is a natural map, and to check that it is an isomorphism we may reduce mod $p$, because our category is $p$-complete. This reduction identifies with the assembly map 
\[
\bigoplus_{\alpha \in \sI} \Hom_{A_n}(\Sigma^\infty_+X\otimes A_n,\cR_{A_n}(M_n)[p_\alpha](q_\alpha)) \otimes \F_p \rightarrow \Hom_{A_{n}}\Big(\Sigma^\infty_+X \otimes  A_{n},\bigoplus_{\alpha \in \sI} \cR_{A_n}(M_n)[p_\alpha] (q_\alpha)\Big) \otimes \F_p.
\]
Up to a shift, tensoring the $\Hom$ with $\F_p$ is equivalent to replacing $\Sigma^\infty_+X\otimes A_n$ with $\Sigma^\infty_+X/p\otimes A_n$. But this object is compact (indeed, $\Sigma^\infty_+X/p$ is already compact in $\MS_k$), so the map 
\[
\bigoplus_{\alpha \in \sI} \Hom_{A_n}(\Sigma^\infty_+X/p\otimes A_n,\cR_{A_n}(M_n)[p_\alpha](q_\alpha))  \rightarrow \Hom_{A_{n}}\Big(\Sigma^\infty_+X /p\otimes A_{n},\bigoplus_{\alpha \in \sI} \cR_{A_n}(M_n)[p_\alpha] (q_\alpha)\Big) 
\]
is an isomorphism. 
Similarly, for the RHS of \eqref{eq:Beck--Chevalley-levelwise} we have 
\begin{align*}
\Hom_{A_{n+1}}(\Sigma^\infty_+X \otimes A_{n+1},\cR_{A_{n+1}}(M_n \otimes_{A_n} A_{n+1})) &\cong \Hom_{A_{n+1}}(\Sigma^\infty_+X \otimes A_{n+1},M_n \otimes_{A_n} A_{n+1})\\ & \cong \Hom_{A_{n}}(\Sigma^\infty_+X \otimes A_{n},M_n \otimes_{A_n} A_{n+1}) \\ &\cong 
\Hom_{A_{n}}\Big(\Sigma^\infty_+X \otimes A_{n},\bigoplus_{\alpha \in \sI} M_n[p_\alpha](q_\alpha)\Big) 
.    
\end{align*}
Moreover, these identifications carry \eqref{eq:Beck--Chevalley-levelwise} to the identity map, as desired.
\end{proof}

\begin{cor}\label{cor:sH-X-lifting}
For all $X\in \Sm_k$ the natural map 
$
\wt\iota^*\sH^X \rightarrow \cH^X \in \cD(\FSyn{})
$
is an isomorphism. 
\end{cor}

\begin{proof}
We have 
\begin{align*}
\wt\iota^*\sH^X&:= \MSZp\otimes_{\psienh(\Sph)}  \sH^X = \MSZp\otimes_{\psienh(\Sph)} \cR_{\mbb{S}}(\hom(\Sigma^\infty_+X,\psienh(\Sph))) \\
&\stackrel{(1)}\cong \cR_{\ZZ_p}(\MSZp \otimes_{\psienh(\Sph)}\hom(\Sigma^\infty_+ X, \psienh(\Sph))) \\
&\stackrel{(2)}\cong \cR_{\ZZ_p}(\hom(\Sigma^\infty_+X,\MSZp)) \cong \cH^X 
\end{align*}
where step (1) is \Cref{prop:right-adjoint-base-change} and step (2) is \Cref{prop:inner-hom-base-change}.
\end{proof}

\begin{remark}
The formulas \eqref{eq:cH^X-formula} for $\cH^{(-)}$ and \eqref{eq:spherical-Fgauge} for $\sH^{(-)}$ make sense for any qcqs scheme $X/k$, and carry cosifted limits to cosifted limits. Therefore, the same formulas define functors from the category $\Sch_k^{\op}$ of qcqs schemes over $k$, which are Kan-extended from $\Sm_k^{\op}$. Thus we obtain more generally a commutative triangle
\begin{equation}\label{eq:motive-lift-qcqs}
\begin{tikzcd}
\Sch_k^{\op} \ar[r, "\sH^{(-)}"] \ar[dr, "\cH^{(-)}"']  &\pFGauge{\mbb{S}} \ar[d, "\wt \iota^*"] \\
& \FGauge{\ZZ_p}
\end{tikzcd}
\end{equation}
\end{remark}

\begin{cor}\label{cor:sph_syn_proper_dualizable}
Let $X$ be smooth and proper over $k$. Then $\sH^X$ is dualizable in $\pFGauge{\mbb{S}}$. 
\end{cor}

\begin{proof}
Since $\pFGauge{\mbb{S}}\cong \lim_{[n]\in \Delta} \FGauge{A_n}$ we can check the dualizability of $\sH^X$ by base-changing to each of the $A_n$'s. Thus it would suffice to show that $\sH^X \otimes A_n$ is dualizable in $\FGauge{A_n}$ for all $[n] \in \Delta$. 

Since $\sH^X \otimes A_n \cong \sH^X\otimes A_0 \otimes_{A_0} A_n$, it would suffice to show that $\sH^X\otimes A_0 = \sH^X\otimes \MSFp$ is dualizable in $\FGauge{\FF_p}$.  Thanks to \Cref{cor:sH-X-lifting}, the equivalence $\FGauge{\FF_p}\cong \cD(\cS)$ intertwines $\sH^X\otimes \MSFp$ with $\ol\cH^X\in \cD(\cS)$. The latter object is dualizable by the syntomic Poincar\'e duality of \cite[Theorem 4.5.3]{Bha22}.
\end{proof}

One motivation for Bhatt--Lurie to construct $\FSyn{}$ was to ``salvage'' the K\"unneth formula, in the sense that $\cH^{(-)}$ is symmetric monoidal. This statement is refined by the symmetric monoidality of $\sH^{(-)}$, which we now prove. 

\begin{cor}\label{cor:sH-symmetric-monoidal}
The functor $\sH^{(-)}\colon \Sm_k^\op \to \pFGauge{\mbb{S}}$ is symmetric monoidal. In particular, for every $X,Y\in \Sm_k$ we have 
\begin{equation}\label{eq:sH-symmetric-monoidal}
\sH^X\otimes \sH^Y\iso \sH^{X\times Y}.
\end{equation}
\end{cor}

\begin{proof}
We have already equipped $\sH^{(-)}$ with a lax symmetric monoidal structure, so it suffices to show that the comparison map \eqref{eq:sH-symmetric-monoidal} is an isomorphism. Since $\iota^*\colon \pFGauge{\mbb{S}} \to \FGauge{\FF_p}$ is conservative and symmetric monoidal, it suffices to show that $\iota^*\sH^{(-)}$ is symmetric monoidal. Now \Cref{cor:sH-X-lifting} identifies this functor with $\ol\cH^{(-)}\colon \Sm_k^\op \to \FGauge{\FF_p}$, which is symmetric monoidal, as explained in the proof of \Cref{prop: classical QCoh}. 
\end{proof}

\section{Prismatization of syntomic Steenrod operations}\label{sec: prismatization of Steenrod}

Prismatization lifts the syntomic cohomology of $X$ to a prismatic $F$-gauge $\ol \cH^X \in \FGauge{\FF_p}$. In this section, we will prismatize the syntomic Steenrod operations. Concretely, this means lifting the syntomic Steenrod algebra $\Asyn^{*,*}$ to an algebra $\sAsyn \in \pFGauge{\mbb{S}}$, and also lifting the action of $\Asyn^{*,*}$ on syntomic cohomology to an action of $\sAsyn$ on $\ol \cH^X$.

\subsection{The spectral prismatization}
We continue to abbreviate 
$\cS := \FSyn{\F_p}$. We constructed in Proposition \ref{prop: classical QCoh} a symmetric monoidal equivalence $\cD(\cS) \cong \FGauge{\FF_p}$ of stable $\infty$-categories. 

We introduce the abbreviation $\cD(\Splus) := \pFGauge{\mbb{S}}$. The notation suggests that $\cD(\Splus)$ should be viewed as the category of quasicoherent sheaves on a spectral stack $\Splus$.\footnote{At this time, we are only willing to assert this for the ``true'' category $\FGauge{\mbb{S}}$, cf. Remark \ref{rem:true-fgauge}. Again, it makes little difference for our purposes in this paper, because we are interested in eventually connective objects.} We will not (in this paper) invoke the geometricity of the stack $\Splus$, so this is mainly a notational device which lends geometric intuition and notation to various functors that we will work with, e.g., the natural symmetric monoidal functor $\iota^* \co \cD(\Splus) \rightarrow \cD(\cS)$ from \Cref{prop:iota_affine} should be imagined as pullback along a map $\iota \co \cS \rightarrow \Splus$. 

Thus we have identifications 
\begin{equation}\label{eq: prismatization 1}
\cD(\Splus) \cong \pFGauge{\mbb{S}} \quad \text{and} \quad \cD(\cS) \cong \FGauge{\FF_p}
\end{equation} 
by tautology in the first instance, and 
\Cref{prop: classical QCoh} in the second. We transport the functors from \Cref{prop:iota_affine} to adjoint functors $\iota^* \co \cD(\Splus) \adj \co  \cD(\cS) \co \iota_*$. 

Recall from \Cref{cor:fgauge-module-category} that $\iota_*\colon \cD(\cS) \to \cD(\Splus)$ can be identified with the forgetful functor from the module category over the commutative algebra $\MSFp$. In particular, it preserves all limits and colimits, and it is linear over $\cD(\Splus)$. Part of \Cref{prop:iota_affine} gives the projection formula  
\begin{equation}\label{eq:projection_formula_for_iota}
\iota_*(\cF \otimes \iota^*\sG)\cong (\iota_*\cF) \otimes \sG  \quad \text{for all $\cF\in \cD(\cS)$ and $ \sG\in \cD(\Splus).$}
\end{equation}
 We will interpret the (dual) syntomic Steenrod algebra in terms of quasicoherent sheaves on the stacks $\cS, \Splus$. 

\subsection{Internal Hom}  
The $\infty$-category 
$\Mod_{\psi(\SphFp)}(\MS_k)$ is presentably symmetric monoidal, being a limit of such categories along symmetric monoidal colimit-preserving functors. The full subcategory $\cD(\Splus) \subseteq \Mod_{\psi(\SphFp)}(\MS_k)$ is a symmetric monoidal subcategory closed under colimits, hence itself presentably symmetric monoidal. It follows from the Adjoint Functor Theorem \cite[Corollary 5.5.2.9]{HTT} that it is closed symmetric monoidal, i.e., admits internal Hom objects such that $\cRHom(X,-)$ is right adjoint to $X\otimes (-)$. 

\begin{defn}
We denote the internal Hom functor of $\pFGauge{\mbb{S}}=: \cD(\Splus)$ by $\cRHom_{\Splus}(-,-)$ and the one of $\FGauge{\FF_p} \cong \cD(\cS)$ by $\cRHom_{\cS}(-,-)$.
\end{defn}

\subsection{Prismatization of the syntomic Steenrod algebra} 

We write $\cO_{\cS} \in \cD(\cS)$ and $\cO_{\Splus} \in \cD(\Splus)$ for the units of these symmetric monoidal categories. In terms of \eqref{eq: prismatization 1}, we have $\cO_{\cS} \leftrightarrow \psienh(\MHFp) \in \FGauge{\FF_p}$ and $\cO_{\Splus} \leftrightarrow \psienh(\Sph) \in \pFGauge{\mbb{S}}$. 

\begin{defn}[Prismatized syntomic Steenrod algebra] Let
\[
\sAsyn := \cRHom_{\Splus}(\iota_* \cO_{\cS}, \iota_* \cO_{ \cS}),
\]
a priori considered as an associative algebra over $\iota_* \cO_{\cS}$ in $\mrm{Alg}(\cD(\Splus))$. 

\end{defn}

\begin{notation}\label{not: global sections}
Note that
\[
\RGamma(\cS; -) := \RHom_{\cD(\cS)}(\cO_{\cS}, -) \cong \RHom_{\FGauge{\FF_p}}(\psienh(\MHFp), -)
\]
as a functor $\cD(\cS) \rightarrow \cD(\F_p)$. Motivated by this, we formally define 
\[
\RGamma(\Splus; -) := \RHom_{\cD(\Splus)}(\cO_{\Splus}, -) = \RHom_{\pFGauge{\mbb{S}}}(\psienh(\Sph), -)
\]
as a functor $\cD(\Splus) \rightarrow \Sp$. When this functor is evaluated on sheaves pushed forward from $\cD(\cS)$ via $\iota_*$, we may regard it as factoring through $\cD(\F_p) \rightarrow \Sp $. 

We define $\rH^*(\cS; -) := \rH^*(\RGamma(\cS; -))$ and $\rH^*(\Splus; -)$ analogously. For $\cF \in \cD(\cS)$ and $b \in \Z$, we write $\rH^{a,b}(\cS; \cF) := \rH^a(\cS; \cF\BK{b})$ and
\[
\rH^{a,*}(\cS; \cF) := \bigoplus_{b \in \Z} \rH^{a}(\cS; \cF \BK{b}) \quad \text{ and } \quad \rH^{*,*}(\cS; \cF) := \bigoplus_{a \in \Z} \rH^{a,*}(\cS; \cF).
\]
Note that if $\cF \in \Perf(\cS)$, then it follows from \cite[Proposition 4.4.3]{Bha22} that $\rH^{*}(\cS; \cF \BK{b}) = 0$ for all but finitely many $b \in \Z$. 

We analogously define $\rH^{*,*}(\Splus; \sF)$ for $\sF \in \cD(\Splus)$. 
\end{notation}

From \eqref{eq: prismatization 1} and \Cref{cor:Ext_F_p_in_F_Gauge}, we obtain an identification 
\[
\Asyn^{*,*} = \rH^{*,*}(\Splus; \sAsyn) 
\]
of bigraded associative algebras over $\Hsynpt = \rH^{*,*}(\cS; \cO_{\cS})$ (later to be upgraded to an identification of bigraded cocommutative Hopf algebras). Thus we view $\sAsyn$ as the prismatization of the syntomic Steenrod algebra.

\subsection{Prismatization of the dual syntomic Steenrod algebra} We translate some of the discussion from \S\ref{ssec: syntomic dual Steenrod algebra} into the prismatized language. 

\begin{prop}\label{prop: coherent psi of tensor} The symmetric monoidal equivalence \eqref{eq: prismatization 1} carries $
 \psienh(\MHFp \otimes \MHFp) \in \pFGauge{\mbb{S}}$ to 
$\iota_* \cO_{\cS} \otimes_{\cO_{\Splus}} \iota_* \cO_{\cS} \in \cD(\Splus)$.
\end{prop}

\begin{proof}
By construction, \eqref{eq: prismatization 1} sends 
\begin{align*}
\psienh(\MHFp)  &\mapsto \iota_* \cO_{\cS} \\
\psienh(\Sph) &\mapsto \cO_{\Splus}
\end{align*}
Hence Proposition \ref{prop:Psi_enh_on_Fp^2} translates into the desired statement.
\end{proof}

From  Proposition \ref{prop: coherent psi of tensor} and \eqref{eq: tensor split}, we obtain a splitting
\begin{equation}\label{eq: coherent tensor split}
\iota_*  \cO_{ \cS} \otimes_{\cO_{\Splus}} \iota_*  \cO_{ \cS} \cong \bigoplus_{\alpha \in \sI} \iota_* \cO_{\cS} \, \xi_\alpha \in \cD(\Splus),
\end{equation} 
where $\xi_\alpha$ has cohomological degree $-p_\alpha$ and twist $q_\alpha$, so that $\cO_{\cS} \, \xi_\alpha \cong  \cO_{\cS} [p_\alpha]\BK{q_\alpha}$.

\begin{defn}[Prismatized dual syntomic Steenrod algebra] Under our identifications, the object $\dsAsyn$ from \Cref{defn:dual-sAsyn} may be viewed as 
\[
\dsAsyn = \iota^* \iota_*  \cO_{ \cS} \in \cD(\cS).
\]
This is a commutative Hopf algebra\footnote{Using \Cref{lem: dual Hopf algebroid is algebra} to see that the a priori Hopf algebroid structure refines to a Hopf algebra structure.} over $\cO_{ \cS}$ in $\cD(\cS)$, with a natural isomorphism
\[
\iota_* \dsAsyn \cong \iota_*  \cO_{ \cS} \otimes_{\cO_{\Splus}} \iota_*  \cO_{ \cS} \cong \bigoplus_{\alpha \in \sI} \iota_* \cO_{\cS} \, \xi_\alpha.
\]
\end{defn}

Note that by adjunction,
\begin{equation}\label{eq: dual steenrod adjunction}
\iota_* \cRHom_{\cS} (\dsAsyn,  \cO_{\cS}) \cong \cRHom_{\Splus} (\iota_* \cO_{\cS}, \iota_*  \cO_{ \cS}) \cong \sAsyn,
\end{equation}
a priori as $\iota_* \cO_{\cS}$-modules, and then as associative $\iota_* \cO_{\cS}$-algebras. More colloquially, the $\cO_{ \cS}$-dual of $ \dsAsyn$ is $\sAsyn$. Using this, we transfer the commutative Hopf algebra structure on $\dsAsyn$ to the dual structure on $\sAsyn$.

Although the global sections functors (cf. \Cref{not: global sections}) are not symmetric monoidal, they are compatible with tensor products specifically on $\sAsyn$ and $\dsAsyn$, as articulated below.

\begin{lemma}\label{lem: sAsyn cup product}
(1) The cup product 
\begin{equation}\label{eq: dsAsyn cup product}
\rH^{*,*}(\cS, \dsAsyn) \otimes_{\Hsynpt} \rH^{*,*}(\cS, \dsAsyn)  \rightarrow  
\rH^{*,*}(\cS, \dsAsyn \otimes_{\cO_{\cS}} \dsAsyn)
\end{equation}
is an isomorphism.  

(2) The cup product 
\begin{equation}\label{eq: sAsyn cup product}
\rH^{*,*}(\Splus, \sAsyn) \otimes_{\Hsynpt} \rH^{*,*}(\Splus, \sAsyn)  \rightarrow  
\rH^{*,*}(\Splus, \sAsyn \otimes_{\iota_* \cO_{\cS}} \sAsyn)
\end{equation} 
is an isomorphism. 
\end{lemma}

\begin{proof}
(1) Using \eqref{eq: coherent tensor split}, we have 
\begin{align*}
\iota_*(\dsAsyn \otimes_{\cO_{\cS}} \dsAsyn) & \cong \iota_* \dsAsyn \otimes_{\iota_* \cO_{\cS}} \iota_* \dsAsyn  \\
& \cong \left(\bigoplus_{\alpha \in \sI} \iota_* \cO_{\cS} \, \xi_\alpha\right) \otimes_{\iota_* \cO_{\cS}} \left(\bigoplus_{\alpha' \in \sI} \iota_* \cO_{\cS} \, \xi_{\alpha'} \right) \\
& \cong \bigoplus_{\alpha, \alpha' \in \sI} \iota_* \cO_{\cS} \, \xi_\alpha  \xi_{\alpha'}.
\end{align*}
Hence we have
\begin{align*}
\rH^{*,*}(\cS; \dsAsyn \otimes_{\cO_{\cS}} \dsAsyn) & \cong \rH^{*,*}(\Splus; \iota_* \dsAsyn \otimes_{\iota_* \cO_{\cS}} \iota_* \dsAsyn)  \\
& \cong \rH^{*,*}(\Splus; \bigoplus_{\alpha, \alpha' \in \sI} \iota_* \cO_{\cS} \, \xi_\alpha \xi_{\alpha'} ).
\end{align*}
Therefore, with respect to the identifications
\[
\rH^{*,*}(\Splus; \iota_* \cO_{\cS}\, \xi_\alpha) \cong \Hsynpt \, \xi_\alpha \quad \text{and} \quad 
\rH^{*,*}(\Splus; \iota_* \cO_{\cS}\, \xi_{\alpha'}) \cong \Hsynpt \, \xi_{\alpha'}
\]
the map \eqref{eq: dsAsyn cup product} reads\footnote{using that only finitely many twists contribute, or alternatively that cohomology commutes with direct sums in general on $\cS$.}
\begin{align*}
\rH^{*,*}(\cS, \dsAsyn) \otimes_{\Hsynpt} \rH^{*,*}(\cS, \dsAsyn) & \cong \bigoplus_{\alpha \in \sI} \Hsynpt \, \xi_\alpha \otimes_{\Hsynpt}  \bigoplus_{\alpha' \in \sI} \Hsynpt \, \xi_{\alpha'} \\
& \rightarrow \bigoplus_{\alpha, \alpha' \in \sI} \Hsynpt \, \xi_\alpha \xi_{\alpha'} \cong \rH^{*,*}(\cS, \dsAsyn \otimes_{\cO_{\cS}} \dsAsyn),
\end{align*}
which is visibly an isomorphism. 

(2) Follows from a similar argument. 
 
\end{proof}

From \eqref{eq: prismatization 1} and \Cref{lem: sAsyn cup product}, we obtain an identification 
\[
\dAsyn_{*,*} = \rH^{-*,-*}(\cS; \dsAsyn).
\]
of bigraded commutative Hopf algebras over $\Hsynpt = \rH^{*,*}(\cS)$. Thus we view $\dsAsyn$ as the prismatization of the dual syntomic Steenrod algebra.

\subsection{Prismatization of Steenrod actions} 
We already know that $\Asyn^{*,*}$ acts on the syntomic cohomology groups of a variety with coefficients in $\F_p^{\syn}(\bu)$. We now explore more refined structure that can be articulated on the prismatization.

\subsubsection{Action on sheaves pulled back via $\iota^*$}\label{subsec:action_on_i_star}
We will construct a tautological action of $\Asyn^{*,*}$ on the syntomic cohomology of any sheaf $\cF \in  \cD(\cS)$ of the form $\iota^* \sF$ for $\sF \in \cD(\Splus)$. In fact, this arises from the more refined structure of an action of $\sAsyn$ on $\cF$ in $\cD(\cS)$. Informally speaking, the slogan is that \emph{``$\sAsyn$ naturally acts on a sheaf that admits a spectral lift''}. This is a special case of the structural pattern discussed in \cite[\S 4.2]{SHA}. 

Let $\sF \in \cD(\Splus)$ and $\cF : = \iota^* \sF \in \cD(\cS)$. Then we may write 
\[
\iota_* \cF \cong \iota_* \iota^* \sF \cong \sF \otimes_{\cO_{\Splus}} \iota_* \cO_{\cS} \in \cD(\Splus)
\]
Via the tautological action of $\sAsyn$ on $\iota_* \cO_{\cS}$, this induces a tautological action of $\sAsyn$ on $\iota_* \cF$. Taking cohomology, we get a (bigraded) action of $\Asyn^{*,*}$ on 
\[
\rH^{*,*}(\Splus; \iota_* \cF) \cong \rH^{*,*}(\cS; \cF). 
\]
By construction, this recovers the action defined earlier in \S \ref{sec: syntomic Steenrod operations}. 

\subsubsection{Action on (symmetric monoidal) duals of pullbacks}\label{subsec:action_on_duals} Recall the notion of \emph{dualizable object} in a symmetric monoidal category $\msf{C}$. For a dualizable object $c \in \msf{C}$, the dual will be denoted $c^\vee$. 

Let $\sF \in \cD(\Splus)$ and $\cF := \iota^* \sF \in \cD(\cS)$. Then 
\begin{equation}\label{eq:dual-asyn-action}
\cRHom_{\cS}(\cF, \cO_{\cS}) \cong \cRHom_{\cS}(\iota^* \sF, \cO_{\cS}) \cong \cRHom_{\Splus}(\sF, \iota_* \cO_{\cS})
\end{equation}
also has a tautological action of $\sAsyn$ (through its action on the target). 

If $\sF$ is dualizable in $\cD(\Splus)$, then $\cF$ is dualizable in $\cD(\cS)$ (since it is the image of a dualizable object under a symmetric monoidal functor) and we have
\[
\iota_* \cF^\vee \cong \iota_* \cRHom_{\cS}(\cF, \cO_{\cS}) \cong \cRHom_{\Splus}(\sF, \iota_* \cO_{\cS}) \cong \sF^\vee \otimes_{\cO_{\Splus}} \iota_* \cO_{\cS}.
\]
This isomorphism is equivariant for the $\sAsyn$-actions, which for the left term was defined via \eqref{eq:dual-asyn-action}, and for the right term is induced by the action on $\iota_*\cO_{\cS}$. Taking global sections, we obtain an action of $\Asyn^{*,*}$ on $\RGamma(\cS; \cF^\vee)$ for any such $\cF$.

\subsection{Coproduct}\label{ssec: coproduct}

The product on the commutative Hopf algebra $\dsAsyn$ is dual to a coproduct
\begin{equation}\label{eq: sAsyn coproduct}
\sAsyn \rightarrow \sAsyn \otimes_{\iota_* \cO_{\cS}} \sAsyn.
\end{equation}
After applying $\rH^{*,*}(-)$ to \eqref{eq: sAsyn coproduct}, and using \Cref{lem: sAsyn cup product}, we obtain the coproduct on $\Asyn^{*,*}$, which is identified explicitly by the Cartan formula of Proposition \ref{prop: Cartan}. 

In particular, the coproduct $\Asyn^{*,*} \rightarrow \Asyn^{*,*} \otimes_{\Hsynpt} \Asyn^{*,*}$ equips the category of $ \Asyn^{*,*}$-modules with a monoidal structure $(-)\otimes_{\Hsynpt}(-)$. Similarly, the coproduct $\rAsyn^{*,*} \rightarrow \rAsyn^{*,*} \otimes_{\F_p} \rAsyn^{*,*}$ equips the category of $\rAsyn^{*,*}$-modules with a monoidal structure $(-)\otimes_{\F_p}(-)$.

For $\cF, \cG \in \cD(\cS)$, we have a cup product map 
\begin{equation}\label{eq: Fsyn cup product}
\rH^{*,*}(\cS;\cF)\otimes_{\Hsynpt} 
\rH^{*,*}(\cS;\cG) \to 
\rH^{*,*}(\cS;\cF\otimes_{\cO_{\cS}} \cG).
\end{equation}
Suppose $\cF= \iota^*\sF$ and $\cG = \iota^*\sG$ for $\sF, \sG \in \cD(\Splus)$. Then the source has a canonical $\Asyn^{*,*}$-action induced by the action on each factor (\S \ref{subsec:action_on_i_star}) and the coproduct, while the target has a canonical $\Asyn^{*,*}$-action via the presentation $\cF \otimes_{\cO_{\cS}} \cG \cong \iota^* (\sF \otimes_{\cO_{\Splus}} \sG)$.

\begin{lemma}\label{Cartan} 
Let $\sF, \sG \in \cD(\Splus)$ and $\cF := \iota^* \sF, \cG := \iota^* \sG \in \cD(\cS)$. Then the cup product \eqref{eq: Fsyn cup product} is $\Asyn^{*,*}$-equivariant with the natural actions described above. 
\end{lemma}

\begin{proof}
We have two actions (in the sense of groupoids) of $\sAsyn$ on $\iota_* (\cF \otimes_{\cO_{\cS}} \cG)$: one is the tautological action coming from the isomorphisms 
\[
\iota_* (\cF \otimes_{\cO_{\cS}} \cG) \cong \iota_* \iota^* (\sF \otimes_{\cO_{\Splus}} \sG) \cong (\sF \otimes_{\cO_{\Splus}} \sG) \otimes_{\cO_{\Splus}} \iota_* \cO_{\cS}
\]
and the other coming from the coproduct \eqref{eq: sAsyn coproduct} composed with the natural action of $\sAsyn \otimes_{\iota_* \cO_{\cS}} \sAsyn$ on $\iota_* (\cF \otimes \cG) \cong \iota_* \cF \otimes_{\iota_* \cO_{\cS}} \iota_* \cG$ obtained by tensoring the tautological actions. We claim that these two actions are canonically identified; from this the Lemma follows immediately. 

Writing 
\begin{align*}
\iota_* \cF \otimes_{\iota_* \cO_{\cS}} \iota_* \cG & \cong (\sF  \otimes_{\cO_{\Splus}} \iota_* \cO_{\cS}) \otimes_{\iota_* \cO_{\cS}} (\sG  \otimes_{\cO_{\Splus}} \iota_* \cO_{\cS}) \\
&\cong \sF  \otimes_{\cO_{\Splus}} \left(\iota_* \cO_{\cS}  \otimes_{\iota_* \cO_{\cS}} \iota_* \cO_{\cS}  \right)  \otimes_{\cO_{\Splus}} \sG 
\end{align*}
we see that the claim follows from the statement that the two actions of $\sAsyn = \cRHom_{\Splus}(\iota_* \cO_{\cS}, \iota_* \cO_{\cS})$ on $(\iota_* \cO_{\cS} \otimes_{\iota_* \cO_{\cS}} \iota_* \cO_{\cS})$ coincide: one via the identification $(\iota_* \cO_{\cS} \otimes_{\iota_* \cO_{\cS}} \iota_* \cO_{\cS}) \cong \iota_* \cO_{\cS}$, and the other by the coproduct composed with the tensor product of the natural actions. This last statement amounts to the tautological compatibility of the coproduct with the unit. 
\end{proof}

\begin{remark}\label{dual cartan}
There is also an analogous compatibility for the cup product
\[
\rH^{*,*}(\cS; \cF^\vee) \otimes_{\Hsynpt} \rH^{*,*}(\cS; \cG^\vee) \rightarrow \rH^{*,*
}(\cS; \cF^\vee \otimes_{\cO_{\cS}} \cG^\vee)
\]
where $\cF = \iota^* \sF $ and $\cG = \iota^* \sG$, so that there are Steenrod actions by \S \ref{subsec:action_on_duals}. We will only invoke the compatibility for $\sF$ and $\sG$ which are dualizable, in which case it follows from \Cref{Cartan}, so we omit the proof of the more general statement.
\end{remark}

\section{Spectral Serre duality and Steenrod equivariance}\label{sec: Steenrod equivariance}

\subsection{Steenrod equivariance for arithmetic duality}\label{ssec: duality} 

The goal of this section is to prove the compatibility statement from \S \ref{intro: spectral prismatization}. We will reformulate it slightly, keeping the notation there. Let $X$ be a smooth, proper, geometrically connected variety over a characteristic $p$ finite field $k$. We dualize the cup product 
\[
\Hsyn^{*,*}(X) \otimes_{\F_p} \Hsyn^{*,*}(X) \rightarrow \Hsyn^{*,*}(X \times_{k} X)
\]
over $\F_p$, and apply Poincar\'e duality to obtain a commutative diagram 
\begin{equation}\label{eq: varphi diagram}
\begin{tikzcd}
\Hsyn^{*,*}(X)^\vee \otimes_{\F_p} \Hsyn^{*,*}(X)^\vee \ar[d, "\wr", "\text{Poincar\'e}"']  & \ar[l] \Hsyn^{*,*}(X \times_{k} X)^\vee \ar[d, "\wr","\text{Poincar\'e}"']  \\
\Hsyn^{*,*}(X) \otimes_{\F_p} \Hsyn^{*,*}(X) & \ar[l, "\varphi_*"] \Hsyn^{*,*}(X \times_{k} X)
\end{tikzcd}
\end{equation} The top horizontal arrow preserves the bi-grading; the other maps do not, but we record that $\varphi_*$ increases the bidegree by $+(1,0)$. 

In the bottom row of \eqref{eq: varphi diagram}, both the source and target have natural actions of $\rAsyn^{*,*}$, the target by the coproduct (recall \S \ref{ssec: coproduct}).

\begin{thm}\label{thm: steenrod push equivariant} The map $\varphi_*$ from \eqref{eq: varphi diagram} is equivariant with respect to the action of $\rAsyn^{*,*}$.
\end{thm}

The proof of Theorem \ref{thm: steenrod push equivariant} will occupy the rest of the section. It will be long, so let us give a high-level overview. \emph{Throughout this section, we continue to use the notation $\cS := \FSyn{\F_p}$, and 
\[
\iota^* \co \cD(\Splus) \rightleftarrows \cD( \cS) \co \iota_* 
\]
as in \S \ref{sec: prismatization of Steenrod}.}
\begin{enumerate}
\item Firstly, in \S \ref{ssec: prismatization of varphi} we localize \Cref{thm: steenrod push equivariant} onto $\cS$: we formulate a ``prismatization of $\varphi_*$'' as a map $\varphi^{\prism}$ of sheaves on $\cS$, whose compatibility with the prismatized Steenrod action recovers \Cref{thm: steenrod push equivariant} upon taking global sections. 

\item The main input to defining $\varphi^{\prism}$ is Serre duality on $\cS$. Since $\sAsyn$ is defined in terms of endomorphisms of $\cS$ over $\Splus$, the desired compatibility will ultimately come from the fact that Serre duality itself lifts to $\cD(\Splus)$. 

\item In \S \ref{ssec: spectral Serre duality}, we formulate and prove a form of coherent duality for $\cD(\Splus)$ that we call \emph{spectral Serre duality}. The key idea is that this should implement Brown--Comenetz duality on spectral syntomic cohomology. 

\item Finally, in \S \ref{ssec: 11.4} and \S \ref{ssec: 11.5} we study the interaction between spectral Serre duality and $\sAsyn$, and prove the prismatized version of Theorem \ref{thm: steenrod push equivariant}.
\end{enumerate}


\subsection{Prismatization of $\varphi$}\label{ssec: prismatization of varphi}

The map
\begin{equation}\label{eq: varphi}
\varphi_* \colon \Hsyn^{*,*}(X \times_{k} X) \rightarrow \Hsyn^{*,*}(X) \otimes_{\F_p} \Hsyn^{*,*}(X)
\end{equation}
was constructed above in \eqref{eq: varphi diagram} using Poincar\'e duality for syntomic cohomology. However, it admits another description in terms of prismatization, based on a prismatization of Poincar\'e duality. Indeed, in the work of Bhatt--Lurie \cite{BL22a}, Poincar\'e duality for syntomic cohomology is reproven as a combination of ``geometric Poincar\'e duality'' for prismatic $F$-gauges (due to Longke Tang \cite{Tang22}), and Serre duality on $\FSyn{}$. We shall see that we can construct $\varphi$ only from the latter Serre duality. 

Let $\cF,\cG \in \Perf(\cS)$, the $\infty$-category of perfect complexes on $\cS$.  We will construct a map 
\begin{equation}\label{eq: varphi'}
\varphi^{\prism}\colon \RGamma(\cS;\cF\otimes_{\cO_{\cS}} \cG) \to \RGamma(\cS;\cF) \otimes_{\F_p} \RGamma(\cS;\cG)[1]. 
\end{equation}
and then show that $\rH^*(\varphi^{\prism})$ recovers $\varphi_*$ in a suitable special case. 

\subsubsection{Serre duality on $\cS$}
We remind that $(-)^{\vee}$ denotes the formation of symmetric monoidal duals, which can be calculated by inner hom to the unit. The duality results of \cite[\S 4.5]{Bha22} can be interpreted as follows. The stack $\cS$ enjoys a form of Serre duality:
\begin{equation}\label{eq: Serre duality}
\RGamma(\cS;\cK)^\vee \cong \RGamma(\cS;\DD_{\cS} \cK) \in \cD(\F_p),
\end{equation}
where $\DD_{\cS}  \cK  := \cRHom(\cK, \omega_{\cS})$ is the internal Hom into the dualizing sheaf $\omega_{\cS}$, and furthermore there is an isomorphism 
\begin{equation}\label{eq: dualizing on S}
\omega_{\cS} \cong \cO_{\cS}[1].
\end{equation}
The isomorphisms \eqref{eq: Serre duality} and \eqref{eq: dualizing on S} come from \cite[Theorem 4.5.2]{Bha22}, which says that we have natural isomorphisms of functors $\Perf(\cS)^\op \rightarrow \cD(\F_p)$,
\begin{equation}\label{eq: varphi' eq 0}
\Serre \co \RGamma(\cS;(-))^\vee \xrightarrow{\eqref{eq: Serre duality}}  \RGamma(\cS;\DD_{\cS} (-)) \xrightarrow{\eqref{eq: dualizing on S}}  \RGamma(\cS;(-)^\vee[1]) 
\end{equation}
whose composite we call $\Serre$. 

\subsubsection{Construction of $\varphi^{\prism}$}
Let $\cF,\cG \in \Perf(\cS)$. We have a natural isomorphism in $\Perf(\cS)$,
\begin{equation}\label{eq: varphi' eq 1}
\cF^\vee \otimes_{\cO_{\cS}} \cG^\vee \cong \cRHom_{\cS}(\cF,\cO_{\cS}) \otimes_{\cO_{\cS}} \cRHom_{\cS}(\cG,\cO_{\cS}) \xrightarrow{\sim} \cRHom_{\cS}(\cF\otimes_{\cO_{\cS}} \cG,\cO_{\cS}) \cong (\cF\otimes_{\cO_{\cS}} \cG)^\vee.
\end{equation}
Applying the lax symmetric monoidal functor $\RGamma \colon \Perf(\cS) \to \cD(\F_p)$ to \eqref{eq: varphi' eq 1} and composing with the cup product gives a map 
\begin{equation}\label{eq: varphi' eq 2}
\RGamma(\cS;\cF^\vee)\otimes_{\F_p} \RGamma(\cS;\cG^\vee) \to \RGamma(\cS;\cF^\vee \otimes_{\cO_{\cS}} \cG^\vee) \cong  \RGamma(\cS;(\cF\otimes_{\cO_{\cS}} \cG)^\vee). 
\end{equation}
Taking $\F_p$-linear duals and using that $\RGamma(\cS;-)$ takes dualizable objects to dualizable objects \cite[Proposition 4.5.1]{Bha22}, we obtain a sequence of maps
\begin{equation}\label{diag: varphi'}\begin{tikzcd}
   \RGamma(\cS;(\cF\otimes_{\cO_{\cS}} \cG)^\vee)^\vee \ar[r, "\eqref{eq: varphi' eq 2}^\vee"] \ar[dd, "\wr"', "\Serre^{-1}"]  &  (\RGamma(\cS;\cF^\vee)\otimes_{\F_p} \RGamma(\cS;\cG^\vee))^\vee \ar[d, "\wr"]  \\     
   & \ar[d, "\wr", "\Serre^{-1}"'] \RGamma(\cS;\cF^\vee)^\vee\otimes_{\F_p} \RGamma(\cS;\cG^\vee)^\vee  \\
\RGamma(\cS; \cF\otimes_{\cO_{\cS}} \cG)[1]  \ar[r, dashed] &  \RGamma(\cS;\cF)[1]\otimes_{\F_p} \RGamma(\cS;\cG)[1] 
\end{tikzcd}
\end{equation}
Note the resemblance between this construction and \eqref{eq: varphi diagram}. Shifting the dashed map by $-1$ gives the desired map 
\begin{equation}\label{eq: varphi' eq 3}
\varphi^{\prism}\colon \RGamma(\cS; \cF\otimes_{\cO_{\cS}} \cG) \to \RGamma(\cS;\cF)\otimes_{\F_p} \RGamma(\cS;\cG)[1]. 
\end{equation}
In particular, upon taking cohomology we obtain a map 
\begin{equation}\label{eq: varphi' eq 3.5}
\varphi^{\prism}_* \colon \rH^{*,*}(\cS; \cF\otimes_{\cO_{\cS}} \cG) \to \rH^{*,*}(\cS;\cF)\otimes_{\F_p} \rH^{*,*}(\cS;\cG)[1]. 
\end{equation}

\subsubsection{Comparison with $\varphi_*$}\label{sssec: comparison varphi}
Recall that for smooth and proper $f \co X \rightarrow \Spec k$, we have a perfect complex $\cH^X := Rf^{\Syn}_* (\cO_{X^{\Syn}})$, whose mod $p$ reduction is $\ol \cH^X \in \Perf(\cS)$.\footnote{More generally, this formula defines $\cH^X$ for any quasicompact smooth scheme $f\co X \rightarrow \Spec k$, and if $X$ is proper then $\cH^X$ is perfect, as explained in \cite[Remark 4.2.3]{Bha22}.} In \S \ref{ssec:spectral-cohomological-motives}, we lifted this construction to $\sH^X \in \cD(\Splus)$, and equipped it with a natural isomorphism $\iota^* \sH^X = \ol\cH^X$.

Taking $\cF = \cG  := \ol\cH^X$ in \eqref{eq: varphi' eq 3.5}, and using the K\"{u}nneth isomorphism $\ol\cH^X\otimes \ol\cH^X \cong \ol\cH^{X\times_{k} X}$, we obtain a map
\begin{equation}\label{eq: varphi' eq 4}
\varphi^{\prism}_X \colon \Hsyn^{*,*}(X\times_{k} X) \to \Hsyn^{*,*}(X)\otimes_{\F_p} \Hsyn^{*,*}(X)[1].
\end{equation}

\begin{lemma}\label{lem: the maps agree}
Let $X$ be a smooth, proper, geometrically connected variety over $k$ of dimension $d$. Then the map $\varphi^{\prism}_X$ from \eqref{eq: varphi' eq 4} agrees with the map $\varphi_*$ from \eqref{eq: varphi diagram}.
\end{lemma}

\begin{proof}
As explained in \cite[\S 4.5]{Bha22}, the Poincar\'e duality isomorphism $\Hsyn^{*,*}(X)^\vee \cong 
 \Hsyn^{2d+1-*,d-*}(X)$ is a combination of the Serre duality isomorphism 
\[
\Serre \colon \rH^{*,*}(\cS; (\ol \cH^X)^\vee[1] ) \xleftarrow{\sim} \rH^{*,*}(\cS;\ol \cH^X)^\vee 
\]
and the ``geometric Poincar\'e duality'' isomorphism \cite{Tang22}
\[
\PD_X\colon \ol \cH^X[2d]\BK{d} \iso (\ol \cH^X)^\vee  \in \Perf(\cS).
\]

Moreover, the geometric Poincar\'e duality isomorphism is compatible with products, in the sense that the diagram
\[
\xymatrix{
\ol \cH^{X\times_k X'}[2d+2d']\BK{d+d'} \ar^\sim[r]\ar_\wr^{\PD_{X\times X'}}[d] & \ol \cH^X[2d]\BK{d} \otimes \ol\cH^{X'}[2d']\BK{d'}\ar_\wr^{\PD_{X}\otimes \PD_{X'}}[d]\\
(\ol \cH^{X\times_k X'})^\vee\ar^\sim[r]  & (\ol \cH^X)^\vee\otimes (\ol \cH^{X'})^\vee 
}
\]
commutes.\footnote{Inspecting the construction of Poincar\'e duality in \cite{Tang22}, this is a consequence of the additivity property for Thom classes in \cite[Theorem 4.2]{Tang22}.}

Consider the following diagram of cohomology groups of sheaves on $\cS$, which we omit from the notation for ease of reading:
\[
\adjustbox{scale=1,center}{\begin{tikzcd}
    {\rH^{*,*}(\ol\cH^X)^\vee \otimes_{\F_p} \rH^{*,*}(\ol\cH^X)^\vee} \arrow[d, "\wr"', "\Serre"] &
    {\rH^{*,*}(\ol\cH^X\otimes_{\cO_{\cS}} \ol\cH^X)^\vee} \arrow[d, "\wr"', "\Serre"] \arrow[l, "\mathrm{\lax^\vee}"] &
    {\rH^{*,*}(\ol\cH^{X\times_{k} X})^\vee}  \arrow[l, "\sim"', "\text{K\"unneth}"] \arrow[d, "\wr"', "\Serre"]  \\
    {\rH^{*,*}((\ol\cH^X)^\vee[1])\otimes_{\F_p} \rH^{*,*}((\ol\cH^X)^\vee[1])}  &
    {\rH^{*,*}((\ol\cH^X)^\vee \otimes_{\cO_{\cS}} (\ol\cH^X)^\vee[1])} \arrow[l, "\varphi^{\prism}_*"] &
    {\rH^{*,*}((\ol\cH^{X\times_k X})^\vee[1])} \arrow[l, "\sim"', "\text{K\"unneth}"] \\
    {\rH^{*,*}(\ol\cH^X)\otimes_{\F_p} \rH^{*,*}(\ol\cH^X)[4d+2]} \arrow[u, "\PD_X\otimes  \PD_X"',"\wr"] &
    {\rH^{*,*}(\ol\cH^X\otimes_{\cO_{\cS}}  \ol\cH^X[4d+1])} \arrow[l, "\varphi^{\prism}_*"] \arrow[u, "\wr", "\PD_X\otimes \PD_X"'] &
    {\rH^{*,*}(\ol\cH^{X\times_{k} X}[4d+1])} \arrow[l, "\sim"', "\text{K\"unneth}"] \arrow[u, "\wr"', "\PD_{X\times X}"]
\end{tikzcd}}
\]
We claim that each square in the diagram commutes. Indeed:
\begin{itemize}
\item The upper left square commutes by definition of $\varphi^{\prism}_*$. 
\item The upper right square commutes by naturality of Serre duality. 
\item The lower left square commutes by naturality of $\varphi^{\prism}_*$. 
\item The lower right square commutes by the aforementioned compatibility of Poincar\'e duality with products.
\end{itemize}
Hence the entire diagram commutes. 

Now, the map $\varphi_*$ from \eqref{eq: varphi diagram} is the composite map in the bottom row of the diagram, while the map $\varphi^{\prism}_X$ \eqref{eq: varphi' eq 4} is the composition of the other three faces of the outer square, so they agree. 
\end{proof}

\subsection{Spectral Serre duality}\label{ssec: spectral Serre duality} We will establish an incarnation of Serre duality on $\Splus$ that ``lifts'' Serre duality on $\cS$ in an appropriate sense.

Following classical coherent duality, we might try to start by defining a dualizing sheaf on $\Splus$ as a right adjoint to the global sections functor. Unfortunately, in $\pFGauge{\mbb{S}}$ the unit $\psienh(\Sph) = \psi(\SphFp)$ is not $p$-completely compact.\footnote{The unit of the ``true'' prismatic stable homotopy category, $\FGauge{\mbb{S}}$, is supposed to be $p$-completely compact.} As a result, the global sections functor $\cD(\Splus) \to \Sp$ is not colimit-preserving, hence cannot have a right adjoint. 

\subsubsection{Compatibility of Ind-completion and adjunctions} There is a general categorical fix for this situation: passing to \emph{$\Ind$-completion}.\footnote{Applying $\Ind$ to big categories causes set-theoretic issues. These can be easily fixed using a choice of big cardinal $\kappa$ and applying $\cC\mapsto \Ind(\cC^\kappa)$ instead of $\Ind(\cC)$; we leave it for the reader to carry this modification.} 
Indeed, if $F\colon \cC \to \cD$ is an exact functor between stable $\infty$-categories, the functor $\Ind(F)\colon \Ind(\cC) \to \Ind(\cD)$ (that we will usually denote simply by $F$) admits a right adjoint $F^R \colon \Ind(\cD) \to \Ind(\cC)$. 

We will need to see that passing to Ind-completion does not ruin the good properties enjoyed by the functors between $\FGauge{\MSFp}$, $\pFGauge{\mbb{S}}$, and $\Sp$ that we have been considering. 
We now collect the facts about the construction $\Ind(-)$ that will ensure this. 

\begin{prop}\label{prop:ind-completion-right-adjoint}
Let $F\colon \cC \to \cD$ be an exact functor of $\infty$-categories with finite colimits, which admits an exact right adjoint $G\colon \cD \to \cC$.\footnote{If $\cC$ and $\cD$ are stable, which will be the case when we apply this result, then the adjoint is automatically exact if it exists.} Then $\Ind(G)$ is canonically right adjoint to $\Ind(F)$.
\end{prop}

\begin{proof}
The functor $\Ind$ is a functor of $(\infty,2)$-categories from the category of categories with finite colimits and exact functors between them to the category of presentable categories. Therefore, it carries adjunctions to adjunctions.  
\end{proof}

\begin{prop}\label{prop:ind_closed_symmetric_monoidal}
Let $\cC$ be a symmetric monoidal $\infty$-category with biexact tensor product. Then: 
\begin{itemize}
\item The fully faithful embedding $i\colon \cC \to \Ind(\cC)$ is symmetric monoidal; in particular,
\[
i(c \otimes c') \cong i(c)\otimes i(c') \in \Ind(\cC) \quad \text{ for all } c,c'\in \cC.
\]
\item If $\cC$ is closed symmetric monoidal, so that we have an internal Hom-functor $\hom_\cC$, then $i$ is closed; in particular,
\[
i(\hom_{\cC}(c,c')) \cong \hom_{\Ind(\cC)}(i(c),i(c')) \in \Ind(\cC) \text{ for all } c,c'\in \cC.
\]
\end{itemize}
\end{prop}

\begin{proof}
The first point follows immediately from the construction of the symmetric monoidal structure on the functor $\Ind$. The second point follows from the first using the compatibility of $\Ind$ with the formation of right adjoints (Proposition \ref{prop:ind-completion-right-adjoint}). To spell this out: for every $c\in \cC$ the functor $\hom_\cC(c,-)$ is right adjoint to $c\otimes (-)$. Hence the ind-completion of $\hom_{\cC}(c,-)$ is right adjoint to the ind-completion of the functor $c\otimes (-)$, i.e., to $i(c)\otimes (-)$. We deduce that $\hom_{\Ind(\cC)}(i(c),-) \cong i(\hom_{\cC}(c,-))$. Restricting this to the essential image of $i$, we obtain the result.     
\end{proof}

\subsubsection{Compatibility of Ind-completion and projection formula}
Given an exact, symmetric monoidal functor $f^*\colon \cD \to \cC$ with lax symmetric monoidal right adjoint $f_*\colon \cC \to \cD$, we have a natural projection map 
\[
(f_*c)\otimes d \to
f_*(c\otimes f^*d),  \quad c \in \cC, d \in \cD.
\]
If this map is an isomorphism, so that $f_*$ is $\cD$-linear, let us say that the adjunction $f^*\dashv f_*$ satisfies the \emph{projection formula}.

\begin{prop}\label{prop:ind-completion-projection-formula}
Let $f^*\dashv f_*$ be a symmetric monoidal adjunction as above. If it satisfies the projection formula, then so does $\Ind(f^*) \dashv \Ind(f_*)$. 
\end{prop}

\begin{proof}
The functors 
\begin{equation}\label{eq:projection-formula-ind}
(c,d) \mapsto \Ind(f_*) (c\otimes \Ind(f^*)d) \quad \text{ and } \quad (c,d) \mapsto (\Ind(f_*)c)\otimes d
\end{equation}
both preserve filtered colimits in the $c$ and the $d$ variables separately. Hence, to check that a natural transformation between them is an isomorphism, it suffices to check that it is an isomorphism at $c,d$ in the essential images of the embeddings $\cC \to \Ind(\cC)$ and $\cD \to \Ind(\cD)$, respectively. For such objects, both functors in \eqref{eq:projection-formula-ind} land in $\cD\subseteq \Ind(\cD)$ and the projection map between them coincides with the projection map of the non-Ind-completed adjunction $f^*\dashv f_*$, which we assumed to be an isomorphism.
\end{proof}

Finally, assume that $\cC$ and $\cD$ are closed symmetric monoidal and that the functor $f_*$ itself has a right adjoint $f^!$. Then, if $f^*\dashv f_*$ satisfies the projection formula, we can pass to the right adjoints to obtain an isomorphism
\[
\hom_{\cD}(f_* c,d) \cong f_*\hom_{\cC}(c,f^!d) 
\]

\begin{cor}\label{cor:ind-complete-iota-functors}
Let $f^*\colon \cC \to \cD$ be an exact, symmetric monoidal functor between stable closed symmetric monoidal $\infty$-categories which admits a right adjoint $f_*$. Assume that $f_*$ admits a further right adjoint $f^!$ and that the adjunction $f^*\dashv f_*$ satisfies the projection formula. Then all these properties remain true after Ind-completion. Namely, the functor $\Ind(f^*)$ admits $\Ind(f_*)$ as a right adjoint, which admits $\Ind(f^!)$ as a further right adjoint. The adjunction $\Ind(f^*) \dashv \Ind(f_*)$ satisfies the projection formula and hence we have 
\[
\hom_{\Ind(\cD)}(\Ind(f_*)X,Y)\cong \Ind(f_*)\hom_{\Ind(\cC)}(X,\Ind(f^!)Y) 
\]
\end{cor}

\begin{convention}
To avoid cumbersome notation, from now on we shall denote the ind-completion of a functor $F\colon \cC \to \cD$ simply by $F \co \Ind(\cC) \rightarrow \Ind(\cD)$. The compatibilities established above imply that this abuse of notation does not affect the validity of assertions involving $F$. In particular, we shall denote by $\hom_\cC(-,-)$ the internal hom of $\Ind(\cC)$ when $\cC$ is closed symmetric monoidal. 
\end{convention}

\subsubsection{Brown--Comenetz duality}
Brown--Comenetz \cite{BC76} introduced the \emph{Brown--Comenetz} spectrum, which represents the generalized cohomology theory associating to a spectrum $E$ the Pontrjagin dual of its homotopy groups, 
\[
\bI^*(E) = \Hom_{\Z}(\pi_{-*} (E), \Q/\Z).
\]

Let $\bI \in \Sp$ be the $p$-completion of the Brown--Comenetz spectrum. Thus, if $E \in \Sp$ is a $p$-complete spectrum of bounded $p$-power torsion (meaning that $p^N\colon E\to E$ factors over the zero map for sufficiently large $N$), then the mapping spectrum $\cRHom_{\Sp}(E,\bI)$ has homotopy groups 
\[
\pi_{i}\cRHom_{\Sp}(E,\bI) \cong \Hom(\pi_{-i}E,\Q/\Z).
\]
Thus (at least for spectra of bounded $p$-torsion) the functor $E \mapsto \cRHom_{\Sp}(E,\bI)=: \bI E$ can be viewed as a generalization of Pontrjagin duality to spectra. For example, there is a natural isomorphism 
\begin{equation}\label{eq: BC of HF_p}
\bI \F_p \cong  \F_p.
\end{equation}

\subsubsection{Spectral dualizing sheaf}\label{sssec:spectral-dualizing-sheaf}

Let $\cC$ be a $p$-complete presentably symmetric monoidal stable $\infty$-category (in particular, a module over $\Sp$ in presentable categories). Let $(\pi_\cC)^* \colon\Sp \to  \cC$ be the unit functor and $(\pi_{\cC})_*$ its right adjoint. \Cref{cor:ind-complete-iota-functors} implies that after $\Ind$-completion, we have a further right adjoint $\pi_\cC^! \colon \Ind(\Sp) \to \Ind(\cC)$. 
\begin{defn} 
We define the \emph{dualizing object} of $\cC$ to be the ind-object $\omega_\cC:= \pi_\cC^!\bI \in \Ind(\cC)$, and the corresponding \emph{Serre duality functor} to be 
\[
\DD_{\cC}:= \hom_{\Ind(\cC)}(-,\omega_\cC)\colon \Ind(\cC) \to \Ind(\cC)^\op.
\]
\end{defn}

\begin{remark}
This definition gives a reasonable notion of dualizing object only in a ``sufficiently $p$-torsion'' setup (which is the only situation in which we will apply it). 
\end{remark}

More generally, if $f^*\colon \cC \to \cD$ is a colimit-preserving symmetric monoidal functor, then by the same procedure we produce a functor $f^!\colon \Ind(\cC) \to \Ind(\cD)$.

\begin{example}\label{ex:F_p_linear_dualizing}
Let $\cD(\F_p)$ be the derived $\infty$-category of $\F_p$-vector spaces, and let $\phi \colon \cD(\F_p) \to \Sp$ be the forgetful functor, whose left adjoint is given by tensoring with $\F_p$ over $\sph$. Then $\phi^! = \hom_{\Sp}(\F_p,-)$. Since $\hom_{\Sp}(\F_p,\BC) \cong \F_p$ by \eqref{eq: BC of HF_p}, we deduce that $\omega_{\F_p} := \phi^!\bI \cong \F_p$. 

More generally, if $\cC$ is a stable presentably symmetric monoidal category which is linear over $\F_p$, then the global sections (i.e., RHom from the unit object) functor $\cC \to \Sp$ factors as $\cC \oto{\pi'} \cD(\F_p) \oto{\phi} \Sp$, so that 
\[
\omega_{\cC} = \pi_{\cC}^!\BC = (\pi')^!\omega_{\F_p} \cong (\pi')^!\F_p.
\]
In particular, if $\cC = \QCoh(X)$ for a scheme $X/\F_p$, then $\omega_{\cC} = \omega_X$ is the usual dualizing (ind-)sheaf in the theory of coherent duality.   
\end{example}
	
\subsubsection{Spectral Serre duality}
Applying the discussion of \S\ref{sssec:spectral-dualizing-sheaf} to $\cC := \cD(\Splus)$, we can now define an ind-object 
\[
\omega_\Splus \in \Ind(\cD(\Splus)),
\]
and a corresponding duality functor 
\[
\DD_\Splus \colon \Ind(\cD(\Splus)) \to \Ind(\cD(\Splus))^\op.
\]
Similarly, we have the dualizing object $\omega_\cS = \iota^!\omega_{\Splus} \in \Ind(\cD(\cS))$ and a corresponding Serre duality functor $\DD_\cS$. In fact, since the unit of $\cD(\cS) \cong \FGauge{\FF_p}$ is $p$-completely compact, the global sections functor already has a right adjoint. Therefore, the object $\omega_\cS$ belongs to $\cD(\cS) \subset \Ind(\cD(\cS))$, hence the functor $\DD_\cS$ sends ind-constant objects to ind-constant objects by \Cref{prop:ind_closed_symmetric_monoidal}. We abuse notation and denote the resulting functor $\cD(\cS) \to \cD(\cS)^\op$
 again by $\DD_\cS$.

 \begin{prop}\label{prop: iota dual commute}
There is a canonical isomorphism 
\begin{equation}\label{eq: iota dual commute}
\iota_*\DD_\cS \cong \DD_{\Splus}\iota_*\colon \Ind(\cD(\cS)) \to \Ind(\cD(\Splus))^\op.  
\end{equation}
In particular, the functor $\DD_{\Splus}$ carries the essential image of $\iota_*\colon \cD(\cS) \to \cD(\Splus)$ into $\cD(\Splus)\subseteq \Ind(\cD(\Splus))$.  
 \end{prop}

 \begin{proof}
The adjunction $\iota^*\dashv \iota_*$ between the presentably symmetric monoidal stable $\infty$-categories $\cD(\cS)$ and $\cD(\Splus)$ satisfies the projection formula by \Cref{prop:iota_affine}. By \Cref{prop:ind-completion-projection-formula}, the induced adjunction between the ind-completions also satisfies the projection formula. Again by \Cref{prop:iota_affine}, $\iota_*$ preserves colimits hence admits a right adjoint $\iota^!$ \cite[Corollary 5.5.2.9]{HTT}, so we have a canonical isomorphism of the form 
\[
\iota_*\hom_{\cS}(\cF,\iota^!\sG)\cong \hom_{\Splus}(\iota_*\cF,\sG)
\]
for all $\cF\in \Ind(\cD(\cS))$ and $\sG\in \Ind(\cD(\Splus))$.
The natural isomorphism \eqref{eq: iota dual commute} follows by taking $\sG = \omega_{\Splus}\in \Ind(\cD(\Splus))$. The ``in particular'' part follows because $\iota_*$ and $\DD_\cS$ carry ind-constant objects to ind-constant objects: for $\iota_*$ this is clear, and for $\DD_{\cS}$ it follows from \Cref{prop:ind_closed_symmetric_monoidal}.
 \end{proof}

\subsubsection{Duality involution on the prismatized Steenrod algebra} Combining the identification $\omega_{\cS} \cong \cO_{\cS}[1]$ from \eqref{eq: dualizing on S} with \eqref{eq: iota dual commute} yields isomorphisms
\begin{equation}\label{eq:iota-O-dual-identity}
\iota_* \cO_{\cS}[1]  \cong \iota_* \DD_{\cS}(\cO_{\cS})  \cong \DD_{\Splus}(\iota_* \cO_{\cS}).
\end{equation}
From this we get a sequence of isomorphisms
\begin{align}\label{eq: sigma}
\sAsyn & = \cRHom_{\Splus}(\iota_* \cO_{\cS}, \iota_* \cO_{\cS}) \xrightarrow{\sim} \cRHom_{\Splus}(\DD_{\Splus} \iota_* \cO_{\cS}, \DD_{\Splus} \iota_* \cO_{\cS})^{\op} \nonumber \\
& \stackrel{\eqref{eq:iota-O-dual-identity}}\cong \cRHom_{\Splus}(\iota_* \cO_{\cS}[1], \iota_* \cO_{\cS}[1])^{\op}  \cong \cRHom_{\Splus}(\iota_* \cO_{\cS}, \iota_* \cO_{\cS})^{\op} = \sAsyn^{\op}.
\end{align}

\begin{defn}[The involution $\sigma$]
We let $\sigma \co \sAsyn \rightarrow \sAsyn^{\op}$ be the composition of the maps in \eqref{eq: sigma}. Noting that $\sigma$ is an involution, we also write $\sigma \co \sAsyn^{\op} \rightarrow \sAsyn$ for the opposite map. 

\end{defn}

The algebra $\sAsyn$ has a tautological action on $\iota_*\cO_{\cS}$, while $\sAsyn^\op$ has a tautological action on the functor $\cRHom_{\Splus}(\iota_*\cO_{\cS},-)$. By construction, $\sigma$ is characterized by the following property: the composition of isomorphisms
\begin{equation}\label{eq: sigma prop 1}
\iota_*(\cO_{\cS}[1])\cong \iota_*\omega_{\cS} \cong \iota_* \iota^!\omega_{\Splus} \cong \cRHom_{\Splus}(\iota_*\cO_{\cS},\omega_{\Splus}) 
\end{equation}
is $\sigma$-semilinear when equipped with the tautological $\sAsyn$-action on the left and the tautological $\sAsyn^{\op}$-action on the right, both via the tautological $\sAsyn$-action on $\iota_* \cO_{\cS}$. In other words, \eqref{eq: sigma prop 1} promotes to an isomorphism 
\begin{equation}\label{eq: sigma characterization}
\sigma^*\iota_*(\cO_{\cS}[1]) \cong \cRHom_{\Splus}(\iota_*\cO_{\cS},\omega_{\Splus}) 
\end{equation}
of $\sAsyn^\op$-modules in $\cD(\Splus)$.

\subsubsection{Serre duality and Steenrod action} We will see that, informally speaking, ``the involution $\sigma$ intertwines the syntomic Steenrod action with Serre duality''. Using \Cref{ex:F_p_linear_dualizing}, we have an identification
\begin{equation}\label{eq: dualizing compatibility}
\omega_{\cS} \cong  \iota^! \pi_{\cD(\Splus)}^! \BC \cong \iota^! \omega_{\Splus} \in \cD(\cS).
\end{equation}

\begin{defn}[Steenrod action on Serre dual]\label{def:action_on_serre_duals}
Let $\sF \in \cD(\Splus)$ and $\cF = \iota^* \sF \in \cD(\cS)$. Then by adjunction, the projection formula, and the identification \eqref{eq: dualizing compatibility}, we have a chain of isomorphisms
\begin{equation}\label{eq:action_on_serre_duals}
\iota_*\DD_{\cS}(\cF)\cong \iota_*\cRHom_{\cS}(\cF,\iota^!\omega_{\Splus})\cong \cRHom_{\Splus}(\iota_*\cF,\omega_{\Splus}) \cong \cRHom_{\Splus}(\sF\otimes_{\cO_{\Splus}} \iota_*\cO_{\cS},\omega_{\Splus}), 
\end{equation}
in which all objects a priori belong to $\Ind(\cD(\Splus))$, but are in fact contained in the full subcategory $\cD(\Splus)$ (by \Cref{prop: iota dual commute} and the discussion preceding it). The tautological $\sAsyn$-action on $\iota_*\cO_{\cS}$ induces a tautological $\sAsyn^\op$-action on $\cRHom_{\Splus}(\sF\otimes_{\cO_{\Splus}} \iota_*\cO_{\cS},\omega_{\Splus})$. We then equip $\iota_*\DD_{\cS}(\cF)$ with the $\sAsyn^\op$-action induced by transport along \eqref{eq:action_on_serre_duals}.
\end{defn}

Now suppose that $\sF \in \Perf(\Splus)$ and $\cF:= \iota^*\sF \in \Perf(\cS)$. Since $\sF$ is dualizable and $\iota^*$ is symmetric monoidal, $\cF$ is dualizable and there is a natural isomorphism $\iota^*(\sF^\vee) \cong \cF^\vee$, equipping $\iota_* \cF^\vee[1]$ with a natural $\sAsyn$-action (cf. \S \ref{subsec:action_on_duals}). On the other hand, the identification of the dualizing sheaf on $\cS$ in \eqref{eq: dualizing on S} supplies an isomorphism $\DD_{\cS}(\cF) \cong \cF^\vee[1]$, which equips
\begin{equation}\label{eq: Serre dual of F}
\iota_* \DD_{\cS}(\cF) \cong \iota_* \cF^\vee[1]
\end{equation}
with a natural $\sAsyn^\op$-action by Definition \ref{def:action_on_serre_duals}. We will relate these two actions.

\begin{prop}\label{prop:Serre_Steenrod}
Let $\sF \in \Perf(\Splus)$ and $\cF  := \iota^* \sF \in \Perf(\cS)$. Equipping each side of \eqref{eq: Serre dual of F} with the tautological $\sAsyn^{(\op)}$-action described above, it promotes to an $\sAsyn$-equivariant isomorphism
\begin{equation}\label{eq:Serre_Steenrod}
\sigma^*\iota_*\DD_{\cS}(\cF) \cong \iota_*\cF^\vee[1] \in \cD(\Splus).
\end{equation}
\end{prop}
\begin{proof}
The $\sAsyn^{\op}$-action on the left side of \eqref{eq: Serre dual of F} is obtained from the tautological $\sAsyn$-action on $\iota_* \cO_{\cS}$ in \eqref{eq:action_on_serre_duals}. We can rewrite this using tensor-Hom adjunction as
\[
\begin{tikzcd}
\iota_*\DD_{\cS}(\cF) \ar[r, "\sim", "\eqref{eq:action_on_serre_duals}"'] &  \cRHom_{\Splus}(\sF \otimes \iota_* \cO_{\cS}, \omega_{\Splus}) \ar[r, "\sim"] &   \cRHom_{\Splus}(\sF,\cRHom_{\Splus}(\iota_*\cO_{\cS},\omega_{\Splus}))
\end{tikzcd}
\]

The $\sAsyn$-action on the right side of \eqref{eq: Serre dual of F} is obtained from the tautological $\sAsyn$-action on $\iota_* \cO_{\cS}$ in 
\[
\iota_*\cF^\vee[1] \cong \cRHom_{\Splus}(\sF,\iota_*\cO_{\cS}[1]).
\]

In both cases, the algebra $\sAsyn$ acts via its action on the target of the mapping objects, so the $\sAsyn$-equivariance of \eqref{eq:Serre_Steenrod} reduces to the $\sAsyn$-equivariance of the isomorphism 
\[
\sigma^*\iota_*\cO_{\cS}[1]\cong \cRHom_{\Splus}(\iota_*\cO_{\cS},\omega_{\Splus}),
\] 
which holds by the definition of $\sigma$, as discussed around \eqref{eq: sigma characterization}. 
\end{proof}

\begin{cor}\label{cor: serre dual steenrod} Let $\sF \in \Perf(\Splus)$ and $\cF = \iota^* \sF \in \cD(\cS)$. Then the identification 
\[
\rH^{*,*}(\cS; \DD_{\cS} (\cF))  \cong \rH^{*,*}(\cS; \cF^\vee[1]),
\]
coming from \eqref{eq: dualizing on S}, promotes to a $(\Asyn^{*,*})^{\op}$-equivariant isomorphism 
\[
\sigma^* \rH^{*,*}(\cS; \DD_{\cS} (\cF))  \cong \rH^{*,*}(\cS; \cF^\vee[1]).
\]
\end{cor}

\begin{proof}
Apply Proposition \ref{prop:Serre_Steenrod} to $\sF\BK{b}$ for $b \in \Z$, then take global sections and then direct sum over all $b$. 
\end{proof}

\subsection{$\sigma$-equivariance of the coproduct}\label{ssec: 11.4} We will need the following compatibility of the duality involution $\sigma$ with the coproduct $\Delta \co \rAsyn^{*,*} \rightarrow \rAsyn^{*,*} \otimes_{\F_p} \rAsyn^{*,*}$. 

\begin{prop}\label{prop: inv Delta commute}
The anti-involution $\sigma\colon \Asyn^{*,*} \to (\Asyn^{*,*})^\op$ preserves the subalgebra $\rAsyn^{*,*}$, and defines a map of Hopf algebras $\rAsyn^{*,*}  \to (\rAsyn^{*,*} )^{\op}$. In particular, the diagram 
\[
\begin{tikzcd}
\rAsyn^{*,*} \ar[d, "\sigma"] \ar[r, "\Delta"] & \rAsyn^{*,*} \otimes_{\F_p} \rAsyn^{*,*} \ar[d, "\sigma \otimes \sigma"] \\
(\rAsyn^{*,*})^{\op} \ar[r, "\Delta^{\op}"] & (\rAsyn^{*,*})^{\op} \otimes_{\F_p} (\rAsyn^{*,*})^{\op}
\end{tikzcd}
\]
commutes.
\end{prop}

The rest of this subsection is devoted to the proof of Proposition \ref{prop: inv Delta commute}. 

\subsubsection{Initial reformulations} The dual syntomic Steenrod algebra $\dAsyn_{*,*}$ is a commutative Hopf algebra over $\Hsynpt$. Hence it has an antipode $\mf{s} \co \dAsyn_{*,*} \rightarrow (\dAsyn_{*,*})^{\op}$, which corresponds geometrically to inversion on the corresponding group scheme. Write $\mf{s}^* \co \Asyn^{*,*} \to (\Asyn^{*,*})^\op$ for the (opposite) map induced by duality over $\Hsynpt$ of $\mf{s}$. The antipode $\mf{s}$ is a Hopf $\Hsynpt$-algebra involution, so its dual $\mf{s}^*$ is a Hopf $\Hsynpt$-algebra involution; in particular, it is compatible with $\Delta$, while we want to prove that $\sigma$ is compatible with $\Delta$. So in order to prove Proposition \ref{prop: inv Delta commute}, it suffices to identify 
\begin{equation}\label{eq: sigma = dual antipode}
\mf{s}^* = \sigma \co \Asyn^{*,*} \to (\Asyn^{*,*})^\op.
\end{equation}

\subsubsection{Relating $\sigma$ to the swap map} Recall that $\sigma$ comes from the composition 
\begin{align*}
\sAsyn & = \cRHom_{\Splus}(\iota_* \cO_{\cS}, \iota_* \cO_{\cS}) \xrightarrow{\sim} \cRHom_{\Splus}(\DD_{\Splus} \iota_* \cO_{\cS}, \DD_{\Splus} \iota_* \cO_{\cS})^{\op} \nonumber \\
& \cong \cRHom_{\Splus}(\iota_* \cO_{\cS}[1], \iota_* \cO_{\cS}[1])^{\op}  \cong \cRHom_{\Splus}(\iota_* \cO_{\cS}, \iota_* \cO_{\cS})^{\op} = \sAsyn^{\op}
\end{align*}
where the passage to the second line used \eqref{eq:iota-O-dual-identity}. On the other hand, let $\eta$ be the composite of the isomorphisms of $\iota_* \cO_{\cS}$-modules below,
\begin{align*}
\eta \colon \cRHom_{\Splus}(\iota_*\cO_{\cS},\iota_*\cO_{\cS}) &= \cRHom_{\Splus}(\iota_*\cO_{\cS},\iota_*\cO_{\cS}[1])[-1] 
 \stackrel{\eqref{eq:iota-O-dual-identity}}\cong  \cRHom_{\Splus}(\iota_*\cO_{\cS},\DD_{\Splus} \iota_*\cO_{\cS})[-1] \\
 &= \cRHom_{\Splus}(\iota_*\cO_{\cS},\cRHom_{\Splus}(\iota_*\cO_{\cS},\omega_{\Splus}))[-1]  \cong \cRHom_{\Splus}(\iota_*\cO_{\cS} \otimes_{\cO_{\Splus}} \iota_*\cO_{\cS} ,\omega_{\Splus})[-1].
\end{align*}

\begin{lemma}\label{lem: swap sigma} Let $\sw \in \End_{\iota_* \cO_{\cS}}( \iota_*\cO_{\cS} \otimes_{\cO_{\Splus}} \iota_*\cO_{\cS})$ be the map swapping the two tensor factors. Then the square 
\begin{equation}\label{diag: swap 0}
\xymatrix{
\cRHom_{\Splus}(\iota_*\cO_{\cS},\iota_*\cO_{\cS})\ar^\wr_\eta[d]\ar^{\sigma}[r] & \cRHom_{\Splus}(\iota_*\cO_{\cS},\iota_*\cO_{\cS})^\op \ar^\wr_\eta[d]
\\
 \cRHom_{\Splus}(\iota_*\cO_{\cS}\otimes_{\cO_{\Splus}} \iota_*\cO_{\cS} ,\omega_{\Splus})[-1]\ar^{\sw^*}[r] & \cRHom_{\Splus}(\iota_*\cO_{\cS}\otimes_{\cO_{\Splus}} \iota_*\cO_{\cS} ,\omega_{\Splus})[-1]
}
\end{equation}
commutes.
\end{lemma}

\begin{proof}
We begin with generalities. For $\sF,\sG \in \cD(\Splus)$, we have by Hom-tensor adjunction natural isomorphisms 
\[
\cRHom_{\Splus}(\sF\otimes_{\cO_{\Splus}} \sG,\omega_{\Splus}) \cong \cRHom_{\Splus}(\sF, \cRHom_{\Splus}(\sG,\omega_{\Splus})) = \cRHom_{\Splus}(\sF,\DD_{\Splus}\sG).
\]
Moreover, the composite identification fits into a commutative diagram
\begin{equation}\label{diag: swap 1}
\xymatrix{
\cRHom_{\Splus}(\sF \otimes \sG,\omega_{\Splus})\ar^{\sw^*}[dd]\ar^\sim[r]
& \cRHom_{\Splus}(\sF,\DD_{\Splus}\sG)\ar^{\DD_{\Splus}}_\wr[d] \\ 
& \cRHom_{\Splus}(\DD_{\Splus}^2\sG,\DD_{\Splus}\sF)\ar_\wr[d]
\\
\cRHom_{\Splus}(\sG\otimes \sF,\omega_{\Splus}) \ar^\sim[r]& \cRHom_{\Splus}(\sG,\DD_{\Splus}\sF) 
}
\end{equation}
Apply this to $\sF := \iota_*\cO_{\cS}$ and $\sG := \iota_*\omega_{\cS}$, and embed it as the middle rectangle in the following large diagram. (We abbreviate $\DD := \DD_{\Splus}$ and $\otimes$ for $\otimes_{\cO_{\Splus}}$ for ease of notation.)
\[
\adjustbox{scale = 0.75, center}{
\begin{tikzcd}[column sep=small, row sep=huge]
\cRHom_{\Splus}(\iota_*\cO_{\cS} \otimes \iota_*\cO_{\cS}, \omega_{\Splus})[-1] 
  \arrow[dd, "\sw^*"']
  \arrow[r, "\sim", "\eqref{eq:iota-O-dual-identity}"']
  \arrow[rrr, bend left=30, "\eta^{-1}"'] &
\cRHom_{\Splus}(\iota_*\cO_{\cS} \otimes \iota_*\omega_{\cS}, \omega_{\Splus}) 
  \arrow[dd, "\sw^*"']
  \arrow[r, "\sim"] &
\cRHom_{\Splus}(\iota_*\cO_{\cS}, \DD \iota_*\omega_{\cS}) 
  \arrow[d, "\DD"] 
  \arrow[r, "\sim"] &
\cRHom_{\Splus}(\iota_*\cO_{\cS}, \iota_*\cO_{\cS}) 
  \arrow[d, "\DD"] 
  \arrow[ddd, bend left=60, "\sigma"] \\
& &
\cRHom_{\Splus}(\DD^2\iota_*\omega_{\cS}, \DD \iota_*\cO_{\cS}) 
  \arrow[d, "\wr"]
  \arrow[r, "\sim"] &
\cRHom_{\Splus}(\DD \iota_*\cO_{\cS}, \DD \iota_*\cO_{\cS}) 
  \arrow[d, "\wr"] \\
\cRHom_{\Splus}(\iota_*\cO_{\cS} \otimes \iota_*\cO_{\cS}, \omega_{\Splus})[-1] 
  \arrow[r, "\sim", "\eqref{eq:iota-O-dual-identity}"']
  \arrow[rrrd, bend right=15, "\eta^{-1}"] &
\cRHom_{\Splus}(\iota_*\omega_{\cS} \otimes \iota_*\cO_{\cS}, \omega_{\Splus}) 
  \arrow[r, "\sim"] &
\cRHom_{\Splus}(\iota_*\omega_{\cS}, \DD \iota_*\cO_{\cS}) 
  \arrow[r, "\sim"] &
\cRHom_{\Splus}(\iota_*\omega_{\cS}, \iota_*\omega_{\cS}) 
  \arrow[d, "\wr"] \\
& & &
\cRHom_{\Splus}(\iota_*\cO_{\cS}, \iota_*\cO_{\cS})
\end{tikzcd}}
\]
We claim that the entire diagram commutes. 
\begin{itemize}
\item The middle diagram commutes, as a special case of \eqref{diag: swap 1}. 
\item The left square commutes by naturality of $\sw^*$.
\item The right upper rectangle commutes by the naturality of $\DD$.
\item The right lower rectangle commutes by the ``exchange law'' (alias ``interchange law''). 
\item All the boundary curved 
triangles commute by the definitions of $\eta$ and $\sigma$. 
\end{itemize}
This establishes the claim. In particular, the outer diagram commutes. But this outer diagram is, up to inverse, the diagram \eqref{diag: swap 0} which we wanted to show was commutative, so we are done. 

\end{proof}

\subsubsection{Relating $\mf{s}$ to the swap map} 

By Serre duality, we have 
\begin{equation}\label{eq: antipode eq 2}
(\dAsyn_{-*,-*})^\vee \cong  \rH^{*,*}(\Splus; \iota_* \dsAsyn )^\vee \cong   \Ext^{*,*}_{\Splus}(\iota_*\cO_{\cS} \otimes_{\cO_{\Splus}} \iota_*\cO_{\cS},\omega_{\Splus}). 
\end{equation}
(Recall that $(-)^\vee$ means dual over $\F_p$ here.) By definition, the swap map induces the antipode of the Hopf $\iota_*\cO_{\cS}$-algebroid $\iota_* \dsAsyn =  \iota_*\cO_{\cS} \otimes_{\cO_{\Splus}} \iota_*\cO_{\cS} \in \cD(\Splus)$. Applying $\rH^{*,*}(\Splus; -)$ then shows that it induces also the antipode of the Hopf $\Hsynpt$-algebra $\dAsyn_{*,*}$. Hence \eqref{eq: antipode eq 2} carries the swap map $\sw^*$ on the RHS to the dual of the antipode of $\dAsyn_{-*,-*}$ on the left side. 

 A special case of Serre duality gives $\Hsynpt^\vee \cong \Hsynpt[1]$. Hence we have (using \eqref{eq:Asyn-dual-of-dAsyn} in the first step)
\begin{align}\label{eq: antipode eq 3}
\Asyn^{*,*}  & \cong \Hom_{\Hsynpt}(\dAsyn_{*,*}, \Hsynpt) \cong \Hom_{\F_p}(\drAsyn_{*,*}, \Hsynpt) \nonumber \\
& \cong  (\drAsyn_{*,*})^\vee \otimes_{\F_p} \Hsynpt \cong  (\drAsyn_{*,*})^\vee \otimes_{\F_p} \Hsynpt^\vee[-1] \cong (\dAsyn_{*,*})^\vee [-1].
\end{align} 
This is related to the earlier identifications in the following way. 

\begin{lemma}\label{lem: second duality} The composition of \eqref{eq: antipode eq 2} with $\eta^{-1}$ is the composite identification $(\dAsyn_{-*,-*})^\vee [-1] \cong \Asyn^{*,*}$ from \eqref{eq: antipode eq 3}. 
\end{lemma}

\begin{proof}
This is immediate upon comparing the definitions. 
\end{proof}

\subsubsection{Completion of the proof} The proof of \eqref{eq: sigma = dual antipode}, hence also of Proposition \ref{prop: inv Delta commute}, will be completed by the following Proposition.

\begin{prop}\label{prop: antipode sigma dual}
The isomorphism 
\begin{equation}\label{eq: antipode sigma dual}
\Hom_{\Hsynpt}(\dAsyn_{*,*}, \Hsynpt)  \cong  \Asyn^{*,*}
\end{equation}
intertwines the anti-involution $\mf{s}^*$ on the left side with the anti-involution $\sigma$ on the right side. 
\end{prop}

\begin{proof} Applying \Cref{lem: second duality} identifies the isomorphism \eqref{eq: antipode sigma dual} with the composite isomorphism 
\[
\begin{tikzcd}
(\dAsyn_{*,*})^\vee[-1] \ar[r, "{\sim}"', "\eqref{eq: antipode eq 2}"] &  \Ext^{*,*}_{\Splus}(\iota_*\cO_{\cS} \otimes_{\cO_{\Splus}} \iota_*\cO_{\cS},\omega_{\Splus})[-1] \ar[r, "\eta", "\sim"'] &  \Asyn^{*,*}.
\end{tikzcd}
\]
As discussed above, the first isomorphism carries $\mf{s}^*$, the dual of the antipode on $\dAsyn_{*,*}$, to $\sw^*$ on the middle term. Then by \Cref{lem: swap sigma}, the second isomorphism carries $\sw^*$ to $\sigma$. 

\end{proof}

\subsection{Proof of Theorem \ref{thm: steenrod push equivariant}}\label{ssec: 11.5} In this subsection, we will (finally!) combine the preceding ingredients to prove Theorem \ref{thm: steenrod push equivariant}. We do so by proving more generally: 

Let $\sF \in \Perf(\Splus)$ and $\cF := \iota^* \sF$. From now on, we adopt the following conventions.
\begin{itemize}
\item We regard $\iota_*\cF \cong \iota_*\iota^*\sF$ as an $\sAsyn$-module via the construction from \S \ref{subsec:action_on_i_star}. This induces an $\Asyn^{*,*}$-module structure on $\rH^{*,*}(\cS;\cF)$. 
 
\item We regard $\iota_*\DD_{\cS}(\cF)$ as an $\sAsyn^\op$-module via \Cref{def:action_on_serre_duals}, which induces an $(\Asyn^{*,*})^\op$-module structure on $\rH^{*,*}(\Splus; \iota_*\DD_{\cS}(\cF))$. 

\item If $M$ is an $\Asyn^{*,*}$-module, then we regard $M^\vee:=\Hom_{\F_p}(M,\F_p)$ as an $(\Asyn^{*,*})^\op$-module via the action on the source. 

\item If $M,N$ are $\rAsyn^{*,*}$-modules, then we regard $M\otimes_{\F_p} N$ as an $\rAsyn^{*,*}$-module by restriction along the coproduct $\Delta \co \rAsyn^{*,*} \rightarrow \rAsyn^{*,*} \otimes_{\F_p} \rAsyn^{*,*}$, as in \S \ref{ssec: coproduct}. 
\end{itemize}

\begin{prop}\label{dual action}
The Serre duality isomorphism $\rH^{*,*}(\cS; \cF)^\vee \cong \rH^{*,*}(\cS; \DD_{\cS} (\cF))$ is  $(\Asyn^{*,*})^{\op}$-equivariant, with respect to the actions defined above. 
\end{prop}

\begin{proof} 
We have commutative triangles\footnote{The functor $\pi_{\cS *} \co \cD(\cS) \rightarrow \Sp$ factors over $\cD(\F_p)$, but we are forgetting that structure here.}
\[
\begin{tikzcd}
\cD(\cS) \ar[r, "\iota_*"] \ar[dr, "\pi_{\cS *}"']  & \cD(\Splus) \ar[d, "\pi_{\Splus *}"] \\ 
  & \Sp
\end{tikzcd}
\quad \text{and} \quad 
\begin{tikzcd}
\cD(\cS)    & \cD(\Splus) \ar[l, "\iota^*"'] \\ 
  & \ar[ul, "\pi_{\cS}^*"] \Sp \ar[u, "\pi_{\Splus}^*"']
\end{tikzcd}
\]
Recall that $\bI$ is the Brown--Comenetz spectrum. In these terms, the Serre duality isomorphism assumes the form 
\[
\RGamma(\cS;\cF)^\vee  \cong \RGamma(\Splus; \iota_* \cF)^\vee \cong \RHom_{\Sp}(\pi_{\Splus *} \iota_*\iota^*\sF,\bI)\cong \RHom_{\Splus}(\iota_*\iota^*\sF,\pi_{\Splus}^!\bI)
\]
Chasing through these identifications, we see that the $(\Asyn^{*,*})^{\op}$-action on $
\rH^{*,*}(\cS; \cF)^\vee$ is obtained from the tautological $\sAsyn$-action of $\iota_* \iota^* \sF \cong \sF \otimes \iota_* \cO_{\cS}$.  According to \Cref{prop:Serre_Steenrod}, this is identified with the tautological $\sAsyn^{\op}$-action on $\DD_{\cS}(\cF)$ under the further sequence of identifications
\[
\cRHom_{\Splus}(\iota_*\iota^*\sF,\pi_{\Splus}^!\bI)  \cong  \cRHom_{\Splus}(\iota_*\iota^*\sF,\omega_{\Splus})  \cong \iota_* \DD_{\cS}(\cF).
\]
Passing to global sections on $\Splus$ completes the proof. 

\end{proof}

The following Lemma is elementary, but we record it for convenience. 

\begin{lemma}\label{lem: dual tensor}
The isomorphism 
\[
\left(\rH^{*,*}(\cS; \cF) \otimes_{\F_p} \rH^{*,*}(\cS; \cG) \right)^\vee \cong \rH^{*,*}(\cS; \cF)^\vee \otimes_{\F_p} \rH^{*,*}(\cS; \cG)^\vee
\]
is equivariant for the actions of $\rAsyn^{\op}$. 
\end{lemma}

\begin{proof}
More generally, if $R$ is any $\F_p$-algebra and $M,N$ are $R$-modules whose underlying $\F_p$-module is dualizable, then we have an identification 
\[
M^\vee \otimes_{\F_p} N^\vee \cong(M\otimes_{\F_p} N)^\vee
\] 
of $R^\op\otimes_{\F_p} R^\op$-modules. Applying this to $R := \rAsyn^{*,*}$ and restricting along $\Delta^{\op}$ yields the result.  
\end{proof}

By the discussion at the beginning of \S \ref{sssec: comparison varphi}, and the comparison of Lemma \ref{lem: the maps agree}, Theorem \ref{thm: steenrod push equivariant} is a special case of the following more general result.

\begin{thm}
Let $\sF, \sG \in \Perf(\Splus)$, and let $\cF=\iota^*\sF, \cG=\iota^*\sG \in \Perf(\cS)$. 
Then the map 
\[
\varphi^{\prism}_* \colon \rH^{*,*}(\cS;\cF\otimes_{\cO_{\cS}} \cG)  \to \rH^{*,*}(\cS;\cF)\otimes_{\F_p} \rH^{*,*}(\cS;\cG) [1]
\]
from \eqref{eq: varphi' eq 3.5} is $\rAsyn^{*,*}$-equivariant.
\end{thm}

\begin{proof}
It follows from \Cref{Cartan} and \Cref{dual cartan} that the cup product 
\begin{equation}\label{eq: final 1}
\rH^{*,*}(\cS; \cF^\vee) \otimes_{\F_p} \rH^{*,*}(\cS; \cG^\vee) \rightarrow \rH^{*,*}(\cS; \cF^\vee \otimes_{\cO_{\cS}} \cG^\vee) 
\end{equation}
is $\rAsyn^{*,*}$-equivariant. Dualizing \eqref{eq: final 1} and using \Cref{lem: dual tensor}, we learn that the dual map 
\begin{equation}\label{eq: final 2}
\rH^{*,*}(\cS; \cF^\vee)^\vee \otimes_{\F_p} \rH^{*,*}(\cS; \cG^\vee)^\vee \leftarrow \rH^{*,*}(\cS; \cF^\vee \otimes_{\cO_{\cS}} \cG^\vee)^\vee
\end{equation}
is $(\rAsyn^{*,*})^{\op}$-equivariant. 

By \Cref{dual action}, \eqref{eq: final 2} is identified $(\rAsyn^{*,*})^{\op}$-equivariantly with the map 
\begin{equation}\label{eq: final 3}
\rH^{*,*}(\cS; \DD_{\cS} (\cF^\vee)) \otimes_{\F_p} \rH^{*,*}(\cS; \DD_{\cS}( \cG^\vee)) \leftarrow \rH^{*,*}(\cS; \DD_{\cS}( \cF^\vee \otimes_{\cO_{\cS}} \cG^\vee)).
\end{equation}

By \Cref{cor: serre dual steenrod}, \eqref{eq: final 3} is identified $(\rAsyn^{*,*})^{\op}$-equivariantly with the map
\begin{equation}\label{eq: final 4}
\sigma^*\rH^{*,*}(\cS; \cF[1]) \otimes_{\F_p} \sigma^*\rH^{*,*}(\cS; \cG[1]) \leftarrow \sigma^*\rH^{*,*}(\cS; \cF \otimes_{\cO_{\cS}} \cG[1]),
\end{equation}
where the underlying map of graded $\F_p$-vector spaces is the dashed map in \eqref{diag: varphi'}, which is $\varphi^{\prism}$ up to shifting by 1. We still have to track that it is equipped with the prescribed Steenrod action. By the identification $(\sigma \otimes \sigma)\circ\Delta = \Delta \circ \sigma$ of \Cref{prop: inv Delta commute} we have an $(\rAsyn^{*,*})^{\op}$-equivariant isomorphism
\begin{equation}\label{eq: final 5}
\sigma^*\rH^{*,*}(\cS; \cF[1]) \otimes_{\F_p} \sigma^*\rH^{*,*}(\cS; \cG[1]) \cong 
\sigma^*(\rH^{*,*}(\cS; \cF[1]) \otimes_{\F_p} \rH^{*,*}(\cS; \cG[1])).
\end{equation}
Combining \eqref{eq: final 4} and \eqref{eq: final 5}, we obtain an $(\rAsyn^{*,*})^{\op}$-equivariant isomorphism 
\[
\sigma^*\rH^{*,*}(\cS; \cF \otimes_{\cO_{\cS}} \cG[1]) \to \sigma^*(\rH^{*,*}(\cS; \cF[1]) \otimes_{\F_p} \rH^{*,*}(\cS; \cG[1])).
\]
Finally, we conclude by applying $(\sigma^*)^{-1}$ and shifting. 
\end{proof}

\part{Characteristic classes}\label{part: characteristic classes}

In this Part, we develop a theory of mod $2$ characteristic classes in syntomic cohomology for $p=2$. The purpose of this theory is to facilitate explicit calculation of syntomic Steenrod operations. Our path is inspired by classical results in the algebraic topology of manifolds (as can be found in \cite{MS74}, for example): 
\begin{itemize}
\item Steenrod operations to the top degree are controlled by ``Wu classes''. 
\item (Wu formula) Wu classes can be expressed in terms of Stiefel--Whitney classes. 
\item For complex vector bundles, Stiefel--Whitney classes are reductions of Chern classes.
\end{itemize}
We will establish parallel results which will ultimately be used to prove the vanishing of certain syntomic Steenrod operations.

\section{Syntomic Stiefel--Whitney classes}

Let $X$ be a scheme over a field $k$ of characteristic $p=2$ and $E \rightarrow X$ be a vector bundle. In this section we construct ``syntomic Stiefel--Whitney classes'' $w^{\syn}_i \in \Hsyn^{i,\lfloor i/2 \rfloor}(X)$, following the approach of \cite{Feng20}, which in turn is based on ideas of Thom. A posteriori, these syntomic Stiefel--Whitney classes will turn out to be merely the reductions modulo $2$ of Chern classes, but this is a non-trivial calculation which is essential to the eventual applications.

\subsection{Cohomology with supports}\label{subsec: steenrod on relative coh}

Let $i \co Z \hookrightarrow X$ be a closed subscheme and $j \co U \inj X$ the complementary open embedding. For an \'{e}tale sheaf $\cF$ on $X$, we have the unit map $\cF \rightarrow j_* j^* \cF$, which induces
\begin{equation}\label{eq: open restriction}
\RGamma(X, \cF) \rightarrow \RGamma(X, j_* j^* \cF) \cong \RGamma(U, j^* \cF).
\end{equation}
The cohomology of $X$ with supports in $Z$ is defined as the derived kernel of the restriction map \eqref{eq: open restriction}, so that we have an exact triangle 
\[
\RGamma_Z(X; \cF) \rightarrow \RGamma(X; \cF) \rightarrow \RGamma(X; j_* j^*\cF).
\]
For any $\cF = \F_p(b)^{\syn}_X$, $b \in \Z$, we have from \S \ref{part: syntomic Steenrod algebra} an action of the syntomic Steenrod algebra on $\RGamma(X; \cF) \rightarrow \RGamma(U; j^* \cF)$, hence also on $\RGamma_Z(X; \cF)$. 

\subsection{Syntomic Steenrod operations}\label{ssec: syntomic Steenrod summary}
We summarize some of the general theory of the syntomic Steenrod algebra from \S \ref{sec: syntomic Steenrod operations} and \S \ref{sec: comparing operations}, specialized to the case $p=2$. We have explicit cohomology operations:
\begin{enumerate}
\item ``$\beta$'', which acts as the Bockstein differential 
\[
\beta \co \Hsyn^{a,b}(-) \rightarrow \Hsyn^{a+1,b}(-)
\]
induced by the exact triangle $\Z/2(b)^{\syn} \rightarrow \Z/4(b)^{\syn}  \rightarrow \Z/2(b)^{\syn} $,  and 
\item ``$\Sqs^{2i} = \Ps^i $'' for each $i \geq 0$, acting as 
\[
\Hsyn^{a,b}(-) \xrightarrow{\Sqs^{2i}} \Hsyn^{a+2i,b+i}(-).
\]
\end{enumerate}

\begin{example}\label{ex: steenrod operation examples}
The operation $\Sqs^0$ acts as the identity. 

By Corollary \ref{cor: Sqm special cases}, the operation 
\[
\Sqs^{2i} \co \Hsyn^{2i,i}(-) \rightarrow \Hsyn^{4i,2i}(-)
\]
is given by squaring.
\end{example}

Set $\Sqs^{2i+1} := \beta \circ \Sqs^{2i}$. The comultiplication on $\Asyn^{*,*}$ takes the form
\[
\Sqs^{2i}(u \cdot v) = \sum_{j=0}^i \Sqs^{2j}(u) \cdot \Sqs^{2i-2j}(v) 
\]
and 
\[
\Sqs^{2i+1}(u \cdot v) = \sum_{j=0}^i \left( \Sqs^{2j+1}(u) \cdot \Sqs^{2i-2j}(v) + \Sqs^{2j} (u)  \cdot \Sqs^{2i-2j+1} (v) \right) 
\]

\subsection{Construction of syntomic Stiefel--Whitney classes}\label{subsec: SW construction}

Let $i \colon Z \hookrightarrow Y$ be a regular embedding of pure codimension $d$ over $k$. Then we have a \emph{cycle class} $s_{Z/Y} \in \rH^{2d}_Z(Y; \SZp(d))$. In the greater generality of syntomic cohomology (in mixed characteristic settings), this has been constructed by Longke Tang \cite{Tang22}. In the case of smooth varieties over $k$, it follows from results of Milne \cite{Milne86} and Gros \cite{Gros85}.

We are going to apply this with $Y$ being the total space of a vector bundle $\pi \co E \rightarrow X$ of rank $d$, $Z=X$, and $i \co X \hookrightarrow E$ being the zero section. This gives a cycle class $s_{X/E} \in \rH^{2d}_X(E; \SZp(d))$. 

We will define a pushforward map 
\begin{equation}\label{eq: push vb 1}
\pi_* \co \rH^{a}_X(E; \SZp(b)) \rightarrow \Hsyn^{a-2d}(X; \Z_p(b-d)).
\end{equation}
Let $\ol \pi \co  \ol{E} := \bP(E \oplus \cO_X) \rightarrow X$, a compactification of $\pi \co E \rightarrow X$. Note that Mayer--Vietoris gluing implies a derived Cartesian square
\[
\begin{tikzcd}
\RGamma_{\syn}(\ol E; -) \ar[r] \ar[d] & \RGamma_{\syn}(\ol E \setminus X; -) \ar[d] \\ \RGamma_{\syn}(E; -) \ar[r] & \RGamma_{\syn}(E \setminus X; -)
\end{tikzcd}
\]
which induces a canonical isomorphism $\rH^{a}_X(E; \SZp(b)) \cong  \rH^a_X(\ol{E}; \SZp(b))$. Since the map $\ol{\pi} \co \ol{E} \rightarrow X$ is smooth and proper, we have by the construction in \cite[Chapitre II, \S 1]{Gros85} a map 
\begin{equation}\label{eq: push vb 2}
\ol{\pi}_* \SZp(b)_{\ol{E}}  \rightarrow \SZp(b-d)_X[-2d].
\end{equation}
Taking cohomology of \eqref{eq: push vb 2}, and passing through the identification $\rH^{a}_X(E; \SZp(b)) \cong  \rH^a_X(\ol{E}; \SZp(b))$, gives \eqref{eq: push vb 1}. 

Similarly we get a pushforward on mod $p$ syntomic cohomology,
\begin{equation}\label{eq: push vb 3}
\pi_* \co \rH^{a,b}_X(E) \rightarrow \Hsyn^{a-2d,b-d}(X).
\end{equation}

\begin{defn}\label{def: SW class} Suppose $k$ has characteristic $p=2$. Let $\pi \co  E \rightarrow X$ be a vector bundle of rank $d$. Let $\ol{s}_{X/E} \in \rH^{2d,d}_X(E)$ be the reduction modulo 2 of the cycle class of the zero section. For $j \geq 0$, we define the \emph{$j$th syntomic Stiefel--Whitney class} of $E$ to be
\begin{equation}\label{eq: defn w_i}
w_j^{\syn} (E) = \pi_* (\Sqs^j (\ol{s}_{X/E}) ) \in \Hsyn^{j,\lfloor j/2 \rfloor}(X).
\end{equation}
Define the \emph{total Stiefel--Whitney class} to be $w^{\syn}(E) := \sum_j w_j^{\syn}(E)$. If no vector bundle is mentioned, then by default we set $w_j := w_j^{\syn}(TX)$ for the tangent bundle $TX$, and $w^{\syn} := \sum_j w_j^{\syn}$. 
\end{defn}

\begin{remark}It is crucial that in Definition \ref{def: SW class} we use the \emph{syntomic} Steenrod operations $\Sqs^i$ instead of the $\EE_\infty$ Steenrod operations $\Sqe^i$; the classes that would come out of using the latter operations would not be well-behaved (as can be seen just by considering weights -- see Remark \ref{rem: E SW classes}). 
\end{remark}

\subsection{Properties of the syntomic Stiefel--Whitney classes}\label{ssec: SWproperties}

We now record that the syntomic Stiefel--Whitney classes, as constructed in \S \ref{subsec: SW construction}, enjoy the usual properties of characteristic classes. We continue to assume that $p=2$ through the rest of the section. 

\begin{enumerate}
\item If $E$ is a vector bundle on $X$, then we have $w^{\syn}_0(E)=1$ and $w_i^{\syn}(E)=0$ for $i> 2 \rank E$. 

\item ({\sc naturality}) If $f \colon X' \rightarrow X$ is a morphism and $E$ is a vector bundle on $X$, then we have
\[
f^* w^{\syn}_i(E) = w^{\syn}_i(f^* E).
\]
\item ({\sc Whitney sum formula}) If $E,E'$ are vector bundles on $X$, then we have
\begin{equation}\label{eq:whitney-sum-even}
w^{\syn}_{2i}(E \oplus E') = \sum_{j=0}^i w^{\syn}_{2j}(E) \cdot w^{\syn}_{2i-2j}(E')
\end{equation}
and
\begin{equation}\label{eq:whitney-sum-odd}
w^{\syn}_{2i+1}(E \oplus E') = \sum_{j=0}^i \Big(w^{\syn}_{2j+1}(E) \cdot w^{\syn}_{2i-2j}(E') + w^{\syn}_{2j}(E) \cdot w^{\syn}_{2i-2j+1}(E')\Big).
\end{equation}

\end{enumerate}

These follow formally from the properties of Steenrod operations in \S \ref{ssec: syntomic Steenrod summary}, as in \cite[\S 5.4]{Feng20}. 

\begin{lemma}
If 
\[
0 \rightarrow E' \rightarrow E \rightarrow E'' \rightarrow 0
\]
is a short exact sequence of vector bundles on $X$, then we have
\begin{equation}
w^{\syn}_{2i}(E) = \sum_{j=0}^i w^{\syn}_{2j}(E') \cdot w^{\syn}_{2i-2j}(E'')
\end{equation}
and
\begin{equation}
w^{\syn}_{2i+1}(E) = \sum_{j=0}^i \Big(w^{\syn}_{2j+1}(E') \cdot w^{\syn}_{2i-2j}(E'') + w^{\syn}_{2j}(E') \cdot w^{\syn}_{2i-2j+1}(E'')\Big).
\end{equation}
\end{lemma}

\begin{proof}
Applying \cite[Proposition 9.2.9]{BL22a} in the same way as in the proof of \cite[Theorem 9.2.7]{BL22a}, we may reduce to the case where $E \cong E' \oplus E''$, which was handled above.
\end{proof}

\subsection{Relation to Chern classes} Chern classes in syntomic cohomology are defined in \cite[\S 9]{BL22a}. For a vector bundle $E \rightarrow X$, there are defined in \cite[Construction 9.2.1]{BL22a} \emph{syntomic Chern classes} $c_j^{\syn}(E) \in \Hsyn^{2j}(X; \Z_p(j))$. 

Here we prove that our syntomic Stiefel--Whitney classes are simply reductions of Chern classes. This fact will be of significance later in \S \ref{sec: final}, where we will want to know that our syntomic Stiefel--Whitney classes lift to integral cohomology. The reason that we have defined them by this complicated construction involving syntomic Steenrod operations, rather than simply defining them to be the reduction of Chern classes, is that our definition will interface well with later constructions.

\begin{prop}\label{prop: chern classes}
Let $X$ be a smooth proper variety over $\F_2$ and $E$ a vector bundle on $X$ of rank $r$. For each $j \in \N$, let $c^{\syn}_j(E) \in \Hsyn^{2j}(X; \Z_2(j))$ be the $j^{\mrm{th}}$ Chern class of $E$ and $\ol{c}^{\syn}_j(E) \in  \Hsyn^{2j,j}(X)$ be its reduction modulo $2$. Then we have:
\begin{equation}\label{eq: SW formula}
w_i^{\syn}(E) := \begin{cases} \ol{c}^{\syn}_{i/2}(E) & \text{$i$ even}, \\ 0 & \text{$i$ odd}.\end{cases}
\end{equation}
\end{prop}

\begin{remark}
The form of Proposition \ref{prop: chern classes} is actually \emph{simpler} than the analogous statement when $\ell \neq p$ \cite[Theorem 5.10]{Feng20}. This can be ultimately traced to the difference between $\Sqs^1$ and the Bockstein operation used in that case, which differ by $\rho$. 
\end{remark}

\begin{proof}
By standard reductions for characteristic classes, it suffices to check that the formula above is correct for all line bundles.  

Let $\pi \co L \rightarrow X$ be a line bundle on $X$. We view $X$ as embedded in $L$ via the zero section. We then have the associated cycle class $s_{X/L} \in \rH^2_X(L; \Z_2^{\syn}(1))$.\\

\noindent \textbf{Calculation of $w^{\syn}_1$.} To show that $w^{\syn}_1(L) =0$, it suffices to show that $\Sqs^1 (\ol{s}_{X/L}) = 0$. We have $[\ol{s}_{X/L}] \in \rH^{2,1}_X(L)$ so the operation $\Sqs^1$ is the Bockstein map for the exact triangle
\[
\Z/2^{\syn}(1) \rightarrow \Z/4^{\syn} (1)\rightarrow \Z/2^{\syn}(1)
\]
But the cycle class $\ol{s}_{X/L}$ even lifts to $\rH^2_X(L; \Z_2^{\syn}(1))$, so this Bockstein map vanishes. \\

\noindent \textbf{Calculation of $w^{\syn}_2$.} We claim that $\Sqs^2 (\ol s_{X/L}) = \pi^* (\ol{c}^{\syn}_1(L)) \cdot \ol{s}_{X/L}$. Granting this claim, we deduce that 
\begin{align*}
w^{\syn}_2(L)  &= \pi_* (\Sqs^2 (\ol{s}_{X/L}) ) = \pi_* ( \pi^* (\ol{c}_1(L)) \cdot \ol{s}_{X/L}) = \ol{c}_1^{\syn}(L) 
\end{align*}
using the projection formula in the last equality.

It remains to prove the claim. By Example \ref{ex: steenrod operation examples}, we have that $\Sqs^2 (\ol{s}_{X/L}) = \ol{s}_{X/L} \cdot \ol{s}_{X/L}$. Hence it suffices to show that the map $\rH^{2,1}_X(L) \rightarrow \Hsyn^{2,1}(L)$ sends $\ol{s}_{X/L}$ to $\pi^* \ol{c}^{\syn}_1(L)$. 

Consider the commutative diagram 
\[
\begin{tikzcd}
\Pic(X)  = \Het^1(X; \G_m) \ar[r] \ar[d, "\pi^*"] & \Hsyn^{2,1}(X) \ar[d, "\pi^*"] \\
\Pic(L) = \Het^1(L; \G_m) \ar[r] & \Hsyn^{2,1}(L)
\end{tikzcd}
\]
where the horizontal arrow is the association of the first Chern class. The line bundle $L \rightarrow X$ pulls back to $\cO_L(X)$ on $L$, i.e., the line bundle associated to the divisor of the zero-section in $L$. Hence the cycle class of the zero section in $L$ coincides with $c^{\syn}_1(\pi^* L) = \pi^* c^{\syn}_1(L)$. This completes the proof.
\end{proof}

\begin{remark}\label{rem: E SW classes}
We could have constructed ``$\EE_{\infty}$ Stiefel--Whitney classes'' $w^{\EE}_i$ in an analogous way, using the $\Sqe^i$ instead of $\Sqs^i$. Elementary weight considerations reveal that if $E \rightarrow X$ has rank $r$ then $w^{\EE}_i(E) \in \Hsyn^{i,r}(X)$. In particular, for $i \neq 2r$ they cannot lie on the ``motivic line'' $\Hsyn^{2*,*}$. Moreover, if $r$ is large compared to $\dim X$, then they all vanish (even for $i=0$!), which indicates that their behavior differs from what would be expected of characteristic classes. 
\end{remark}

\section{Arithmetic Wu formula}\label{sec: Wu}

In this section we prove an \emph{arithmetic Wu formula}, which calculates certain syntomic Steenrod operations in terms of the Stiefel--Whitney classes just defined.

The adjective ``arithmetic'' refers to the absolute nature of our cohomology theory, which has a non-trivial contribution from the base field. We note that \emph{geometric} Wu formulas have been proved in \cite{Urabe96, SS24, Ben24} (for $\ell$-adic \'etale cohomology, $\ell \neq p$) and \cite{Pri20, AE} (for mod $p$ motivic cohomology in characteristic $p$).
 Although the final formulations look the same, our arithmetic Wu formula is much subtler to prove. The fact that even a point has non-trivial cohomology, while being geometrically trivial, is the root of all our problems, and in fact the entire purpose of Part \ref{part: spectral prismatization} was to provide a technical ingredient to surmount this difficulty.\footnote{The analogous problem for \emph{arithmetic} $\ell$-adic \'etale cohomology was solved in \cite{Feng20}, using some tricks with \'etale homotopy theory, but those do not apply here.}

\subsection{The syntomic Wu classes} Throughout this section, we let $X$ be a smooth, proper, and geometrically connected variety over a finite field $k$ of characteristic $p$, of dimension $d$. For $p=2$, consider the syntomic Steenrod operation 
\[
\Sqs^i \co \Hsyn^{2d+1-i, d-\lfloor i/2 \rfloor}(X) \rightarrow \Hsyn^{2d+1,d}(X).
\]

Recall that the cup product 
\[
\Hsyn^{2d+1-i, d-\lfloor i/2 \rfloor}(X) \otimes \Hsyn^{i,\lfloor i/2 \rfloor}(X) \rightarrow \Hsyn^{2d+1,d}(X)
\]
induces a perfect pairing. Therefore, there exists a unique $v^{\syn}_i \in \Hsyn^{i,\lfloor i/2 \rfloor}(X)$ such that 
\[
\Sqs^i (\alpha) = v^{\syn}_i \cdot \alpha \text{ for all $\alpha \in \Hsyn^{2d+1-i,d-\lfloor i/2 \rfloor}(X)$}.
\]
We call $v^{\syn}_i$ the \emph{$i^{\mrm{th}}$ syntomic Wu class} and we call 
\[
v^{\syn} := \sum_i v^{\syn}_i \in \Hsyn^{*,*}(X) 
\]
the \emph{total syntomic Wu class}.

\begin{example}
If $i = 0$, then $\Sqs^i $ is the identity map, so $v^{\syn}_i = 1$. 

If $i > d+1$, then \Cref{cor: Sqm special cases} implies that $\Sqs^i(\alpha) =0$ on all $\alpha \in \Hsyn^{2d+1-i,d-\lfloor i/2 \rfloor}(X)$, so $v^{\syn}_i =0$. 
\end{example}

The purpose of this section is to prove the following theorem. 

\begin{thm}[Arithmetic Wu formula]\label{thm: syntomic Wu theorem} Let $X$ be a smooth, proper, geometrically connected variety over a finite field of characteristic $2$. Then we have 
\begin{equation}\label{eq:arithmetic-wu-even}
w^{\syn}_{2i}(TX) = \sum_{j=0}^i \Sqs^{2j} (v^{\syn}_{2i-2j}) \in \Hsyn^{*,*}(X). 
\end{equation}
and
\begin{equation}\label{eq:arithmetic-wu-odd}
w^{\syn}_{2i+1}(TX) = \sum_{j=0}^i \Big(\Sqs^{2j+1} (v^{\syn}_{2i-2j}) + \Sqs^{2j}(v^{\syn}_{2i-2j+1}) \Big)\in \Hsyn^{*,*}(X). 
\end{equation}
\end{thm}

\begin{remark}
Since $\Sqs^0 = \Id$ and $\Sqs^i$ is obviously nilpotent for $i>0$ whenever the cohomological dimension of $X$ is finite, equations \eqref{eq:arithmetic-wu-even} and \eqref{eq:arithmetic-wu-odd} can be inverted to solve for $\{v^{\syn}_j\}$ in terms of the $\{w_j^{\syn}\}$.

The analogous definition for the $\EE_{\infty}$ Steenrod operations would lead to an ``$\EE_\infty$ Wu class'' $v^{\EE}$ which is different from $v^{\syn}$, as can be seen from weight considerations. 
Moreover, if $X$ has dimension $d>0$ then $w^{\EE}_0 (TX)\in \Hsyn^{0,d}(X)$ would vanish, so that the analogous ``arithmetic Wu formula'' could not be inverted to calculate the $\{v^{\EE}_j\}$ in terms of the $\{w^{\EE}_j(TX)\}$. This is another reason we use the syntomic Steenrod operations instead of the $\EE$ Steenrod operations for the purpose of defining characteristic classes.  

\end{remark}

\subsection{Deformation to the diagonal} In this subsection, $p = \mrm{char}(k)$ can be an arbitrary prime. 

\subsubsection{More on pushforward maps}\label{sssec: pushforwards} 
Let $f \co X \rightarrow Y$ be a proper map of connected smooth varieties over a finite field $k$. Set $d := \dim X - \dim Y$. Then there is a pushforward map \cite[Chapitre II, \S 1]{Gros85}
\begin{equation}\label{eq: pushforward}
f_* \SZp(b)_X \rightarrow \SZp(b-d)_Y[-2d].
\end{equation}
It is dual to the pullback map $\SZp(b)_Y \rightarrow f_* \SZp(b)_X$. We also write 
\begin{equation}
f_* \co \Hsyn^a(X; \Z_p(b)) \rightarrow \Hsyn^{a-2d}(Y; \Z_p(b-d))
\end{equation}
for the induced map on cohomology, which is dual to the pullback map on cohomology. 


It is elementary to show (see for example \cite[Chapitre II, \S 2]{Gros85}) that the map \eqref{eq: pushforward} enjoys the following properties: 
\begin{itemize}
\item It is functorial: if $f \co X \rightarrow Y$ is a proper map of smooth varieties over $k$, and $g \co Y \rightarrow Z$ is another proper map of smooth varieties over $k$, then 
\[
(gf)_* = g_* \circ f_*.
\] 
\item We have the \emph{projection formula} for all $\alpha  \in \Hsyn^{*,*}(X; \Z_p)$ and $\beta \in \Hsyn^{*,*}(Y; \Z_p)$,
\begin{equation}\label{eq: product formula}
f_*(\alpha \cdot f^* \beta) = (f_* \alpha) \cdot \beta.
\end{equation}
\item If $f \co X \hookrightarrow Y$ is a regular embedding of codimension $a$, let $\mrm{cl}_Y(X)$ be the image of $s_{X/Y} \in \rH_X^{2a}(Y; \Z_p^{\syn}(a))$ under the ``forget supports'' map to $\rH^{2a}(Y; \Z_p^{\syn}(a))$. Then $f_* f^*$ is multiplication by $\mrm{cl}_Y(X)$. In particular, $f_* (1) = \mrm{cl}_Y(X)$ where the input element $1$ is regarded in $\Hsyn^0(X; \Z_p) \cong \Z_p$.
\end{itemize}

By descent, if $f \co \cX \rightarrow \cY$ is a (representable) smooth and proper map of stacks, then there exists a pushforward map 
\begin{equation}\label{eq: pushforward on cohomology stacks}
f_* \co \Hsyn^a(\cX; \Z_p(b)) \rightarrow \Hsyn^{a-2d}(\cY; \Z_p(b-d))
\end{equation}
satisfying the same properties.

\subsubsection{Weighted deformation to the normal cone}\label{sssec:weighted-deformation-normal-cone}
In $\A^1$-invariant cohomology theory, Lemma \ref{lem: step 1} may be proved by deformation to the normal cone (using Morel--Voevodsky purity). However, syntomic cohomology theory is \emph{not} $\A^1$-invariant, which causes substantial technical complications for us. However, it enjoys a weaker property called \emph{weighted homotopy invariance} that turns out to be sufficient for our purposes.\footnote{We learned of this from \cite[\S 3]{Tang22}, which credits the idea to course notes of Dustin Clausen, who credits it to Burt Totaro.}

More generally, let $\iota \co X \inj Y $ be a regular embedding. The \emph{weighted deformation to the normal cone} \cite[Definition 5.11]{Tang22} of $\iota$ is a flat family $\cY \rightarrow [\A^1/\G_m]$ (the stack quotient for the standard scaling action) whose restriction to $\pt = [\G_m/\G_m] \inj [\A^1/\G_m]$ is $Y \rightarrow \pt$ and whose restriction to $[0/\G_m] \inj [\A^1/\G_m]$ is the stack quotient of the normal bundle $N_{\iota}$ of $\iota$ by the inverse scaling action of $\G_m$. Moreover, $\cY$ is equipped with a closed embedding from $\cX := X \times [\A^1/\G_m]$, which restricts over $\pt = [\G_m/\G_m] \inj [\A^1/\G_m]$ to the given $\iota \co X \inj Y$ and over $[0/\G_m] \inj [\A^1/\G_m]$ to the zero section of $X$ in $N_\iota$. We write $\cX_0, \cY_0$ for the special fibers of $\cX, \cY$, respectively, over $[0/\G_m] \inj [\A^1/\G_m]$. This is summarized in the diagram below. 
\begin{equation}\label{eq: weighted deformation}
\begin{tikzcd}
\cX_0 = X/\G_m \ar[d, hook] \ar[r, hook] & \cX = X \times \A^1/\G_m \ar[d, hook] & X  \ar[l, hook'] \ar[d,hook, "\iota"] \\
\cY_0 = N_{\iota}/\G_m \ar[d] \ar[r, hook] & \cY  \ar[d]  & Y  \ar[l, hook'] \ar[d] \\
0/\G_m \ar[r, hook] & \A^1/\G_m
 & \G_m/\G_m \ar[l, hook'] 
\end{tikzcd}
\end{equation}
Then it follows from \cite[Theorem 5.20]{Tang22} that we have the following ``weighted homotopy invariance'': the obvious pullback maps induce a commutative diagram with indicated isomorphisms, 
\begin{equation}\label{eq: weighted deformation cohog}
\begin{tikzcd}
\rH^{*,*}_{\cX_0}(\cY_0) \ar[d]  & \rH^{*,*}_{\cX}(\cY) \ar[d] \ar[l,"\sim"] \ar[r] &  \rH_X^{*,*}(Y) \ar[d] \\
\Hsyn^{*,*}(\cY_0)  & \Hsyn^{*,*}(\cY) \ar[l,"\sim"] \ar[r] & \Hsyn^{*,*}(Y)
\end{tikzcd}
\end{equation}
Furthermore, the construction of cycle classes is arranged so that \eqref{eq: weighted deformation cohog} has the following effect on cycle classes \cite[\S 5]{Tang22}:
\begin{equation}\label{eq: deform cycle classes}
\begin{tikzcd}[row sep = tiny]
\rH^{2*,*}_{\cX_0}(\cY_0)  & \rH^{2*,*}_{\cX}(\cY)\ar[l,"\sim"] \ar[r] & \rH^{2*,*}_X(Y) \\
s_{\cX_0/\cY_0} & s_{\cX/\cY} \ar[l] \ar[r] & s_{X/Y} 
\end{tikzcd}
\end{equation}

\subsubsection{The case of the diagonal embedding}\label{sssec:deformation-diagonal-embedding} Now assume that $X$ is smooth and proper of (equi)dimension $d$ over $k$. Let $Y := X \times_k X$ and $\iota \co X \inj Y$ be the diagonal embedding. Note that we have a commutative diagram 
\begin{equation}\label{eq: mod diagram 1}
\begin{tikzcd}
\cY_0 = N_{\iota}/\G_m \ar[d, "\pi"] \ar[r, hook] & \cY  \ar[d, "\wt{\pr}"]  & Y  \ar[l, hook'] \ar[d, "\pr_1"] \\ 
\cX_0 = X/\G_m \ar[r, hook] & \cX & X  \ar[l, hook'] 
\end{tikzcd}
\end{equation}
which provides a retraction to the upper row of vertical arrows in \eqref{eq: weighted deformation}. We will construct a commutative diagram of pushforward maps 
\begin{equation}\label{diag: pushforward}
\begin{tikzcd}
\rH^{*,b}_{\cX_0}(\cY_0) \ar[d, "\pi_*"]  & \rH^{*,b}_{\cX}(\cY) \ar[d, "\wt{\pr}_*"] \ar[l,"\sim"] \ar[r] & \rH^{*,b}_X(Y) \ar[d, "\pr_{1*}"] \\
\Hsyn^{*-2d,b-d}(\cX_0)  & \Hsyn^{*-2d,b-d}(\cX) \ar[l,"\sim"] \ar[r]  & \Hsyn^{*-2d,b-d}(X)
\end{tikzcd}
\end{equation}

\begin{itemize}
\item Since $\pr_{1*}$ is smooth and proper, we have on general grounds (\S \ref{sssec: pushforwards}) a pushforward map $\pr_{1*} \co \Hsyn^{*,b}(Y) \rightarrow \Hsyn^{*-2d,b-d}(X)$. Abusing notation, we define the vertical map $\pr_{1*}$ in the right column of \eqref{diag: pushforward} to be its composition with the ``forget support'' map $\rH^{*,b}_{X}(Y) \rightarrow \Hsyn^{*,b}(Y)$. 

\item For $C := \bP_{\cX}(N_{X/Y} \oplus \cO)$, let $\cC$ be the quotient stack $C/\G_m$, where $\G_m$ acts on $N_{X/Y} \oplus \cO$ via inverse scaling. As $\pi$ is a vector bundle, the vertical map $\pi_*$ in the left column of \eqref{diag: pushforward} is defined by the procedure of \S \ref{subsec: SW construction}: letting $\ol{\pi} \co \cC \rightarrow \cX_0$ be the projection for the compactification, we define $\pi_*$ to be the isomorphism $\rH^{*,b}_{\cX_0}(\cY_0) \cong \rH^{*,b}_{\cX_0}(\cC)$ composed with the ``forget supports'' map $\rH^{*,b}_{\cX_0}(\cC) \rightarrow \Hsyn^{*,b}(\cC)$ followed by the pushforward map $\Hsyn^{*,b}(\cC) \rightarrow \Hsyn^{*-2d,b-d}(\cX_0)$.

\item Then the vertical map $\wt{\pr}_*$ in the middle column of \eqref{diag: pushforward} is determined uniquely by the commutativity of the left square. 

\end{itemize}

\begin{lemma}\label{lem: pushforwards compatible} With the definitions above, diagram \eqref{diag: pushforward} commutes.
\end{lemma}

\begin{proof} 
Although the Thom isomorphism fails for syntomic cohomology since it is not $\A^1$-invariant, there is a ``weighted Thom isomorphism'' \cite[Theorem 4.2]{Tang22} which says that $\rH^{*,*}_{\cX_0}(\cY_0)$ is free of rank one over $\Hsyn^{*,*}(\cX_0)$, with distinguished generator being the Thom class $\Th_{TX} \in \rH^{2d,d}_{\cX_0}(\cY_0)$. All cohomology groups in question are modules over $\Hsyn^{*,*}(\cX_0)$ under pullback and cup product, and the horizontal maps (being pullbacks) are obviously linear over $\Hsyn^{*,*}(\cX_0)$. Moreover, the pushforward maps $\pi_*$ and $\pr_{1*}$ are also linear over $\Hsyn^{*,*}(\cX_0)$: this is a reformulation of the projection formula. So it suffices to check the commutativity on the Thom class $\Th_{TX}$. 

By the very construction of cycle classes in \cite[\S 5]{Tang22}, the diagonal cycle class in $\mrm{cl}_X(X \times X)  \in \Hsyn^{2d,d}(X \times X)$ is the image of $\Th_{TX}$ from the top left to the bottom right term\footnote{meaning the inverse of the leftwards restriction, composed with the rightwards restriction, followed by the ``forget supports'' map} of \eqref{eq: weighted deformation cohog}. By the uniqueness properties of cycle classes, this definition of $\mrm{cl}_X(X \times X)$ agrees with the pushforward of $1 \in \Hsyn^{0,0}(X)$ under the diagonal map $X \inj X \times X$. Then from the functoriality of the pushforward, we see that 
\[
\pr_{1*} (\mrm{cl}_X(X \times X))  = 1 \in \Hsyn^{0,0}(X).
\]
On the other hand, by the general property of Thom classes we have $\pi_* (\Th_{TX}) = 1 \in \Hsyn^{0,0}(\cX_0) $, which maps to $1 \in \Hsyn^{0,0}(X)$ under the composition from left to right (inverting the middle isomorphism)
\[
 \Hsyn^{*,*}(\cX_0) = \Hsyn^{*,*}(X \times \rB\G_m) \xleftarrow{\sim} \Hsyn^{*,*}(X \times \A^1/\G_m) \rightarrow \Hsyn^{*,*}(X),
\]
since $\Hsyn^{*,*}(X \times \rB\G_m) \cong \Hsyn^{*,*}(X)[c^{\syn}_1]$ by \cite[Lemma 9.3.2]{BL22a}, with the restriction map to $\Hsyn^{*,*}(X)$ being the quotient by $c^{\syn}_1$. Thus we see the desired commutativity. 
\end{proof}

\subsection{Calculation of syntomic Wu classes}  
For the rest of this section, $p=2$. Recall that $X$ is a smooth, proper, geometrically connected variety of dimension $d$ over a finite field $k$ of characteristic $2$. 

\subsubsection{Stiefel--Whitney classes of the tangent bundle} The normal bundle of $X$ in its diagonal embedding into $X \times X$ is isomorphic to the tangent bundle $TX$. Comparing this to the definition of the syntomic Stiefel--Whitney classes in \S \ref{subsec: SW construction} motivates the following. 

\begin{lemma}\label{lem: step 1}
Let $\ol{s}_{X/X\times X} \in \rH^{2d,d}_X(X \times X)$ be the mod $2$ cycle class of the diagonal. Then we have 
\begin{equation}\label{WuProof1}
\pr_{1*} \Sqs^i (\ol{s}_{X/X\times X}) =  w^{\syn}_i \in \Hsyn^{i,\lfloor i/2 \rfloor}(X),
\end{equation}
where $\mrm{pr}_1 \co X \times X\rightarrow X$ denotes projection to the first factor.
\end{lemma}

\begin{proof} We will use the notation of \S\ref{sssec:weighted-deformation-normal-cone} and \S \ref{sssec:deformation-diagonal-embedding}. By naturality of Steenrod operations under pullback, applying $\Sqs^i$ to \eqref{eq: deform cycle classes} shows that the maps
\[
\begin{tikzcd}
\rH^{2d+i,d+ \lfloor i/2 \rfloor}_{\cX_0}(\cY_0)  & \rH^{2d+i,d + \lfloor i/2 \rfloor}_{\cX}(\cY) \ar[l,"\sim"] \ar[r] & \rH^{2d+i,d + \lfloor i/2 \rfloor}_X(Y)
\end{tikzcd}
\]
carry
\[
\begin{tikzcd}
\Sqs^i (\ol s_{\cX_0/\cY_0}) & \Sqs^i (\ol s_{\cX/\cY}) \ar[l] \ar[r] & \Sqs^i (\ol s_{X/Y} )
\end{tikzcd}
\]
Then by Lemma \ref{lem: pushforwards compatible}, we see that applying the pushforward maps from \eqref{diag: pushforward} carry cohomology classes as in the commutative diagram below
\begin{equation}\label{eq: def diagram 1}
\begin{tikzcd}
\Sqs^i (\ol s_{\cX_0/\cY_0}) \ar[d, "\pi_*"]  & \Sqs^i (\ol s_{\cX/\cY}) \ar[l] \ar[r] \ar[d, "\wt{\pr}_*"]  & \Sqs^i (\ol s_{X/Y} ) \ar[d, "\pr_{1*}"] \\
\pi_* \Sqs^i (\ol s_{\cX_0/\cY_0}) & \wt{\pr}_* \Sqs^i (\ol s_{\cX/\cY}) \ar[l, "\mrm{inc}_0^*"'] \ar[r, "\mrm{inc}_{\G_m}^*"] & \pr_{1*}\Sqs^i (\ol s_{X/Y} ) 
\end{tikzcd}
\end{equation}
where the maps $\mrm{inc}_0^*$ and $\mrm{inc}_{\G_m}^*$ are the pullbacks
\[
\begin{tikzcd}
\Hsyn^{i,\lfloor i/2 \rfloor}(X/\G_m)  & \ar[l, "\mrm{inc}_0^*"']  \Hsyn^{i,\lfloor i/2 \rfloor}(X \times \A^1/\G_m) \ar[r, "\mrm{inc}_{\G_m}^*"]  &  \Hsyn^{i, \lfloor i/2 \rfloor}(X)
\end{tikzcd}
\]

Consider the Cartesian square 
\[
\begin{tikzcd}
TX \ar[d, "\pi"] \ar[r, "q'"] & \cY_0 \ar[d, "\pi_0"] \\
X \ar[r, "q"] & \cX_0 
\end{tikzcd}
\]
where the maps $q$ are the quotient by $\G_m$. We have $\Sqs^i(\ol s_{\cX_0/\cY_0}) = \pi_0^*(w_i(\cY_0)) s_{\cX_0/\cY_0}$, hence by naturality of Thom classes and Stiefel--Whitney classes, we obtain that
\begin{equation}\label{eq: def diagram 2}
q^* \pi_{0*} \Sqs^i(\ol s_{\cX_0/\cY_0})  = q^* (w_i^{\syn}(\cY_0)) =   w_i^{\syn}(TX). 
\end{equation}

We explained earlier that $\mrm{inc}_0^*$ is an isomorphism. Comparing \eqref{eq: def diagram 1} and \eqref{eq: def diagram 2}, we see that it suffices to show that $\mrm{inc}_{\G_m}^* \circ (\mrm{inc}_0^*)^{-1} = q^*$. Since the projection map $\pr_1 \co X \times \A^1/\G_m \rightarrow X / \G_m$ satisfies $\pr_1 \circ \mrm{inc}_0 = \Id$, and $\mrm{inc}_0^*$ is an isomorphism, we must have $\pr_1^* \cong (\mrm{inc}_0^*)^{-1}$. This implies that
\[
\mrm{inc}_{\G_m}^* \circ (\mrm{inc}_0^*)^{-1} = \mrm{inc}_{\G_m}^* \circ \pr_1^* = (\pr_1 \circ \mrm{inc}_{\G_m})^* = q^*
\]
as desired.   
\end{proof}

\subsubsection{The map $\varphi_*$} Let 
\begin{equation}\label{eq: push tensor}
\varphi_* \colon  \Hsyn^{*,*}(X \times_k X) \rightarrow \Hsyn^{*,*}(X) \otimes_{\F_p} \Hsyn^{*,*}(X)
\end{equation}
be the map studied at the beginning of \S \ref{ssec: duality}. 

We will now discuss features specific to characteristic $p = 2$. We have the syntomic Steenrod operations $\Sqs$ on $\Hsyn^{*,*}(X \times_k X)$, and also on $\Hsyn^{*,*}(X)$. Although $\Hsyn^{*,*}(X) \otimes_{\F_p} \Hsyn^{*,*}(X)$ is not the syntomic cohomology of a variety over $k$, it admits a natural $\Asyn^{*,*}$-module structure because of the Hopf algebra structure of $\Asyn^{*,*}$, as discussed in \S \ref{ssec: coproduct}. Concretely, for $x \otimes y \in \Hsyn^{*,*}(X) \otimes_{\F_p} \Hsyn^{*,*}(X)$, we have 
\[
\Sqs^{2i}(x \otimes y) = \sum_{j=0}^i \Sqs^{2j}(x) \otimes \Sqs^{2i-2j}(y) 
\]
and 
\[
\Sqs^{2i+1}(x \otimes y) = \sum_{j=0}^i \left( \Sqs^{2j+1}(x) \otimes \Sqs^{2i-2j}(y) + \Sqs^{2j} (x)  \otimes \Sqs^{2i-2j+1} (y) \right).
\]

From Theorem \ref{thm: steenrod push equivariant}, we have that 
\begin{equation}\label{eq: commutes with Sqm}
\Sqs \varphi_* = \varphi_* \Sqs.
\end{equation}

We give another perspective on the map $\varphi_*$ from \eqref{eq: push tensor}. Note that $\Hsyn^{*,*}(X \times_k X)$ acts by correspondences on $\Hsyn^{*,*}(X)$, inducing an $\F_p$-vector space map  
\begin{equation}\label{eq: end}
\Hsyn^{*,*}(X \times_k X) \rightarrow \End_{\F_p} \left( \Hsyn^{*,*}(X) \right)
\end{equation}
which sends $\alpha \in \Hsyn^{*,*}(X \times_k X)$ to the endomorphism
\begin{equation}\label{eq: correspondence}
\Hsyn^{*,*}(X ) \ni u \mapsto (\pr_1)_* ( \alpha \cdot \pr_2^* u ) \in \Hsyn^{*,*}(X ).
\end{equation}
Using Poincar\'e duality, we may identify 
\[
\End_{\F_p}  \left( \Hsyn^{*,*}(X) \right) \cong \Hsyn^{*,*}(X)^\vee \otimes_{\F_p} \Hsyn^{*,*}(X) \cong \Hsyn^{*,*}(X) \otimes_{\F_p} \Hsyn^{*,*}(X).
\]
After unraveling the definitions, one sees that under this identification,  \eqref{eq: end} agrees with \eqref{eq: push tensor}.

\begin{lemma}\label{diag_is_identity}
Let $\Delta := \cla_{X \times X}(X) \in \Hsyn^{2d,d}(X \times_k X)$. Then the map \eqref{eq: end} sends $\Delta \mapsto \Id$. 
\end{lemma}

\begin{proof}
Let $f \co X\hookrightarrow X\times_k X$ denote the diagonal embedding. Then we may view $\Delta = f_*(1)$ as explained in \S \ref{sssec: pushforwards}. Taking $\alpha = \Delta$ in \eqref{eq: correspondence}, we see that $\Delta$ corresponds to the endomorphism 
\[
\gamma \mapsto (\pr_1)_* (f_*(1) \cdot \pr_2^* \gamma)  = (\pr_1)_* f_* (1 \cdot f^* \pr_2^* \gamma) 
\]
by the projection formula. But since $\pr_1 \circ f = \pr_2 \circ f = \Id_X$, this last expression is just $\gamma$ again. 
\end{proof}

\subsubsection{Decomposition of the diagonal} We will now compute $\Sq(\varphi_* \Delta)$ explicitly and use it to prove \Cref{thm: syntomic Wu theorem}.

\begin{lemma}\label{lem: Delta_basis}
Let $X$ be a smooth, proper variety over a finite field (of any characteristic). Let 
\begin{itemize}
\item $\{e_m\}$ be a basis for $ \Hsyn^{*,*}(X)$ which is homogeneous in bi-degree, and 
\item $\{f_m\}$ be the dual basis of $\Hsyn^{*,*}(X)$ under Poincar\'e duality. 
\end{itemize}
For a cohomology class $x \in \Hsyn^{i,*}(X)$ let $|x| = i$. Then, letting $\Delta$ be as in Lemma \ref{diag_is_identity}, we have
\begin{equation}\label{eq: push of delta}
\varphi_* \Delta = \sum_m (-1)^{|e_m|} e_m \otimes  f_m \in \Hsyn^{*,*}(X) \otimes \Hsyn^{*,*}(X),
\end{equation}
where $\varphi_*$ is as in \eqref{eq: push tensor}.
\end{lemma}

\begin{proof}
Lemma \ref{diag_is_identity} says that the action of $\Delta$ induced on $\Hsyn^*(X)$ by \eqref{eq: correspondence} is just the identity map. Therefore, it suffices to show that the  right hand side of \eqref{eq: push of delta} acts as the identity on $\Hsyn^*(X)$, which is a straightforward linear algebra exercise about dual bases in graded vector spaces (cf. \cite[Theorem 11.11]{MS74}). 
\end{proof}

Let $(\pr_1')_*$ and $(\pr_2')_*$ denote the ``pushforward'' maps 
\[
\Hsyn^{*,*}(X)  \xleftarrow{(\pr_1')_*}  \Hsyn^{*,*}(X) \otimes_{\F_p} \Hsyn^{*,*}(X)  \xrightarrow{(\pr_2')_*}   \Hsyn^{*,*}(X)
\]
which are dual to the obvious ``pullbacks''
\[
\Hsyn^{*,*}(X)  \xrightarrow{(\pr_1')^*}  \Hsyn^{*,*}(X) \otimes_{\F_p} \Hsyn^{*,*}(X)  \xleftarrow{(\pr_2')^*} \Hsyn^{*,*}(X)  .
\]
Note that $(\pr_1')_*$ and $(\pr_1')^*$ are not induced by a map of varieties $\pr_i'$, since indeed $\Hsyn^{*,*}(X) \otimes_{\F_p} \Hsyn^{*,*}(X)$ is not the cohomology of a variety; we are just introducing formal notation. Concretely, $(\pr_1')_*$ and $(\pr_2')_*$ are characterized by the identities 
\begin{align}\label{eq:pr'-proj-formula}
(\pr_1')_* (\pr_1'^* x \cdot \pr_2'^* y) &= x \cdot \left( \int_X y \right)  \quad \text{ for all $x, y \in \Hsyn^{*,*}(X)$}, \\
(\pr_2')_* (\pr_1'^* x \cdot \pr_2'^* y) &= \left( \int_X x \right) \cdot y  \quad \text{ for all $x, y \in \Hsyn^{*,*}(X)$},
\end{align}
where $\int_X$ is the projection to $\Hsyn^{2d+1, d}$ followed by the trace map.

\begin{proof}[Proof of Theorem \ref{thm: syntomic Wu theorem}]
Since the pullback $\Hsyn^{*,*}(X) \xrightarrow{\pr_1^*} \Hsyn^{*,*}(X \times_k X)$ obviously factors as 
\[
\Hsyn^{*,*}(X) \xrightarrow{(\pr_1')^*} \Hsyn^{*,*}(X) \otimes_{\F_p} \Hsyn^{*,*}(X)  \xrightarrow{\pr_1^* \cdot  \pr_2^*}  \Hsyn^{*,*}(X \times X)
\]
(morally, ``$\pr_1 =   \pr_1' \circ \varphi$'') we have a corresponding factorization of the pushforward maps as
\begin{equation}\label{eq: step 3 factor push}
(\pr_1)_* = (\pr_1')_* \varphi_*. 
\end{equation}
Combining Lemma \ref{lem: step 1} and  \eqref{eq: step 3 factor push}, we find that
\begin{equation}\label{eq: step 3 eqn2}
w^{\syn} = (\pr_1)_* \Sqs (\Delta) = (\pr_1')_* \varphi_* \Sqs (\Delta) \in \Hsyn^{*,*}(X)
\end{equation}
where $\Sqs := \Sqs^0 + \Sqs^1 + \ldots$ is the total syntomic Steenrod operation. We saw in \eqref{eq: commutes with Sqm} that $\varphi_* \Sqs = \Sqs \varphi_*$. Let $\{e_m\}$ and $\{f_m\}$ be dual bases for $\Hsyn^{*,*}(X)$. Then by Lemma \ref{lem: Delta_basis}, we may rewrite \eqref{eq: step 3 eqn2} as\footnote{Omitting signs for simplicity, since we are working in characteristic $2$.s}
\begin{equation}\label{eq: step 3 eqn3}
w^{\syn} = (\pr_1')_* \Sqs \left(\sum_m  (\pr_1')^* e_m \cdot  (\pr_2')^* f_m \right).
\end{equation}
By the Cartan formula for $\Sqs$ and the ``projection formula for $(\pr_1')_*$'' \eqref{eq:pr'-proj-formula}, we have 
\begin{align}\label{eq: step 3 eqn4}
(\pr_1')_* \Sqs^{2i} \left(\sum_m  (\pr_1')^* e_m \cdot  (\pr_2')^* f_m \right) &=  \sum_{j=0}^i  \sum_m   (\pr_1')_* \left( (\pr_1')^* \Sqs^{2j} (e_m) \cdot (\pr_2')^* \Sqs^{2i-2j}(f_m) \right) \nonumber\\
&= \sum_{j=0}^i \sum_m  \left( \Sqs^{2j} (e_m)  \int_X \Sqs^{2i-2j}(f_m) \right) .
\end{align}
Combining \eqref{eq: step 3 eqn3} and \eqref{eq: step 3 eqn4} and using that $\int_X \Sqs^{2i-2j}(f_m) = \int_X (v^{\syn}_{2i-2j}  f_m)$ by definition of $v^{\syn}$, we find that
\begin{align*}
w^{\syn}_{2i}  = \sum_{j=0}^i \sum_m \left(\Sqs^{2j} (e_m)  \int_X (v^{\syn}_{2i-2j}  f_m) \right)  = \sum_{j=0}^i \Sqs^{2j}  \left(\sum_m e_m  \int_X v^{\syn}_{2i-2j} f_m \right) = \sum_{j=0}^i \Sqs^{2j} (v^{\syn}_{2i-2j}),
\end{align*}
with the last equality using that $\{e_m\}$ and $\{f_m\}$ are dual bases. This proves \eqref{eq:arithmetic-wu-even}; the proof of \eqref{eq:arithmetic-wu-odd} is similar. 
\end{proof} 

\part{Applications to Brauer groups}\label{part: Brauer groups}

In this final Part, we will assemble the preceding theory for applications to arithmetic duality on Brauer groups. In \S \ref{sec: MAT pairing}, we define the Milne--Artin--Tate pairing $\langle -, - \rangle_{\MAT}$ for surfaces over $\F_p$ and higher dimensional generalizations. We prove that $\langle -, - \rangle_{\MAT}$ is always skew-symmetric and non-degenerate. Our goal is to show that it is \emph{symplectic}, so the interesting case left is $p=2$. Thanks to the skew-symmetry, the map 
\[
u \mapsto \langle u, u \rangle_{\MAT}
\]
defines a linear functional on $\Br(X)_{\nd}[p^\infty]$, which we want to explicate. In \S \ref{sec: E_infty and pairing}, we will prove a formula for this functional in terms of $\EE_\infty$ Steenrod operations, of the form 
\[
``\langle u, u \rangle_{\MAT} = \int_X [2^{n-1}] \circ \Pe^{d}(\ol{\beta_n (u)})."
\]
See Theorem \ref{thm: MAT form} for explanation of the terms. Crucially, this is one of the ``edge'' cases of \Cref{cor: operations agree} where the $\EE_\infty$ power operation $\Pe^d$ coincides with the syntomic power operation $\Ps^d$. That allows us to transfer the formula to one involving syntomic Steenrod operations instead, which we can then calculate in terms of characteristic classes using the results of Part \ref{part: characteristic classes}. In \S \ref{sec: final}, we carry out this computation to finally show that the pairing is symplectic. 

\section{The Milne--Artin--Tate pairing}\label{sec: MAT pairing}

In this section, we define various pairings of interest, including the Milne--Artin--Tate pairing on the Brauer group of a surface and its higher dimensional generalizations, and prove their skew-symmetry. 

\subsection{The Brauer group}\label{ssec: Brauer}

Let $X$ be a scheme. The \emph{(cohomological) Brauer group} \cite[\S 6.6.1]{Poo17} of $X$ is $\rH^2_{\et}(X; \G_m) \cong \rH^2_{\fppf}(X; \G_m)$. If $X$ is a smooth projective surface over a finite field $k$ of characteristic $p$, then $\Br(X)$ is a torsion abelian group, which is conjecturally finite. We will express its $p$-primary part in terms of syntomic cohomology.

The short exact sequence of fppf sheaves
\[
0 \rightarrow \mu_{p^n} \rightarrow \G_m \xrightarrow{p^n} \G_m \rightarrow 0
\]
induces a long exact sequence in cohomology 
\begin{equation*}
\begin{tikzcd}[row sep = tiny]
 & \ldots \ar[r] & \Hf^1(X; \G_m) \ar[r, "p^n"] & \Hf^1(X; \G_m)   \\
\ar[r] & \Hf^2(X; \mu_{p^n}) \ar[r] & \Hf^2(X; \G_m) \ar[r, "p^n"] & \Hf^2(X; \G_m) \\
\ar[r] & \Hf^3(X; \mu_{p^n}) \ar[r] & \Hf^3(X; \G_m) \ar[r, "p^n"] & \ldots 
\end{tikzcd}
\end{equation*}
from which we extract a short exact sequence 
\begin{equation}\label{eq:MAT-defn-1}
0 \rightarrow \frac{\Hf^1(X; \G_m) }{p^n \Hf^1(X; \G_m)} \rightarrow \Hf^2(X; \mu_{p^n}) \rightarrow \Hf^2(X; \G_m)[p^n] \rightarrow 0.
\end{equation}
Taking the colimit over $n \in \N$ in \eqref{eq:MAT-defn-1} along the multiplication-by-$p$ map, we get a short exact sequence
\begin{equation}\label{eq: brauer ses 1}
0 \rightarrow \Hf^1(X; \G_m) \otimes \frac{\Q_p}{\Z_p}  \rightarrow \Hf^2(X; \mu_{p^{\infty}}) \rightarrow \Br(X)[p^{\infty}] \rightarrow 0.
\end{equation}

For an abelian group $G$, we write $G_{\nd}$ for the non-divisible quotient of $G$ (i.e., the quotient of $G$ by its subgroup of divisible elements), and we write $G_{\tors}$ for the torsion subgroup of $G$. Since $\Q_p/\Z_p$ is divisible, \eqref{eq: brauer ses 1} induces an isomorphism  
\[
\Hf^2(X; \mu_{p^{\infty}})_{\nd} \xrightarrow{\sim} \Br(X)[p^{\infty}]_{\nd} \cong \Br(X)_{\nd}[p^{\infty}].
\]
On the other hand, we have by Example \ref{ex: nu sheaves examples} a canonical identification $\Hf^2(X; \mu_{p^{\infty}}) \cong \Hsyn^2(X; \Q_p/\Z_p(1))$, which when combined with the preceding discussion yields
\[
\Hsyn^2(X; \Q_p/\Z_p(1))_{\nd} \cong \Br(X)_{\nd}[p^{\infty}].
\]

\subsection{Duality on the Brauer group}
Now suppose that $X$ is a smooth, proper, geometrically connected surface over a finite field $k$ of characteristic $p$. By Poincar\'e duality, there is a trace map 
\[
\int_X \co \Hsyn^{5,2}(X; \Z/p^n ) \xrightarrow{\sim} \Z/p^n
\]
and the resulting pairing 
\[
\Hsyn^{a,b}(X; \Z/p^n) \times \Hsyn^{5-a,2-b}(X; \Z/p^n) \rightarrow \Hsyn^{5,2}(X; \Z/p^n) \xrightarrow{\int_X} \Z/p^n 
\]
is perfect for every $a,b \in \Z$, and $n \in \N$. Taking colimits in $n$ in one factor and limits in $n$ in the other, we obtain a non-degenerate pairing 
\[
\Hsyn^{2,1}(X; \Q_p/\Z_p) \times \Hsyn^{3,1}(X; \Z_p) \rightarrow \Hsyn^{5,2}(X; \Q_p/\Z_p) \xrightarrow{\sim} \Q_p/\Z_p.
\]
Since the divisible subgroup is the annihilator of the torsion subgroup, this induces a non-degenerate pairing 
\[
\Hsyn^{2,1}(X; \Q_p/\Z_p)_{\nd} \times \Hsyn^{3,1}(X; \Z_p)_{\tors} \rightarrow \Hsyn^{5,2}(X; \Q_p/\Z_p) \xrightarrow{\sim} \Q_p/\Z_p.
\]
On the other hand, the short exact sequence
\[
0 \rightarrow \Z_p(1) \rightarrow \Q_p(1) \rightarrow \Q_p/\Z_p(1) \rightarrow 0
\]
yields a boundary map on cohomology $\wt{\delta} \co \Hsyn^{2,1}(X; \Q_p/\Z_p)  \rightarrow \Hsyn^{3,1}(X; \Z_p)$, and from the long exact sequence we see that it induces an isomorphism 
\[
\wt{\delta} \co \Hsyn^{2,1}(X; \Q_p/\Z_p)_{\nd}  \xrightarrow{\sim} \Hsyn^{3,1}(X; \Z_p)_{\tors}.
\]

\begin{defn}[Milne--Artin--Tate pairing]
Let $X$ be a smooth, proper, geometrically connected surface over a finite field $k$ of characteristic $p$. The \emph{Milne--Artin--Tate pairing} on $\Br(X)_{\nd}[p^{\infty}]$ sends $u,v \in \Br(X)_{\nd}[p^{\infty}] \cong \Hsyn^{2,1}(X; \Q_p/\Z_p)_{\nd}$ to 
\[
\tw{u,v}_{\MAT} := \int_X u \cdot \wt{\delta} v \in \Q_p/\Z_p
\]
\end{defn}

\subsection{Higher Brauer groups}\label{ssec: higher brauer}
We next give a higher-dimensional generalization of the Milne--Artin--Tate pairing. Following Jahn \cite{Jahn15}, for an integer $r \geq 1$ we define the \emph{$r$th higher Brauer group} of a smooth quasiprojective variety $X$ over a field to be the \emph{\'etale-motivic cohomology} 
\[
\Br^r(X) := \rH_{\et}^{2r+1}(X; \Z(r)).
\]
These groups are torsion, since $\rH_{\et}^{2r+1}(X; \Q(r)) \cong \rH_{\mot}^{2r+1}(X; \Q(r)) = 0$. 

\begin{example}
For $r=1$, $\Z(1)[1]  \cong \G_m$ so that $\Br^1(X) = \Br(X)$ recovers the usual cohomological Brauer group. 
\end{example}

\begin{lemma} Assume that $X$ is a smooth quasiprojective variety over a finite field $k$ of characteristic $p$. Then there is a natural isomorphism
\begin{equation}\label{eq: Brd as syntomic}
\Br^r(X)_{\nd}[p^\infty] \cong \Hsyn^{2r,r}(X; \Q_p/\Z_p)_{\nd}.
\end{equation}
\end{lemma}

\begin{proof}
Consider the exact triangles of \'etale-motivic sheaves on $X$,
\[
\Z(r)^{\et} \xrightarrow{p^n}  \Z(r)^{\et} \rightarrow \Z/p^n(r)^{\et}.
\]
As explained in \S \ref{sssec: motivic}, the results of Geisser--Levine \cite{GL00} identify $\Z/p^n(r)^{\et}$ with Milne's logarithmic de Rham--Witt sheaves, which in turn is identified with $\Z/p^n(r)^{\syn}$ by Bhatt--Morrow--Scholze \cite[Corollary 8.21 and Remark 8.22]{BMS2}. (Recall from \S \ref{ssec:syntomic-cohomology} that we regard $\Z/p^n(r)^{\syn}$ as an \'etale sheaf.) In particular, we have natural isomorphisms
\[
\Het^*(X; \Z/p^n(r)^{\et}) \cong \Het^*(X; \Z/p^n(r)^{\syn}).
\]
From these identifications and the associated long exact sequence in cohomology, we obtain exact sequences
\[
\Het^{2r}(X; \Z(r)^{\et}) \otimes \Z/p^n \Z \rightarrow \Hsyn^{2r}(X; \Z/p^n(r)) \rightarrow \Het^{2r+1}(X; \Z(r)^{\et})[p^n] \rightarrow 0.
\]
Taking the direct limit over $n$, the leftmost term becomes  $\Het^{2r+1}(X; \Z(r)) \otimes \Q_p/\Z_p$, which is divisible. Therefore, after passing to non-divisible quotients we obtain the claimed isomorphism. 
\end{proof}

Now suppose that $X$ is furthermore proper and geometrically connected over $k$, of dimension $2d$. For $\ell \neq p$, Jahn defines a pairing on $\Br^d(X)[\ell^\infty]$; this is recalled in \cite[\S 2]{Feng20} and called the \emph{Artin--Tate pairing}. We will now extend the pairing to the full $\Br^d(X)$ by defining it on $
\Br^d(X)[p^\infty]$, and we will also call this extension the \emph{Milne--Artin--Tate pairing}. 

Let 
\begin{equation}\label{eq: wtdelta}
\wt{\delta} \co \Hsyn^{2d,d}(X; \Q_p/\Z_p)_{\nd} \rightarrow \Hsyn^{2d+1,d}(X; \Z_{p})_{\tors}
\end{equation}
be the map induced by the boundary map for the exact triangle of syntomic sheaves
\begin{equation}\label{eq: SES}
\Z_p(d)^{\syn} \rightarrow \Q_p(d)^{\syn}\rightarrow  \Q_p/\Z_p(d)^{\syn}
\end{equation}

\begin{lemma}
The map \eqref{eq: wtdelta} is an isomorphism.
\end{lemma}

\begin{proof} The long exact sequence associated to \eqref{eq: SES} reads
\[
\begin{tikzcd}
& \Hsyn^{2d,d}(X; \Q_p) \ar[r] & \Hsyn^{2d,d}(X; \Q_p/\Z_p)  \ar[dll] \\  
 \Hsyn^{2d+1,d}(X; \Z_p) \ar[r] &  
\Hsyn^{2d+1,d}(X; \Q_p)  
\end{tikzcd}
\]
From this, we see that the image of the boundary homomorphism is $\Hsyn^{2d+1,d}(X; \Z_p)_{\tors}$, and its kernel is divisible. Therefore, it factors over an isomorphism $\wt{\delta}$, as claimed. 
\end{proof}

\begin{defn}\label{def: MAT} Let $X$ be a smooth, proper, geometrically connected variety of dimension $2d$ over a finite field $k$ of characteristic $p$. For $u,v \in \Hsyn^{2d,d}(X; \Q_p/\Z_p)_{\nd}$, we define 
\[
\langle u,v \rangle_{\MAT} := \int_X u \cdot \wt{\delta} v.
\]
By \eqref{eq: Brd as syntomic}, we may view this as a pairing on $\Br^d(X)_{\nd}[p^\infty]$. This extends the Artin--Tate pairing to a pairing on all of $\Br^d(X)_{\nd}$, which we call the \emph{Milne--Artin--Tate pairing}.
\end{defn}

From Poincar\'{e} duality and the fact that \eqref{eq: wtdelta} is an isomorphism, it is evident that this pairing is non-degenerate on $\Br^d(X)_{\nd}$. We prove below that it is skew-symmetric, i.e., 
\[
\langle u,v \rangle_{\MAT}  = - \langle v,u \rangle_{\MAT} \text{ for all } u,v \in \Br^d(X)_{\nd}.
\]
Recall that this is weaker than \emph{symplectic}, or \emph{alternating}, which would say 
\[
\langle u,u \rangle_{\MAT}  = 0 \quad \text{ for all } u \in \Br^d(X)_{\nd}.
\]
We will even prove the stronger property that the pairing is symplectic. 

\begin{thm}\label{thm: MAT higher dimension}
Let $X$ be a smooth, proper, geometrically connected variety of dimension $2d$ over a finite field $k$ of characteristic $p$. Then the Milne--Artin--Tate pairing on $\Br^d(X)_{\nd}$ is symplectic. 
\end{thm}

The proof of \Cref{thm: MAT higher dimension} will be completed in \S\ref{sec: final}. 

\subsection{Skew-symmetry}
Let $X$ be a smooth, proper, geometrically connected variety of dimension $2d$ over a finite field $k$ of characteristic $p$. We will define an auxiliary pairing on the group $\Hsyn^{2d,d}(X; \Z/p^n)$.

\begin{defn}\label{def: pairing mod 2^n} Consider the exact triangle of syntomic sheaves on $X$:
\[
\Z/p^n(d)^{\syn} \rightarrow  \Z/p^{2n} (d)^{\syn} \rightarrow \Z/p^n(d)^{\syn}
\]
and call the induced boundary map 
\begin{equation}\label{eq: bockstein_n}
\beta_n \colon \Hsyn^{i,d}(X; \Z/p^n) \rightarrow \Hsyn^{i+1,d}(X; \Z/p^n).
\end{equation}
We define the pairing 
\[
\langle \cdot , \cdot \rangle_{n} \colon \Hsyn^{2d,d}(X; \Z/p^n) \times \Hsyn^{2d,d}(X; \Z/p^n) \rightarrow \Z/p^n
\]
by 
\[
\langle u,v \rangle_n := \int_X (u \cdot \beta_n v).
\]
\end{defn}

\begin{prop}\label{prop: skew-symmetric}
The pairing $\langle \cdot , \cdot \rangle_n$ is skew-symmetric. 
\end{prop} 

\begin{proof}
The assertion is equivalent to 
\[
x \cdot (\beta_n y)  + y \cdot (\beta_n x)  = 0.
\]
Since $\beta_n$ is a derivation, we have $ x \cdot (\beta_n y)  + y \cdot (\beta_n  x)  = \beta_n(x \cdot y)$. Then the result follows from the next Lemma. 
\end{proof}

\begin{lemma}\label{lem: top boundary vanishes}
The boundary map $\beta_n  \colon \Hsyn^{4d,2d}(X; \Z/p^n) \rightarrow  \Hsyn^{4d+1,2d}(X; \Z/p^n)$ vanishes. 
\end{lemma}

\begin{proof}
By the obvious long exact sequence, the image of $\beta_n$ is the kernel of 
\[
[p^{n}] \co \Hsyn^{4d+1,2d}(X; \Z/p^n) \rightarrow \Hsyn^{4d+1,2d}(X; \Z/p^{2n})
\]
which is identified with the inclusion $p^n \Z/p^{2n}\Z \hookrightarrow \Z/p^{2n}\Z$ by Poincar\'{e} duality. 
\end{proof}

\begin{prop}\label{prop: compatibility}
The boundary map 
\begin{equation}\label{eq: integral bockstein}
\Hsyn^{2d,d}(X; \Z/p^n) \rightarrow \Hsyn^{2d+1,d}(X; \Z_p)
\end{equation}
induced by the exact triangle of syntomic complexes
\[
\Z_p(d)^{\syn}  \xrightarrow{p^n} \Z_p(d)^{\syn}  \rightarrow \Z/p^n(d)^{\syn}
\]
surjects onto $\Hsyn^{2d+1,d}(X; \Z_p)[p^n]$. Moreover, it is compatible for the pairings $\langle \cdot , \cdot \rangle_n$ and $\langle \cdot , \cdot \rangle_{\MAT}$ in the sense that the following diagram commutes
\[
\xymatrix @C=0pc{
\Hsyn^{2d,d}(X; \Z/p^n) \ar@{->>}[d]_{\eqref{eq: integral bockstein}} & \times &  \Hsyn^{2d,d}(X; \Z/p^n)  \ar@{->>}[d]_{\eqref{eq: integral bockstein}}  \ar[rrrrrrr]^{\langle \cdot , \cdot \rangle_n} &&&&&&& \Hsyn^{4d+1,2d}(X; \Z/p^n) \ar[d]_{\wr}  \\
\Hsyn^{2d+1,d}(X; \Z_p)[p^n] & \times & \Hsyn^{2d+1,d}(X; \Z_p)[p^n]  \ar[rrrrrrr]_{\langle \cdot , \cdot \rangle_{\MAT}} &&&&&&& \Hsyn^{4d+1,2d}(X; \Q_p/\Z_p) [p^n]
	}
\]
\end{prop}

\begin{proof} This follows formally from a diagram chase, exactly as in \cite[Proposition 2.5]{Feng20}. 
\end{proof}

\begin{cor}
The Milne--Artin--Tate pairing is skew-symmetric.
\end{cor}

\begin{proof}Combine Proposition \ref{prop: skew-symmetric} and Proposition \ref{prop: compatibility}.
\end{proof}

\begin{cor}\label{cor: alternating_criterion}
If the pairing $\langle \cdot , \cdot \rangle_n$ on $\Hsyn^{2d,d}(X; \Z/p^n )$ is alternating, then so is the pairing $\langle \cdot, \cdot \rangle_{\MAT}$ on $\Br^d(X)_{\nd}[p^n]$. 
\end{cor}

\begin{proof}
This follows immediately from Proposition \ref{prop: compatibility}.
\end{proof}

Therefore, to prove \Cref{thm: MAT higher dimension}, it suffices to establish: 

\begin{thm}\label{thm: pseudo_main}
The pairing $\langle \cdot , \cdot \rangle_n$  is alternating for all $n$. 
\end{thm}

The proof of Theorem \ref{thm: pseudo_main} will be the focus of the rest of the paper. Thanks to Corollary \ref{cor: alternating_criterion}, the only non-trivial case is $p=2$.

\section{Arithmetic duality and $\EE_\infty$ Steenrod operations}\label{sec: E_infty and pairing}

The goal of this section is to prove the following result, which relates the pairing of \Cref{def: pairing mod 2^n} with $\EE_\infty$ Steenrod operations. 

\begin{thm}\label{thm: MAT form}
Let $X$ be a smooth, proper, geometrically connected variety of dimension $2d$ over a finite field $k$ of characteristic $2$. For all $u \in \Hsyn^{2d,d}(X; \Z/2^n)$, we have the identity
\[
u \cdot \beta_n (u) = [2^{n-1}] \circ \Pe^{d}(\ol{\beta_n (u)})
\]
where: 
\begin{itemize}
\item $\beta_n \co \Hsyn^{2d,d}(X; \Z/2^n) \rightarrow \Hsyn^{2d+1,d}(X; \Z/2^n)$ is the Bockstein homomorphism \eqref{eq: bockstein_n}.
\item $\ol{\beta_n (u)}$ denotes the reduction of $\beta_n (u)$ mod $2$. 
\item $[2^{n-1}] \co \Hsyn^{4d+1,2d}(X; \Z/2) \rightarrow \Hsyn^{4d+1,2d}(X; \Z/2^n)$ is the cohomology map induced by the map $\Z/2(2d)^{\syn}_X \rightarrow \Z/2^n(2d)^{\syn}_X$, given by multiplication by $2^{n-1}$. 
\end{itemize}
\end{thm}

\subsection{The Bockstein spectral sequence}
The key ingredient behind \Cref{thm: MAT form} is a more general formula relating the $\EE_\infty$-power operations and the differentials in the Bockstein spectral sequence, which we will now formulate.

For $M\in \cD(\Z)$, we can form the derived $2^n$-adic filtration
\[
M/2^\bullet:\ldots  M/2^{3n} \to M/2^{2n} \to M/2^n \to 0
\]
where all quotients are formed in the derived sense, that is, $M/2^n$ is the derived cofiber of multiplication by $2^n$ on $M$. The limit of this filtration is the derived $2$-adic completion $M^\wedge_2$.  

The filtration $M/2^\bullet$ is obtained by tensoring $M$ with the filtered commutative algebra $\Z/2^\bullet$. Let $t$ be the class of $2^n$ in the first graded piece $\Gr_1(\Z/2^\bullet)\cong 2^n\Z/{2^{2n}\Z}[-1]$. Then the associated graded algebra of $\Z/2^\bullet$ is given by $\Gr(\Z/2^\bullet)\cong \Z/2^n[t]$ (where $t$ has degree $1$). Consequently, the associated graded of $\Gr(M/2^\bullet)$ is 
\[
\Gr(M/2^\bullet) \cong M\otimes \Gr(\Z/2^\bullet)\cong M/2^n[t]. 
\]
The spectral sequence of cohomology groups associated with this filtration is the familiar
\emph{Bockstein spectral sequence}, which therefore assumes the form 
\begin{equation}\label{eq: bockstein SS}
E^1_{k,l} = \rH^{l}(M/2^n) \cdot t^k \Rightarrow \rH^l(M^\wedge_2)
\end{equation}
(where $t$ should be treated as an indeterminate polynomial variable with grading degree $1$ and cohomological degree $0$).

The $d_1$-differentials in the Bockstein spectral sequence are given by $d_1(x) = \beta_n(x)t$ for $x\in \rH^*(M/2^n)$. The $d_2$ differentials are divisible by $t^2$, and may be written as 
\[
d_2(x) = \beta_{n}^{(2)}(x)t^2
\] 
for a partially defined map $\beta_{n}^{(2)}\colon \rH^*(M/2^n) \to \rH^{*+1}(M/2^n)$ that we refer to as the \emph{secondary Bockstein homomorphism}. It is defined on the kernel of $\beta_n$ and modulo the image of $\beta_n$ in the appropriate degrees.

If $A = \bigoplus_{b \in \Z} A_b$ is a $\Z$-graded $\EE_\infty$-algebra over $\Z$, we write
\[
\rH^{a,b}(A) := \rH^a(A_b).
\]
Although both sides of the equation in Theorem \ref{thm: MAT form} can be defined for general such $\EE_\infty$-algebras, the equality does not hold in this generality. The abstract ingredient in the proof of \Cref{thm: MAT form} that does hold for general graded $\EE_\infty$-algebras is the following. 

\begin{prop} \label{prop:MAT form general}
Let $A$ be a $\Z$-graded $\EE_\infty$-algebra over $\Z$. Let $\beta_{n}^{(2)}$ be the secondary Bockstein homomorphism for $A$.
Then, with similar notations as in \Cref{thm: MAT form}, for all $u\in \rH^{2a,b}(A/2^n)$ we have
\[
\beta_n(2^{n-1}u^2) = 0
\]
and 
\[
\beta_{n}^{(2)}(2^{n-1}u^2) =  u \cdot \beta_n(u) -[2^{n-1}] \cdot \Pe^a(\ol{\beta_n(u)}).
\]
\end{prop}

\subsection{Calculation of differentials}

In this subsection, we will prove \Cref{prop:MAT form general}.

 Let $\cD^{\gr}(\Z)$ be the $\infty$-category of graded $\Z$-module spectra, or equivalently the graded derived $\infty$-category of $\Z$.  We write $\rH^{a,b}(M)$ for the $a$th cohomology group of the $b$th graded piece of $M$. 

For $M\in \cD(\Z)$, we suggestively denote by\footnote{Note that this notation clashes with the Tate-twist notation for motivic spectra. However, they are compatible in our application: $\MS_S$ is linear over filtered, and hence graded, spectra, in such a way that the Tate twist agrees with tensoring with the graded spectrum $\mathbb{S}(n)$.} $M(a)\in \cD^\gr(\Z)$ the module $M$ concentrated in degree $a$.

\subsubsection{Universal models for cohomology classes} We will develop some explicit models to calculate power operations on cohomology classes. 

\begin{lemma} \label{lem:rep_for_cohomology_mod_a}
Given a graded $\EE_\infty$-algebra $A$ over $\Z$, the data of a class $u\in \rH^{a,b}(A/2^n)$ agrees with that of a homotopy class of maps 
\[
\varphi_u\colon \Z/2^n(-b)[-a-1] \to A.
\]
The base change $\ol{\varphi}_u$ of $\varphi_u$ along the map $\Z \rightarrow \Z/2^n$ classifies the pair $(\beta_n(u) ,u )$, under the identifications of the commutative diagram below.
\[
\begin{tikzcd}
\Z/2^n(-b)[-a-1] \oplus \Z/2^n(-b)[-a] \ar[d, "\wr"] \ar[dr, "\ol{\varphi}_u"] \\
 \Z/2^n(-b)[-a-1]\otimes_\Z \Z/2^n \ar[r, "\varphi_u \otimes_{\Z} \Z/2^n"'] &  A/2^n
\end{tikzcd}
\]
\end{lemma}

\begin{proof}We have
\begin{align*}
\rH^{a,b}(A/2^n) &=  \Ext^{a,b}(\Z, A/2^n)  \cong \Ext^0(\Z(-b)[-a],A\otimes \Z/2^n)\\ &\cong \Ext^0(\Z(-b)[-a] \otimes (\Z/2^n)^\vee, A) \cong \Ext^0(\Z/2^{n}(-b)[-a-1],A).
\end{align*}
where in the last step we used that the dual of $\Z/2^n$ in $\cD^{\gr}(\Z)$ is $\Z/2^n[-1]$.

Explicitly, we can represent $\Z/2^{n}(-b)[-a-1]$ by the complex 
\begin{equation}\label{eq:M-a-b-complex}
\dots 0 \too \Z(-b) \oto{2^n} \Z(-b) \too 0 \dots
\end{equation}
where the non-zero terms are in cohomological degrees $a$ and $a+1$. Hence if $A$ is represented by a graded chain complex of $\Z$-module with levelwise torsion-free terms, then the map $\varphi_u$ can be chosen to send the summand in degree $a$ to an integral lift $\tilde{u}$ of $u$ and the summand in degree $a+1$ to $d(\tilde{u})/2^n$.
From this description we immediately get the desired description of the map $\ol{\varphi}_u$ by reducing modulo $2^n$.   
\end{proof}

We can now use this discussion to compute power operations acting on $u$ and $\beta_n(u)$. 
To ease notation, from now on we abbreviate $M_{a+1,b} := \Z/2^n(-b)[-a-1]$, for which we take the explicit chain complex \eqref{eq:M-a-b-complex} as a representative.

Recall that $C_2$ is the cyclic group of order $2$. To compute $u\cdot \beta_n(u)$, we shall consider the composition 
\begin{equation}\label{eq:psi_u}
\psi_u\colon (M_{a+1,b}^{\otimes 2})_{hC_2} \oto{\varphi_u^{\otimes 2}} (A^{\otimes 2})_{hC_2} \oto{\mult} A,
\end{equation}
in which $C_2$ acts on the tensor factors by swapping them, and $(-)_{hC_2}$ is the homotopy quotient by the action of $C_2$. 
Let 
\begin{equation}\label{eq:ol-psi_u}
\ol{\psi}_u  = \psi_u \otimes_{\Z} \Z/2^n \co (M_{a+1,b}^{\otimes 2})_{hC_2} \otimes_{\Z} \Z/2^n \rightarrow A/2^n
\end{equation}
be the base change of $\psi_u$ along $\Z \rightarrow \Z/2^n$.

\begin{prop}\label{prop:distributivity_psi}
We have 
\[
(M_{a+1,b}^{\otimes 2})_{hC_2} \otimes_\Z \Z/2^n \cong (M_{2a+2,2b} / 2^n)_{hC_2} \oplus M_{2a+1,2b} / 2^n \oplus (M_{2a,2b} / 2^n)_{hC_2},
\]
where $C_2$ acts on the summands according to the Koszul sign rule. Via this identification, the map $\ol{\psi}_u$ from \eqref{eq:ol-psi_u} is given by the triple of cohomology classes
\[
(\PP^2(\beta_n(u)) , u \cdot \beta_n(u), \PP^2(u) ),
\]
where $\PP^2(u)\colon (M_{2a,2b} / 2^n)_{hC_2} \to A/2^n$ is the ``total square`` (in the sense of \Cref{def:total_power}, for $p=2$) of $u$, and similarly for $\PP^2(\beta_n(u))$.
\end{prop}

\begin{proof}
Since the functor $(-)\otimes_\Z \Z/2^n\colon \cD(\Z)^{\gr} \to \cD(\Z/2^n)^{\gr}$ is colimit-preserving and symmetric monoidal, it commutes with the operations in the definition of $\psi_u$ \eqref{eq:psi_u}. Interchanging the order of operations, the map $\ol{\psi}_u$ can be described as the composition 
\[
(M_{a+1,b} / 2^n)^{\otimes 2}_{hC_2} \oto{\ol{\varphi}_u^{\otimes 2}} (A/2^n)^{\otimes 2}_{hC_2} \oto{\mult} A/2^n
\]
where the tensor products are now over $\Z/2^n$. Using \Cref{lem:rep_for_cohomology_mod_a}, this identifies with the composition 
\begin{equation}\label{eq:dist-psi-eq1}
(M_{a+1,b}/2^n   \oplus M_{a,b}/2^n )^{\otimes 2}_{hC_2} \oto{(\beta_n(u) \oplus u  )^{\otimes 2}} (A/2^n)^{\otimes 2}_{hC_2} \oto{\mult} A/2^n.
\end{equation}
The result now follows from the ($C_2$-equivariant) distributivity of the tensor product. More precisely, by distributivity the first map in \eqref{eq:dist-psi-eq1} identifies with the direct sum of the maps $\beta_n(u)^{\otimes 2}$, $u\otimes \beta_n(u) \oplus \beta_n(u)\otimes u$, and $u^{\otimes 2}$, which give the desired map formula after composing with the multiplication map by the definition of the total square. 
\end{proof}

\subsubsection{Calculation in the universal case} By functoriality, the preceding results allow us to focus on the universal case of $M_{a+1,b}$. To carry out the computation in this case, we shall work with the explicit presentation \eqref{eq:M-a-b-complex} of $M_{a+1,b}$ used in the proof of  \Cref{lem:rep_for_cohomology_mod_a}. We let $x \in M_{a+1,b}$ be the generator in cohomological degree $a$ and $y$ be the generator in degree $a+1$, so that $d(x) = 2^n y$, and both $x$ and $y$ have grading degree $b$. We also assume from now that $a$ is even, though a similar analysis applies to the odd case (giving a different formula that we shall not need). 

We will build an explicit model for the group homology chains $C_*(C_2;M)$. Let $\sigma$ be a generator for $C_2$. For the trivial module $\Z$ of $C_2$, there is the standard resolution 
\begin{equation}\label{eq:standard-C_2-resolution}
\ldots \xrightarrow{1-\sigma}  \Z[C_2]e_2  \xrightarrow{1+\sigma} \Z[C_2]e_1 \xrightarrow{1-\sigma} \Z[C_2]e_0 
\end{equation}
which in degree $i$ is the free rank one $\Z[C_2]$-module on a generator $e_i$, satisfying 
\[
d(e_i) = (1+ (-1)^i \sigma )e_{i-1}.
\]
If $M$ is a complex of modules over $C_2$, then tensoring $M$ with \eqref{eq:standard-C_2-resolution} over $\Z[C_2]$ gives a model for $C_*(C_2;M)$ by the complex 
\begin{equation}
\ldots \xrightarrow{d} e_2 \otimes M \xrightarrow{d} e_1 \otimes M \xrightarrow{d} e_0 \otimes M
\end{equation}
whose differentials are described by 
\begin{equation}\label{eq:diff-C_2-cohomology}
d(e_i \otimes m) = e_{i-1} \otimes (m + (-1)^i \sigma(m)) + (-1)^i e_i \otimes d(m).
\end{equation}
Accordingly, the object 
$
(M_{a,b}^{\otimes 2})_{hC_2} 
$
is representable by a graded complex of abelian groups with generators 
\[
e_i \otimes x^{\otimes 2}, \quad e_i \otimes x\otimes y, \quad e_i \otimes y\otimes x,\quad e_i \otimes y^{\otimes 2} \quad i \ge 0. 
\]
The differential is determined by \eqref{eq:diff-C_2-cohomology} and the $\sigma$-action (coming from the Koszul sign rule)
\[
\sigma(x^{\otimes 2}) = x^{\otimes 2}, \quad \sigma(x\otimes y) = y\otimes x,\quad \sigma(y^{\otimes 2}) = -y^{\otimes 2}. 
\]
Let $\bar{x}$ (resp. $\bar{y}$) denote the reductions of $x$ (resp. $y$) mod $2^n$. 
We are now ready for the computation in our universal example. 

\begin{prop} \label{prop:second_bock_symm_power}
Let $a \in 2 \Z$ and $n \in \Z_{\geq 1}$. Then in $(M_{a+1,b}^{\otimes 2})_{hC_2} / 2^n$, the cochain $2^{n-1}\bar{x}^{\otimes 2}$ is a cocycle, whose cohomology class $[2^{n-1}\bar{x}^{\otimes 2}]$ satisfies 
\[
\beta_n([2^{n-1}\bar{x}^{\otimes 2}]) = 0
\]
and
\[
\beta_{n}^{(2)}([2^{n-1}\bar{x}^{\otimes 2}]) \equiv [\bar{x}\otimes \bar{y}] -  [2^{n-1} e_1 \otimes \bar{y}^{\otimes 2}]
\]
modulo the image of $\beta_n$. 
\end{prop}

\begin{proof}
The fact that $2^{n-1}\bar{x}^{\otimes 2}$ is a cocycle is seen by calculating directly that $d(x^{\otimes 2})$ is divisible by $2$.

To calculate the Bockstein operations, consider the integral lift 
\[
z := 2^{n-1} x\otimes x + 2^{2n -1} e_1 \otimes x\otimes y 
\]
of $2^{n-1}\ol{x}\otimes \ol{x}$. From the formulas above, we have 
\begin{align}\label{eq:bock-lift-1}
d(z) &= 2^{2n-1}(x\otimes y + y\otimes x) + 2^{2n-1}(x\otimes y - y\otimes x) - 2^{3n-1} e_1 \otimes y\otimes y \nonumber \\
& = 2^{2n}\left(x\otimes y - 2^{n-1}e_1\otimes y \otimes y \right).  
\end{align}
Since \eqref{eq:bock-lift-1} is divisible by $2^{2n}$, we learn that $\beta_n([2^{n-1}\bar{x}^{\otimes 2}])=0$, and the secondary Bockstein is
\begin{equation}\label{eq: cochain rep secondary bock}
\beta_{n}^{(2)}(2^{n-1}\bar{x}^{2}) = \ol{\frac{d(z)}{2^{2n}}} =  \ol{x}\otimes \ol{y} - 2^{n-1} e_1 \otimes \ol{y} \otimes \ol{y}, 
\end{equation}
yielding the result. 
\end{proof}

\subsubsection{Reduction to the universal case} To complete the computation of the Bockstein images of the classes of the form $2^{n-1} u^2$, it remains to identify the image of the classes appearing in \Cref{prop:second_bock_symm_power} under the map $\ol{\psi}_u \colon  (M_{a+1,b}^{\otimes 2})_{hC_2} / 2^n \to A/2^n$, which is given by our explicit chain complex presentations as above. 

\begin{prop}\label{prop:image_classes_under_psi_u}
The map $\ol{\psi}_u$ satisfies 
\begin{equation}\label{eq:image_classes_under_psi_u 1}
\ol{\psi}_u(x\otimes y) = u \cdot \beta_n(u)
\end{equation}
and 
\begin{equation}\label{eq:image_classes_under_psi_u 2}
\ol{\psi}_u(e_1 \otimes y \otimes y) \bmod{2} =  \Pe^{\frac{a}{2}}(\ol{\beta_n(u)}). 
\end{equation}
\end{prop}

\begin{proof}
Recall from \Cref{prop:distributivity_psi} that under the identification
\[
(M_{a+1,b}^{\otimes 2})_{hC_2} \otimes_\Z \Z/2^n \cong (M_{2a+2,2b} / 2^n)_{hC_2} \oplus M_{2a+1,2b} / 2^n \oplus (M_{2a,2b} / 2^n)_{hC_2},
\]
the map $\ol{\psi}_u$ is the sum of the maps  $\PP^2(u), u\cdot \beta_n(u)$, and $\PP^2(\beta_n(u))$. Explicitly, the first summand is the restriction to the subcomplex spanned by the tensors $e_i \otimes x \otimes x$, the second to $e_i \otimes x\otimes y$, and the third to $e_i \otimes y \otimes y$, where the first and third complexes give the standard presentations of $(M_{2a+2,2b}/2^n)_{hC_2}$ and $(M_{2a,2b}/2^n)_{hC_2}$ as chain complexes (by tensoring with the resolution \eqref{eq:standard-C_2-resolution} over $\Z[C_2]$). 

It follows that $\ol{\psi}_u(\ol{x}\otimes \ol{y}) = u\cdot \beta_n(u)$, and that $\ol{\psi}_{u}(e_1\otimes y\otimes y)$, after reduction modulo $2$, is given by the composition 
\[
M_{2a+2,2b}/2[1] \oto{t_1} M_{2a+2,2b}/2 \otimes \rB C_2 \cong (M_{2a+2,2b}/2)_{hC_2} \oto{\PP^2(\ol{\beta_n(u)})} A/2   
\]
where $t_1$ is the homology class of $e_1$ modulo 2. Comparing with \Cref{defn:norm-E-infty-ops} (which is related to our other definition of the $\EE_\infty$ power operation by \Cref{prop:compare_Pe_defs}), we see that this composition is precisely $\Pe^{\frac{a}{2}}(\ol{\beta_n(u)})$. 
\end{proof}

\begin{proof}[Proof of \Cref{prop:MAT form general}]
Let $A$ be a $\Z$-graded $\EE_\infty$ algebra over $\Z$ and let $u\in \rH^{2a,b}(A/2^n)$. By \Cref{lem:rep_for_cohomology_mod_a}, $u$ corresponds to a map $\varphi_u\colon M_{2a+1,b} \to A$ which after tensoring with $\Z/2^n$ classifies the classes $(\beta_n(u),u)$ in $\rH^{*,*}(A/2^n)$. Taking the $C_2$-equivariant square, we obtain a map 
\[
\psi_u\colon (M_{2a+1,b}^{\otimes 2})_{hC_2} \to A
\]
whose reduction $\ol{\psi}_u$ modulo $2^n$ classifies the triple $(\PP^2(\beta_n (u)),u \cdot \beta_n(u),\PP^2(u))$, by \Cref{prop:distributivity_psi}. Now, by \Cref{prop:image_classes_under_psi_u} the map $\ol{\psi}_u$ satisfies 
\begin{equation}\label{eq:psi-of-classes}
\ol{\psi}_u(\ol{x}\otimes \ol{y}) = u \cdot \beta_n(u), \quad \ol{\psi}_u(e_1\otimes \ol{y}\otimes \ol{y}) \bmod{2} = \Pe^a(\beta_n(u))
\end{equation}
Note that the last equality promotes to an equality mod $2^n$,
\begin{equation}\label{eq:psi-of-classes-2}
\ol{\psi}_u(2^{n-1} e_1\otimes \ol{y}\otimes \ol{y})  = [2^{n-1}] \Pe^a(\beta_n(u)).
\end{equation}

Since $\ol{\psi}_u$ lifts to an integral map, it is compatible with the Bockstein differentials. Then from \Cref{prop:second_bock_symm_power} we get the identities: 
\[
\beta_n(2^{n-1}u^2) = \beta_n(\ol \psi_u(2^{n-1}\ol{x}^{\otimes 2})) = \ol \psi_u(\beta_n(2^{n-1}\ol{x}^{\otimes 2})) = 0
\]
and, using \eqref{eq:psi-of-classes} and \eqref{eq:psi-of-classes-2},
\begin{align*}
\beta_n^{(2)}(2^{n-1}u^2) &= \beta_n^{(2)}(\ol \psi_u(2^{n-1}\ol{x}^{\otimes 2})) = \ol \psi_u(\beta_n^{(2)}(2^{n-1}\ol{x}^{\otimes 2})) \\
&= \ol{\psi}_{u}(\ol{x}\otimes \ol{y} - 2^{n-1} e_1\otimes \ol{y}^{\otimes 2}) = u \cdot \beta_n(u) - [2^{n-1}] \Pe^a(\beta_n(u)),
\end{align*}
as desired. 
\end{proof}

\subsection{Proof of \Cref{thm: MAT form}}
We now prove \Cref{thm: MAT form}. Applying \Cref{prop:MAT form general} to the graded $\EE_\infty$-algebra $A := \RGamma_{\syn}^{*,*}(X; \Z_2)$, we deduce that 
\[
\beta_n^{(2)}(2^{n-1} u^2) = u \cdot \beta_n (u) - [2^{n-1}] \Pe^d(\ol{\beta_n (u)}) \in \frac{\Hsyn^{4d+1,2d}(X; \Z/2^n)}{\Ima(\beta_n)}.
\]

It then clearly suffices to show the following two items:
\begin{enumerate}
\item $\Ima(\beta_n) \subset \Hsyn^{4d+1,2d}(X; \Z/2^n)$ is the zero subspace, and then
\item $\beta_n^{(2)}(2^{n-1} u^2)$ vanishes in   $\Hsyn^{4d+1,2d}(X; \Z/2^n)$. 
\end{enumerate}
Thanks to the interpretations of $\beta_n$ and $\beta_n^{(2)}$ as differentials in the Bockstein spectral sequence, both vanishing statements follow from the torsion-freeness of $\Hsyn^{4d+1,2d}(X; \Z_2)$. We can also spell this out explicitly as follows. 

By the commutative diagram of syntomic complexes on $X$,
\[
\begin{tikzcd}
\Z_2(2d)^{\syn} \ar[r, "2^n"] \ar[d]  & \Z_2(2d)^{\syn} \ar[r] \ar[d] & \Z/2^n(2d)^{\syn} \ar[d, equals] \\
\Z/2^n(2d)^{\syn}  \ar[r, "2^n"] & \Z/4^n(2d)^{\syn} \ar[r] & \Z/2^n(2d)^{\syn}
\end{tikzcd}
\]
we see that $\beta_n$ is the reduction modulo $2^n$ of the connecting map $\wt{\beta}_n$ for the long exact sequence associated to the upper row. But the long exact sequence shows that the image of $\wt{\beta}_n$ is contained $\Hsyn^{4d+1,2d}(X; \Z_2)[2^n]$, which is torsion-free thanks to Poincar\'e duality, so we obtain the desired vanishing of $\Ima(\beta_n) \subset \Hsyn^{4d+1,2d}(X; \Z/2^n)$.  

Similarly, equation \eqref{eq: cochain rep secondary bock} exhibits $\beta_n^{(2)}(2^{n-1} u^2)$ as the reduction of a $2^{2n}$-torsion class in $\Hsyn^{4d+1,2d}(X; \Z_2)$, which then necessarily vanishes. 
\qed

\section{Symplectic structure on Brauer groups}\label{sec: final}

Let $X$ be a smooth, proper, geometrically connected variety of dimension $2d$ over a finite field of characteristic $2$. We will complete the proof of Theorem \ref{thm: pseudo_main}, which (as already noted in Corollary \ref{cor: alternating_criterion}) implies \Cref{thm: MAT higher dimension}.

\subsection{Reductions to characteristic classes}\label{ssec:symplecticity-initial-reductions}
We wish to show that 
\[
u \cdot \beta_n (u)  = 0  \text{ for all } u \in \Hsyn^{2d,d}(X;\Z/2^n). 
\]
From Theorem \ref{thm: MAT form}, we have that
\begin{equation}\label{eq: final comp 1}
u \cdot \beta_n (u) = [2^{n-1}] \circ \Sqe^{2d} (\ol{\beta_n (u)}).
\end{equation}

Thanks to the position of the weight with respect to the degree of the Steenrod operation, \Cref{cor: operations agree} implies that
\begin{equation}
\Sqe^{2d} (\ol{\beta_n (u)}) = \Sqs^{2d} (\ol{\beta_n (u)}) \text{ for all } u \in \Hsyn^{2d,d}(X; \Z/2^n).
\end{equation}
Inserting this into \eqref{eq: final comp 1}, we have to show that 
\begin{equation}\label{eq: final comp 2.5}
[2^{n-1}] \circ \Sqs^{2d} (\ol{\beta_n (u)})  = 0  \quad \text{for all $u \in \Hsyn^{2d,d}(X;\Z/2^n)$}.
\end{equation}
From the definition of the syntomic Wu classes, we have 
\begin{equation}\label{eq: final comp 3}
\Sqs^{2d} (\ol{\beta_n (u)})   = v^{\syn}_{2d} \cdot \ol{\beta_n(u)}.
\end{equation}
It is immediate from the definition of $[2^{n-1}]$ that 
\begin{equation}\label{eq: final comp 4}
[2^{n-1}] (v^{\syn}_{2d} \cdot \ol{\beta_n(u)}) = ([2^{n-1}] v^{\syn}_{2d}) \cdot \beta_n(u) \text{ for all } u \in \Hsyn^{2d,d}(X;\Z/2^n). 
\end{equation}
Since $\beta_n$ is a derivation, we have 
\begin{equation}\label{eq: final comp 5}
([2^{n-1}]  v^{\syn}_{2d}) \cdot \beta_n (u)   = \beta_n(([2^{n-1}] v^{\syn}_{2d})  \cdot u )  - \beta_n( [2^{n-1}] v^{\syn}_{2d}  ) \cdot u. 
\end{equation}
Stringing together \eqref{eq: final comp 2.5}, \eqref{eq: final comp 3}, and \eqref{eq: final comp 4}, we see that it suffices to show that \eqref{eq: final comp 5} vanishes for all $u$. Since $([2^{n-1}] v^{\syn}_{2d})  \cdot u  \in \Hsyn^{4d,2d}(X; \Z/2^n)$, \Cref{lem: top boundary vanishes} implies that $\beta_n( ([2^{n-1}] v^{\syn}_{2d}  ) \cdot u) = 0$. Hence we have reduced to showing that 
\begin{equation}\label{eq: final comp 6}
\beta_n( [2^{n-1}] v^{\syn}_{2d}  )  = 0 \text{ for all } n \geq 1.
\end{equation}

\begin{lemma}
For all $v \in \Hsyn^*(X; \Z/2(b))$ we have 
\begin{equation}\label{eq: bocksteins equal}
\beta_n([2^{n-1}]v) =  \beta_{2,2^n}(v)
\end{equation}
where $\beta_{2,2^n}$ is the boundary map for the exact triangle of syntomic sheaves
\begin{equation}\label{eq: beta_{2,2^n}}
\Z/2^n (b)^{\syn}_X \xrightarrow{2}  \Z/2^{n+1}(b)^{\syn}_X  \rightarrow \Z/2(b)^{\syn}_X.
\end{equation}
\end{lemma}

\begin{proof}Consider the commutative diagram of exact triangles,
\[
\begin{tikzcd}
 \Z/2^n (b)^{\syn}_X  \ar[r, "2"] \ar[d, equals] &   \Z/2^{n+1} (b)^{\syn}_X  \ar[r] \ar[d, "2^{n-1}"] &  \Z/2 (b)^{\syn}_X   \ar[d, "2^{n-1}"]   \\
 \Z/2^n  (b)^{\syn}_X  \ar[r, "2^{n}"] & \Z/2^{2n}  (b)^{\syn}_X   \ar[r] & \Z/2^n (b)^{\syn}_X  
\end{tikzcd}
\]
Then \eqref{eq: bocksteins equal} results from comparing the two induced maps on cohomology from the top right to bottom left sheaves. 
\end{proof}

Applying \eqref{eq: bocksteins equal} in \eqref{eq: final comp 6}, it suffices to show that $\beta_{2,2^n}(v^{\syn}_{2d}) = 0$ for all $n$. This will be done in the next subsection. 

\subsection{Calculations with characteristic classes} So far, we have merely reformulated the question in terms of syntomic Wu classes. Now, we will draw on the earlier computations that express syntomic Wu classes in terms of characteristic classes.

\begin{lemma}\label{lem: squares of SW classes}
For any $w_j^{\syn} \in \Hsyn^{j,\lfloor j/2 \rfloor}(X)$, its syntomic Steenrod square $\Sqs^i(w_j^{\syn})$ can be expressed as a polynomial in the syntomic Stiefel--Whitney classes $\{w^{\syn}_{j'}\}$ with coefficients in $\F_p$. 
\end{lemma}

\begin{proof}Using the Adem relations from Proposition \ref{prop: adem p=2} and the Cartan formula from Proposition \ref{prop: Cartan}, this follows from the same inductive argument as in the proof of \cite[Lemma 5.4]{Feng20}.
\end{proof}

\begin{lemma}\label{lem: wu in terms of SW}
Every syntomic Wu class $v^{\syn}_j \in \Hsyn^{j,\lfloor j/2 \rfloor}(X)$ can be expressed as a polynomial in the syntomic Stiefel--Whitney classes $\{w^{\syn}_{j'}\}$. 
\end{lemma}

\begin{proof}
We induct on $j$. The base case is $v^{\syn}_0 = w^{\syn}_0 =1 \in \rH^{0,0}(X)$. Suppose $j=2i$ is even. Since $\Sqs^0 = \Id$, from \eqref{eq:arithmetic-wu-even} we have 
\begin{equation}\label{eq:w-v-relation}
v^{\syn}_{2i} + \Sqs^2 (v^{\syn}_{2i-2}) + \ldots  = w_{2i}^{\syn}.
\end{equation}
By the induction hypothesis, for $i \geq 1$ each term $v^{\syn}_{2i-2l}$ is a polynomial in the $\{w^{\syn}_{j'}\}_{j' \leq 2i-2l}$ with coefficients in $\F_p$, so by the Cartan formula and Lemma \ref{lem: squares of SW classes}, each $\Sqs^i (v^{\syn}_{2i-2l})$ is a polynomial in the $\{w^{\syn}_{j'}\}_{j' \leq 2i-2l}$ with coefficients in $\F_p$. Then solving for $v^{\syn}_{2i}$ in \eqref{eq:w-v-relation} completes the induction step if $j=2i$ is even. If $j$ is odd, then a similar argument works, using \eqref{eq:arithmetic-wu-odd} instead.
\end{proof}

\begin{cor}\label{cor: wu class lifts} For every $j \in \Z$, the syntomic Wu class $v^{\syn}_j \in \Hsyn^{j,\lfloor j/2 \rfloor}(X)$ is the reduction of a class in $\Hsyn^{j, \lfloor j/2 \rfloor}(X;\Z_2)$.
\end{cor}

\begin{proof}
Combine Lemma \ref{lem: wu in terms of SW} and Proposition \ref{prop: chern classes}. 
\end{proof}

\begin{proof}[Completion of the proof of Theorem \ref{thm: pseudo_main}] At the end of \S \ref{ssec:symplecticity-initial-reductions}, we reduced to showing that $\beta_{2,2^n}(v^{\syn}_{2d}) = 0$ for all $n$. Inspecting the long exact sequence associated to \eqref{eq: beta_{2,2^n}} whose boundary map is $\beta_{2,2^n}$, we see that this vanishing is equivalent to the property that $v^{\syn}_{2d}$ lifts mod $2^{n+1}$ for all $n$, which is guaranteed by Corollary \ref{cor: wu class lifts}. 
\end{proof}

\bibliographystyle{amsalpha}

\bibliography{Bibliography}

\providecommand{\bysame}{\leavevmode\hbox to3em{\hrulefill}\thinspace}
\providecommand{\MR}{\relax\ifhmode\unskip\space\fi MR }
\providecommand{\MRhref}[2]{%
  \href{http://www.ams.org/mathscinet-getitem?mr=#1}{#2}
}
\providecommand{\href}[2]{#2}
\begin{thebibliography}{AMMN22}

\bibitem[ACS19]{ACS19}
Edo Arad, Shachar Carmeli, and Tomer~M. Schlank, \emph{\'{E}tale {H}omotopy
  {O}bstructions of {A}rithmetic {S}pheres}, 2019.

\bibitem[AE25]{AE}
Toni Annala and Elden Elmanto, \emph{Motivic power operations at the
  characteristic via inﬁnite ramiﬁcation}, 2025.

\bibitem[AHI24]{AHI2}
Toni Annala, Marc Hoyois, and Ryomei Iwasa, \emph{Atiyah duality for motivic
  spectra}, 2024.

\bibitem[AHI25]{AHI1}
Toni Annala, Marc Hoyois, and Ryomei Iwasa, \emph{Algebraic cobordism and a
  {C}onner-{F}loyd isomorphism for algebraic {K}-theory}, J. Amer. Math. Soc.
  \textbf{38} (2025), no.~1, 243--289.

\bibitem[AI22a]{Kuniversal}
Toni Annala and Ryomei Iwasa, \emph{Motivic spectra and universality of $ k
  $-theory}, arXiv preprint arXiv:2204.03434 (2022).

\bibitem[AI22b]{MotSpec}
\bysame, \emph{Motivic spectra and universality of $ k $-theory}, arXiv
  preprint arXiv:2204.03434 (2022).

\bibitem[AMM22]{AMM22}
Benjamin Antieau, Akhil Mathew, and Matthew Morrow, \emph{The {K}-theory of
  perfectoid rings}, Doc. Math. \textbf{27} (2022), 1923--1952.

\bibitem[AMMN22]{AMMN22}
Benjamin Antieau, Akhil Mathew, Matthew Morrow, and Thomas Nikolaus, \emph{On
  the {B}eilinson fiber square}, Duke Math. J. \textbf{171} (2022), no.~18,
  3707--3806.

\bibitem[BC76]{BC76}
Edgar~H. Brown, Jr. and Michael Comenetz, \emph{Pontrjagin duality for
  generalized homology and cohomology theories}, Amer. J. Math. \textbf{98}
  (1976), no.~1, 1--27.

\bibitem[BCSY24]{ChromFourier}
Tobias Barthel, Shachar Carmeli, Tomer~M. Schlank, and Lior Yanovski, \emph{The
  chromatic {F}ourier transform}, Forum Math. Pi \textbf{12} (2024), Paper No.
  e8, 96.

\bibitem[Bei87]{Bei87}
A.~A. Beilinson, \emph{On the derived category of perverse sheaves},
  {$K$}-theory, arithmetic and geometry ({M}oscow, 1984--1986), Lecture Notes
  in Math., vol. 1289, Springer, Berlin, 1987, pp.~27--41.

\bibitem[BEM]{BEM}
Tom Bachmann, Elden Elmanto, and Matthew Morrow, \emph{$\mathbf{A}^1$-invariant
  motivic cohomology of schemes}.

\bibitem[Ben25]{Ben24}
Olivier Benoist, \emph{Steenrod operations and algebraic classes}, Tunis. J.
  Math. \textbf{7} (2025), no.~1, 53--89.

\bibitem[BH25]{BH25}
Tom Bachmann and Michael Hopkins, \emph{Stable operations in motivic homotopy
  theory}, 2025, Available at \url{https://tom-bachmann.com/ops.pdf}.

\bibitem[Bha22]{Bha22}
Bhargav Bhatt, \emph{Prismatic ${F}$-gauges}, 2022, Lecture notes, available at
  \url{https://www.math.ias.edu/~bhatt/teaching/mat549f22/lectures.pdf}.

\bibitem[BJ15]{BJ15}
Patrick Brosnan and Roy Joshua, \emph{Comparison of motivic and simplicial
  operations in mod-{$l$}-motivic and \'{e}tale cohomology}, Feynman
  amplitudes, periods and motives, Contemp. Math., vol. 648, Amer. Math. Soc.,
  Providence, RI, 2015, pp.~29--55.

\bibitem[BK25]{BK25}
Tess Bouis and Arnab Kundu, \emph{Beilinson--{L}ichtenbaum phenomenon for
  motivic cohomology}, 2025.

\bibitem[BL22]{BL22a}
Bhargav Bhatt and Jacob Lurie, \emph{Absolute prismatic cohomology}, 2022.

\bibitem[BM23]{BM23}
Bhargav Bhatt and Akhil Mathew, \emph{Syntomic complexes and {$p$}-adic
  \'{e}tale {T}ate twists}, Forum Math. Pi \textbf{11} (2023), Paper No. e1,
  26.

\bibitem[BMS19]{BMS2}
Bhargav Bhatt, Matthew Morrow, and Peter Scholze, \emph{Topological
  {H}ochschild homology and integral {$p$}-adic {H}odge theory}, Publ. Math.
  Inst. Hautes \'{E}tudes Sci. \textbf{129} (2019), 199--310.

\bibitem[Bou24]{Bou24}
Tess Bouis, \emph{Motivic cohomology of mixed characteristic schemes}, 2024.

\bibitem[BS22]{BS22}
Bhargav Bhatt and Peter Scholze, \emph{Prisms and prismatic cohomology}, Ann.
  of Math. (2) \textbf{196} (2022), no.~3, 1135--1275.

\bibitem[Car23]{carmeli2023strict}
Shachar Carmeli, \emph{On the strict picard spectrum of commutative ring
  spectra}, Compositio Mathematica \textbf{159} (2023), no.~9, 1872--1897.

\bibitem[CN17]{CN17}
Pierre Colmez and Wies\l~awa Nizio\l, \emph{Syntomic complexes and {$p$}-adic
  nearby cycles}, Invent. Math. \textbf{208} (2017), no.~1, 1--108.

\bibitem[Dri24]{Drin24}
Vladimir Drinfeld, \emph{Prismatization}, 2024.

\bibitem[Fen20a]{Feng20}
Tony Feng, \emph{\'{E}tale {S}teenrod operations and the {A}rtin-{T}ate
  pairing}, Compos. Math. \textbf{156} (2020), no.~7, 1476--1515.

\bibitem[Fen20b]{SHA}
Tony Feng, \emph{The {S}pectral {H}ecke {A}lgebra}, 2020.

\bibitem[FM87]{FM87}
Jean-Marc Fontaine and William Messing, \emph{{$p$}-adic periods and {$p$}-adic
  \'{e}tale cohomology}, Current trends in arithmetical algebraic geometry
  ({A}rcata, {C}alif., 1985), Contemp. Math., vol.~67, Amer. Math. Soc.,
  Providence, RI, 1987, pp.~179--207.

\bibitem[FS18]{FS18}
Martin Frankland and Markus Spitzweck, \emph{Towards the dual motivic
  {S}teenrod algebra in positive characteristic}, 2018.

\bibitem[Gei05]{Gei05}
Thomas Geisser, \emph{Motivic cohomology, {$K$}-theory and topological cyclic
  homology}, Handbook of {$K$}-theory. {V}ol. 1, 2, Springer, Berlin, 2005,
  pp.~193--234.

\bibitem[GL00]{GL00}
Thomas Geisser and Marc Levine, \emph{The {$K$}-theory of fields in
  characteristic {$p$}}, Invent. Math. \textbf{139} (2000), no.~3, 459--493.

\bibitem[GL01]{GL01}
\bysame, \emph{The {B}loch-{K}ato conjecture and a theorem of
  {S}uslin-{V}oevodsky}, J. Reine Angew. Math. \textbf{530} (2001), 55--103.

\bibitem[Gro85]{Gros85}
Michel Gros, \emph{Classes de {C}hern et classes de cycles en cohomologie de
  {H}odge-{W}itt logarithmique}, M\'{e}m. Soc. Math. France (N.S.) (1985),
  no.~21, 87.

\bibitem[HKOsr17]{HKO}
Marc Hoyois, Shane Kelly, and Paul~Arne \O~stv\ae r, \emph{The motivic
  {S}teenrod algebra in positive characteristic}, J. Eur. Math. Soc. (JEMS)
  \textbf{19} (2017), no.~12, 3813--3849.

\bibitem[HW19]{HW19}
Christian Haesemeyer and Charles~A. Weibel, \emph{The norm residue theorem in
  motivic cohomology}, Annals of Mathematics Studies, vol. 200, Princeton
  University Press, Princeton, NJ, 2019.

\bibitem[Ill79]{Ill79}
Luc Illusie, \emph{Complexe de de {R}ham-{W}itt et cohomologie cristalline},
  Ann. Sci. \'{E}cole Norm. Sup. (4) \textbf{12} (1979), no.~4, 501--661.

\bibitem[Jah15]{Jahn15}
Thomas Jahn, \emph{The order of higher {B}rauer groups}, Math. Ann.
  \textbf{362} (2015), no.~1-2, 43--54.

\bibitem[Kat87]{Kat87}
Kazuya Kato, \emph{On {$p$}-adic vanishing cycles (application of ideas of
  {F}ontaine-{M}essing)}, Algebraic geometry, {S}endai, 1985, Adv. Stud. Pure
  Math., vol.~10, North-Holland, Amsterdam, 1987, pp.~207--251.

\bibitem[KT03]{KT}
Kazuya Kato and Fabien Trihan, \emph{On the conjectures of {B}irch and
  {S}winnerton-{D}yer in characteristic {$p>0$}}, Invent. Math. \textbf{153}
  (2003), no.~3, 537--592.

\bibitem[LLR05]{LLR05}
Qing Liu, Dino Lorenzini, and Michel Raynaud, \emph{On the {B}rauer group of a
  surface}, Invent. Math. \textbf{159} (2005), no.~3, 673--676.

\bibitem[LM23]{LM23}
Morten L\"{u}ders and Matthew Morrow, \emph{Milnor {$K$}-theory of {$p$}-adic
  rings}, J. Reine Angew. Math. \textbf{796} (2023), 69--116.

\bibitem[Lur09]{HTT}
Jacob Lurie, \emph{Higher topos theory}, Princeton University Press, 2009.

\bibitem[Lur17]{HA}
\bysame, \emph{Higher {A}lgebra}, 2017, Version of Sep. 18, 2017, available at
  \url{https://www.math.ias.edu/~lurie/papers/HA.pdf}.

\bibitem[Lur18]{DAGXIII}
\bysame, \emph{Derived {A}lgebraic {G}eometry {X}{I}{I}{I}: {R}ational and
  p-adic {H}omotopy {T}heory}, 2018, Version of Dec. 15, 2011, available at
  \url{https://people.math.harvard.edu/~lurie/papers/DAG-XIII.pdf}.

\bibitem[Lur21]{DAGXI}
\bysame, \emph{Derived algebraic geometry xi: Descent theorems. 2011},
  Preprint, available at http://www. math. harvard. edu/lurie/papers/DAG-XI.
  pdf \textbf{4} (2021).

\bibitem[Man67]{Manin67}
Ju.~I. Manin, \emph{Rational surfaces over perfect fields. {II}}, Mat. Sb.
  (N.S.) \textbf{72 (114)} (1967), 161--192.

\bibitem[Man86]{Manin86}
Yu.~I. Manin, \emph{Cubic forms}, second ed., North-Holland Mathematical
  Library, vol.~4, North-Holland Publishing Co., Amsterdam, 1986, Algebra,
  geometry, arithmetic, Translated from the Russian by M. Hazewinkel.

\bibitem[May70]{May70}
J.~Peter May, \emph{A general algebraic approach to {S}teenrod operations}, The
  {S}teenrod {A}lgebra and its {A}pplications ({P}roc. {C}onf. to {C}elebrate
  {N}. {E}. {S}teenrod's {S}ixtieth {B}irthday, {B}attelle {M}emorial {I}nst.,
  {C}olumbus, {O}hio, 1970), Lecture Notes in Mathematics, Vol. 168, Springer,
  Berlin, 1970, pp.~153--231.

\bibitem[Mil75]{Milne75}
J.~S. Milne, \emph{On a conjecture of {A}rtin and {T}ate}, Ann. of Math. (2)
  \textbf{102} (1975), no.~3, 517--533.

\bibitem[Mil76]{Milne76}
\bysame, \emph{Duality in the flat cohomology of a surface}, Ann. Sci.
  \'{E}cole Norm. Sup. (4) \textbf{9} (1976), no.~2, 171--201.

\bibitem[Mil86]{Milne86}
\bysame, \emph{Values of zeta functions of varieties over finite fields}, Amer.
  J. Math. \textbf{108} (1986), no.~2, 297--360.

\bibitem[MS74]{MS74}
John~W. Milnor and James~D. Stasheff, \emph{Characteristic classes}, Princeton
  University Press, Princeton, N. J.; University of Tokyo Press, Tokyo, 1974,
  Annals of Mathematics Studies, No. 76.

\bibitem[Niz98]{Niz98}
Wies\l~awa Nizio\l, \emph{Crystalline conjecture via {$K$}-theory}, Ann. Sci.
  \'{E}cole Norm. Sup. (4) \textbf{31} (1998), no.~5, 659--681.

\bibitem[NS18]{NS18}
Thomas Nikolaus and Peter Scholze, \emph{On topological cyclic homology}, Acta
  Math. \textbf{221} (2018), no.~2, 203--409.

\bibitem[Poo17]{Poo17}
Bjorn Poonen, \emph{Rational points on varieties}, Graduate Studies in
  Mathematics, vol. 186, American Mathematical Society, Providence, RI, 2017.

\bibitem[Pri20]{Pri20}
Eric Primozic, \emph{Motivic {S}teenrod operations in characteristic {$p$}},
  Forum Math. Sigma \textbf{8} (2020), Paper No. e52, 25.

\bibitem[PS99]{PS99}
Bjorn Poonen and Michael Stoll, \emph{The {C}assels-{T}ate pairing on polarized
  abelian varieties}, Ann. of Math. (2) \textbf{150} (1999), no.~3, 1109--1149.

\bibitem[Pst23]{pstrkagowski2023synthetic}
Piotr Pstr{\k{a}}gowski, \emph{Synthetic spectra and the cellular motivic
  category}, Inventiones mathematicae \textbf{232} (2023), no.~2, 553--681.

\bibitem[Sch12]{Sch12}
Peter Scholze, \emph{Perfectoid spaces}, Publ. Math. Inst. Hautes \'{E}tudes
  Sci. \textbf{116} (2012), 245--313.

\bibitem[SGA73]{SGA4}
\emph{Th\'{e}orie des topos et cohomologie \'{e}tale des sch\'{e}mas. {T}ome
  3}, Lecture Notes in Mathematics, Vol. 305, Springer-Verlag, Berlin-New York,
  1973, S\'{e}minaire de G\'{e}om\'{e}trie Alg\'{e}brique du Bois-Marie
  1963--1964 (SGA 4), Dirig\'{e} par M. Artin, A. Grothendieck et J. L.
  Verdier. Avec la collaboration de P. Deligne et B. Saint-Donat.

\bibitem[SS24]{SS24}
Federico Scavia and Fumiaki Suzuki, \emph{Two coniveau filtrations and
  algebraic equivalence over finite fields}, 2024.

\bibitem[{Sta}24]{stacks-project}
The {Stacks project authors}, \emph{The stacks project},
  \url{https://stacks.math.columbia.edu}, 2024.

\bibitem[SV00]{SV98}
Andrei Suslin and Vladimir Voevodsky, \emph{Bloch-{K}ato conjecture and motivic
  cohomology with finite coefficients}, The arithmetic and geometry of
  algebraic cycles ({B}anff, {AB}, 1998), NATO Sci. Ser. C Math. Phys. Sci.,
  vol. 548, Kluwer Acad. Publ., Dordrecht, 2000, pp.~117--189.

\bibitem[Tan24a]{Tang24}
Longke Tang, \emph{Slicing criterion for ind-smooth ring maps}, 2024.

\bibitem[Tan24b]{Tang22}
Longke Tang, \emph{Syntomic cycle classes and prismatic {P}oincar\'{e}
  duality}, Compos. Math. \textbf{160} (2024), no.~10, 2322--2365.

\bibitem[Tat63]{Tate62}
John Tate, \emph{Duality theorems in {G}alois cohomology over number fields},
  Proc. {I}nternat. {C}ongr. {M}athematicians ({S}tockholm, 1962), Inst.
  Mittag-Leffler, Djursholm, 1963, pp.~288--295.

\bibitem[Tat95]{Tate66}
\bysame, \emph{On the conjectures of {B}irch and {S}winnerton-{D}yer and a
  geometric analog}, S\'eminaire {B}ourbaki, {V}ol.\ 9, Soc. Math. France,
  Paris, 1995, pp.~Exp.\ No.\ 306, 415--440.

\bibitem[Tsu99]{Tsu99}
Takeshi Tsuji, \emph{{$p$}-adic \'{e}tale cohomology and crystalline cohomology
  in the semi-stable reduction case}, Invent. Math. \textbf{137} (1999), no.~2,
  233--411.

\bibitem[Ura96]{Urabe96}
Tohsuke Urabe, \emph{The bilinear form of the {B}rauer group of a surface},
  Invent. Math. \textbf{125} (1996), no.~3, 557--585.

\bibitem[Voe02]{Voe02}
Vladimir Voevodsky, \emph{Open problems in the motivic stable homotopy theory.
  {I}}, Motives, polylogarithms and {H}odge theory, {P}art {I} ({I}rvine, {CA},
  1998), Int. Press Lect. Ser., vol.~3, Int. Press, Somerville, MA, 2002,
  pp.~3--34.

\bibitem[Voe03a]{voevodsky2003zeroslicespherespectrum}
Vladimir Voevodsky, \emph{On the zero slice of the sphere spectrum}, 2003.

\bibitem[Voe03b]{Voe03}
Vladimir Voevodsky, \emph{Reduced power operations in motivic cohomology},
  Publ. Math. Inst. Hautes \'{E}tudes Sci. (2003), no.~98, 1--57.

\bibitem[Voe10]{Voe10}
\bysame, \emph{Motivic {E}ilenberg-{M}ac{L}ane spaces}, Publ. Math. Inst.
  Hautes \'{E}tudes Sci. (2010), no.~112, 1--99.

\end{thebibliography}

\end{document}